\documentclass[12pt]{amsart}
\usepackage{amssymb}
\usepackage{amsmath}
\usepackage{amsopn}
\usepackage{amscd}
\catcode`\@=11
\def\marginnote#1{}

\def\appname{Appendix}
\newcounter{app}
\def\theapp{\Alph{app}}
\def\app{\par
   \addvspace{4ex}
   \@afterindentfalse
  \secdef\@app\@dapp}
\def\@app[#1]#2{\ifnum \c@secnumdepth >\m@ne
        \refstepcounter{app}
        \addcontentsline{toc}{app}{\theapp
        \hspace{1em}#1}\else
      \addcontentsline{toc}{app}{ #1}\fi
   {\parindent \z@ \raggedright
    \Large \bf \appname~\theapp .
   \Large  \bf \hspace{1em}    #2}\nobreak
   \vskip -2ex   \noindent
\setcounter{equation}{0}
\def\theequation{\Alph{app}.\arabic{equation}}}
\def\@dapp#1{%
{\parindent \z@ \raggedright  \bf #1}\par\nobreak}
\def\l@app#1#2{\addpenalty{\@secpenalty}%
   \addvspace{1em plus\p@}%
   \begingroup
   \@tempdima 3em
     \parindent \z@ \rightskip \@pnumwidth
     \parfillskip -\@pnumwidth
     { \bf
     \leavevmode
     #1\hfil \hbox to\@pnumwidth{\hss #2}}\par
     \nobreak
   \endgroup}
\def\mylangle{\langle}
\def\myrangle{\rangle}
\def\lform{(}
\def\rform{)}

\DeclareMathOperator*{\ooplus}{{\oplus}}
\def\j{\mbox{${\scriptstyle{\mathbb J}}$}}

\def\UG{{\mathcal W\,}}
\newcommand{\MD}{\mbox{${\mathcal M}$}}
\def\DBG{D^{b}}
\def\DB{D^{b}}

\def\FF{{\mathcal F}}
\def\OS{{\mathcal A}}
\newcommand{\xxi}[1]{\vec{\xi}_{#1}\ }
\newcommand{\f}[1]{\vec{f}_{#1} }
\newcommand{\ff}[1]{\overline{f}_{#1}}
\newcommand{\QG}{\mbox{$Qui_{\Gamma,\{G_\b,\widetilde{G}_\b,G_{\a,\b}\}}$}}
\newcommand{\QGn}{\mbox{$Qui_{\Gamma^n,\{ G_\b,\widetilde{G}_\b,G_{\a,\b}\}}$}}

\newcommand{\rf}[1]{(\ref{#1})}
\newcommand{\FN}{\mbox{$F_{}\,$}}
\newcommand{\TN}{\mbox{$T_{}\,$}}
\newcommand{\T}[1]{\mbox{$T_{#1}$\,}}

\newcommand{\Ta}{\mbox{$T_{\alpha}\,$}}
\newcommand{\Fa}{\mbox{$F_{\alpha}\,$}}
\newcommand{\TTa}{\mbox{$T'_{\alpha}\,$}}
\newcommand{\FFa}{\mbox{$F'_{\alpha}\,$}}
\newcommand{\TTb}{\mbox{$T'_{\beta}\,$}}
\newcommand{\FFb}{\mbox{$F'_{\beta}\,$}}

\newcommand{\TTT}[1]{\mbox{${\mathcal T}_{#1}\,$}}

\newcommand{\Tb}{\mbox{$T_{\beta}$\,}}
\newcommand{\Fb}{\mbox{$F_{\beta}$\,}}

\def\Om{{\Omega}}
\def\om{\overline{\omega}}
\def\Oman{{\Omega}^{an}}
\def\Oan{\mathcal{O}^{an}}
\def\Dan{\mathcal{D}^{an}}

\newcommand{\Oa}{\mbox{$\overline{\Omega}_{\alpha}\,$}}

\newcommand{\OOa}{\mbox{$\widehat{\Omega}_{\alpha}\,$}}
\newcommand{\Ob}{\mbox{$\overline{\Omega}_{\beta}$\,}}

\newtheorem{proposition}{Proposition}[section]
\newtheorem{theorem}[proposition]{Theorem}
\newtheorem{lemma}[proposition]{Lemma}
\newtheorem{corollary}[proposition]{Corollary}
\renewcommand{\theequation}{{\thesection}.{\arabic{equation}}}
\def\Hom {{\rm Hom}}
\def\a{\alpha}                        \def\d{\delta} \def\dd{{d^*}}
\def\b{\beta}                         \def\D{{\mathcal D}}
\def\g{\gamma}
\def\l{\lambda}                       
\def\L{\Lambda}                       \def\o{\omega}
\def\ii {i}

\def\ho{\widehat{\omega}}                 \def\ot{\otimes}

\newcommand\der[1]{\mbox{$\frac{\partial}{\partial_{#1}}$}}
\let\cal=\mathcal
 
\def\ve{\varepsilon}
\def\h{{\rm\hbar}}

\def\A{{\cal C}}              
\def\X{\overline{X}}
\def\cXa{\overline{X}_\a}
\def\cXb{\overline{X}_\b}
\def\Xa{{\overline{X}_\a}}
\def\Xb{{\overline{X}_\b}}

\def\codim{{\rm codim\ }}  \def\dim{{\rm dim\ }}

\def\G{\Gamma}
\def\Ga{\Gamma}

\def\ee{\emptyset}

\def\VA{{\cal V}}
\def\VG{\VA_\G}
\def\VGn{\VA_{\G^n}}

\def\EVA{E\VA}
\def\QVA{{\cal Q}\VA}

\def\Qui{Qui}
\newcommand{\CC}{{\mathbb C}}
\newcommand{\NN}{{\mathbb N}}

\newcommand{\ZZ}{{\mathbb Z}}
\let\C=\CC
\def\R{{\mathcal R}}
\def\R{\overline{{\mathcal R}}}

\def\v{I}
\def\e{E}
\def\H{H}
\def\U{X}
\def\B{B}
\def\risom{\stackrel{\thicksim}{\to}}
\def\R{{\mathcal R}}
\def\Rlong{\overline{{\mathcal R}}}
\def\GR{\widetilde{{\mathcal R}}}
\def\GGR{\widehat{{\mathcal R}}}

\newcommand{\GM}[1]{\widetilde{#1}}
\newcommand{\GGM}[1]{\widehat{#1}}

\def\ifa{I(F_\alpha)}
\def\ita{I(T_\alpha)}

\def\V{\overline{V}}

\newcommand{\TG}[1]{L_{#1}(g)}
\newcommand{\TXG}[2]{L_{#1}({#2})}

\def\forget{\mathcal N}
\def\DD{D}
\def\mycircle{{\mbox{\begin{picture}(3.5,1.5)\put(1.5,0.5){\circle*{1.5}}
\end{picture}}}}
\def\upcircle{{\mbox{\begin{picture}(3.5,1.5)\put(1.5,2.0){\circle*{1.5}}
\end{picture}}}}
\def\my!{{\,!\,}}
\def\inverse{\mycircle}
\def\ggg{{\mathfrak g}}
\def\h{{\mathbf h}}

\def\Abb1{{\cal A}} \def\Ab2{{\cal B}}
\def\hh{{\mathfrak h}}
\def\nn{{\mathfrak n}}

\def\nnn{\overline{\nn}}
\def\gggg{\overline{\ggg}}
\def\ds{\displaystyle}
\def\Ker{{\text{\rm Ker}\,}}
\def\<{\langle}
\def\>{\rangle}
\def\gr{{\rm gr}\,}
\def\Gr{{\rm gr}\,}
\def\CG{\overline{G}}
\def\Tprime{\bar{T}}

\newcommand{\tE}{\tilde E}
\newcommand{\tVA}{\tilde {\VA}}

\begin{document}
\begin{center}
\bigskip\bigskip
{\Large\bf Quiver $D$-Modules and Homology of Local
Systems
over an Arrangement of Hyperplanes}\\
\bigskip
\bigskip
{\bf S. Khoroshkin$^{\circ}$ \ \ and \ A. Varchenko$^{\star}$
}\bigskip\\
$^\circ${\it Institute of Theoretical and Experimental Physics, 117259
Moscow, Russia}\\
$^\star${\it Department of Mathematics, University of North Carolina at
Chapel Hill, NC 27599-3250, USA}\\
\end{center}
\begin{center}
{\bf Abstract}
\medskip
\end{center}

{\small{Let $\A$ be an arrangement of hyperplanes in $\CC^N$, $D$
the ring of algebraic differential operators on $\CC^N$. We
define a category of quivers associated with $\A$. A quiver is a collection
of vector spaces, attached to strata of the arrangement,
 and suitable linear maps
between the  spaces . To a quiver we assign a $D$-module on $\CC^N$,
called a quiver $D$-module. We describe
 basic operations for 
$D$-modules in terms of linear algebra of
quivers. We give an explicit construction of a free resolution of a quiver
$D$-module and use the construction  to describe
the associated perverse sheaf.
As an application, we calculate
the cohomology of $\CC^N$ with coefficients in
the quiver perverse sheaf (under certain assumptions).

}}
{\small \tableofcontents

}

\section{Introduction}
 {\bf 1}. It is known that  nontrivial homological characteristics
of a space can
 be expressed in terms of linear algebra, if the space is
an affine space with an arrangement of hyperplanes. An
example is the description by Arnold, Brieskorn, Orlik and Solomon of
 the cohomology groups of the complement to the arrangement in terms of
 combinatorics of
intersections of hyperplanes \cite{Ar,Brs,OS}.
On the other hand, holonomic
$D$-modules with regular singularities and constructible perverse sheaves
on algebraic manifolds also admit combinatorial description,
 according to  Beilinson, Kashiwara, Malgrange \cite{Be1,Kash2,M}.
In this paper we  combine the two points of view and describe
a certain category of holonomic $D$-modules and associated perverse
sheaves, attached to an arrangement of hyperplanes in an affine space.

\bigskip

{\bf 2}.
 Let $\A$ be an arrangement of  hyperplanes in a
complex affine space $X$. The intersections of hyperplanes
 define in $X$ the structure of a stratified space. The strata
 and their adjacencies are described by the graph $\Ga$ of the arrangement.
The vertices $\a\in \v(\G)$  correspond to strata $X_\a$.
 If the closure $\overline{X}_\a$ of a
stratum $X_\a$ contains a stratum $X_\b$ and
$\codim_{\overline{X}_\a}X_\b =1$, then the
vertices $\a$ and $\b$ are connected by an edge and we write
$\a\succ\b$.
To each vertex $\a\in\v(\G)$  we attach a nonnegative integer
$l(\a)$, equal to the codimension of $X_\a$ in $X$.

We assign to the graph $\Ga$  quivers with relations,
which form an abelian category $Qui_{\Ga}$. A quiver $\VA\in Qui_\G$
consists of a collection of finite-dimensional
vector spaces $V_\a$, $\a\in\v(\G)$, and linear
maps  $A_{\a,\b}: V_\b\to V_\a,$ $\a,\b\in\v(\G)$,
such that
\begin{enumerate}
\item[(a)]   $A_{\a,\b}=0$ if $\a$ and $\b$ are  not adjacent;
\item[(b)]
$ \sum_\b A_{\a,\b}A_{\b,\g} = 0 $
 if $|l(\a)-l(\g)|=2$ ;
\item[(c)]
$\sum_\b A_{\a,\b}A_{\b,\g} = 0$ if
$ l(\a)=l(\g),\;\a\not=\g $, and there exists $ \d\in\v(\G)$,
  adjacent to $\a$ and $\g$, such that $l(\d)=l(\a)+1$.
\end{enumerate}

\medskip
{\it Example 1.1}. Let $X$ be $\CC$ with coordinate $z$. Let $\A=\{z=0\}$.
The space $X$ has two
strata: $X_\ee=\{z\not= 0\}$ and $X_\a=\{0\}$. The graph $\Ga$ has two
vertices, $\ee$ and $\a$, connected by an edge, $l(\a)=1, \ l(\ee)=0$.
 A quiver $\VA\in Qui_\Ga$
is a collection of two vector spaces, $V_\ee$ and $V_\a$, and two linear
maps, $A_{\ee,\a}: V_\a\to V_\ee$ and $A_{\a,\ee}: V_\ee\to V_\a$.

\medskip

Let $D_X$ denote the ring of algebraic differential operators
on $X$ and  $\D_X$ the corresponding sheaf of rings on  $X$.
 To a quiver $\VA\in Qui_\Ga$ we attach a holonomic
$D_X$-module $E\VA$. This $D_X$-module is given by generators and
relations. The generators are elements of the vector space
$\oplus_{\a\in\v(\G)} {\overline{V}}_\a$.
 Here  ${\overline{V}}_\a={\Oa}\ot V_\a$, and  ${\Oa}$ denotes the
one-dimensional vector space of top exterior forms on $\overline{X}_\a$,
invariant with respect to affine translations along $\overline{X}_\a$.
The relations in $E\VA$ are:
\begin{enumerate}
\item[(a)] for any generator $\om_\a\ot v_\a\in\overline{ V}_\a$
and any vector field $\xi$ with constant coefficients, which is parallel to
$X_\a$, the element  $\xi(\om_\a\ot v_\a)$ equals to a certain
linear combination
of the vectors $\om_\b\ot A_{\b,\a}(v_\a)$, $\b\in\v(\G),\ \a\succ\b$,
\item[(b)] for any
 affine function
$f$,  vanishing on $X_\a$, the element $f(\om_\a\ot v_\a)$
equals to a certain  linear combination
of the vectors $\om_\b\ot A_{\b,\a}(v_\a)$, $\b\in\v(\G),\ \b\succ \a$,
\end{enumerate}
 see formulas \rf{2.3}, \rf{2.4}.

\medskip \noindent
{\it Example 1.2}.
In  Example 1.1,  the $D_{X}$-module $E\VA$, assigned to the quiver
$\,\VA=\{V_\ee, V_\a,$ $ A_{\ee, \a},$ $ A_{\a,\ee} \}$
is generated by  vectors
$dz\ot v_\ee\in {\overline{V}}_\ee$, where $v_\ee\in V_\ee$, and
$1\ot v_\a\in {\overline{V}}_\a$, where $v_\a\in V_\a$.
The relations are: $\frac d {d z}\, (dz\ot v_\ee)=1\ot A_{\a,\ee}(v_\ee)$
 for any $v_\ee\in V_\ee$ and
$z\, (1\ot v_\a)=dz\ot A_{\ee,\a}(v_\a)$ for any $v_\a\in V_\a$.

\medskip
The assignment $\VA\mapsto E\VA$ defines a functor $E$ from the category
$Qui_\G$ to the category $\MD^{hol}_X$ of holonomic $\D_X$-modules.
We describe  a full subcategory in $Qui_\G$ of non-resonant quivers, such
that the restriction of $E$ to this subcategory  is full and faithful.

 For $n=0,...,N$, let $X^n$ denote the complement in $X$ to the union of all
strata of codimension greater than $n$. This set $X^n$ is
 called the  principal open subset.
 The principal subset $X^0$ is the complement to the union of hyperplanes of the
arrangement. All $X^n$ are stratified spaces. The adjacency structure
of the stratification of $X^n$ is described by the  graph  $\G^n$,
 which is a
suitable truncation of $\Ga$. We introduce an abelian category
 $Qui_{\G^n}$
of quivers of $\Ga^n$. The quivers of $Qui_{\G^n}$ are also called level
$n$ quivers. To  each level $n$ quiver
$\VA_{\G^n}\in Qui_{\G^n}$ we assign a sheaf of $\D_{X^n}$-modules, denoted by
$E^n\VA_{\G^n}\in  \MD^{hol}_{X^n}$.

\medskip \noindent
{\it Example 1.3}.
In  Example 1.1,  we have $X^0=X_\ee=\{z\not=0\}$, $X^1=X=\CC$, $\G^1=\G$.
The graph $\G^0$ has one vertex $\ee$ and one edge, which is a
loop. A quiver $\VA_{\G^0}\in Qui_{\G^0}$ consists of a vector space $V_\ee$
and a linear map $A_\ee^\ee: V_\ee\to V_\ee$. The $D_{X^0}$-module
$E^0\VA_{\G^0}$, assigned to this quiver, is generated by the elements
$dz\ot v_\ee \in {\overline{V}}_\ee$, subject to the relations
$z \frac d {d z}\,(dz\ot v_\ee) = dz\ot A_\ee^\ee(v_\ee)$.

\medskip
For $k\leq l$ we have inclusions
$j_{l,k}: X^{k}\to X^l$ of principal open subsets and related functors
$j_{l,k}^*$ and $j_{l,k,*}$ of inverse and direct images of $\D$-modules.
We define quiver inverse and direct image functors
 $\j_{l,k}^* : \Qui_{\G^l} \to \Qui_{\G^k}$ and
 $\j_{l,k,*} : \Qui_{\G^k} \to \Qui_{\G^l}$, respectively. The functor
$\j_{l,k}^*$ is left adjoint to the functor $\j_{l,k,*}$. We have an
isomorphism of $\D_{X^k}$-modules
 $E^k\j_{l,k}^*(\VA_{\Ga^l}) \approx j_{l,k}^*E^l(\VA_{\Ga^l})$
for any level $l$ quiver $\VA_{\Ga^l}$. We construct an isomorphism of
$\D_{X^l}$-modules
  $j_{l,k,*}E^k (\VA_{\G^k})$ and
$ E^l \j_{l,k,*}(\VA_{\G^k})$  under certain
non-resonance conditions on the level $k$ quiver $\VA_{\G^k}$.

\medskip \noindent
{\it Example 1.4}.
In  Examples 1.1 and 1.3, the level zero quiver $\j^*_{1,0}(\VA)$, associated
to a level one quiver $\VA = \{V_\ee, V_\a, A_{\ee, \a}, A_{\a,\ee} \}$,
is  $ \{U_\ee, B_\ee^\ee\}$, where
$U_\ee=V_\ee$, and $ B_\ee^\ee=A_{\ee,\a}A_{\a,\ee}$. The level one quiver
$\j_{1,0,*}(\VA_{\Ga^0})$, associated to a level zero quiver
$\VA_{\G^0} = \{V_\ee,  A_{\ee}^{\ee} \}$, is
 $\{W_\ee, W_\a, C_{\ee, \a}, C_{\a,\ee} \}$,
where $W_\ee=V_\ee$, $W_\a=V_\ee$, $C_{\ee, \a}={\rm Id}_{V_\ee}$,
$C_{\a,\ee}=A_\ee^\ee$. The quiver $\j_{1,0,*}(\VA_{\Ga^0})$ is non-resonant,
if any two eigenvalues of the operator $A_\ee^\ee$ do not differ by a nonzero
 integer, and the only possible integer eigenvalue of $A_\ee^\ee$ is zero.

\bigskip

There is a {\it specialization} construction in the theory of $\D$-modules due to
Kashi\-wa\-ra \cite{Kash2}.
Let $Y$ be a smooth complex algebraic variety, $Z\subset Y$ a
subvariety,  $T_ZY \to Z$  the normal bundle.
Given a $\D_{Y}$-module $M$ and some additional data,
the specialization construction produces a  $\D$-module on  the total
space $T_ZY$ of the normal bundle. This $\D_{TY_Z}$-module is called
{\it the specialization} of $M$ to $Z$ and denoted by $Sp_Z(M)$.

The specialization construction gives a collection of functors between quiver
$\D$-modules.

We consider the case when $Y$ is the affine
space $X=\CC^N$ with a central arrangement $\A$ and
$M$ is the quiver $\D$-module $\EVA$ associated with a quiver $\VA$ of $\A$.
Then it turns out that
for any closed stratum $\X_\a$ of $X$, the specialization
$Sp_{\X_\a}(\EVA)$ is isomorphic to a suitable quiver $\D$-module on
$T_{\X_\a}X$.

More precisely, we construct an arrangement of hyperplanes $\A_\a$ in
$T_{\X_\a}X$ and a quiver $Sp_\a(\VA)$ of the arrangement $\A_\a$ such that
the specialization $Sp_{\X_\a}(\EVA)$ is isomorphic to the quiver $\D$-module
$ESp_\a(\VA) $ associated with the quiver  $Sp_\a(\VA)$, if the linear maps of the
initial quiver $\VA$ satisfy certain non-resonance conditions.

\medskip
\noindent
{\it Example 1.5}.
In Examples 1.1
and 1.2, we identify the space $T_{\X_\a}X$ with $X=\CC$. Then
the arrangement $\A_\a$ is identified with $\A$, the quiver
$Sp_\a(\VA)$ is identified with $\VA$. Introduce the linear operator
$$
S = A_{\emptyset,\a}A_{\a,\emptyset}+
A_{\a,\emptyset}A_{\emptyset,\a}\ :\
V_\emptyset\oplus V_\a\
\to\ V_\emptyset\oplus V_\a\ .
$$
If any two eigenvalues of $S$ do not differ by a nonzero integer and if
$S$ has no nonzero integer eigenvalues,
then the  $\D_{T_{\X_\a}}$-modules
$Sp_{\X_\a}(\EVA)$ and $ ESp_\a(\VA)$ are isomorphic.

\bigskip
 For any
quiver $\VA_{\Ga}\in Qui$, we construct a free $\D_X$-module resolution
$\R E\VA_{\G}$ of the quiver $\D_X$-module $E\VA_{\G}$.

\medskip \noindent
{\it Example 1.6}. In Example 1.2, the free $\D_X$-module resolution $\R E\VA$
of the quiver $\D_X$-module $E\VA$ is the complex
$$
\R E\VA\ :\  0\to\R^{-1} E\VA\ \stackrel{d}{\to}\  \R^{0}
 E\VA\ \stackrel{\nu}{\to}\  E\VA\ \to\ 0\ .
$$
 Here
$\R^{0} E\VA$ is the free $\D_X$-module
$$
\D_X\ot {\overline{V}}_\ee
 \
\oplus\
 \D_X\ot {\overline{V}}_\a\ .
$$
The term $\R^{-1} E\VA$ is the
 free $\D_X$-module
$$
\D_X\ot T_\ee\ot {\overline{V}}_\ee\
 \oplus\
 \D_X\ot F_\a\ot {\overline{V}}_\a\ ,
$$
where $T_\ee$ is the
 one-dimensional space $\CC\cdot \frac d {d z}$, and
 $F_\a$ is the
one-dimensional space $\CC\cdot  z$ of affine functions, vanishing on
$X_\a$. The differentials $d$, $\nu$ are given by the formulas:
\begin{align*}
&d(D\ot \frac d {d z}\ot dz \ot v_\ee)=D\cdot \frac d {d z}
 \ot dz\ot v_\ee -D\ot 1\ot A_{\a,\ee}(v_\ee)\, ,&&\\
&d(D\ot   z\ot 1 \ot v_\a)=D\cdot  z
 \ot 1\ot v_\a -D\ot dz\ot A_{\ee,\a}(v_\a)\,,&&\\
&\nu(D\ot 1 \ot v_\a)=D\cdot(1\ot v_\a), \qquad
\nu(D\ot dz \ot v_\ee)=D\cdot(dz\ot v_\ee)\, ,&&
\end{align*}
 for any $D\in D_X$, $v_\ee\in V_\ee$, $v_\a\in V_\a$. Here $D\cdot(1\ot v_\a)$
and $D\cdot(dz\ot v_\ee)$ denote the application of $D$ to elements of
$E\VA$.

\bigskip
The de Rham functor transforms holonomic $\D_X$-modules to perverse
sheaves \cite{Kash,Meb}, see also \cite{B}.
 The application
of the de Rham functor to the complex $\R E\VA_{\G}$ gives an explicit
description of the associated perverse sheaf $\QVA$,  as a complex
of free ${\mathcal O}_X$-modules.

\medskip \noindent
{\it Example 1.7}. In Example 1.2,
the perverse sheaf $\QVA$, associated with
$E\VA$, is given by the complex
$$
\QVA\ :\  0\ \to\ \QVA^{-1}\ \stackrel{d}{\to}\ \QVA^{0}\ \to\ 0\ .
$$
Here
$$
\QVA^{-1}\ =\ \Om_X^{an}\ot T_\ee\ot {\overline{V}}_\ee\
\oplus\
\Om_X^{an}\ot F_\a\ot {\overline{V}}_\a \ ,
$$
and
 $$
\QVA^{0}\ =\ \Om_X^{an}\ot {\overline{V}}_\ee\
\oplus\
\Om_X^{an}\ot {\overline{V}}_\a \ ,
$$
where $\Om_X^{an}$ denotes the space of holomorphic
one-forms on $X$. The differential $d$ is given by the formulas:
\begin{align*}
d(\omega\ot  \frac d {d z}\ot dz \ot v_\ee)&=
-Lie_{\frac d {d z}}(\omega)
 \ot dz\ot v_\ee -\omega\ot 1\ot A_{\a,\ee}(v_\ee) ,\\
d(\omega\ot   z\ot 1 \ot v_\a)&=z\omega
 \ot 1\ot v_\a -\omega\ot dz\ot A_{\ee,\a}(v_\a),
\end{align*}
 for any $\omega\in\Om_X^{an}$, $v_\ee\in V_\ee$, $v_\a\in V_\a$.
Here $Lie_{\frac d {d z}}(\omega)$ denotes the Lie derivative of
$\omega$ in the direction of the vector field $\frac d {d z}$.

\bigskip
Let
$\Delta_X:\MD_X^{hol}\to \MD_X^{hol}$ and
$\Delta_{X^0}:\MD_{X^0}^{hol}\to \MD_{X^0}^{hol}$ be the duality functors
on holonomic $\D$-modules \cite{Kash3, KK}, see also \cite{Kash4,B}. For any $k=0,1...,N$,
 we define an anti-automorphism $\tau_k$
of the category $Qui_{\Ga^k}$. We establish an isomorphism
$\Delta_{X^0} (  E^0(\VA_{\G^0}))\, \approx \, E^0(\tau_0 (\VA_{\G^0}))$
of $\,\D_{X^0}$-modules for any level zero quiver $\VA_{\G^0}$ and an
isomorphism
$\Delta_{X} (  E(\VA_{\G}))$ $ \approx \ E(\tau_N (\VA_{\G}))$
of $\,\D_{X}$-modules for any level $N$ quiver $\VA_{\G}$.

\medskip \noindent
{\it Example 1.8}. In  Example 1.3, for any level zero quiver
 $\VA_{\G^0} = \{V_\ee,  A_{\ee}^{\ee} \}$, the quiver $\tau_0 \VA_{\G^0}$
 is $ \{U_\ee,  B_{\ee}^{\ee}
\}$, where $U_\ee=V_\ee^*$ is the dual space  and $B_{\ee}^{\ee}=\left(
 A_{\ee}^{\ee}\right)^t$ is the transpose operator.

\bigskip

In the theory of $\D$-modules, there are 
two functors of direct
image, $j_{l,k,!}: \MD_{X^k}^{hol}\to \MD_{X^l}^{hol}$, and
$j_{l,k,*}: \MD_{X^k}^{hol}\to \MD_{X^l}^{hol}$, associated to the
inclusion $j_{l,k}: X^k\to X^l$, $k<l$. In addition to the functor
$\j_{l,k,*}: Qui_{\G^k}\to Qui_{\G^l}$,
we introduce the second functor
of quiver direct image $\j_{l,k,!}: Qui_{\G^k}\to Qui_{\G^l}$, and construct
an isomorphisms of $\D_X$ modules
$\j_{N,0!}\VA_{\G^0}$ and $\tau_N \j_{N,0,*}\tau_0\VA_{\G^0}$ for a
 non-resonant quiver $\j_{N,0,!}\VA_{\G^0}$.

\medskip
\noindent
{\it Example 1.9}.
In  Example 1.3, for any level zero quiver
 $\VA_{\G^0} = \{V_\ee,  A_{\ee}^{\ee} \}$,  the quiver
$\j_{1,0,!}(\VA_{\Ga^0})$,  is
 $\{W'_\ee, W'_\a, C'_{\ee, \a}, C'_{\a,\ee} \}$,
where $W'_\ee=V_\ee$, $W'_\a=V_\a$, $C'_{\a, \ee}={\rm Id}_{V_\ee}$,
$C'_{\ee,\a}=A_\ee^\ee$. It is non-resonant if and only if the quiver
$\j_{1,0,*}(\VA_{\Ga^0})$ is non-resonant.

\medskip

Let $z_1,...,z_N$ be linear coordinates on $X=\CC^N$ and $\xi_1,...,
\xi_N$ the dual coordinates on the dual space $X^*$.
The assignment
$
\xi_i \,\mapsto\, -\frac{\partial}{\partial{z_i}}\ ,
\
\frac{\partial}{\partial{\xi_i}}\mapsto  z_i\
$
for  $i=1,\dots,N$,  defines an isomorphism
 of the rings $D_{X^*}$ and $D_{X}$ called {\it the Fourier transform.}
The Fourier transform defines a functor from the category $\MD_X$ to the category
$\MD_{X^*}$.

We consider a  central arrangement
 $\A $ of hyperplanes in $X$, a quiver $\VA$ of $\A$ and the
associated quiver $\D_X$-module $\EVA$.
We describe the Fourier transform of $E\VA$.  It
turns out that the Fourier transform of $E\VA$ is the quiver
$\D_{X^*}$-module associated with a suitable arrangment in $X^*$.

\medskip
{\it Example 1.10}.
In Example 1.1, the space $X^*$ has two strata:
 $X^*_\ee = \{\xi=0\}$ and $X^*_\a= \{\xi\neq  0\}$.
The Fourier transform of the
$\D_X$-module $\EVA$ from Example 1.2 is generated by vectors $1\ot v_\ee
$, where $v_\ee\in V_\ee$, and $d\xi\ot v_\a  $, where $v_\a\in V_\a$.
The relations are:
$\xi\,(1\ot v_\ee)=-d\xi\ot A_{\a,\ee}(v_\ee)$ for any $v_\ee\in
V_\ee$ and $\frac{d}{d\xi}\, (d\xi\ot v_\a)=1\ot A_{\ee,\a}(v_\a)$ for
any $v_\a\in V_\a$.

\bigskip
In the theory of perverse sheaves there is a notion of the
MacPherson extension \cite{BBD, GMc}.
In our situation it is a morphism from the category of perverse sheaves
on $X^0$ to the category of perverse sheaves on $X$. A perverse sheaf on
 $X^0$ is a local system on $X^0$. Thus the  MacPherson extension is an
extension of a local system on $X^0$ to a suitable perverse sheaf on $X$.

The MacPherson extension of perverse sheaves corresponds to a morphism
$j_{N,0,!*}: \MD^{hol}_{X^0}\to \MD^{hol}_{X}$  of holonomic
$\D$-modules. The morphism $j_{N,0,!*}$ can be defined in two steps.
First, one defines a morphism $s_0^{(D)}: j_{N,0,!}(M)\to j_{N,0,*}(M)$ for any
$\D_X$-module $M$ and then one defines $j_{N,0,!*}(M)$ as the image of
$j_{N,0,!}(M)$ in $j_{N,0,*}(M)$ under the morphism $s_0^{(D)}$. The
MacPherson extension of the perverse sheaf associated to the
$\D_{X^0}$-module $M$, is the perverse sheaf, associated to the
$\D_{X}$-module $j_{N,0,!*}(M)$.

In this paper we construct a morphism $s_0: \j_{N,0,!}\to \j_{N,0,*}$ of
functors of quiver categories. We
introduce the quiver MacPherson extension of a level zero quiver
$\VA_{\G^0}$ as
 the image of
$\j_{N,0,!}(\VA_{\G_0})$ under the map $s_0$. We prove that for a non-resonant
level zero quiver $\VA_{\G_0}$, the $\D_X$-module  morphisms
$Es_0: E\j_{N,0,!}(\VA_{\G_0})\to E\j_{N,0,*}(\VA_{\G_0})$ and
$s_0^{(D)}E^0: j_{N,0,!}(E^0\VA_{\G_0})\to j_{N,0,*}(E^0\VA_{\G_0})$ are
equal
under the identifications of $j_{0,N,!}(E^0 \VA_{\G^0})$ with
$ E(\j_{0,N,!}\VA_{\G^0})$, and $j_{0,N,*}(E^0 \VA_{\G^0})$ with
$ E(\j_{0,N,*}\VA_{\G^0})$, mentioned above.
 Thus we have  the isomorphism of $\D_X$-modules:
 $j_{N,0,!*}(E\VA_{\G_0})\approx E\j_{N,0,!*}(\VA_{\G_0})$, which relates
the geometric and quiver MacPherson extensions.

\medskip \noindent
{\it Example 1.11}. In  Examples 1.4 and 1.9, the
morphism $s_0$ of the level one quiver $\j_{1,0,!}\VA_{\G^0}$ to the
level one
quiver $\j_{1,0,*}\VA_{\G^0}$ is the collection of two linear maps,
$(s_0)_\ee: W'_\ee\to W_\ee$ and $(s_0)_\a: W'_\a\to W_\a$. They are
$(s_0)_\ee=Id_{V_\ee}$, and $(s_0)_\a=A_\ee^\ee$. The level one quiver
$\j_{1,0,!*}\VA_{\G^0}$ is  $\{W''_\ee, W''_\a,
C''_{\ee, \a}, C''_{\a,\ee} \}$,
where $W''_\ee=V_\ee$, $W''_\a={\rm Im}\, A^\ee_\ee\subset V_\a$,
$C''_{\a, \ee}=A_\ee^\ee$, and
$C''_{\ee,\a}$ is the inclusion of ${\rm Im}\, A^\ee_\ee$ to $V_\ee$.

\medskip
 As an application of quiver constructions, we calculate the
cohomology groups of $X$ (and $X^0$) with coefficients in a quiver perverse sheaf
if the arrangement  $\A$ is  central
and all linear maps of the quiver  are close to zero.
Namely, for a quiver $\VA = \{ V_\a, A_{\a,\b}\}\in \Qui_\G$,
we construct a finite-dimensional complex in terms of the spaces $V_\a$
and linear maps $A_{\a,\b}$ and show that this complex
calculates the cohomology groups of $X$ with  coefficients in
the quiver  perverse sheaf associated with $\VA$.

\medskip \noindent
{\it Example 1.12}. In Examples 1.1 and 1.2,
consider the perverse sheaf $\QVA$ on $X$, associated with the $\D_{X}$-module
$E\VA$. If the operators
$A_{\a,\ee}$ and $A_{\ee,\a}$ are close to zero,
 then the cohomology groups of $X$ with coefficients
 in the perverse sheaf $\QVA$  are calculated by the
complex
$$0\ \to\ V_\ee\ \stackrel{A_{\a,\ee}}{\to}\   V_\a\ \to\ 0\ ,$$
where $\deg V_\ee=-1$, and $\deg V_\a=0$.

\medskip

\medskip \noindent
{\it Example 1.13}.  In Example 1.3, consider the perverse sheaf
 ${\mathcal Q}^0\VA^0$ on $X^0$, associated with the $\D_{X^0}$-module
$E^0\VA_{\G^0}$. If the operator
$A_{\ee}^\ee$ is close to zero,  then the cohomology groups
of $X^0$ with coefficients  in ${\mathcal Q}^0\VA^0$,  are calculated by the
complex
$$0\ \to\ V_\ee\ \stackrel{A_{\ee}^\ee}{\to}\   V_\ee\ \to\ 0\ .$$

\medskip \noindent
{\it Example 1.14}. In Example 1.11, consider the perverse sheaf
 ${\mathcal Q}^0\VA^0$ on $X^0$. If the  operator
$A_{\ee}^\ee$ is close to zero,  then the cohomology groups
of $X$ with coefficients in  the MacPherson extension of
 ${\mathcal Q}^0\VA^0$ are calculated by the
complex
$$0\ \to\ V_\ee\ \stackrel{A_{\ee}^\ee}{\to}\  {\rm Im}\, A_{\ee}^\ee
 \ \to\ 0\ .$$

\bigskip
{\bf 3}. All results of the paper admit an equivariant extension to
the case when the arrangement $\A$ admits a finite symmetry group.

As an application of equivariant constructions, we consider in Section
\ref{section5.6} the discriminantal arrangements and
one-dimensional local systems appearing in hypergeometric solutions of
the KZ equations \cite{SV}.  We calculate the equivariant intersection
cohomology groups of the discriminantal arrangements and formulate the result
in terms of representation theory in Corollary \ref{cor5.9}. The statement of
Corollary \ref{cor5.9} is one of the main motivations of this paper.

\bigskip
 {\bf 4}. In some of the statements of the paper we impose two types
of non-resonance conditions on a quiver $\VA_\G=\{V_\a,
A_{\a,\b}\}$.
For every vertex $\a\in \G$, we construct a finite collection of
auxiliary linear operators of the form $\sum A_{\a,\b}A_{\b,\a}:
V_\a\to V_\a$ with suitable range of summation over $\b$.  In the
first type of non-resonance conditions, we require that the auxiliary
operators have no eigenvalues which differ by an integer (the weak
non-resonance condition).  In the second type of non-resonance
conditions, we require that the auxiliary operators have no
eigenvalues which differ by an integer and zero is the only admissible
integer eigenvalue (the non-resonance condition).  The importance of
non-resonance conditions is illustrated in the Example 1.13 below.

\medskip \noindent {\it Example 1.13}. Under conditions of Example
1.3, let $\VA_{\G^0}=\{ V_\ee, A_\ee^\ee\}$ and ${\mathcal
  W}_{\G^0}=\{W_\ee, B_\ee^\ee\}$ be two level zero quivers and
$E^0\VA_{\G^0}$ and $E^0{\mathcal W}_{\G^0}$ the associated
$D_{X^0}$-modules.  Assume that $G\subset\CC$ is a subset, such that
$x-y\notin\ZZ\setminus\{0\}$ for all $x,y \in G$.  Assume that
eigenvalues of $A_\ee^\ee$ and $B_\ee^\ee$ belong to $G$. Under these
assumptions we will construct an isomorphism of vector spaces
$\Hom_{D_{X^0}}(E^0\VA_{\G^0},E^0{\mathcal W}_{\G^0})$ and
$\Hom_{Qui_{\G^0}}(\VA_{\G^0},{\mathcal W}_{\G^0})$.  Indeed, let
$\Phi\, :\, E^0\VA_{\G^0}\to E^0{\mathcal W}_{\G^0}$ be a morphism of
$D_{X^0}$-modules. The morphism $\Phi$ is uniquely determined by its
restriction $\varphi=\Phi|_{{\overline V}_\ee}$ to the space
${\overline V}_\ee\subset E^0\VA_{\G^0}$. The map $\varphi$ is a
regular rational function on $X^0$ with values in $\Hom \,({\overline
  V}_\ee, {\overline W}_\ee)$.  Write its Laurent series,
$$
\varphi=(1\ot \varphi_{k}){z^{k}}+(1\ot \varphi_{k+1}){z^{k+1}}+
\cdots\, ,
$$ where $\varphi_l:V_\ee\to W_\ee$ are linear maps and $\varphi_k\neq
0$. The condition $[\Phi,\frac d {d z}]=0$ implies
$\varphi_k A_\ee^\a= (B^\a_\ee+k)\varphi_k$. Since all eigenvalues of
$A_\ee^\a$ and $B_\ee^\a$ are in $G$, we conclude that $k=0$. The
leading term $\varphi_0$ is a morphism of quivers, $\VA_{\G^0}\to
{\mathcal W}_{\G^0}$.  Let $E^0(\varphi_0)\,:\,
E^0\VA_{\G^0}\to E^0{\mathcal W}_{\G^0}$ be the $D_{X^0}$-module morphism
associated with $\varphi_0$. The restriction of $E^0(\varphi_0)$ to
$\overline{V}_\ee$ is $1\ot \varphi_0$. The leading term
$\widetilde{\varphi}_k$ of the $D_{X^0}$-module morphism
$\widetilde{\Phi}=\Phi-E^0(\varphi_0)$ has $k>0$. Hence $\widetilde{\Phi}=0$
and $\Phi=E^0(\varphi_0)$.

\bigskip
{\bf 5}. In the theory of arrangements, two complexes play an
important role. The first is the Aomoto complex \cite{A} on the space
of the Orlik-Solomon algebra and the second is the flag complex on the
dual space to the Orlik-Solomon algebra \cite{SV}.  The Shapovalov
form of \cite{SV} defines a morphism of the flag complex to the Aomoto
complex. The image of the Shapovalov form is called the complex of
flag forms.  The cohomology groups of the flag complex give the
cohomology groups of $X^0$ with compact support.  The cohomology
groups of the Aomoto complex give the cohomology groups of $X^0$ with
coefficients in a suitable local system on $X^0$. An important problem
is to give a topological meaning to the complex of flag forms.

 The main quiver constructions of this paper are the constructions of
the two quiver direct images $\j_{N,0,*}$ and $\j_{N,0,!}.$ The
constructions of $\j_{N,0,*}$ and $\j_{N,0,!}$ are respective
generalizations of the constructions of the Aomoto complex and the
flag complex.  The morphism of functors $s_0: \j_{N,0,!}\to
\j_{N,0,*}$ is a generalization of the Shapovalov morphism
of the flag complex to the Aomoto complex, see Section \ref{secdir}.
The quiver MacPherson extension $\j_{N,0,!*}$ is a generalization of
the complex of flag forms in \cite{SV}.

From this point of view, our calculation of the cohomology groups of $X$ with
coefficients in the MacPherson extension of a quiver perverse sheaf on
$X^0$ gives a topological meaning to the complex of  flag forms.

To find a topological meaning of the complex of flag forms was one of the
main goals of this paper.

\bigskip
{\bf 6}. The plan of the paper is as follows.  In Section 2 we
introduce graphs and stratified spaces, associated with an arrangement
of hyperplanes.
In Section 3 we study quivers, assigned to the
arrangement, and operations on quivers.  Section 4 is devoted to
quiver $\D$-modules. In Section 5 we describe perverse sheaves related
to quiver $\D$-modules and calculate the cohomology groups of the
affine space with coefficients in these sheaves. In Section 6 we study
the equivariant version of quiver constructions.
 All proofs are collected in Section 7.

\bigskip
 {\bf 7}. We thank V. Schechtman, who inspired and initiated the study of $\D$-modules
associated with an arrangement of hyperplanes \cite{Kh1, KhS1} and who
suggested an idea of the construction of a free $\D$-module resolution of
a quiver $\D$-module, see an announcement in \cite{KhS}.

We thank the referee for many helpful suggestions.

The work of S.Kh. was supported  by CRDF grant RM1-2334, the grants
INTAS-OPEN-03-51-3350,
 RFBR grant 04-01-00642 and RFBR grant
to support scientific schools  NSh-1999.2003.2.
The work of A.V. was supported by NSF grant
DMS-0555327.

\section{Arrangements}

\subsection{Graph of an arrangement}\label{Sec.1.1}
Let $\A =\{H_j\}$, $j\in J(\A)$,  be an arrangement of finitely many
hyperplanes in the complex
 affine space $X=\C^N$.  Hyperplanes $H_j$ define in $X$ the structure
 of a stratified space which we denote by the same letter $X$ or by $X^\A$.
 The closure of a stratum $X_\a\subset \C^N$ is the intersection of some
 hyperplanes  $H_j$, $j\in J_\a\subset J(\A)$,
 and the stratum $X_\a$ itself is the complement in $\Xa$ to the union
of the closures of all smaller strata contained in $\Xa$:
$$
\cXa= \bigcap_{j\in J_\a\subset J(\A)}H_j,\qquad
X_\a=\cXa\ \setminus\ (\bigcup_{\cXb\subset \cXa,
\cXb\not=\cXa}\cXb)\ .
$$
We denote by $X_\ee$  the complement to the union of all
hyperplanes $H_j$ which is the unique stratum of $X$ open in $X$.

The combinatorial structure of the adjacency of
strata can be described by a non-oriented graph $\Ga=\Ga(\A)$.
The vertices $\a \in \v(\Ga)$ of the graph are in one-to-one
correspondence with the strata $X_\alpha$ of $X$.
The edges $(\a,\b)\in \e(\Ga)$
 describe adjacency in codimension one:
 $$
(\a,\b)\in \e(\Ga) \Leftrightarrow
\left\{
\begin{array}{ll} \Xa\supset\Xb,\,& \codim_\Xa\Xb=1,
\quad {\rm or}\\
\Xb\supset\Xa,\,& \codim_\Xb\Xa=1.
\end{array}\right.
$$
Two vertices $\a$ and $\b$  will be called {\it adjacent},
if $(\a,\b)\in\e(\Ga)$.

To any vertex $\a$  we attach a
nonnegative integer $l(\a)$,
$$
l(\a)=\codim_{\C^N}X_\a.
$$
We write $\a\succ\b$
for vertices $\a,\b$ if $(\a,\b)$ is an edge and $l(\a)<
l(\b)$. In other words,
$$
\a\succ\b\quad \Leftrightarrow\quad \Xa\supset\Xb\ , \quad
\codim_\Xa\Xb=1\ . 
$$
We set
$$\varepsilon(\a,\b)=\left\{\begin{array}{cc}1\,&
{\rm if}\quad \a\succ\b,\\
-1\, &{\rm if}\quad \b\succ\a,\\
0\, & {\rm otherwise}.
\end{array}\right. $$
{} For an integer $n\geq 0$ denote by $\v_n(\Ga)\subset \v(\Ga)$ the subset of all vertices $\a$
such that $l(\a)=n$. We have $\v_0(\G)=\ee$, $\v_1(\G)=J(\A)$.

We write $\a>\b$ for $\a,\b\in I(\Ga)$ if there exists a chain
$\a_1\succ\a_2\succ...\succ\a_n$ with $\a_1=\a$ and $\a_n=\b$.
We write $\a\geq b$ if $\a>\b$ or $\a=\b$.
Clearly,
$$\a\geq\b\quad\Leftrightarrow\quad \X_\a\supset\X_\b.$$
The labeled graph $\Ga=\{ \v(\Ga), \e(\Ga),\ {} l:\, \v(\Ga)\to \ZZ_{\geq 0} \}$
has the following properties.
\begin{enumerate}

\item[(i)]
 There exists the unique vertex $\emptyset\in\v(\Ga)$ such that
$l(\emptyset)=0$ and $\emptyset\geq\b$ for any $\b\in\v(\Ga)$;
\item[(ii)]
 If $(\a,\b)\in \e(\Ga)$, then $ |l(\a)-l(\b)|=1$;
\item[(iii)]
 {}For any $\a,\b\in\v(\Ga)$ the set $\L(\a,\b)=\{\g\in\v(\Ga):\a\geq\g,\
\b\geq\g\}$ is either empty or has the unique element $\gamma$ maximal
with respect to the partial order $\geq$. Moreover, we have
$l(\g)\leq l(\a)+l(\b)$. This element is called {\it the intersection of
$\alpha$ and $\beta$} and denoted $\a\wedge\b$.

\item[(iv)]
 {}For any $\a,\b\in\v(\Ga)$ with non-empty
$\L(\a,\b)$, denote $U(\a,\b)=\{\d\in\v(\Ga):\delta\geq\a,
\
\d\geq\b\}$. Then for any $\d\in U(\a,\b)$ we have $l(\d)\leq l(\a)+l(\b)-l(\a\wedge\b)$.
Moreover, if there is $\delta \in U(\a,\b)$ such that
$l(\d)=l(\a)+l(\b)-l(\a\wedge\b)$, then $\d$ is the unique
minimal element of $U(\a,\b)$ with respect to the partial order $\geq$.
\end{enumerate}

One can see that if $\beta\wedge\gamma$ and
$\a\wedge(\beta\wedge\gamma)$ are defined, then
$\a\wedge\b$ and $(\a\wedge\b)\wedge\gamma$ are defined,
and
$\a\wedge(\beta\wedge\gamma)=(\a\wedge\b)\wedge\gamma$.

Any labeled graph with properties (i-iv) will be called {\it admissible}.

\bigskip

\subsection{Framings of an arrangement}\label{Framing}
Let $\A \subset \C^N$ be an arrangement.  To each stratum $X_\a\subset X^\A$
we attach certain spaces of affine functions, constant vector
fields, and the top exterior forms on $\Xa$ with constant coefficients.

Denote by $F$ the space of affine functions on $\CC^N$,
$$
  f\in\FN \quad \Leftrightarrow \quad f=a_0+\sum_{i=1}^{N}
a_iz_i, \qquad a_i\in\CC,
  $$
where $z_1, \dots , z_N$ are affine coordinates on $\CC^N$.

Denote
by \TN\ the space of vector fields $\xi$ on $\C^N$ with constant
coefficients,
 $$\xi\in\TN \quad\Leftrightarrow \quad \xi=\sum b_i
  \frac{\partial}{\partial{z_i}},\qquad
 b_i\in\CC.$$
{} For a stratum $X_\a$ denote by  $\Fa\,\subset F$ the subspace
of all functions   vanishing on $X_\a$
  and by $\Ta\subset T$ the subspace of all vector fields
 parallel to $X_\a$, that is, $\xi \in \Ta$ if
$\mylangle \xi,f\myrangle=0$ for any $f\in\Fa$.
Here $\mylangle \xi,f\myrangle=i_\xi(df)$ denotes
the substitution of the vector field $\xi$ to the differential form $df$.

A function $f\in \Fa$ is called {\it generic},
 if for any stratum $X_\b$  the intersection
$$
\{ f=0\}\cap X_\b
$$
has the minimal possible dimension among all $f\in\Fa$.

 Denote  by $\TTT{\a}$ the $N$-dimensional complex vector space,
$$
\TTT{\a}=
 \Ta \oplus \Fa\ .
$$

\bigskip

Let  $X_\a$ and $X_\b$ be  strata, such that
$\X_\a\supset X_\b$. Then there are tautological
  inclusions
\begin{equation}\label{inclusions}
\mu_{\b,\a}:\ \Fa\hookrightarrow\Fb,\qquad \overline{\mu}_{\a,\b}:\
 \Tb\hookrightarrow\Ta\ ,
\end{equation}
We  denote
$$
\mu_{\b,\a}:\ \bigwedge^p\Fa\hookrightarrow\bigwedge^p\Fb ,
\qquad
\overline{\mu}_{\a,\b}:\ \bigwedge^p\Tb\hookrightarrow\bigwedge^p\Ta ,
\qquad p > 0,
$$
the induced inclusions.

For $\a\succ\b$  and $f\in \Fb$, introduce the linear maps
$$
i_f^{\b,\a}\ :\ \bigwedge^p\Ta\ \to\
\bigwedge^{p-1}\Tb ,
\qquad
\vec{\xi} \ \mapsto \ i_{df}(\vec{\xi})\  .
$$
Here $i_{df}(\vec{\xi})$ is the substitution of  $df$
to the $p$-vector $\vec{\xi}\in\bigwedge^p\Ta$ with the result
considered as a
$(p-1)$-vector in $\bigwedge^{p-1}\Tb$.
More precisely, if $\xxi{}=\xi_1\wedge...\wedge\xi_p$
and  $\xi_2, \dots , \xi_p \in \Tb$, then
\begin{equation}\label{if}
i_f^{\b,\a}(\xxi{})\ =\ \mylangle \xi_1,f\myrangle\,\xi_2\wedge...\wedge\xi_p \ .
\end{equation}

For  $\a\succ\b$ and $\xi\in\Ta$ introduce the linear map
\begin{equation}
\ii_\xi^{\a,\b}\ :\ \bigwedge^p \Fb\ \to\ \bigwedge^{p-1}\Fa\ ,
\qquad
\f{} \ \mapsto i_\xi (\f{}) \ .
\notag
\end{equation}
More precisely, if $\f{} = f_1\wedge...\wedge f_p$  with
 $f_2, \dots , f_p \in \Fa$, then
\begin{equation}
\ii_\xi^{\a,\b}(\f{})\ =\ \mylangle \xi,f_1\myrangle\, f_2\wedge...\wedge f_p \ .
\label{ixi}
\end{equation}

\medskip

Denote by $\overline{\Omega}^p$ the complex vector space   of
holomorphic  exterior $p$-forms on $\C^N$,
 invariant with respect to affine translations. Each  form in
$\overline{\Omega}^p$ is a linear combination  of
the forms $df_1\wedge...\wedge df_p$ with $f_1, \dots , f_p \in \FN$.

For a stratum $X_\a$, denote by \Oa  the one-dimensional vector space of the top
exterior holomorphic forms on $X_\a$ invariant with respect to translations
along $X_\a$.
The space \Oa is isomorphic to the quotient of the space
$\overline{\Omega}^{N-l(\a)}$ over the subspace generated by the forms
$df_1\wedge...\wedge df_{N-l(\a)}$, where $f_1, \dots , f_{N-l(\a)}
\in \FN$ and at least one of them belongs to  \Fa.

  It is technically convenient to
  fix a vector space structure on $\C^N$
  and fix  a non-degenerate bilinear form $\lform \,, \rform$
in this vector space  $\C^N$.
 In this setting to each stratum $X_\a$ we  also attach
the following spaces $\TTa$ and $\FFa$.

The space \TTa consists of constant vector fields orthogonal to
\Ta.
Let $z_1,...,z_N$ be linear coordinates in $\C^N$ with respect to an
orthonormal basis.
Let $\xi'=\sum b_i \frac{\partial}{\partial{z_i}}$ be a constant vector
field. Then $\xi'$ belongs to \TTa,
if $b_1a_1+...+b_Na_N=0$ for any
 $\xi=\sum a_i \frac{\partial}{\partial{z_i}} \in \Ta$.
  The space
 \FFa consists of all linear functions $f'=\sum c_iz_i$, which are
 annihilated by vector fields in \TTa, i.e. $\mylangle \xi',f'\myrangle=0$ for any $\xi' \in
 \TTa$.

\bigskip


We use two types of framings of the arrangement $\A$: edge framings and
vertex framings.

{\it An edge framing} is  a choice for any edge
 $(\a,\b)\in E(\G)$, $l(\a)+1=l(\b)$,
 of a function $f_{\b,\a}\in \Fb\setminus F_a$
 and a vector field $\xi_{\a,\b} \in
\Ta\setminus T_\beta$. For any edge  $(\a,\b)\in E(\G)$, the edge framing
determines the linear map $\pi_{\b,\a}:\Oa\to\Ob$, such that for
$\overline{\o}_\a\in\Oa$ and $\overline{\o}_\b\in\Ob$ we have
\begin{equation}
\label{pi1}
\begin{array}{ccccc}
\pi_{\b,\a}(\overline{\o}_\a)=\overline{\o}_\b&\Leftrightarrow&
df_{\b,\a}\wedge\overline{\o}_\b=\om_\a,&{\rm if}& \a\succ\b,\\
\pi_{\b,\a}(\overline{\o}_\a)=\overline{\o}_\b&\Leftrightarrow&
i_{\xi_{\b,\a}}(\overline{\o}_\b)=\om_\a,&{\rm if}& \b\succ\a.
\end{array}
\end{equation}

{\it A vertex framing}
of  $\A$ is a choice
 of a  function $\ff{\a}\in F_\a$ for any stratum $X_\a$.
A vertex framing  ${\mathcal F}$ is called {\it generic} if for any
$\a\in \v(\Ga)$ the function $\ff{\a}$ is generic.

For a nonnegative integer $n$,
{\it a level $n$ vertex framing}
is a choice  of a  function $\ff{\a}\in F_\a$
for any $\a\in \v_{n+1}(\Ga)$. It is generic if every chosen $\ff{\a}$ is generic.

\bigskip

\subsection{Principal open subsets $X^n$}\label{Sec.Xn}
{} For a nonnegative integer $n$ and an arrangement
$\A$, define the open
subset $X^n\subset X$ as the complement to the closure of the
union  of all codimension $n+1$ strata:
$$
X^n\ =\ \CC^N\setminus \bigcup\limits_{\a,\ l(\a)=n+1}\X_\a\ .
$$
We call $X^n$ {\it the principal level $n$ open subset of $X$}.
In particular, $X^0$ coincides with $X_\ee$, the open stratum in $X$,
and $X^N=X$.
The sets $X^n$ determine a filtration of $X$:
$$
X^0\ \subset \ X^1\ \subset\ ...\ \subset\ X^N=X\ ,
$$

For $k\leq l$ denote by $j_{l,k}$ the inclusion
$$j_{l,k}: X^{k}\to X^l.$$
We have $j_{k_3,k_1}=j_{k_3,k_2}j_{k_2,k_1}$ for $k_1\leq k_2\leq k_3$
 and we use the notation $j_0$ for the inclusion $j_{N,0}:X^0 \to X=X^N$.

 Each level $n$ vertex framing ${\mathcal F}_{n}=\{\ff{\a}\,|\ \a\in\v_{n+1}(\Ga)\}$
determines an open subset
 $U_{\FF_n}\subset X^n$,
\begin{equation}
\label{UFn}
U_{\FF_n}=X\setminus \bigcup\limits_{\a \in\v_{n+1}(\Ga)}\{\ff{\a}=0\} .
\end{equation}
The open subsets  $U_{\FF_n}$ for different choices of framings
$\FF_n$ form a cover of the space $X^n$ by open affine subsets.

\subsection{Truncated graphs}
The sets $X^n$ are stratified spaces; their strata are the strata $X_\a$
of $X$ with $l(\a)\leq n$.
We define the graph (with loops)  of the stratified space
$\U^n$ as follows.   The  vertices of this graph
$\Ga^n$ are the vertices $\a\in\v(\Ga)$ such that $l(\a)\leq
 n$.  The edges $(\a,\b)\in\e(\Ga^n)$, which are not loops, coincide with
the edges $(\a,\b)$ of  $\Ga$, connecting the vertices of
 $\Ga^n$. For any pair $\a\succ\b$ with $l(\a) =n$, the graph
$\Ga^n$ has a loop  at the vertex $\a$, which we denote  by $(\a,\a)^\b$.

 The graph $\G^n$ will be called
{\it the truncated graph at level $n$} associated to the admissible graph
$\G$.

{\em Examples.} The graph $\G^0$ has one vertex and $|J(\A)|$ loops.
 Every loop corresponds to a hyperplane of the initial
arrangement $\A$. The graph $\G^N$ coincides with $\G$.

\bigskip

\subsection{Orlik-Solomon algebra and flag complex}\label{Orlik}
Here we recall the construction of the Orlik-Solomon algebra and
the flag complex
 of an arrangement, see
 \cite{SV}.
 For an arrangement $\A$ in $ \C^N$ define the complex
vector spaces $\OS^p(\A)$, $p = 0,  \dots, N$.
 For $p=0$ set $\OS^p(\A)=\C$. For  $p \geq 1$ the vector space
 $\OS^p(\A)$   is generated by symbols
$(H_{j_1},...,H_{j_p})$, where $H_{j_k}\in\A$, such that
\begin{enumerate}
\item[(i)]
$$
(H_{j_1},...,H_{j_p})=0
$$
if $H_{j_1}$,...,$H_{j_p}$ are not in general position, that is if their
intersection $H_{j_1}\cap ... \cap H_{j_p}$ is empty or
 its codimension is
 less than $p$;
\item[(ii)]
$$
(H_{j_{\sigma(1)}},...,H_{j_{\sigma(p)}})={\rm sgn}(\sigma)\cdot
(H_{j_1},...,H_{j_p})
$$
for any permutation $\sigma\in S_p$;
\item[(iii)]
\begin{equation}
\sum\limits_{k=1}^{p+1}(-1)^k(H_{j_1},...,\widehat{H_{j_k}},...,H_{j_{p+1}})=0
\label{OSrelation}
\end{equation}
for any $(p+1)$-tuple $H_{j_1},...,H_{j_{p+1}}$ of hyperplanes
in $\A$ which are
not in general position and such that $H_{j_1}\cap...\cap H_{j_{p+1}}\not=\ee$.
\end{enumerate}

The direct sum $\OS(\A)=\ooplus\limits_{p=1}^{N}\OS^p(\A)$ is a
differential graded skew commutative algebra with respect to
 multiplication
 $$(H_{j_1},...,H_{j_p})\cdot(H_{j_{p+1}},...,H_{j_{p+q}}) =
 (H_{j_1},...,H_{j_p},H_{j_{p+1}},...,H_{j_{p+q}})$$
 and the differential $\delta_{\OS}$:
 $\OS^p(\A)\to \OS^{p-1}(\A)$,
 \begin{equation}
 \delta_{\OS}:\ (H_{j_1},...,H_{j_p}) \ \mapsto \
 \sum\limits_{k=0}^p(-1)^{k-1}(H_{j_1},...,\widehat{H_{j_k}},...,H_{j_p}).
 \label{dOS}
 \end{equation}
 The algebra is  called {\it the Orlik-Solomon algebra} of the arrangement $\A$.

The Orlik-Solomon algebra is determined by the admissible graph $\Ga$
 only. For an admissible graph $\G$, its Orlik-Solomon
algebra $\OS(\Ga)$  can be defined as the algebra
 generated by the skew-symmetric symbols
$(H_{j_1},...,H_{j_k})$, where $j_k\in\v_1(\Ga)$. We set
$(H_{j_1},...,H_{j_p})=0$ if $H_{j_1},...,H_{j_k}$ are not in general position,
and
 $H_{j_1},...,H_{j_k}$ are not in general position if the vertex
$j_1\wedge...\wedge j_k$ does not exist or its codimension
is less than $k$.
Relation \rf{OSrelation} has to take place for any $(p+1)$-tuple
of vertices  $j_i \in \v_1(\Ga)$ for which the intersection
$j_1\wedge...\wedge j_{p+1}$ exists and its codimension is less than
$p+1$.
\medskip

{} For $\a\in I(\Ga)$, {\it a flag starting at $X_\a$} is a sequence
$$
\X_{\a_0} \supset ... \supset \X_{\a_p}=\X_\a
$$
of closed strata such that
$ \codim_{\C^N}X_{\a_k} = j$ for $j = 0, \dots , p$.

For $\a\in I(\Ga)$, we define $\widetilde{\FF}_\a$  to be  the
complex vector space  with basis vectors
$\widetilde{F}_{\a_0,...,\a_p=\a}$ la\-bel\-ed by the elements of
the set of all flags  starting at $X_\a$.

 Define  $\FF_\a$ as the quotient of
$\widetilde{\FF}_\a$ over the subspace generated by the elements
\begin{equation}
\label{flagrelation} \sum\limits_{\b,\ {}\
\a_{k-1}\succ\b\succ\a_{k+1}}
\widetilde{F}_{\a_0,...,\a_{k-1},\b,\a_{k+1},...,\a_p=\a}\ ,
\end{equation}
where $k = 1, \dots , p-1$ and
$\X_{\a_0}\supset...\supset \X_{\a_{k-1}} \supset \X_{\a_{k+1}}\supset...\supset
\X_{\a_p}=\X_\a$ is any incomplete flag with  $l(\a_j)=j$.

Denote  by ${F}_{\a_0,...,\a_p}$ the image in $\FF_\a$ of the basis vector
$\widetilde{F}_{\a_0,...,\a_p}$. Set
$$
{\FF}^p(\A)\ =\ \oplus_{\a\in\Ga,\, l(\a)=p}\ {\FF}_\a\ .
$$
The direct sum
 $$
{\FF}(\A)\ =\ \ooplus\limits_{p=0}^N\,{\FF}^p(\A)
$$
is a complex
with respect to the differential
 \begin{equation}
d_{\FF} :  \FF^p(\A)\to\FF^{p+1}(\A),
\qquad
{F}_{\a_0,...,\a_p} \mapsto (-1)^p
\sum\limits_{\a_{p+1},\ \a_p\succ
\a_{p+1}}{F}_{\a_0,...,\a_p,\a_{p+1}}\ .
\label{df}
\end{equation}
The  complex $(\FF^\mycircle(\A), d_{\FF})$   is called
{\it the flag complex of the arrangement $\A$}.

As well as the Orlik-Solomon
algebra, it is determined by the graph $\Ga=\Ga(\A)$ only and will be
denoted in this context by $(\FF^\mycircle(\Ga), d_{\FF})$.

The vector spaces $\OS^p(\A)$ and $\FF^p(\A)$ are dual
\cite{SV}. The pairing $ \OS^p(\A)\otimes\FF^p(\A) \to \C$ is defined as follows.
{}For $H_{j_1},...,H_{j_p}$ in
general position, set
$$
F(H_{j_1},...,H_{j_p})=F_{\a_0,\a_1,...,\a_p}
$$
where
$$
\X_{\a_0}=\C^N,\quad \X_{\a_1}=H_{j_1},\quad
\X_{\a_2}=H_{j_1}\cap H_{j_2},\quad...\quad,\ {} \ {}
\X_{\a_p}=
H_{j_1}\cap H_{j_2} \cap...\cap H_{j_p}.
$$
Then set
$$
\langle (H_{j_1},...,H_{j_p}), F \rangle = \text{sign}\,(\sigma)\ ,
$$
if $F = F(H_{j_{\sigma(1)}},...,H_{j_{\sigma(p)}})$ for some $\sigma \in S_p$,
and
$$
\langle (H_{j_1},...,H_{j_p}), F \rangle = 0
$$
otherwise.

The differentials $\d_{\OS}$ and $d_{\FF}$ are adjoint with
respect to this duality and the map dual to the multiplication map in
$\OS(\A)$ turns  $\FF(\A)$ into a skew cocommutative
 differential graded coalgebra.

Let $a: \A\to \C$ be a map which assigns to each hyperplane $H$
a complex number $a(H)$ called {\it the exponent} of the hyperplane.
 Set
 $$
\o(a)\ =\ \sum_{H\in\A}\,a(H)\,H\ {}\ \in \ {}\ \OS^1(\A)\ .
$$
 The multiplication by $\o(a)$ defines a differential
 $$
d_a\ :\   \OS^p(\A)\ \to\ \OS^{p+1}(\A) ,
\qquad
 x \ \mapsto\ \o(a)\cdot x \ ,
$$
 in the  vector space of the
Orlik-Solomon algebra $\OS(\A)$. The complex
$(\OS^\mycircle(\A), d_a)$ is called
{\it the Aomoto complex}.

 The collection of exponents determines {\em the Shapovalov map}
 $$
S_a\ :\ \FF(\A)\ \to\ \OS(\A)\ ,
$$
 $$
  {F}_{\a_0,\a_1,...,\a_p}\ \mapsto \
\sum \ a(H_{j_1})\, a(H_{j_2}) \,\dots \,a(H_{j_p})\ (H_{j_1}, \dots , H_{j_p})\ ,
$$
 where the sum is taken over all $p$-tuples $(H_{j_1},...,H_{j_p})$ such that
$$
H_{j_1}\supset \X_{a_1},\ {}\ .\ .\ .\ {} , \ {} H_{j_p}\supset \X_{\a_p}\ .
$$
 The image  $S_a({F}_{\a_0,\a_1,...,\a_p}) \in \OS^p$ of a flag
 ${F}_{\a_0,\a_1,...,\a_p}$ is called {\it a flag form}.

According to \cite{SV}, the map
$S_a$ is a morphism of the complex $(\FF^\mycircle(\A), d_{\FF})$  to
the complex $(\OS^{\mycircle}(\A), d_a)$.
The image of
 the flag complex in the Aomoto complex is called {\it
the complex of  flag forms} of $\A$. We denote it by
$(\overline{\FF}^\mycircle_a(\A), d_a)$.

Identifying $\OS(\A)$ with $\FF^{*}(\A)$, we can consider
 the map $S_a$ as a bilinear form on the vector space $\FF(\A)$.
 The bilinear form is symmetric and  called {\it the Shapovalov form}.

\setcounter{equation}{0}

\bigskip

\section{Quivers}
\label{SEC quivers}

\subsection{Quivers of an admissible  graph}\label{Sec.2.1}
\quad Let  $\G=\{ \v(\G),\, \e(\G),$ $
l: $ $\v(\G) $ $\to $ $\ZZ_{\geq 0}\}$ be an admissible graph.

{\it A quiver} $\VG=\{V_\a, A_{\a,\b}\}$ of
$\G$ is a collection of finite dimensional complex
 vector spaces $V_\a$, $\a\in\v(\G)$, and a collection of
 linear maps $A_{\a,\b}: V_\b\to V_\a,$ $\a,\b\in\v(\G)$,
such that
\begin{enumerate}
\item[(a)]   $A_{\a,\b}=0$ if $\a$ and $\b$ are  not adjacent;
\item[(b)]
$ \sum_\b A_{\a,\b}A_{\b,\g} = 0 $
 if $|l(\a)-l(\g)|=2$ ;
\item[(c)]
$\sum_\b A_{\a,\b}A_{\b,\g} = 0$ if
$ l(\a)=l(\g),\;\a\not=\g $, and there exists $ \d\in\v(\G)$,
such that $ \a\succ\d, \g\succ \d$ .

\end{enumerate}

  The possibly nonzero
operators $A_{\a,\b}$ of a quiver $\VG$ correspond to oriented edges
of the graph $\G$.
 There may be only one
 or two summands in the left hand side of any relation of type (c)
due to properties (iii) and (iv) of  $\G$ in Section \ref{Sec.1.1}

Let $\{V_\a, A_{\a,\b}\},\,\{V'_\a, A'_{\a,\b}\}$ be two quivers of $\Ga$.
{\it A morphism of quivers} is a collection of maps $\{f_\a : V_\a \to V'_\a\}$
such that $ f_\a A_{\a,\b} = A'_{\a,\b} f_\b$ for all $\a, \b$.

Denote by $\Qui_\G $ the category of all quivers  of $\Ga$.

A quiver  can be considered as a
module over the path algebra $\B_\G$.

{\it The path algebra} of  $\G$ is the associative
algebra  $\B_\G$  with  unit, generated by two sets of generators
$i_\a$, $\a\in\v(\G)$,  and $a_{\a,\b}$, where $(\a,\b)$ runs through the set of
pairs of adjacent vertices of $\G$.
The generators are subject to the relations:
\begin{enumerate}
\item[(i)]  $i_\a^2=i_\a$ for all $\a$;
\item[(ii)]  $i_\a i_\b = 0$ for all $\a \neq \b$;
\item[(iii)] $\sum_\a i_\a=1$;
\item[(iv)] Let  $\delta_{\b,\g}$ be the Kronecker symbol, then
$a_{\a,\b}\,i_\g \,=\, \d_{\b,\g}\,a_{\a,\b},$\ and \
$i_\g \,a_{\a,\b}\,=\,\d_{\g,\a}\,a_{\a,\b}$  for all $\a, \b, \g$.

\item[(v)]
$ \sum_\b a_{\a,\b}a_{\b,\g} = 0 $ for all $\a,  \g$ with
 $ |l(\a)-l(\g)|=2$;
\item[(vi)]
$\sum_\b a_{\a,\b}a_{\b,\g} = 0 $ \ for
 all $\a,  \g$, such that $\a \neq \g$, $l(\a)=l(\g)$,
and there exists $\d$ such that   $\a\succ\d, \g\succ \d$.
\end{enumerate}

A quiver in $ \Qui_\G$ has a natural structure of a $\B_\G$-module.
The category $\Qui_\G$ is equivalent to the category of
finite-dimensional $\B_\G$-modules.

\medskip

For $\a,\b \in \v(\Ga)$, let $\varepsilon(\b,\a)$ be defined as in Section \ref{Sec.1.1}.

Let $\VG=\{V_\a, A_{\a,\b}\}$ be a quiver of $\Ga$. Then the collection
${\mathcal W}_\G  = \{W_\a, B_{\a,\b}\}$, where $W_\a = V_\a^*$ and
$B_{\a,\b} = \varepsilon(\b,\a)\, A_{\b,\a}^*\ :\ W_\b\,\to\, W_ \a$ determines a
quiver $\tau(\VG)$,
which we call {\it the dual} quiver.
The assignment $\VG \mapsto \tau(\VG)$
determines an anti-automorphism of the category $\Qui_\Ga$ of the fourth
order. Moreover, for any quiver $\VG$ the quiver $\tau^2(\VG)$ is
isomorphic to $\VG$. The isomorphism $\phi:\VG\to \tau^2(\VG)$ is given by
the relation $\phi(v_\a)=(-1)^{l(\a)}v_\a$.

The linear maps $A_{\a,\b}$ of a quiver $\VG$ can be collected into two
differentials, acting in the total space of the quiver,
 $C(\VG)=\oplus_{\a\in\v(\G)}V_\a$.
 Namely, let
\begin{equation}
\label{complex} C_k(\VG)\ =\ \ooplus\limits_{\a,\,{}\
l(\a)=k}V_\a\ ,\qquad C(\VG)\ =\ \ooplus_kC_k(\VG)\ .
\end{equation}
For $ v\in V_\a$, set
\begin{equation}
\label{differentials}
d(v)\ =\ \sum_{\b,\,{} \ \a\succ\b}A_{\b,\a}(v)\ ,
\qquad{\rm and}
\qquad
\partial(v)\ =\ \sum_{\b,\,{}\ \b\succ\a}A_{\b,\a}(v) \ .
\end{equation}
We have
$$
d (C_k(\VG)) \subset C_{k+1}(\VG) ,
\qquad
d^2 = 0 ,
$$
$$
\partial( C_k(\VG))\subset C_{k-1}(\VG)) ,
\qquad
\partial^2 = 0 .
$$
Denote the complex $(C(\VG),d)$ by $C_+(\VG)$ and the complex
$(C(\VG),\partial)$ by $C_-(\VG)$.

{}For a vertex $\b$ and an edge $(\a,\b)$,
denote  by $A_\b^\a$ the operator
\begin{equation}\label{Aab}
A_\b^\a\ =\ A_{\b,\a}A_{\a,\b}\ :\ {} \ V_\b\ \to\  V_\b\ .
\end{equation}
Let ${\cal O}_\b(\VG)$ be the $\C$-algebra of endomorphisms of $V_\b$,
generated by operators $A_\b^\a$, where $\a$ runs through the set of all
vertices adjacent to $\b$.
  We call ${\cal O}_\b(\VG)$ {\it the local algebra} of the quiver
at the vertex $\b$.

 To a vertex $\b$, $\beta \neq \emptyset$, we attach three
operators $S_\b$, $T_\b$ and $\Tprime_\b$,
 \begin{equation}
 \label{Tb0}
 S_\b  = \sum_{\a,\, \a\succ\b} A_\b^\a\ : \ {}\ V_\b\ \to\  V_\b\ ,
\end{equation}
\begin{equation}      
\label{2.22}
T_\b = \sum_{\a, \a'} A_{\a,\b}A_{\b,\a'}\ :\ {}\ \ooplus\limits_{\a, \
\a\succ\b}V_\a\ \to\ \ooplus\limits_{\a,\  \a\succ\b}V_\a \ ,
\end{equation}
where the summation is over all $\a, \a'$ such that
$ \a\succ\b,\ \a'\succ\b$ and
\begin{equation}      
\label{2.22a}
\Tprime_\b = \sum_{\a, \a'\!, \d} A_{\a,\d}A_{\d,\a'}\ :\ {}\ \ooplus\limits_{\a, \
\a\succ\b}V_\a\ \to\ \ooplus\limits_{\a,\  \a\succ\b}V_\a \ ,
\end{equation}
where the summation is over all $\a, \a',\d$ such that
$ \a\succ\b,\ \a'\succ\b$, and $ \d\succ\a,\ \d\succ\a'$.

{}For a vertex $\beta$, $l(\b)\not=N$, we set
\begin{equation}\label{2.22b}
\widetilde{S}_\b=\sum_{\a,\, \b\succ\a} A_\b^\a\ : \ {}\ V_\b\ \to\  V_\b\ .
\end{equation}
We also set $S_\ee=0: V_\ee\to V_\ee$ and $\widetilde{S}_\b=0$: $V_\b\to
V_\b$, if $l(\b)=N$. We define the operator $S: \oplus_{\a \in \v(\G)}V_\a
\to \oplus_{\a_\in \v(\G)}V_\a$ as
\begin{equation}\label{S}
S=\sum_{\a\in\v(\Ga)}(S_\a+\widetilde{S}_\a)=\sum_{\a,\b\in\v(\G)}
A_{\a,\b}A_{\b,\a}\ .
\end {equation}
\smallskip

 \begin{proposition}
 \label{lemmaTb}
${}\ \ {}$

\begin{enumerate}
\item[(i)] The operator $S_\b$ is a central element in the local
 algebra ${\cal O}_\b(\VA)$.
\item[(ii)]
 Any nonzero eigenvalue of\,  $T_\b$
 is equal to an eigenvalue of the operator 

\noindent $S_\a+A_\a^\b$  for some $\a$, $\a\succ\b$.
\item[(iii)]
 Any nonzero eigenvalue of  $\Tprime_\b$
 is equal to an eigenvalue of the operator 

\noindent $S_\a+A_\a^\b$  for some $\a$, $\a\succ\b$.
 \end{enumerate}
\end{proposition}
\begin{proposition}
\label{propcc}
Let $\A$ be a central arrangement (that is, all hyperplanes of $\A$

\noindent  contain
$0\in\CC^N$). Then
\begin{itemize}
\item[(i)]  The operator $\widetilde{S}_\b$ is a central
element in the local
 algebra ${\cal O}_\b(\VA)$.
\item[(ii)] For any $\a,\b\in\v(\G)$, we have
\begin{equation}\label{commutc}
(S_\a+\widetilde{S}_\a) A_{\a,\b}= A_{\a,\b}(S_\b+\widetilde{S}_\b)\ .
\end{equation}
In particular, the element $S$, defined in \rf{S},
is central in the path algebra $B_\G$.
\end{itemize}
\end{proposition}
 Propositions \ref{lemmaTb} and \ref{propcc} are proved in Section \ref{proofs}.

\bigskip

A subset $G\subset \C$ is called {\it non-resonant}
if  it satisfies the  following two conditions:
\begin{enumerate}
\item[(i)] The intersection $G\cap\ZZ$ is either  $\{0\}$ or empty.
\item[(ii)] For any $x,y$ in $G$, we have $x-y \in \C \setminus \ZZ$.
\end{enumerate}
A subset $G\subset \C$ is called {\it weakly non-resonant} if
it satisfies condition (ii).

\smallskip

 Suppose that for every vertex $\b\in \v(\Ga)$, $\b\neq \emptyset$, a
  non-resonant subset $G_\b\subset\C$ is chosen.  Denote by $Qui_{\G,
  \{G_\b\}}$ the category of all quivers $\{V_\a, A_{\a,\b}\}
  \in\Qui_\G$ such that
\begin{itemize}
\item[(a)]
 for
any $\b\in I(\G)$, \,$\b\neq \emptyset$, the
eigenvalues of the operator $T_\b$ belong to $G_\b$,
\item[(b)]
 for any $\b\in I(\G)$,\,$\b\neq \emptyset$, none of the
eigenvalues of the operator $\Tprime_\b$ is a positive integer.
\end{itemize}
 Suppose we have
\begin{itemize}
\item[]
  a non-resonant subset $G_\b\subset\C$
for every vertex $\b\in \v(\Ga)$, $\b\neq \emptyset$,
\item[]  a weakly non-resonant set $\widetilde{G}_\b$
for every vertex $\b\in \v(\Ga)$,
\item[] a weakly non-resonant
subset $G_{\a,\b}$  for any edge $(\a,\b)\in E(\Ga)$, $\a\succ\b$.
\end{itemize}
Denote by
 $\QG$ \, the  subcategory of all quivers
 $\{V_\a, A_{\a,\b}\}
\in \Qui_\G$, such that
\begin{itemize}
\item[(a)]
 for
any $\b\in I(\G)$, \,$\b\neq \emptyset$, the
eigenvalues of the operator $T_\b$ belong to $G_\b$ and
\item[(c)]
 for any $\b\in I(\G)$,\, 
   the
eigenvalues of the operator $\widetilde{S}_\b$ belong to
 $\widetilde{G}_\b$ and
\item[(d)]   for any edge $(\a,\b)\in E(\Ga)$, $\a\succ\b$, the eigenvalues
 of the operator $A_\a^\b$ belong to $G_{\a,\b}$.
\end{itemize}

 The categories $\QG$ \, and  ${\Qui}_{\G,\{G_\b\}}$
 are full abelian subcategories of the category
$\Qui_\G$.  The categories \QG \, and
${\Qui}_{\G,\{G_\b\}}$ are  closed under extensions.

We call a quiver  {\it non-resonant} if it belongs to the
subcategory ${\Qui}_{\G,\{G_\b\}}$   for some collection
 $\{G_\b\}$.
We call a quiver  {\it strongly non-resonant} if it belongs to the
subcategory $\QG$
for some collections  $\{G_\b\}$, $\{\widetilde{G}_\b\}$, and $\{G_{\a,\b}\}$.

\bigskip

\subsection{Quivers of level $n$} \label{Seclevn}
Suppose that
$\G$ is an admissible graph, $n\in\{0, \dots , N\}$, and $\G^n$
is the
corresponding truncated graph.

{\it A quiver $\VGn$ of} $\G^n$
 (or  {\it a level $n$  quiver})
consists of a collection of finite dimensional vector spaces $V_\a$,
 $\a\in\v(\G^n)$, a collection of
 linear maps
$A_{\a,\b}: V_\b\to V_\a,$\ $\a,\b\in\v(\G^n)$,
and a collection of
linear maps $A_\a^\b:\ V_\a\to V_\a$, labeled by loops $(\a,\a)^\b\in\e(\G^n)$.
The operators are subject to the following relations:
\begin{enumerate}
\item[(i)]
$A_{\a,\b}=0 $ if $\a, \b$ are not adjacent in $\G$;
\item[(ii)]
$\sum_{\b\in\v(\G^n)} A_{\a,\b}A_{\b,\g}=0 $ if
$ |l(\a)-l(\g)| = 2 $;
\item[(iii)]
$\sum_{\b\in\v(\G^n)} A_{\a,\b}A_{\b,\g} = 0$ if
$ l(\a)=l(\g),\;\a\not=\g $, and there exists $\d\in\v(\G^n)$,
such that $ \a\succ\d,$ 
$ \g\succ \d $ ;
\item[(iv)]
If $l(\a)=n$ and $ \d\succ\a\succ\b$ \; in $\G^{n+1}$, then
\begin{gather}\nonumber
 A_\a^\b A_{\a,\d}\ =\ A_{\a,\d}
\left(\sum_\g A_{\d,\g}A_{\g,\d}\right) \ ,\\ \nonumber
 A_{\d,\a}A_\a^\b \ =\ \
\left(\sum_\g A_{\d,\g}A_{\g,\d}\right) A_{\d,\a}\
\end{gather}
where the summation is over
$\g\in\v(\G^n)$ such that  $\g\not=\a,\
 \d\succ\g\succ\b$.
\item[(v)]  If $l(\a)=n$ and $\a\succ\b\succ\g$ in  $\G^{n+2}$, then
$$
\left[A_\a^\b\,,\,\sum_\d A_\a^{\d}\right] = 0
$$
where the summation is over
$\d\in\v(\G^n) ,\ \a\succ\d\succ\g$.

\end{enumerate}
A morphism of two level $n$ quivers is defined in the standard way.

The level $n$ quivers of $\Gamma^n$  form an abelian category  denoted by
$\Qui_{\G^n}$.
It is equivalent to the category of finite-dimensional modules over the
suitable path algebra.

Note that
 the structure of the defining  relations (i-v)
for a level $n$ quiver is completely determined by the adjacency matrix of the
level
$n+2$ truncated graph $\G^{n+2}$.

For $\a,\b \in \v(\Ga)$, let $\varepsilon(\b,\a)$ be defined as in Section
\ref{Sec.1.1}.

Let $\VA_{\G^n}$ $=$ $\{V_\a, A_{\a,\b}, A_\a^\b\}$ be a level $n$ quiver. Then
the collection ${\mathcal W}_{\G^n}$ $=$ $\{W_\a, B_{\a,\b}, B_\a^\b\}$, where
 $W_\a=V_\a^*$, $B_{\a,\b}= \varepsilon(\b,\a)A_{\b,\a}^*:\ W_\b\to W_ \a$,
 $B_{\a}^\b=- \left(A_\a^\b\right)^*:\ W_\a\to W_ \a$, determines a level $n$ quiver
$\tau_n(\VA_{\G^n})$ called {\it the dual} quiver. The assignment
$\VA_{\G^n} \mapsto \tau_n(\VA_{\G_n})$ determines an anti-automorphism of the category
$\Qui_{\G^n}$ of the fourth order, such that  for any level $n$
quiver $\VA_{\G^n}$ the quiver $\tau_n^2(\VA_{\G^n})$ is
isomorphic to $\VA_{\G^n}$.  The isomorphism
$\phi:\VA_{\G^n}\to \tau^2(\VA_{\G^n})$ is given by
the relation $\phi(v_\a)=(-1)^{l(\a)}v_\a$.

{}For a vertex $\b$  and an edge $(\a,\b)$ of $\G^n$,
denote  by $A_\b^\a$ the operator
$$
A_\b^\a\ =\ A_{\b,\a}A_{\a,\b}\ :\ {} \ V_\b\ \to\  V_\b\ .
$$

For a quiver $\VGn$ of level $n$ and
$\b\in I(\G^{n})$, \ $\b\neq \emptyset$, define  the {\it monodromy operators}
$$
T_\b:
\ooplus\limits_{\a,\ \a\succ\b}V_\a\to \ooplus\limits_{\a,\ \a\succ\b}V_\a,
\qquad
\Tprime_\b:
\ooplus\limits_{\a,\ \a\succ\b}V_\a\to \ooplus\limits_{\a,\ \a\succ\b}V_\a,
$$
by formulas \rf{2.22} and \rf{2.22a}
and for any $\b\in I(\G^{n})$ define the operator
$$
\widetilde{S}_\b : \ {}\ V_\b\ \to\  V_\b\ ,
$$
by formula \rf{2.22b}.

 Suppose that for every vertex $\b\in \v(\Ga^{n})$, $\b\neq
\emptyset$, a non-resonant subset $G_\b\subset\C$ is chosen.  Denote
by $Qui_{\G^n,\{G_\b\}}$ the category of all quivers $\{V_\a,
A_{\a,\b}, A_\a^\b\} \in\Qui_{\G^n}$ such that
\begin{itemize}
\item[(a)]
 for
any $\b\in I(\G^{n})$, \,$\b\neq \emptyset$, the
eigenvalues of the operator $T_\b$ belong to $G_\b$,
\item[(b)]
 for any $\b\in I(\G^{n})$,\, $\b\neq \emptyset$, none of the
eigenvalues of the operator $\Tprime_\b$ is a positive integer.
\end{itemize}

 Suppose we have
\begin{itemize}
\item[]
  a non-resonant subset $G_\b\subset\C$
for every vertex $\b\in \v(\Ga^{n})$, $\b\neq \emptyset$,
\item[]  a weakly non-resonant set $\widetilde{G}_\b$
for every vertex $\b\in \v(\Ga^n)$,
\item[] a weakly non-resonant
subset $G_{\a,\b}$  for every edge $(\a,\b)\in E(\Ga^n)$, $\a\succ\b$, and
loop $(\a,\a)^\b$ of $\G^n$.
\end{itemize}
Denote by
 $\QGn$ \, the  subcategory of all quivers
 $\{V_\a, A_{\a,\b},A_\a^\b\} \in Qui_{\G^n}$, such that
\begin{itemize}
\item[(a)]
 for
any $\b\in I(\G^n)$, \,$\b\neq \emptyset$, the
eigenvalues of the operator $T_\b$ belong to $G_\b$,
\item[(c)]
 for any $\b\in I(\G^n)$,\, 
   the
eigenvalues of the operator $\widetilde{S}_\b$ belong to
 $\widetilde{G}_\b$ and
\item[(d)]   for any edge $(\a,\b)\in E(\Ga^n)$, $\a\succ\b$, and
loop $(\a,\a)^\b$ the eigenvalues
 of the operator $A_\a^\b$ belong to $G_{\a,\b}$.
\end{itemize}

 The categories $\QGn$ \, and  ${\Qui}_{\G^n,\{G_\b\}}$
 are full abelian subcategories of the category
$\Qui_{\G^n}$.  The categories \QGn \, and
${\Qui}_{\G^n,\{G_\b\}}$ are  closed under extensions.

We call a level $n$ quiver  {\it non-resonant} if it belongs to the
subcategory ${\Qui}_{\G^n,\{G_\b\}}$   for some collection $\{G_\b\}$.
We call a level $n$ quiver  {\it strongly non-resonant} if it belongs to the
subcategory  $\QGn$
for some collections  $\{G_\b\}$,  $\{\widetilde{G}_\b\}$, and $\{G_{\a,\b}\}$.

\bigskip

\subsection{Restrictions  of quivers}\label{Rdiq}
 {}For integers $k,l$,  $0\leq k<l\leq N$, we define a functor
 $\j_{l,k}^* : \Qui_{\G^l} \to \Qui_{\G^k}$. If
$\VA_{\G^l}=\{V_\a, A_{\a,\b}\}$ is a quiver from $\Qui_{\G^l}$, then the
quiver $\j_{l,k}^*(\VA_{\G^l})=\{W_\a, B_{\a,\b},  B_{\a}^\b\}\in
\Qui_{\G^{k}}$ is defined as follows.

For every $\a\in I(\G^{k})$, the  space $W_\a$
 coincides with the space $V_\a$.
For every $\a,\b\in I(\Gamma^{k}),\, \a\neq\b$,
the map $B_{\a,\b}:W_\b\to W_\a$
coincides with the map $A_{\a,\b}:V_\b\to V_\a$.
For every loop $(\a,\a)^\b\in E(\Gamma^{k})$,
the map $B_\a^\b$
is equal to the product $A_{\a,\b}A_{\b,\a}$.

One can check that  relations (i-v) in Section 2.2
for the level $l$ quiver $\VA_{\G^l}$ imply
relations (i-v) in Section 2.2 for the level $k$  quiver
$\j_{l,k}^*(\VA_{\G^l})$.
If $\VA_{\G^l}=\{V_\a, A_{\a,\b}\}$ and
$\VA'_{\G^l}=\{V'_\a, A'_{\a,\b}\}$ are two level $l$ quivers,
and $\varphi=\{\varphi_\a:V_\a\to V'_\a\ |\ {\a\in\v(\G^l)}\}$ a morphism of
these quivers, then the collection
$\{\varphi_\a:V_\a\to V'_\a \ |\ {\a\in\v(\G^k)}\}$ defines a morphism
$\j_{l,k}^*(\varphi)$ of level $k$ quivers $\j_{l,k}^*(\VA_{\G^l})$ and
$\j_{l,k}^*(\VA'_{\G^l})$.

One can check that  the functor $\j_{l,k}^*$ is exact and it maps the subcategory
of non-resonant $\G^l$ quivers  to the subcategory of
non-resonant $\G^{k}$ quivers. For $0\leq k_1<k_2<k_3\leq N$, we
have
\begin{equation}\label{compj1}
 \j_{k_3,k_1}^*\ =\ \j_{k_2,k_1}^* \cdot \j_{k_3,k_2}^*\ ,
\end{equation}
and for any $0\leq k<l\leq N$
\begin{equation}\label{jtau}
 \j_{l,k}^*\cdot\tau_l\ =\ \tau_k\cdot\j_{l,k}^* \ .
\end{equation}
We will use the notation $\j_0^*$ for the functor
$\j_{N,0}^*\ :\ \Qui_{\G}=\Qui_{\G^N}\ \to\ \Qui_{\G^0}$.

The quiver
 $\j_{l,k}^*(\VA_{\G^l})$ is called {\it the  restriction at level $k$
of the quiver $\VA_{\G^l}$}.
The functors $\j_{l,k}^*$ are called {\it the restriction
functors}.
\medskip

\subsection{Direct images of quivers}\label{section2.4}
Now we recall the definition of the left and right adjoint functors.

Let ${\mathcal A}$ and ${\mathcal B}$ be  additive categories; let
$F:{\mathcal A}\to {\mathcal B}$   and
$G:{\mathcal B}\to {\mathcal A}$   be additive functors. The functor
$F$ is called {\it a left adjoint} to $G$ and the functor $G$ is  called
{\it a right adjoint}
 to $F$, if there exists a collection  of
isomorphisms
\begin{equation}
\label{radjoint}
\a_{X,Y}:\,\Hom_{{\mathcal A}}\left(X,G (Y)\right)\risom
\Hom_{{\mathcal B}}\left(F(X), Y\right),\qquad X\in {\mathcal A},\ {} \ Y\in
{\mathcal B} ,
\end{equation}
such that
\begin{enumerate}
\item [$\bullet$]
for any $X,X'\in {\mathcal A}$, $Y\in{\mathcal B}$,
$\phi\in \Hom_{{\mathcal A}}\left(X,G (Y)\right)$, $\gamma\in
\Hom_{{\mathcal A}}\left(X',X\right)$, we have\
 $\a_{X',Y}(\phi\gamma)=\a_{X,Y}(\phi)F(\g)$,

\item [$\bullet$]
 for any $X,\in {\mathcal A}$, $Y,Y'\in{\mathcal B}$,
$\phi\in \Hom_{{\mathcal B}}\left(F(X),Y)\right)$, $\delta\in
\Hom_{{\mathcal B}}\left(Y,Y'\right)$, we have
\ $\a_{X,Y'}(G(\delta)\phi)=\delta\cdot\a_{X,Y}(\phi)$.
\end{enumerate}

For a given additive functor
$F:{\mathcal A}\to {\mathcal B}$ its left  adjoint functor
$G: {\mathcal B}\to {\mathcal A}$ is unique up to an isomorphism, if exists.
Analogously, for a given additive functor
$G:{\mathcal B}\to {\mathcal A}$ its right  adjoint functor
$F: {\mathcal A}\to {\mathcal B}$ is unique up to an isomorphism, if exists.

The isomorphisms \rf{radjoint} determine
 collections of canonical morphisms
\begin{equation}\label{stau}
\sigma_{F,G}:\ X\to GF(X),\qquad
\tau_{F,G}:\  FG(Y)\to Y,\qquad X\in{\mathcal A},\ Y\in {\mathcal B} ,
\end{equation}
where
$\sigma_{F,G}=\a^{-1}_{X,F(X)}(\ {\rm id}_{F(X)})$ and
$\tau_{F,G}=\a^{}_{G(Y),Y}(\ {\rm id}_{G(Y)})$.
 The collections determine
 natural transformations  of the functors
 $\sigma_{F,G}:\ Id_{\mathcal A}\to GF$  and
 $\tau_{F,G}:\  FG\to Id_{\mathcal B}$, that is,
\begin{enumerate}
\item [$\bullet$]
for any $X,X'\in {\mathcal A}$, $\gamma\in
\Hom_{{\mathcal A}}\left(X',X\right)$ we have the equality
$\sigma_{F,G}\gamma=GF(\gamma)\sigma_{F,G}$ in
$\Hom_{{\mathcal A}}\left(X',GF(X)\right)$,
\item [$\bullet$]
for any $Y,Y'\in{\mathcal B}$,
 $\delta\in
\Hom_{{\mathcal B}}\left(Y,Y'\right)$  we have the equality
$\delta\tau_{F,G}=\tau_{F,G}FG(\delta)$ in
$\Hom_{{\mathcal B}}\left(FG(Y),Y'\right)$.
\end{enumerate}

\bigskip


Our next goal is to describe the functor
$\j_{n,n-1,*} : \Qui_{\G^{n-1}} \to \Qui_{\G^{n}}$, which is the right
adjoint functor to the functor $\j_{n,n-1}^*$.

Let $\VA_{\G^{n-1}}$ $=$ $\{V_\a, A_{\a,\b}, A_\a^\b\}$ be
 a level $n-1$ quiver.
Define a new level $n$ quiver
${\cal W}_{\G^n}$ $=$ $\{W_\a, B_{\a,\b}, B_\a^\b\}$  as follows.

{}For  $\a\in\v(\Ga^{n-1})$, set $W_\a=V_\a$.

{}For $\a\in \v_n(\Ga)$, denote by $\widetilde{W_\a}$ the direct sum
$\widetilde{W_\a}\,=
\,\oplus_{\g\in\v_{n-1}(\Ga),\ \g\succ\a}\,V_\g$.  Set $W_\a$ to be the
subspace of  $\widetilde{W_\a}$, generated by the vectors
$w=\sum_{\g\in\v_{n-1}(\Ga),\ \g\succ\a}v_\g$, $v_\g \in V_\g$, such that
for any $\d \in \v_{n-2}(\Ga)$ with $\d>\a$, the sum 
$\sum_{\g\in\v_{n-1}(\Ga),\ \g\succ\a}A_{\d,\g}(v_\g)$ equals zero.

Set $B_{\a,\b}=A_{\a,\b}$, if $\a$ and $\b$ both belong
to $\v(\Ga^{n-1})$.

Let $\a \in\v_{n-1}(\Ga)$,  $\b \in\v_{n}(\Ga)$, and $\a \succ\b$.
 Set  $B_{\a,\b}\,:\, W_{\b}\,\to\, W_{\a}$ to be the restriction to
the subspace $ W_{\b}$ of the projection of
the space $ \widetilde{W_{\b}}$
to the space $ W_{\a} =  V_{\a}$ along the other direct summands of
$\widetilde{W_{\b}}$. Define the map
 $B_{\b,\a}\,:\, W_{\a}\,\to\, W_{\b}$ by the relation
$$
B_{\b,\a}\ (v_{\a}) =  A_{\a}^{\b}\ (v_{\a})\ -\
\sum\ A_{\a', \d}
A_{\d,\a}\ (v_{\a})\ ,
$$
where $v_\a \in V_\a=W_\a$ and the sum is over all $\a' \in \v_{n-1}(\Ga)$
and $\d \in \v_{n-2}(\Ga)$, such that $\d \succ\a$ and $\d\succ\a'$.

Let $\a \in\v_{n}(\Ga)$,  $\b\in\v_{n+1}(\Ga)$, and $\a\succ\b$.
Define the map
$B_{\a}^{\b}\,:\, W_{\a}\,\to\, W_{\a}$  by the relation
\begin{equation}
\nonumber
B_{\a}^{\b} \,(\!\!\!\!\sum_{\g\in\v_{n-1}(\Ga),\ \g\succ\a}\!\! v_\g)\ =
\sum\limits_{\g,\ \g \succ\a}\ {} \
\sum\limits_{\d,\ \g \succ\d\succ\b,\, \d\not=\a} A_{\g}^{\d}\ (v_{\g}) ,
\end{equation}
where $v_\g\in V_\g$.

\begin{proposition} \label{jnn0}
The collection   ${\cal W}_{\G^n} = \{W_\a, B_{\a,\b}, B_\a^\b\}$
is  a level $n$ quiver.
\hfill{$\square$}
\end{proposition}

The assignment $\VA_{\G^{n-1}} \mapsto {\cal W}_{\G^n}$
determines a functor $\j_{n,n-1,*}: \Qui_{\G^{n-1}}\to \Qui_{\G^n}$.

Let $\VA'_{\G^{n}}$ $=$ $\{V'_\a, {A'}_{\a,\b}, {A'}_\a^\b\}$ be
 a level $n$ quiver and $\phi$ a morphism of $\G_{n-1}$ quivers
$\phi:\ \j^*_{n,n-1}(\VA'_{\G^n})\to \VA_{\G^{n-1}}$, given by a
collection of linear maps $\phi_\a: V'_\a\to V_\a$, $\a\in\v(\G^{n-1})$.

Define a new morphism of level $n$ quivers
$\widetilde{\phi}:  \VA'_{\G^n}\to {\cal W}_{\G^n}$ by the following rules.

Set $\widetilde{\phi}_\a=\phi_\a:{V'}_\a\to W_\a=V_\a$ if $l(\a)<n$.

Let $l(\b)=n$ and $v'_\b\in V'_\b$. Set $\widetilde{\phi}_\b(v'_\b)=
\sum_{\a, \a\succ\b}\phi_\a{A'}_{\a,\b}(v'_\b)$. We have
$\widetilde{\phi}_\b(v'_\b)\subset W_\b$.

\begin{proposition} \label{jnn0a}${}$\\
{\em (i)} The assignment
$
\b_{\VA'_{\G^n},\VA_{\G^{n-1}}} : \phi \mapsto \widetilde{\phi}$
determines an isomorphism
$$
\b_{\VA'_{\G^n},\VA_{\G^{n-1}}}:\Hom_{\Qui_{\G^{n-1}}}
(\j_{n,n-1}^*(\VA'_{\G^n}),\VA_{\G^{n-1}})\risom
\Hom_{\Qui_{\G^{n}}}(\VA'_{\G^n},\j_{n,n-1,*}(\VA_{\G^{n-1}})) \ .
$$
{\rm (ii)}
The functors $\j_{n,n-1}^*$ and $\j_{n,n-1,*}$ form a pair of
left and right adjoint functors
with respect to the collection of isomorphisms
 $$
\a_{\VA'_{\G^n},\VA_{\G^{n-1}}}:\
\Hom_{\Qui_{\G^{n}}}(\VA'_{\G^n},\j_{n,n-1,*}(\VA_{\G^{n-1}}))\risom
\Hom_{\Qui_{\G^{n-1}}}
(\j_{n,n-1}^*(\VA'_{\G^n}),\VA_{\G^{n-1}})\ ,$$ where
  $\a_{\VA'_{\G^n},\VA_{\G^{n-1}}}$ $=$
$\b_{\VA'_{\G^n},\VA_{\G^{n-1}}}^{-1}$.
\hfill{ $\square$}
\end{proposition}

{}For  $0\leq k<l\leq N$,  define the functors
$$
\j_{l,k,*}: \ \Qui_{\G^k} \to \Qui_{\G^l}\ ,
\qquad
 \j_{l,k,!}: \ \Qui_{\G^k} \to \Qui_{\G^l}
$$
by the formulas
\begin{equation}
\label{compj2}
\j_{l,k,*} \ =\ \j_{l,l-1,*} \cdots \j_{k+1,k,*}(\VA_{\G^k})\ ,
\qquad
\j_{l,k,!}\ =\ \tau^{-1}_l \cdot \j_{l,k,*} \cdot \tau_k^{}
 \ .
\end{equation}
The functors $\j_{l,k,*}$ and
 $\j_{l,k,!}$ are called {\em the quiver direct
image functors.}

For $0\leq k_1<k_2<k_3\leq N$, we have
$$
 \j_{k_3,k_1,*}\ =\ \j_{k_3,k_2,*}\ \cdot\ \j_{k_2,k_1,*} \ ,
 \qquad
\j_{k_3,k_1,!}\ =\ \j_{k_3,k_2,!}\ \cdot\ \j_{k_2,k_1,!}\ .
$$
We will use the notation $\j_{0,*}$ and $\j_{0,!}$, respectively,
 for the functors
$\j_{N,0,*}\ :\ \Qui_{\G^0}\ \to\ \Qui_{\G^N}=\Qui_{\G}$ and
$\j_{N,0,!}\ :\ \Qui_{\G^0}\ \to\ \Qui_{\G^N}=\Qui_{\G}$.

\bigskip

The functor $\j_{l,k,*}$
is the right adjoint to the functor $\j_{l,k}^*$. For any
$\VA'_{\G^l}$ and $\VA_{\G^k}$ the isomorphism
$\a_{\VA'_{\G^l},\VA_{\G^{k}}}:\
\Hom_{\Qui_{\G^{l}}}(\VA'_{\G^l},\j_{l,k,*}(\VA_{\G^{k}}))\risom
\Hom_{\Qui_{\G^{k}}}
(\j_{l,k}^*(\VA'_{\G^l}),\VA_{\G^{k}})\ $ is the composition
$$\a_{\VA'_{\G^l},\VA_{\G^{k}}}\ =\
\a_{\j_{l,k+1}^*(\VA'_{\G^l}),\VA_{\G^{k}}}\cdots
\a_{\j_{l,l-1}^*(\VA'_{\G^l}),\j_{l-2,k,*}(\VA_{\G^{k}})}
\a_{\VA'_{\G^l},\j_{l-1,k,*}(\VA_{\G^{k}})}\ .$$

The functor $\j_{l,k,!}$ is left adjoint to the functor
$\j_{l,k}^*$. For any $\VA'_{\G^l}$ and $\VA_{\G^k}$ the
corresponding isomorphism
$\a'_{\VA_{\G^k},\VA'_{\G^{l}}}:\
\Hom_{\Qui_{\G^{k}}}
(\VA_{\G^k},\j_{l,k}^*(\VA'_{\G^l}))\risom\Hom_{\Qui_{\G^{l}}}(\j_{l,k,!}(\VA_{\G^k}),\VA'_{\G^{l}})
\ $ is defined by the relation
$$\a'_{\VA_{\G^k},\VA'_{\G^{l}}}(\phi)=
 \tau_l^{-1}\cdot \a_{\tau_l(\VA'_{\G^l}),\tau_k(\VA_{\G^{k}})}
\cdot\tau_k(\phi)\ .$$

Thus $\j_{l,k}^*:\Qui_{\G^l}\to \Qui_{\G^k}$ and
  $\j_{l,k,*}:\Qui_{\G^k}\to \Qui_{\G^l}$ form a left and right
adjoint pair, as well as
$\j_{l,k,!}:\Qui_{\G^k}\to \Qui_{\G^l}$ and
$\j_{l,k}^*:\Qui_{\G^l}\to \Qui_{\G^k}$  form a left and right
adjoint pair.

\bigskip
\noindent
{\bf Remark.} Below we give a direct description of the functor
$\j_{n,n-1,!}$.

Let $\VA_{\G^{n-1}}= \{V_\a, A_{\a,\b}, A_\a^\b\}$ be a level $n-1$ quiver.
Define a new level $n$ quiver
${\cal U}_{\G^n}$ $=$ $\{U_\a, C_{\a,\b}, C_\a^\b\}$ as follows.

{}For  $\a\in\v(\Ga^{n-1})$, we set $U_\a=V_\a$.

{}For $\a\in \v_n(\Ga)$, we set
$\widetilde{U_\a}\,=
\,\oplus_{\g\in\v_{n-1}(\Ga),\ \g\succ\a}\,V_\g$. We set $U_\a$ to
 be the quotient of $\widetilde{U_\a}$ over the subspace,
generated by  the vectors $\sum_{\a,\ \d\succ\a\succ\b } A_{\a,\d} (v_\d)$,
where  $\d\in\v_{n-2}(\G)$,  $\d>\b$, and $v_\d \in V_\d$.
Let $\pi : \widetilde{U_\a} \to U_\a$ be the natural projection.

We set $C_{\a,\b}=A_{\a,\b}$, if $\a$ and $\b$ both belong
to $\v(\Ga^{n-1})$.

Let $\a \in\v_{n-1}(\Ga)$,  $\b \in\v_{n}(\Ga)$, and $\a \succ\b$.
 We set
  $C_{\b,\a}\,:\, U_{\a} \to U_{\b}$ to be the composition of the
inclusion
$U_{\a} = V_{\a} \hookrightarrow \widetilde{U}_{\b} = \,
\oplus_{\g\in\v_{n-1}(\Ga),\ \g\succ\a}\,V_\g$
and  the projection  $\pi: \widetilde{U}_{\b}\to U_\b$.
Define the map   $C_{\a,\b}\,:\, U_{\b}\,\to\, U_{\a}$, \
$u \mapsto C_{\a,\b}(u)$\ by the following rules:

\begin{enumerate}
\item[$\bullet$]
Assume that there exists $v_\a \in V_\a$ such that
$\pi (v_\a) = u$, then we set  $C_{\a,\b}\, (u) \ =\
\pi (A_\a^{\b}\, (v_\a))\ $.
\item[$\bullet$]
Assume that there exists $\g \in I_{n-1}(\G)$, $\g \neq \a$,
and $v_\g \in V_\g$ such that $\pi (v_\g) = u$.
Assume that there exists
$\d \in I_{n-2}(\G)$ such that $\d \succ\a, \d\succ\g$.
Then we set
 $C_{\a,\b}\, (u) \,=\,     - \pi (A_{\a,\d} A_{\d,\g}\, (v_\g))\ .$
\item[$\bullet$]
Assume that there exists $\g \in I_{n-1}(\G)$, $\g \neq \a$,
and $v_\g \in V_\g$ such that $\pi (v_\g) = u$.
Assume that there is no
$\d \in I_{n-2}(\G)$ such that $\d \succ\a, \d\succ\g$.
Then we set  $C_{\a,\b}\, (u) \,=\,0$.
\end{enumerate}

Let $\a \in\v_{n}(\Ga)$,  $\b\in\v_{n+1}(\Ga)$, and $\a\succ\b$.
Define the  map
$C_{\a}^{\b}\,:\, U_{\a}\,\to\, U_{\a}$  by the relation
\begin{equation}
\nonumber
C_{\a}^{\b} \,\pi\,(\!\!\!\!\sum_{\g\in\v_{n-1}(\Ga),\ \g\succ\a}\!\! v_\g)\ =
\sum\limits_{\g,\ \g \succ\a}\ {} \
\sum\limits_{\d,\ \g \succ\d\succ\b,\, \d\not=\a} \pi \,(A_{\g}^{\d}\ (v_{\g})) ,
\end{equation}
where $v_\g\in V_\g$.

\begin{proposition} \label{jnn01}
The collection   ${\cal U}_{\G^n} = \{U_\a, C_{\a,\b}, C_\a^\b\}$
is  a level $n$ quiver. This level $n$ quiver is isomorphic to
$\j_{n,n-1,!}(\VA_{\G^{n-1}})$ .
\hfill{$\square$}
\end{proposition}

\bigskip

\begin{proposition}\label{jnn} For any  $k<l$ and any level $k$ quiver
 $\VA_{\G^{k}}=\{V_\a, A_{\a,\b}, A_\a^\b \}$,
 the canonical morphisms
\begin{equation}\label{sqj}\begin{split}
\sigma_{\j_{l,k,!}, \j_{l,k}^*}
:\ \,\VA_{\G^{k}}\to
\j_{l,k}^*\j_{l,k,!}(\VA_{\G^{k}}) , &\quad{\rm }\\ 
\tau_{\j_{l,k}^*, \j_{l,k,*}}:\ \,
\j_{l,k}^*\j_{l,k,*}(\VA_{\G^{k}})\to
\VA_{\G^{k}}&\end{split}
\end{equation}
are isomorphisms.
\hfill
$\square$
\end{proposition}

The proofs of Propositions \ref{jnn0}--\ref{jnn} are direct
linear algebra calculations.

\bigskip

Define the morphism
$$
s_{l,k}:\ \, \j_{l,k,!}(\VA_{\G^k})\ \to\
\j_{l,k,*}(\VA_{\G^k})
$$
by the formula
\begin{equation}
\label{phikl}
s_{l,k}\ =\ \a_{\VA_{\G^k},\j_{l,k,*}(\VA_{\G^k}) }
\left( \tau^{-1}_{\j_{l,k}^*, \j_{l,k,*}} \right)\ ,
\end{equation}
where
$\a_{X,Y}$ are isomorphisms \rf{radjoint}, establishing the adjunction
 property of the functors $\j_{l,k,!}$ and $\j_{l,k}^*$.

{}For  $0\leq k<l\leq N$ and a level $k$ quiver
 $\VA_{\G^{k}}=\{V_\a, A_{\a,\b}, A_\a^\b \}$, denote by
$\j_{l,k,!*}(\VA_{\G^{k}}) \in \Qui_{\G^l}$ the image of $\j_{l,k,!}(\VA_{\G^{k}})$ in
 $\j_{l,k,*}(\VA_{\G^{k}})$ with respect to  $s_{l,k}$.

We use the notation $s_0$ for the morphism $s_{N,0}$ and
$\j_{0,!*}(\VA_{\G^{0}})$ for the quiver
$\j_{0,N,!*}(\VA_{\G^{0}})\in \Qui_{\G}$.

\subsection{Example}\label{example1}
Let $X$ be $\C^2$ with coordinates $z_1,z_2$. Let $\A = \{H_1, H_2, H_3\}$ be the
arrangement of three lines, $H_1 = \{z_1=0\}$, $H_2 = \{z_2=0\}$,
$H_1 = \{z_1=z_2\}$, intersecting at one point. The stratification
of $X$ consists of 5 strata:
\begin{enumerate}
\item[]
The complement to the union of lines,
$X_\emptyset=\{(z_1,z_2)\ | \ z_1\not=0, z_2\not=0, z_1\not=z_2\}$;
\item[]
The lines without the point,
$X_{\a_i}=H_i\setminus \{(0,0)\}$, $i=1,2,3$;
\item[]
The point, $X_\b=\{(0,0)\}$.
\end{enumerate}
The principal open subset $X^0$ is $X_\emptyset$.
The principal open subset $X^1$ is the complement to the point,
$X^1 = X\setminus \{(0,0)\}$. The principal open subset $X^2$  is
$X$.

The graph $\Ga$ of the arrangement has five vertices, $\emptyset,
\a_1,\a_2,\a_3, \b$, such that $l(\emptyset)=0$,
$l(\a_i)=1$,  $l(\b)=2$. The graph has six edges $(\emptyset,\a_i)$,
$(\a_i,\b)$. We have $\emptyset\succ \a_i \succ \b$.

The level two graph $\G^2$ coincides with $\G$.
The level one graph $\G^1$
has four vertices $\emptyset,
\a_1,\a_2,\a_3$, three edges $(\emptyset,\a_i)$
and three loops $(\a_i,\a_i)^\b$.
The level zero graph $\G^0$ consists of one vertex $\emptyset$ and
three loops $(\emptyset,\emptyset)^{\a_i}$.

\begin{enumerate}
\item[(i)]
A quiver $\VA_\G$ (the same as a level two quiver)
is a collection of five vector spaces,
$V_\emptyset$, $V_\b$, $V_{\a_i}$, $i=1,2,3$, and twelve linear
maps,  $A_{\a_i,\emptyset}: V_\emptyset\to V_{\a_i}$,
 $A_{\emptyset,\a_i}:  V_{\a_i}\to V_\emptyset$,
$A_{\a_i,\b}: V_\b\to V_{\a_i}$,
 $A_{\b,\a_i}:  V_{\a_i}\to V_\b$, such that
$$
\sum_{i=1}^3A_{\b,\a_i}A_{\a_i,\emptyset}=0,\qquad
\sum_{i=1}^3A_{\emptyset,\a_i}A_{\a_i,\b}=0,$$
and
$$A_{\a_j,\emptyset}A_{\emptyset,\a_i}+A_{\a_j,\b}A_{\b,\a_i}=0$$
for all pairs $(i,j)$, $i\not=j$.
\item[(ii)]
A level one quiver $\VA_{\G^1}$ is a collection of four vector spaces,
$V_\emptyset$, $V_{\a_i}$, $i=1,2,3$, and nine linear
maps $A_{\a_i,\emptyset}: V_\emptyset\to V_{\a_i}$,
 $A_{\emptyset,\a_i}:  V_{\a_i}\to V_\emptyset$,
$A_{\a_i}^\b: V_{\a_i}\to V_{\a_i}$,
 such that
\begin{equation}\nonumber
\begin{split}
A_{\a_i}^\b A_{\a_i,\emptyset}=A_{\a_i,\emptyset}
\sum_{j\not=i}A_{\emptyset,\a_j}A_{\a_j,\emptyset},\\
A_{\emptyset,\a_i}A_{\a_i}^\b=
\sum_{j\not=i}A_{\emptyset,\a_j}A_{\a_j,\emptyset}
A_{\emptyset,\a_i}.\nonumber
\end{split}
\end{equation}
\item[(iii)]
A level zero quiver $\VA_{\G^0}$ consists of the vector space $V=V_\emptyset$
 and three linear maps $A^i: V\to V$, $i = 1, 2, 3$.
subject to the relations
$$[A^i,A^1+A^2+A^3] = 0 ,
\qquad
i = 1, 2, 3 .
$$
\end{enumerate}

Let $\VA_\G$ be a  quiver as in (i). Then its restriction
$\j^*_{2,1}(\VA_\G)$ to level one is the  level one quiver
$\{W_{\emptyset}, W_{\a_i}, B_{\emptyset, \a_i}, B_{\a_i,\emptyset},
 B_{\a_i}^\b\}$,
where $W_\emptyset=V_\emptyset$, $W_{\a_i}=V_{\a_i}$, $B_{\emptyset,\a_i}=
A_{\emptyset,\a_i}$, $B_{\a_i,\emptyset}=
A_{\a_i,\emptyset}$, and $B_{\a_i}^\b = A_{\a_i,\b}A_{\b,\a_i}$.

\medskip

Let $\VA_{\G^1}$ be a  level one quiver as in (ii).
Its restriction
$\j^*_{1,0}(\VA_{\G^1})$ to level zero is  the  level zero quiver
$\{W_{\emptyset},  B^i,\}$,
where $W_\emptyset=V_\emptyset$, $B^i= A_{\emptyset, \a_i}A_{\a_i,\emptyset}$.

\medskip

Let $\VA_{\G^0}$ be a  level zero quiver as in (iii).
Its direct image $\j_{1,0,*}(\VA_{\G^0})$ is the level one quiver
$\{W_{\emptyset}, W_{\a_i}, B_{\emptyset, \a_i}, B_{\a_i,\emptyset},
 B_{\a_i}^\b\}$,
where $W_\emptyset=V$, $W_{\a_i}=V$, $B_{\emptyset,\a_i}=
{\rm Id}$, $B_{\a_i,\emptyset}=
A^i$, and $B_{\a_i}^\b =\sum_{j\in\{1,2,3\},\, j\not=i} A^j$.

\medskip

The direct image $\j_{1,0,!}(\VA_{\G^0})$  is the level one quiver
$\{U_{\emptyset}, U_{\a_i},  C_{\emptyset,\a_i},  C_{\a_i,\emptyset},
 C_{\a_i}^\b \}$,
where $U_\emptyset=V$, $U_{\a_i}=V$, $C_{\a_i,\emptyset}=
{\rm Id}$, $C_{\emptyset,\a_i}=
A^i$, and $C_{\a_i}^\b =\sum_{j\in\{1,2,3\},\, j\not=i} A^j$.

\medskip

Let $\VA_{\G^1}$ be a  level one quiver as in (ii).
Its direct image $\j_{2,1,*}(\VA_{\G^1})$ is the level two quiver
$\{W_{\emptyset}, W_{\a_i}, W_\b, B_{\emptyset, \a_i}, B_{\a_i,\emptyset},
B_{\b, \a_i}, B_{\a_i,\b} \}$,
where  $W_\emptyset=V_\emptyset$, $W_{\a_i}=V_{\a_i}$,
$W_\b$ is the subspace of the direct sum
$\widetilde{W}_\b=\oplus_{i}V_{\a_i}$, generated by the vectors
$\sum_{i=1}^3v_{\a_i}$,
 such that $\sum_i A_{\emptyset, \a_i}(v_{\a_i})=0$.
The maps $B_{\emptyset,\a_i}$ are equal to $A_{\emptyset,\a_i}$,
The maps $B_{\a_i,\emptyset}$ are equal to $A_{\a_i,\emptyset}$.
The maps $B_{\a_i,\b}$ and $B_{\b,\a_i}$ are defined by the formulas
$B_{\a_i,\b}(\sum_j v_{\a_j})=v_{\a_i}$ and
$B_{\b,\a_i}(v_{\a_i})=A_{\a_i}^\b(v_{\a_i})-\sum_{j\in\{1,2,3\},\, j\not=i}
A_{\a_j,\emptyset}A_{\emptyset,\a_i}(v_{\a_i})$.

\medskip

The direct image $\j_{2,1,!}(\VA_{\G^1})$  is the level two quiver
$\{U_{\emptyset}, U_{\a_i}, U_\b, C_{\emptyset, \a_i}, C_{\a_i,\emptyset}$,
$C_{\b, \a_i}, C_{\a_i,\b}, \}$, where
$U_\emptyset=V_\emptyset$ and $U_{\a_i}=V_{\a_i}$.
The space $U_\b$ is the quotient space of
$\widetilde{U}_\b=\oplus_{i}V_{\a_i}$ over the subspace, generated
 by the vectors
$\sum_{i=1}^3A_{\a_i,\emptyset}(v_\emptyset)$,
 $v_\emptyset\in V_\emptyset$. Let $\pi: \widetilde{U}_\b\to U_\b$
be the canonical projection.
The maps $C_{\emptyset,\a_i}$ are equal to $A_{\emptyset,\a_i}$.
The maps $C_{\a_i,\emptyset}$ are equal to $A_{\a_i,\emptyset}$.
The maps  $C_{\b,\a_i}$ are defined by the formulas
$C_{\b,\a_i}(v_{\a_i}) = \pi (v_{\a_i})$. The
maps  $C_{\a_i,\b}$ are defined by the formulas
$C_{\a_i,\b}(\pi(v_{\a_i})) = \pi(A_{\a_i}^\b(v_{\a_i}))$ and
$C_{\a_i,\b}(\pi(v_{\a_j}))=-
\pi(A_{\a_i,\emptyset}A_{\emptyset,\a_j}(v_{\a_j}))$ if  $j\not=i$.

\subsection{Direct images of level zero quivers} \label{secdir}
The direct images of level zero quivers can be obtained by inductive
application of the construction of the previous section. Here we give a
direct description of the functors $\j_{0,*}$ and $\j_{0,!}$.

A level zero quiver $\UG_{\G^0} \in \Qui_{\G^0}$
 consists of a vector space $W= W_\emptyset$ and
 linear operators $B^i=B_\emptyset^i:W\to W$,
$i\in\v_1(\G)$, such that for each $\a\in\v_2(\G)$, $i\succ\a$,
we have
\begin{equation}\label{3.11}
[ B^i , \sum_{j,\ \, j\succ\a} B^{j} ]\ =\ 0\ .
\end{equation}

Let $\UG_{\G^0}=\{W,B^i\}$
 be a level zero quiver. Define a
quiver ${\VA}_\G = \{V_\a, A_{\a,\b}\} \ \in \Qui_\G$ as follows.

For $\a\in I(\Gamma)$, set $V_\a = \FF_\a\ot W$,
where $\FF_\a$ is the subspace of the flag space $\FF(\A)$,
corresponding to $\a\in\v(\Ga)$.

For $\b\succ\a $, $l(\b)=m$ define the map
$A_{\a,\b} : {V}_\b \to {V}_\a$ by the rule
\begin{equation}
\label{3.13}
A_{\a,\b}\ :\ F_{\b_0,...,\b_{m-1},\b}\otimes w\ \mapsto \
 (-1)^{m}F_{\b_0,...,\b_{m-1},\b, \a}\otimes w \ .
\end{equation}

To define  the maps $A_{\a,\b} : {V}_\b \to {V}_\a$ for $\a\succ\b$
we need the following construction.

Let $F_{\b_0,..., \b_{m-1},\b_m}
\in\FF(\A)$ be a flag and $\a\in\v_{m-1}(\Ga)$,
$\a\succ\b_m$, a vertex. Assume that there exists a flag
$F_{\a_0,...,\a_{m-1}}\in\FF(\A)$ such that
$\a_{m-1}=\a$ and  $\a_{k}\succ\b_{k+1}$ for
$k = 0, \dots , m-1$. It is easy to see such a flag is unique.
We call $F_{\a_0,...,\a_{m-1}}$ {\it the cutoff} of
$F_{\b_0,...,\b_m}$ in the direction of $\a$ and denote by
$F^{\a}_{\b_0,...,\b_m}$. If such a flag does not exists, then we say
that the cutoff is zero and set $F^{\a}_{\b_0,...,\b_m} = 0$.

We define the number
$(\a; \b_0,...,\b_{m-1},\b_{m})$ to be
the maximal codimension in which the flags $F_{\b_0,...,\b_{m-1}}$
and $F_{\a_0,...,\a_{m-1}}$ coincide. More precisely, we define
$(\a; \b_0, ... , \b_{m}) = k$ if
$\b_0=\a_0, ... , \b_k = \a_k$ for $k < m-1$ and $\b_{k+1} \not= \a_{k+1}$.
We define
$(\a; \b_0,...,\b_{m}) = m-1$ if
$\b_k=\a_k$ for $k = 0,  ... , m-1$. We define $(\a; \b_0,...,\b_{m}) = 0$ if
  $F^{\a}_{\b_0,...\b_m=\b} = 0$.

For $\b\succ\a $, define the map
$A_{\a,\b} : {V}_\b \to {V}_\a$ by the rule
\begin{equation}
\label{flag-}
 A_{\a,\b}\ :\ F_{\b_0,...,\b_{m-1},\b_m=\b} \ot w
\ \mapsto\
 (-1)^{k}
F_{\b_0,...,\b_{m-1},\b_m}^\a \ot
\sum  B^i(w)\ ,
 \end{equation}
where  $k = (\a; \b_0,...,\b_m)$ and the sum is
over all
$i \in \v_1(\Ga)$ such that $i\geq\a_{k+1}$ and
$i\not\geq \a_{k}$.
 Another form of formula \rf{flag-} is
$$
 A_{\a,\b} \ :\ F_{\b_0,...,\b_{m-1},\b_m=\b}\ot w\ \mapsto\
 (-1)^{k}
 F_{\b_0,...,\b_k,\a_{k+1},...,\a_{m-1}=\a} \ot
\sum  B^i(w) \ ,
$$
where the sum is over $i\in\v_1(\Ga)$ such that
$i\wedge \b_k=\b_{k+1},\, i\wedge\a_{k+1}=\b_{k+2},\, ...
,   i\wedge \a_{m-1}=\b_{m}$.
 Here $k = (\a; \b_0,...,\b_m)$
  and $ F_{\b_0,...,\b_k,\a_{k+1},...,\a_{m-1}} =
F_{\b_0,...,\b_{m-1},\b}^\a$.

\begin{proposition}\label{5}
${}$

\noindent
{\em (i)}
The spaces $V_\a$ and operators $A_{\a,\b}$ form a quiver
${\VA}_\G \in \Qui_\G$.\\
{\em (ii)} The quiver ${\VA}_\G$ is isomorphic to
$\j_{0,!}(\UG_{\G^0})$.\\
{\em (iii)} For any quiver ${\cal U} =\{U_\a,
C_{\a,\b}\}\in\Qui_\G$ and a morphism
 $\varphi \in \ \Hom_{\Qui_{\G^0}} (\UG_{\G^0}, \j^*_0({\cal U}))$,
the adjoint morphism $\a_{{\cal W}_{\G^0},{\cal U}}(\varphi)$ $\in$
$\Hom_{\Qui} (\j_{0,!}\,( \UG_{\G^0}) , {\cal U} )$ is given
 by the relations
\begin{equation}
\a_{{\cal W}_{\G^0},{\cal U}}(\varphi)\ :
\begin{array}{ccl}\ F_\ee\ot w&
\mapsto\ & \varphi(w)\\
\label{3.15}
  F_{\a_0,...,\a_m}\otimes w&
\mapsto\  &(-1)^{m-1} C_{\a_m\a_{m-1}} \cdots  C_{\a_1\a_0}\
 \varphi ( w ) \ , \quad m\geq 1 \ .
\end{array}
\end{equation}

{\em (iii)} The functor $\j_{0,!} : { \Qui}_{\G^0}\to { \Qui_\G}$ is exact.
\label{prop2}
\hfill{$\square$}
\end{proposition}

\medskip

{}For a level zero quiver  $\UG_{\G^0}=\{W,B^i\}$,
define a quiver ${\widetilde{\VA}}_\G= \{\widetilde V_\a, \widetilde A_{\a,b})$
$\in \Qui_\G$ as follows.

Let $\OS(\G)$ be the Orlik-Solomon algebra of the graph $\Ga$.
{}For $\a\in\v(\Ga)$, define the projection
 operator
$P_\a:\OS(\Ga)\to\OS(\Ga)$ by the formula
 $$
P_\a(H_{j_1},...,H_{j_m})=\left\{
 \begin{array}{ccc}
(H_{j_1},...,H_{j_m}),&{\rm if}& j_1\wedge...\wedge
j_m=\a,\\
0,&{\rm otherwise}&
\end{array}\right .
$$
Set $\widetilde{V}_\a = P_\a(\OS)\ot W$.   Define the linear maps
$\widetilde{A}_{\a,\b} : \widetilde{V}_\b \to \widetilde{V}_\a$  by the formula
 \begin{equation}
\widetilde{A}_{\a,\b} : (H_{j_1},...,H_{j_m})\ot
w  \ \mapsto \
\sum\limits_{k=1}^{m}\ (-1)^{k-1} P_\a((H_{j_1},...,\widehat{H}_{j_k},...,H_{j_m}))\ot
w
\label{3.16}
\end{equation}
if $\a\succ\b$, and by the formula
\begin{equation}
\widetilde{A}_{\a,\b} : (H_{j_1},...,H_{j_m})\ot
w \ \mapsto \ \sum\limits_{k \in \v_1(\Ga)}\  P_\a((H_k,H_{j_1},...,H_{j_m}))\ot
B^k(w).
\label{3.17}
\end{equation}
if $\b\succ\a$.

\begin{proposition}\label{6}
${}$

\noindent
{\em (i)}
The spaces $\widetilde V_\a$ and linear maps $\widetilde A_{\a,\b}$ form a quiver
$\widetilde{\VA}_\G \in \Qui_\G$.\\
{\em (ii)} The quiver $\widetilde{\VA}_{\G}$
 is isomorphic to $\j_{0,*}(\UG_{\G^0})$.\\
{\em (iii)} The morphism $s_0 :\ \j_{0,!}(\UG_{\G^0}) \to
\j_{0,*}(\UG_{\G^0})$ is given by the formula
\begin{equation}\label{s0}
s_0:  F_{\a_0,...,\a_m}\otimes w
\mapsto\ \sum\left(H_{j_1},\ldots,H_{j_m}\right)
\ot B^{j_m}\cdots B^{j_1}(w)\ ,
\end{equation}
where the sum is taken over all tuples $(H_{j_1}, ... ,H_{j_m})$, such that
$H_{j_1}\supset\X_{\a_1}$, $...$,
 $H_{j_m}\supset \X_{\a_m}$.\\
{\em (iv)} The functor $\j_{0,*} : { \Qui}_{\G^0} \to { \Qui_\G}$ is
exact. 
\label{prop3}
\hfill{$\square$}
\end{proposition}

The morphism $s_0 :\ \j_{0,!}(\UG_{\G^0}) \to
\j_{0,*}(\UG_{\G^0})$ is called {\em the quiver Shapovalov map}.
It is a
 matrix generalization of the Shapovalov map $S_a$ from the flag
complex to the Aomoto complex, see Section \ref{Orlik}.

Identifying $\OS(\A)$ with $\FF^*(\A)$, as in Section \ref{Orlik},
 we can consider
 the map $s_0$ as a bilinear form on the vector space $\FF(\A)$ with values
in ${\rm End} \ W$. This {\em  quiver Shapovalov form}
 is given by the formula:
$$
s_0(F_{\a_0,\a_1,...,\a_m}, F_{\a_0,\a'_1,...,\a'_m})=
\sum\limits_{\sigma\in S_m}(-1)^{l(\sigma)}\sum B^{j_m}\cdots B^{j_1}\ ,
$$
where the inner sum is taken over all $j_1,...,j_m\in J(\A)$, such that
$H_{j_1}\supset\X_{\a_1}$, $...$, $H_{j_m}\supset\X_{\a_m}$, and
$H_{j_{\sigma(1)}}\supset\X_{\a'_1}$, $...$,
$H_{j_{\sigma(m)}}\supset\X_{\a'_m}$.

If all operators $B^j:W\to W$, $j\in J(\A)$, commute, then the quiver
Shapovalov form is symmetric.


\bigskip

\subsection{Monodromy operators  for quivers
$\j_{0,!}(\UG_{\G^0})$ and $\j_{0,*}(\UG_{\G^0})$}\label{Spectrum}
Suppose that
each operator $B^i$ of the level zero quiver $\UG_{\G^0}=\{W, B^i\}$
 has a single eigenvalue $\lambda_i$.
 We call the collection $\{\lambda_i\}$, $i\in\v_1(\Ga)$,
{\it the spectrum} of  $\UG_{\G^0}$.
{}For  $\a\in \v(\Ga)$, set
\begin{equation}
\l_\a\ =\ \sum_{i\in\v_1(\Ga),\ {} i\geq\a} \l_i \ .
\label{lambdaa}
\end{equation}
We call the spectrum $\{\l_i\}$ {\em non-resonant} if
\begin{equation}\label{nonres}
\l_\a\not\in\ZZ\setminus \{0\}\qquad {\rm for\ all}\quad \a\in\v(\Ga) .
\end{equation}
Let $\j_{0,!} (\UG_{\G^0}) = \{ V_\a, A_{\a,\b}\}$ and
$\j_{0,*} (\UG_{\G^0}) = \{ \widetilde V_\a, \widetilde A_{\a,\b}\}$.

\begin{proposition}
${}$

\noindent
{\em (i)} For $\a\in\v(\Ga)$ and the quiver $\j_{0,!} (\UG_{\G^0})$,
consider the central element $S_\a$
of the local algebra ${\cal O}_\a( \j_{0,!} (\UG_{\G^0}) )$, see
 \rf{Tb0}. Then $S_\a$  is the scalar operator
with eigenvalue  $\l_\a$.
\\
{\em (ii)} For  $\a\in\v(\Ga)$, $\a\succ \b$, and the quiver
$\j_{0,!} (\UG_{\G^0})$, consider the operator
$A_\a^\b=A_{\a,\b}A_{\b,\a}\, :\, V_\a \,\to\, V_\a$. Then
$A_\a^\b$ is scalar with eigenvalue  $\l_\b-\l_\a$.
\\
{\em (iii)} For $\a\in\v(\Ga)$ and the quiver $\j_{0,!} (\UG_{\G^0})$,
consider the operators $T_\a$ and $\Tprime_\a$, defined in \rf{2.22} and
\rf{2.22a}, respectively.
Then every eigenvalue of $T_\a$ is either $0$ or   $\lambda_\a$ and
every eigenvalue of $\Tprime_\a$ is either $0$ or   $\lambda_\a$.
\label{lemma2}
\end{proposition}

\begin{proposition}\label{8}
${}$

\noindent
{\em (i)} For $\a\in\v(\Ga)$ and the quiver  $\j_{0,*} (\UG_{\G^0})$,
consider the central element $S_\a$
of the local algebra ${\cal O}_\a( \j_{0,*} (\UG_{\G^0}) )$, see
 \rf{Tb0}. Then $S_\a$  is the scalar operator
with eigenvalue  $\l_\a$.
\\
{\em (ii)} For  $\a\in\v(\Ga)$, $\a\succ \b$, and the quiver
 $\j_{0,*} (\UG_{\G^0})$, consider the operator
$\tilde A_\a^\b = \tilde A_{\a,\b}\tilde A_{\b,\a}\, :\,
\tilde V_\a \,\to\, \tilde V_\a$. Then
$\tilde A_\a^\b$ is scalar with eigenvalue  $\l_\b-\l_\a$.
\\
{\em (iii)} For $\a\in\v(\Ga)$ and the quiver $\j_{0,*} (\UG_{\G^0})$,
consider the operators $T_\a$ and $\Tprime_\a$, defined in \rf{2.22} and
\rf{2.22a}, respectively.
Then every eigenvalue of $T_\a$ is either $0$ or   $\lambda_\a$ and
every eigenvalue of $\Tprime_\a$ is either $0$ or   $\lambda_\a$.
\end{proposition}

\begin{proposition}\label{propcmon}
Let $\A$ be a central arrangement (that is, all hyperplanes of $\A$
 contain
$0\in\CC^N$). Then for each of the  quivers
$\j_{0,!} (\UG_{\G^0})$ and
  $\j_{0,*} (\UG_{\G^0})$,
 the operator  $S$, defined in  \rf{S}, has a single eigenvalue
which is equal to $\l_\infty=\sum_{j\in J(\A)}\l_j$.
\end{proposition}

\begin{corollary} \label{corspec} ${}$

\noindent {\em (i)}
Let  $\UG_{\G^0}=\{W, B^i\}$ be a level zero quiver
 with non-resonant spectrum.
 Then the quivers  $\j_{0,!} (\UG_{\G^0})$,
  $\j_{0,*} (\UG_{\G^0})$, and $\j_{0,!*} (\UG_{\G^0})$
 are  non-resonant. The same conclusion is valid if all operators
$B^i$ are close to zero.

\noindent {\em (ii)} Let $\A$ be a central arrangement and
 $\UG_{\G^0}=\{W, B^i\}$  a level zero quiver
 with non-resonant spectrum.
 Then the quivers  $\j_{0,!} (\UG_{\G^0})$,
  $\j_{0,*} (\UG_{\G^0})$ and $\j_{0,!*} (\UG_{\G^0})$
 are  strongly non-resonant. The same conclusion is valid if all operators
$B^i$ are close to zero.
\end{corollary}

Statements of the corollary for level zero quivers with non-resonant spectrum
follow directly from Propositions \ref{lemma2}-\ref{propcmon}. If
 $\UG_{\G^0}=\{W, B^i\}$ is a level zero quiver with
operators $B^i$  close to zero, then all monodromy operators of
quivers  $\j_{0,!} (\UG_{\G^0})$,
  $\j_{0,*} (\UG_{\G^0})$, and $\j_{0,!*} (\UG_{\G^0})$ are close
to zero and their eigenvalues belong to some non-resonant sets.
\hfill{$\square$}

Propositions \ref{lemma2}-- \ref{propcmon}  are
proved in Section \ref{proofs}.

\bigskip
\subsection{Example}\label{example2}
Let $X$ be $\C^2$ and $\A = \{H_1, H_2, H_3\}$ the
arrangement of three lines, considered in Section \ref{example1}.

In notations of Sections \ref{Orlik} and \ref{example1},
the flag complex $\FF(\A)$ of the arrangement has seven  flags
$F_\ee$, $F_{\ee,\a_1}$,  $F_{\ee,\a_2}$,  $F_{\ee,\a_3}$,
$F_{\ee,\a_1,\b}$,  $F_{\ee,\a_2,\b}$,  $F_{\ee,\a_3,\b}$,
subject to the relation
$$
F_{\ee,\a_1,\b}\ +\ F_{\ee,\a_2,\b}\  +\ F_{\ee,\a_3,\b}\ =\ 0\ .
$$
We have $\FF_\ee=\CC\cdot F_\ee\ $,\ $\FF_{\a_i} = \CC \cdot F_{\ee,\a_i}$,\
 $i=1,2,3,\ $ and $\FF_\b$ is generated by the flags
$F_{\ee,\a_i,\b}$, $i=1,2,3$.

The Orlik--Solomon algebra $\OS(\A)$ is generated
by symbols $1$, $(H_i)$ with $i=1,2,3,$ and by symbols$(H_i,H_j)$, where $i,j = 1,2,3$.
We have  $(H_i,H_j)$ $=$ $-(H_j,H_i)$ and
$$
(H_1,H_2)\ -\ (H_1,H_3)\ +\ (H_2,H_3)\ =\ 0\ .
$$

Let $\VA_{\G^0}$ be a level zero quiver as in Section \ref{example1}, (iii).

The direct image $\j_{0,!}(\VA_{\G^2})$ is the 
 quiver
$\{W_{\emptyset}, W_{\a_i}, W_\b, B_{\emptyset, \a_i}, B_{\a_i,\emptyset},
B_{\b, \a_i}, B_{\a_i,\b} \}$ of $\G$,
where  $W_\g=\FF_\g\ot V$ for all $\g\in\v(\G)$, and
$$
B_{\a_i,\emptyset}(F_\ee\ot v)\ =\ F_{\ee,\a_i}\ot v\ ,
\qquad
 B_{\emptyset, \a_i}(F_{\ee,\a_i}\ot v)\ =\ F_\ee\ot A^i(v)\ ,
$$
\begin{equation}
\begin{split}\nonumber
B_{\b, \a_i}(F_{\ee,\a_i}\ot v)&=-F_{\ee,\a_i,\b}\ot v\ ,\\ \nonumber
 B_{\a_i,\b}(F_{\ee,\a_i,\b}\ot v)&=
F_{\ee,\a_i}\ot \sum_{j,\ j\not=i}A^j(v)\ ,\\ \nonumber
 B_{\a_j,\b}(F_{\ee,\a_i,\b}\ot v)&=-F_{\ee,\a_j}\ot A^i(v)\ ,\quad j\not= i\ .
\end{split}
\end{equation}

The direct image $\j_{0,*}(\VA_{\G^0})$  is the 
 quiver
$\{U_{\emptyset}, U_{\a_i}, U_\b, C_{\emptyset, \a_i}, C_{\a_i,\emptyset}$,
$C_{\b, \a_i}, C_{\a_i,\b}, \}$ of $\G$, where $U_\ee=\CC\cdot 1\ot V$,
$U_{\a_i}=\CC\cdot (H_i)\ot V$, and $U_\b$ is generated
by the vectors $(H_i,H_j)\ot v$, $i,j=1,2,3$, $v\in V$.
We have
$$
C_{\a_i,\ee}(1\ot v)\ =\ (H_i)\ot A^i(v)\ ,
\qquad
C_{\ee, \a_i}((H_i)\ot v)\ =\ 1\ot v\ ,
$$
$$
C_{\b, \a_i}((H_i)\ot v)\ =\ \sum_{j,\ j\not=i}(H_j,H_i)\ot A^j(v)\ ,
$$
$$
C_{\a_i,\b}((H_j,H_i)\ot v)\ =\ (H_i)\ot v\ ,
\qquad
j \not = i\ ,
$$
$$
C_{\a_i,\b}((H_j,H_k)\ot v)\ =\ 0\ ,
\qquad i \notin \{j, k\} .
$$

The Shapovalov
map $s_0 :\ \j_{0,!}(\VA_{\G^0}) \to \j_{0,*}(\VA_{\G^0})$ is given
by the formulas:
$$
s_0(F_\ee\ot v)\ =\ v\ ,
\qquad
s_0(F_{\ee,\a_i}\ot v)\ =\ (H_i)\ot A^i(v)\ ,
$$
$$
s_0(F_{\ee,\a_i,\b}\ot v)\ =\ \sum_{j,\ j\not=i}(H_i,H_j)\ot A^jA^i(v)\ .
$$
 The quiver $\j_{0,!*}(\VA_{\G^0})$ is the subquiver of
$\j_{0,*}(\VA_{\G^0})$, generated by the vectors $1\ot v$, $(H_i)\ot A^i(v)$,
$\sum_{j,\ j\not=i}(H_i,H_j)\ot A^jA^i(v)$, where $i=1,2,3$ and $v\in V$.

Suppose that each of the
operators $A^i: V\to V$ of the level $0$ quiver $\VA_{\G^0}$ has
a single eigenvalue $\l_i$. Then the quivers $\j_{0,!}(\VA_{\G^0})$,
$\j_{0,*}(\VA_{\G^0})$ and $\j_{0,!*}(\VA_{\G^0})$ are strongly non-resonant
if $\l_1+\l_2+\l_3\not\in \ZZ\setminus \{0\}$ and
$\l_i\not\in \ZZ\setminus \{0\}$ for all $i=1,2,3$.

\bigskip

\def\MMD{{\widetilde{\mathcal M}}}

\section{Quiver $\D$-modules}

\subsection{Modules over rings of differential operators}
 Here we briefly recall the basic terminology of modules over rings of
 differential operators \cite{Kash3,Kash4,B}.

Let $Y$ be a smooth  algebraic variety over $\C$.
 We denote  by
 $\D_Y$ the sheaf of rings of linear algebraic differential operators over
  $Y$. Let $M$ be a sheaf of modules over the sheaf $\D_Y$. We call
$M$ a $\D_Y$-module. A $\D_Y$-module $M$ is called coherent if the sheaf
$M$ can be locally presented as the cokernel of a morphism of free
$\D_Y$-modules of finite rank. We denote by $\MD_Y$ the category
of coherent   $\D_Y$-modules. If $Y$ is affine, then the space
$\overline{M}=\G(Y,M)$ of global sections of a coherent $\D_Y$-module $M$
is a module over the ring $\DD_M=\G(Y,\D_Y)$ of global algebraic differential
operators over $Y$, and the sheaf $M$ is determined by localizations of
$\overline M$,
that is, for any open affine $U\subset Y$, we have
$$\G(U,M)=\DD_U\ot_{\DD_Y}\overline{M}\ .$$

Let $U$ be an affine open subset in  an algebraic variety $Y$ and
 $T^*U$ the cotangent bundle of $U$.
A coherent $\DD_U$-module $M_U$ has singular support
$ss.M\subset T^*U$,
which is defined as follows.
Denote by $\DD^k_U\subset \D_U$ the subspace of differential operators of
order not greater than $k$.  The spaces $\DD^k_U$ define an increasing
filtration of the ring $\DD_U$ such that  the corresponding graded ring
 $Gr\DD_U$ is a commutative ring isomorphic to the polynomial ring
$\C[T^*U]$.
An increasing filtration $\{M^k_U\}$ of the module $M_U$ is called {\it
 good} if it agrees with the filtration  $\{\DD^k_U\}$ of the ring $\DD_U$:
\begin{enumerate}
\item[$\bullet$]
 $\DD^k_UM^l_U \subset M^{k+l}_U$ for all $k, l$;
\item[$\bullet$]
 $\DD^k_UM^l_U = M^{k+l}_U$ for all $l\geq  l_0, k\geq 0$;
\item[$\bullet$] All $M^k_U$ are finitely generated ${\cal O}_U$-modules.
\end{enumerate}
{}For a good filtration, the graded space $gr\,M_ U$ is a coherent
module over $Gr \DD_U$ and thus defines an algebraic variety
 $ss.M_U\subset T^*U$. More precisely, the variety is defined by
 the radical of the annihilator ideal $Ann(gr M_U)$.
The variety is called {\it the singular support or characteristic
variety} of the $\DD_U$-module $M_U$.

Bernstein's theorem \cite{Be} says that a coherent $\D_U$-module admits a good
filtration and the singular support of $M_U$ does not depend on the choice of
the good filtration. This implies that
 the characteristic variety for a coherent $\D_Y$-module $M_Y$ on a
nonsingular  variety $Y$ is well defined.
The characteristic variety is an involutive subvariety of $T^*Y$
with respect to the canonical symplectic structure on $T^*Y$, thus its
dimension is not less  than the dimension of $Y$. A $\D_Y$-module $M_Y$
is called {\it holonomic} if its characteristic variety is lagrangian.
Denote by $\MD^{hol}_Y\subset \MD_Y$ the full subcategory of sheaves
 of holonomic $\D_Y$-modules.

Let ${\cal O}_Y$ denote the structure sheaf of regular
functions on $Y$,  ${\cal O}_Y^{an}$ the sheaf of holomorphic
functions on $Y$;  $\Omega^k_Y$ the sheaf of $k$-forms on $Y$;
$\Omega^{k,an}_Y$ the sheaf of holomorphic $k$-forms on $Y$.

The simplest examples of holonomic $\D_Y$-modules are
provided by
flat algebraic connections in algebraic vector bundles over $Y$.
Recall that a connection $\nabla$ in a locally free ${\cal O}_Y$-module
${\cal S}$ of local sections of a vector bundle $S$ is a $\CC$-linear
map
$$
\nabla : {\cal S} \to \Omega_{Y}^1\otimes_{{\cal O}_Y}{\cal S} ,
$$
satisfying the condition $\nabla(fs) = df\otimes s + f\nabla(s)$ for
 $f\in{\cal O}_Y$ and  $s\in {\cal S}$.
 The  curvature $\nabla^{(2)}$ of the connection $\nabla$,
$$\nabla^{(2)}:\Omega_Y^1\otimes_{{\cal O}_Y}{\cal S}\to
\Omega_Y^2\otimes_{{\cal O}_Y}{\cal S},$$
is determined by  conditions  $\nabla^{(2)}(\omega\otimes s)=d\omega\otimes s-
\omega\wedge\nabla(s)$ for any $\omega\in{\Omega}_Y^1$ and  $s\in {\cal S}$.
The connection is flat if its curvature is zero.

A flat connection
$(\nabla,{\cal S})$ can be canonically extended to a holomorphic flat
connection, $(\nabla^{an},{\cal O}_Y^{an}\ot_{{\cal O}_Y}{\cal S})$.

A flat connection
$(\nabla,{\cal S})$ provides the ${\cal O}_Y$-module ${\cal S}$ with a structure of
a $\D_Y$-module ${\cal S}_\nabla$, where the action of a vector field $\xi$ is
given by the relation $\xi(s) = i_{\xi}(\nabla(s))$.


Let $Y$ be a stratified  algebraic variety with  smooth algebraic strata
$Y_j$, $j\in J(Y)$. It means that $Y$ is presented as a finite disjoint
 union of strata $Y_j$, such that the closure of each stratum is a union of
strata. Denote by $\MMD^{hol}_Y\subset\MD_Y$ the full subcategory
of holonomic $\D_Y$-modules $M$ whose singular support
is contained in the  union of the conormal bundles of strata of $Y$,
$ ss.M \subset \bigcup_{j\in J(Y)} T^*_{Y_j}Y$. This category is
abelian and is closed under extensions. If $M \in
\MMD^{hol}_Y$,
then we say that $M$ is {\it smooth along the stratification of $Y$}.

\bigskip

\subsection{Quiver  $\D_X$-modules}\label{Sec.3.2}
Let  $\A$ be an arrangement in $X=\C^N$.
Let  $\{f_{\b,\a}, \xi_{\a,\b}\ |\ \a,\b\in \v(\Ga),\a\succ\b\}$
be an edge framing of $\A$.

Let $\VA$ be a quiver of the graph $\G$ of  $\A$.
We define {\em the associated  quiver    $\D_X$-module $E\VA$} as the quotient
of the free  sheaf $\R^0 E\VA$ $=$
$\D_X\otimes_\C(\oplus_{\a\in \v(\G)} \Oa\ot V_\a)$ of  $\D_X$-modules
over its subsheaf ${\R'} E\VA$ of coherent  $\D_X$-modules,
 whose global sections are the following global sections of  $\R^0 E\VA$:
\begin{enumerate}
\item[(i)] the sections
 \begin{equation}
u\,\xi \ot \om_\a\ot v_\a \ -\ \sum_{\b,\a\succ\b}
\mylangle \xi,f_{\b,\a}\myrangle \,u\ot\pi_{\b,\a}(\om_\a) \ot
A_{\b,\a}(v_\a)\ ,
 \label{2.3}
\end{equation}
where  $\a\in\v(\Ga)$,  $u\ot \om_\a\ot v_\a\in \DD_X\otimes_\C\Oa\ot V_\a$
and  $\xi\in\Ta$,
\item[(ii)] the sections
\begin{equation}
u\, f  \ot \om_\a\ot v_\a\ -\ \sum_{\b,\b\succ\a}
 \mylangle \xi_{\b,\a},f\myrangle\, u\ot
\pi_{\b,\a}(\om_\a)\ot
 A_{\b, \a}(v_\a)\ ,
\label{2.4}
\end{equation}
where $\a\in\v(\Ga)$,  $u \ot\om_\a\ot v_\a\in \DD_X\ot_\C\Oa\ot V_\a$ and
$f\in\Fa$.
\end{enumerate}

It is easy to see that the space generated by sections in \rf{2.3}--\rf{2.4} does
not depend on
the choice of the edge framing. Thus the sheaf $E\VA$ is determined by
the quiver $\VA$ only.

 For  any $\a\in\v(\Ga)$,  $\om_\a\in \Oa$,
$ v_\a\in  V_\a$,\ we denote by $\om_\a\ot v_\a$ the image in $\G(X,E\VA)$
of the global section $1\ot\om_\a\ot v_\a$ $\in$ $\G(X,\R^0E\VA)$ with
respect to
the canonical projection $\R^0E\VA\to  E\VA$.
 The vectors  $\om_\a\ot v_\a$, where  $\a\in\v(\Ga)$,  $\om_\a\in \Oa$,
$ v_\a\in  V_\a$, generate $\G(X,E\VA)$ as $\DD_X$-module.

\bigskip

\noindent
{\bf Remark.} In  \cite{Kh1, KhS}
relations \rf{2.3}--\rf{2.4} were written in a different
but equivalent form.  Namely, let \OOa be the one-dimensional
complex vector space  of $l(\a)$-exterior forms $df_1\wedge
df_2\wedge...\wedge df_{l(\a)}$ with  $f_1, \dots , f_{l(\a)}\in \Fa$.
For any $\a\in \v(\Ga)$,  choose an element $\ho_\a\in\OOa$.
Then the relations of \rf{2.3}--\rf{2.4} are equivalent to the next
relations,
 \begin{equation}
\begin{split}
 \xi ( \om_\a\ot v_\a )\ =\ \sum_{\b,\a\succ\b}
\frac{i_\xi(\ho_\b) \wedge \om_\a}{\ho_\b\wedge\om_\b}\ \om_\b\ot A_{\b,
\a}(v_\a) ,
\label{2.1}
\\
 f ( \om_\a\ot v_\a )\ =\ \sum_{\b,\b\succ\a}
\frac{df\wedge\ho_\b \wedge \om_\a}{\ho_\b\wedge \om_\b}\ \om_\b \ot A_{\b,
\a}(v_\a) .
\end{split}
\end{equation}
The right hand side of formulas
\rf{2.1} does not depend on the choice of the
forms $\ho_\a\in\OOa$.

\bigskip

It is easy to see that a morphism of quivers $\phi : \VA \to \VA'$
induces a natural
morphism $E\phi: \EVA \to \EVA'$  of $\D_X$-modules.
Thus the correspondence $E : \VA \mapsto \EVA$, $E: \phi \mapsto E\phi$
 defines a functor
from the category $\Qui_\Ga$ to the category of $\D_X$-modules.

Choose a vector space structure and a non-degenerate bilinear form $<,>$
on $\C^N$. For $\a\in\v(\Ga)$, introduce the spaces \TTa and \FFa as
in Section \ref{Framing}.
Let $S\left(\TTa\oplus\FFa\right)$ be the symmetric tensor algebra over
the space $\TTa\oplus\FFa$. Denote by
$M_\a=M_\a(\VA)$ the vector space
 \begin{equation}
 \label{2.5}
M_\a\ =\ S\left(\TTa\oplus\FFa\right)\otimes \Oa\ot V_\a \ .
\end{equation}
Clearly, there is a natural linear map of the vector space $M_\a$
to the vector space of global section of $\EVA$.

\medskip

\begin{proposition}\label{lemma1}
The vector space of global sections of the $\D_X$-module $\EVA$
is isomorphic to the direct sum of the spaces $M_\a$,
\begin{equation}\label{2.6}
\Gamma(\EVA)\ \simeq\ \oplus_{\a\in\v(\Ga)} M_\a \ .
\end{equation}
\end{proposition}

The proposition is proved in Section \ref{proofs}.

{}For $k=0 ,...,N$,
denote by $\V_k$ the following finite dimensional subspace of
global sections of  the $\D_{X}$-module $E \VA$:
\begin{equation}\label{vk}
\V_k\ =\ \ooplus\limits_{\a\in\v_k(\Ga)}\, 
\Oa\ot V_\a\ .
\end{equation}
For any $m\geq 0$ and $l, \ 0\leq l\leq N$, \ set
$$
 F_p^mE\VA\ =\ \sum\limits_{k=0}^N\ \D_X^{m- k}\V_k ,
\qquad
F_p^{m,l}E\VA\ =\ \sum\limits_{k=l}^N\ \D_X^{m- k}\V_k
$$
where $\D_X^l$ is the sheaf of differential operators over $X$ of order
less than or equal to $l$. We have
$$
\D^l_X  F_p^mE\VA\subset  F_p^{l+m}E\VA
\qquad
{\rm and}
\qquad
 F_p^{m_1}E\VA\subset
 F_p^{m_2}E\VA
$$
 for $m_1<m_2$, so that
the graded factor $\overline{E}\VA=\oplus_{m\geq 0}
 F_p^mE\VA/ F_p^{m-1}E\VA$ is a module over ${\rm gr}\D_X\simeq \CC[T^*X]$.
We call the
filtration
\begin{equation}\label{pfil}
0= F_p^{-1}E\VA\subset  F_p^0E\VA\subset  F_p^1E\VA\subset...\subset
 F_p^mE\VA\subset... \subset E\VA
\end{equation}
{\it the principal filtration} of
the quiver $\D_X$-module $E\VA$. The sheaves $F_p^{m,l}E\VA$ form a
subfiltration
\begin{equation}\label{subtle}
...\subset F_p^{m-1}E\VA=F_p^{m-1,0}E\VA\subset F_p^{m,N}E\VA
\subset...\subset
F_p^{m,0}E\VA=F_p^{m}E\VA\subset...
\end{equation}
of the principal filtration of the quiver $\D_X$-module $E\VA$.
We call it  {\it the subtle filtration} of the quiver $\D_X$-module $E\VA$.
\medskip

\begin{proposition}
\label{prop1}
${}$

\noindent
 {\em (i)} The principal filtration \rf{pfil} of the quiver $\D_X$-module
 $E\VA$ is good.
 \\
 {\em (ii)} The associated graded factor $\overline{E}\VA$ $=$ ${\rm gr}\ {E}\VA$ of the sheaf
$E\VA$ with respect to the principal filtration\ $F_pE\VA$
 \ is isomorphic to
 $$
\ooplus\limits_{\a\in\v(\Ga)}{\mathcal O}_{\X_\a}\ot\left(S(\TTa)\otimes
\Oa\ot V_\a\right)
$$
 as the sheaf of ${\mathcal O}_X$-modules.\\
  {\em (iii)}
The $\D_X$-module $E\VA$ is  holonomic and smooth along the stratification
 $\{X_\a\}$.
\end{proposition}

The proposition is proved in Section \ref{proofs}.

\medskip

\subsection{Example}${}$

\medskip
\noindent
{(i)} Let $X$ be $\CC$ with coordinate $z$. Let $\A=\{H_1,...,H_n\}$,
where $H_j$ $=$
$\{z-t_j=0\}$. Denote the vertices of the graph of the arrangement by the symbols
$\ee$, and $ j$,\ $j=1,...,n$. Let
$\VA \ =\ \{V_\ee,\ V_j,\ A_{\ee,j},\ A_{j,\ee}\}$
\ be a quiver
of  $\G$.

The vector space $\G(X,\EVA)$ of global sections of the sheaf $\EVA$ is generated
over $\DD_X$ by the vectors
$ {\rm d}z \ot v_\ee\in\overline{\Omega}_{\ee}\ot  V_\ee$, and
$1\ot v_j\in\overline{\Omega}_{j}\ot V_j$, subject to the relations:
$$
\partial_z( {\rm d}z\ot v_\ee)=\sum_{j=1}^n1\ot A_{j,\ee}(v_\ee),
\qquad
(z-t_j)(1\ot v_j)= {\rm d}z\ot A_{\ee,j}(v_j)\  ,
$$
where $\partial_z = \frac{d}{dz}$.

The vectors
\begin{equation}
\begin{split}\nonumber
z^k ({\rm d}z \ot v_\ee) \ ,
\qquad
k\geq 0 , \ \ v_\ee\in  V_\ee\ ,
\\
\partial_z^l (1\ot v_j)\ ,
\qquad
l\geq 0 ,\ \ v_j\in V_j \ ,
\end{split}
\end{equation}
form a $\C$-basis of the vector space $\G(X,\EVA)$.

Let $\{F_p^lE\VA\}$  be the principal filtration of the
$\D_X$-module $\EVA$.

The vectors $z^k ({\rm d}z\ot v_\ee)$,  $k\geq 0$, form a $\C$-basis of
the space of  global sections of the sheaf $F_p^0E\VA$.

 For any $l>0$, the vectors $\partial_z^m(1\ot  v_j)$, $0\leq m\leq l$,\
 and \ $z^k({\rm d}z\ot v_\ee)$,  $k\geq 0$, \
form a $\C$-basis of
the space of  global sections of the sheaf $F_p^lE\VA$.

\bigskip
\noindent
{(ii)}
Let $X$ be $\C^2$ with coordinates $z_1,z_2$.
Let $\A = \{H_1,..., H_n\}$ be the
arrangement of $n$ lines, $H_j = \{z_2+k_jz_1=0\}$, $j = 1, \dots , n$,
 intersecting at one point.  Denote the vertices of the graph of
 the arrangement by the symbols
$\ee$ , $ \a_j$, where $j=1,...,n$, and $\b$, such that $X_\b$ is the point
$(0,0)$.

Let $\VA$ $=$ $\{V_\ee, V_{\a_j}, V_\b, A_{\ee,\a_j},
A_{\a_j,\ee},  A_{\b,\a_j},  A_{\a_j,\b}\}$  be a quiver.

The vector space $\G(X,\EVA)$ of global sections of the sheaf $\EVA$ is generated
over $\DD_X$ by the vectors
${\rm d} z_1\wedge {\rm d}z_2\ot v_\ee$, where  $v_\ee\in V_\ee$,
the vectors
${\rm d} z_1\ot v_{\a_j}$, where  $v_{\a_j}\in V_{\a_j}$, and
 the vectors
$1\ot v_\b$, where  $v_\b\in V_\b$,
subject to the following relations:
\begin{equation}
\begin{split}\nonumber
\partial_{z_1}({\rm d} z_1\wedge {\rm d}z_2\ot v_\ee)=
\!\sum_{j=1}^nk_j{\rm d} z_1\ot A_{\a_j,\ee}(v_\ee),\qquad
&\partial_{z_2}({\rm d} z_1\wedge {\rm d}z_2\ot v_\ee)=
\sum_{j=1}^n{\rm d} z_1\ot A_{\a_j,\ee}(v_\ee),\\
(z_2+k_jz_1)({\rm d} z_1\ot v_{\a_j})=
1\ot A_{\ee,\a_j}(v_{\a_j}),\qquad
&(\partial_{z_1}-k_j\partial_{z_2})({\rm d} z_1\ot v_{\a_j})=
1\ot A_{\b,\a_j}(v_{\a_j}),\\
z_1(1\ot v_\b)=\sum_{j=1}^n{\rm d} z_1\ot A_{\b,\a_j}(v_\b),\qquad
&z_2(1\ot v_\b)=-\sum_{j=1}^nk_j{\rm d} z_1\ot A_{\b,\a_j}(v_\b)
\end{split}
\end{equation}

The vectors:
\begin{equation}
\begin{split}\nonumber
\ z_1^{k_1} z_2^{k_2}({\rm d} z_1\wedge {\rm d}z_2\ot v_\ee) \ ,
&
\qquad
 k_1 , k_2  \geq 0\ ,
\\
z_1^{k_1} \partial_{z_2}^{k_2}({\rm d} z_1\ot v_{\a_j}) \ ,
&
\qquad
 k_1 , k_2 \ \geq 0 \ , \ \ j = 1,...,n \ ,
\\
\partial_{z_1}^{k_1}\partial_{z_2}^{k_2}(1\ot v_{\b})\ ,
 &\qquad
k_1 , k_2\ \geq 0 \ ,
\end{split}
\end{equation}
form a $\C$-basis of the vector space $\G(X,\EVA)$.

The vectors
 $\ z_1^{k_1} z_2^{k_2}({\rm d} z_1\wedge {\rm d}z_2\ot v_\ee)$,\
$k_1, k_2\geq 0$,\
form a $\C$-basis of the space of
global sections of the sheaf $F^0_pE\VA$.

The vectors  $\ z_1^{k_1} z_2^{k_2}({\rm d} z_1\wedge {\rm d}z_2\ot v_\ee)$, \
$k_1, k_2 \geq 0$,\
and the vectors ${\rm d} z_1\ot v_{\a_j}\in V_{\a_j}$,
\ $j = 1, \dots , n$,\
form a $\C$-basis of the space of global sections of the
 sheaf $F^1_pE\VA$.

For $l \geq 2$, the vectors
 $\ z_1^{k_1} z_2^{k_2}({\rm d} z_1\wedge {\rm d}z_2\ot v_\ee)$,
\ $k_1, k_2 \geq 0$, \
the vectors  $z_1^{k_1} \partial_{z_2}^{k_2}({\rm d} z_1\ot v_{\a_j})$,
 \ $k_1 \geq 0,$\  $0 \leq k_2 < l$, \ $j = 1, \dots , n,$\   and the vectors
$\partial_{z_1}^{k_1}\partial_{z_2}^{k_2}(1\ot v_{\b})$, \ $k_1+k_2 < l-1$,\
form a $\C$-basis of the space of global sections of the
 sheaf $F^{l}_pE\VA$.

\subsection{Quiver $\D_{X^n}$-modules}\label{section3.3}
Let $X^n$ be a principal open subset of $X$, see Section \ref{Sec.Xn}. Let
  $\{f_{\b,\a}, \xi_{\a,\b}\ |\ \a,\b\in \v(\Ga),
\a\succ\b\}$ be an edge framing of $\A$, and
 $\VA_{\Ga^n} = \{ V_\a, A_{\a,\b},\  A^\b_\a\} \in \Qui_{\G^n}$ a quiver
 of $\G^n$.

We define {\em the associated  quiver  $\D_{X^n}$-module  $E^n\VA_{\G^n}$}
as the quotient
of the free  sheaf $\R^0 E^n\VA_{\G^n}$ $=$
$\D_{X^n}\otimes_\C(\oplus_{\a\in \v(\G^n)} \Oa\ot V_\a)$ of
 $\D_{X^n}$-modules
over its subsheaf $\R' E^n\VA_{\G^n}$ of coherent  $\D_{X^n}$-modules,
which is described  by its restrictions
to open affine subsets  $U_{\FF_n}$, where  
${\mathcal F}_{n}=\{\ff{\a}\ |\ \a\in\v_{n+1}(\Ga)\}$
 is a level $n$ vertex framing.

The  sections $\G(U_{\FF_n}, \R' E^n\VA_{\G^n})$
of the sheaf $\R' E^n\VA_{\G^n}$ over  $U_{\FF_n}$
are the following  sections in  $\G(U_{\FF_n}, \R^0 E^n\VA_{\G^n})$:

\begin{enumerate}
\item[(i)] the sections
 \begin{equation}
u\,\xi \ot   \om_\a\ot v_\a \ -\ \sum_{\b,\a\succ\b}
\mylangle \xi,f_{\b,\a}\myrangle \,u\ot \pi_{\b,\a}(\om_\a) \ot
A_{\b,\a}(v_\a)\ ,
  \label{3.7}
\end{equation}
where  $\a\in\v(\Ga^n)$ with $l(\a)<n$,
$u\ot \om_\a\ot v_\a\in \DD_{U_{\FF_n}}\otimes_\C\Oa\ot V_\a$, and
 $\xi\in\Ta$,
\item[(ii)] the sections
 \begin{equation}
\begin{split}
 u\,\xi \ot \om_\a\ot  v_\a-&\!\sum_{\b,\ \a\succ\b}
u\,\ff{\b}^{-1}\ot \Big( \mylangle \xi_{},\ff{\b}\myrangle
\om_\a\ot A_{ \a}^\b
(v_\a)+
\\
&
+\sum\limits_{\g,\ \g\succ\b,\g\not=\a}
\mylangle \xi_{},f_{\b,\a}\myrangle\mylangle \xi_{\g,\b},\ff{\b}\myrangle
\pi_{\g,\b}\pi_{\b,\a}(\om_\a)\ot A_{\g,\d(\a,\g)}A_{\d(\a,\g),\a}(v_\a)
 \Big)
 \label{3.9a}
\end{split}
\end{equation}
where $\a\in\v(\Ga^n)$ with $l(\a)=n$,
 $u\ot \om_\a\ot v_\a\in \DD_{U_{\FF_n}}\otimes_\C\Oa\ot V_\a$,
and  $\xi\in\Ta$.

Here $\d(\a,\g)$ is the unique (if exists) element of
$\v(\Ga^n)$, such that $\d(\a,\g)\succ \a$ and $\d(\a,\g)\succ \g$.
If it does not exist, the corresponding summand in
\rf{3.9a} is set to be zero.
\item[(iii)] the sections
\begin{equation}
u\, f\ot\om_\a\ot v_\a -\ \sum_{\b,\b\succ\a}
 \mylangle \xi_{\b,\a},f\myrangle u\ot
\pi_{\b,\a}(\om_\a)\ot
 A_{\b, \a}(v_\a)\ ,
\label{3.9}
\end{equation}
\end{enumerate}
 where  $\a\in\v(\Ga^n)$,
$u\ot \om_\a\ot v_\a\in \DD_{U_{\FF_n}}\otimes_\C\Oa\ot V_\a$ ,
and $f\in\Fa$.

One can see that the space generated by sections in \rf{3.7}--\rf{3.9}
does not depend on
the choice of the edge framing of $\A$.
 Thus the sheaf $E^n\VA_{\G^n}$ is determined by
the quiver $\VA_{\G^n}$ only.

 For  any $\a\in\v(\Ga^n)$,  $\om_\a\in \Oa$,
$ v_\a\in  V_\a$,\ we denote by $\om_\a\ot v_\a$ the image in
 $\G(U_{\FF_n},E^n\VA_{\G^n})$
of the  sections $1\ot\om_\a\ot v_\a$ $\in$ $\G(U_{\FF_n},
\R^0E^n\VA_{\G^n})$ with
respect to
the canonical projection $\R^0E^n\VA_{\G^n}\to E^n\VA_{\G^n} $.

 The vectors  $\om_\a\ot v_\a$, where  $\a\in\v(\Ga^n)$,  $\om_\a\in \Oa$,
$ v_\a\in  V_\a$,  generate $\G(U_{\FF_n},
E^n\VA_{\G^n})$ as a $\DD_{U_{\FF_n}}$-module.

All   statements of the previous section have straightforward analogs for
the quiver $\D_{X^n}$-modules on $X^n$. In particular, the principal
filtration $F_pE^n\VA_{\G^n}$ of a quiver  $\D_{X^n}$-module $E^n\VA_{\G^n}$
is defined by the rule
$$ \G(U_{\FF_n},F_p^mE^n\VA_{\G^n})=
\sum\limits_{k=0}^n\DD_{U_{\FF_n}}^{m- k}\V_k,$$
where $\V_k$ is defined in \rf{vk}.
The graded factor $\overline{E}^n\VA_{\G^n}$ $=$ ${\rm gr}\ {E}^n\VA_{\G^n}$
of the sheaf
$E^n\VA_{\G^n}$ with respect to the principal filtration\ $F_pE^n\VA_{\G^n}$
 \ is isomorphic to
$$\ooplus\limits_{\a\in\v(\Ga^n)} {\mathcal O}_{\X_\a\cap X^n}\ot
\left(S(\TTa)\otimes \Oa\ot V_\a\right)\ $$
 as a sheaf of $\ {\mathcal O}_{X^n}$-modules.
This implies that  the  space of global sections $\G(X^n,E^n\VA_{\G^n})$ of
the $\D_{X^n}$-module $E^n\VA_{\G^n}$ is
isomorphic to
 $$\ooplus\limits_{\a\in\v(\Ga^n)} \CC[\X_\a\cap X^n]\ot
\left(S(\TTa)\otimes
\Oa\ot V_\a\right) \ .$$

{}For any $n=0,1,...,N$,
consider the category  $\QGn$
 associated with given collections $\{G_\a\}$, $\{\widetilde{G}_\b\}$, and
 $\{G_{\a,\b}\}$.

\medskip

\begin{proposition}\label{theorem1}${}$

\noindent {\em (i)} The assignment $\VA_{\G^n}\mapsto E^n\VA_{\G^n}$
defines a faithful functor $E^n: Qui_{\G^n}\to \MD_{X^n}^{hol}$.

\noindent {\em (ii)} Let $\A$ be a central arrangement (that is, all
 hyperplanes of $\A$ contain $0\in\CC^N$). Then the restriction of
 $E^n$ to the subcategory $\QGn$ is full and faithful.

  In particular, the restriction of the functor $E$ to the subcategory
\QG  defines
 an equivalence of the category \QG\
with a full subcategory of the category   $\MD_X^{hol}$, closed under
extensions.
\end{proposition}


Proposition \ref{theorem1} is proved in Section \ref{proofs} by induction on $n$.
The following proposition plays the key role in the inductive proof.

\begin{proposition}\label{theorem1a}
Let $\A$ be a central arrangement. Let $\{G_j\ |\  j\in J(\A)\} \cup
\{\widetilde{G}_\ee\}$ be a collection of weakly non-resonant sets.
Let $\VA=\{V,A^j\}$ and ${\mathcal W}=\{W,B^j\}$ be level zero quivers, such
that
\begin{itemize}
\item[(a)] for any $j\in J(\A)$ the eigenvalues of $A^j$ and eigenvalues of
$B^j$ are in $G_j$,
\item[(b)] the eigenvalues of $\sum_{j\in J(\A)}A^j$ and
 the eigenvalues of $\sum_{j\in J(\A)}B^j$ are in $\widehat{G}_\ee$.
\end{itemize}
Then  the vector spaces
${\rm Hom}_{\Qui_{\G^0}} (\VA, {\mathcal W})$ and
$\,{\rm Hom}_{\D_{X^0}} (E^0 \VA, E^0 {\mathcal W})$
are isomorphic.
\end{proposition}
\medskip
\noindent
{\bf Remark.}  Proposition \ref{theorem1} is a generalization of the
corresponding statement in \cite{Kh1}. Note that the condition of the
centrality of the arrangement is missing in \cite{Kh1}. This is the
mistake, which is in the ignorance of Proposition \ref{theorem1a} in
the inductive proof of the main theorem of \cite{Kh1}.  The statement of the main theorem of \cite{Kh1}
is valid for central arrangements.
\medskip

 Proposition \ref{theorem1a} is proved in Section \ref{proofs}.

\medskip
{\bf Example.}
 Let $X^0$ be the complement
to the union of hyperplanes of the arrangement $\A = \{H_j\}, j\in J(\A)$,\ in
$X=\CC^N$. Let $z_1,..., z_N$ be affine coordinates in $X$.
For $j\in J(\A)$,
let $f_{j}$ be an affine function on $\C^N$ defining the hyperplane $H_j$.
Let $\VA_{\G^0}=\{V,A^j\}\in\Qui_{\G^0}$ be a quiver of level 0.
Recall that $A^j : V \to V$ are linear maps
 such that
for each $\a\in\v_2(\G)$, $j\succ\a$, we have
\begin{equation}
\label{flat}
[ A^j , \sum_{i,\ \, i\succ\a} A^{i}]\ =\  0\ .
\end{equation}
The associated $\D_{X_0}$-module   $E^0\VA_{\Ga^0}$
 is generated by the vectors $\o_\ee\ot v$, where  $\o_\ee = {\rm d}z_1\wedge
... \wedge {\rm d}z_N$ and  $v\in V$,
 subject to the relations, labeled by vectors $v\in V$ and constant
vector fields $\xi$ on $\C^N$,
\begin{equation}\label{3.133}
\xi(\o_\ee\ot v)\ =\ \sum_{j} \frac{\mylangle \xi, f_{j}\myrangle}{f_{j}}
\left(\o_\ee\ot \ A^j(v)\right)\ .
\end{equation}
That means that $E^0\VA_{\Ga^0}$ is isomorphic to
the $\D_{X^0}$-module associated to the following flat connection.

Consider the trivial bundle on $X^0$ with fiber $V$, the
End $V$-valued differential one-form
\begin{equation}
 \sum_{j\in J(\A)}\ A^j \ot d\log f_{j}\ ,
\label{3.12}
\end{equation}
and the associated connection on
the trivial bundle. The connection is flat due to \rf{flat}.
The connection is regular singular. The $\D_{X^0}$-module associated with
the connection is isomorphic to $E^0\VA_{\Ga^0}$. As a sheaf of
${\mathcal O}_{X^0}$-modules, it is isomorphic to the free sheaf
${\mathcal O}_{X^0}\ot V$.
\medskip

\subsection{Direct and inverse images of quiver $\D$-modules}
For $0\leq k\leq l\leq N$, let $X^k, \ X^l$ be the principal open subsets of
$X$. Since the canonical embedding  $j_{l,k}:X^k\to X^l$ is open,
the functors
$j_{l,k}^*$: $\MD_{X^l}\to \MD_{X^k}$  and
$j_{l,k,*}: \MD_{X^k}\to \MD_{X^l}$
of the inverse and direct images of $\D$-modules are well defined and can be
described in the sheaf-theoretical language. Namely, for $M\in \MD_{X^l}$,
the sheaf $j_{l,k}^*(M)$ is the sheaf-theoretical inverse image
$j_{l,k}^\mycircle (M)$ of the sheaf $M$ equipped with the natural structure of the
 sheaf of $\D_{X^k}$-modules. Also for $N\in \MD_{X^k}$,
the sheaf $j_{l,k,*}(N)$ is the sheaf-theoretical direct image
$(j_{l,k})_\mycircle (N)$ of the sheaf $N$ equipped with the natural structure of
the sheaf of $\D_{X^l}$-modules.

More precisely, let $M$ be a coherent $\D_{X^l}$-module and $U\subset X^k$
an open subset. Then $U$ is an open subset of $X^l$ as well and
\begin{equation}\label{gjM}\G(U,j_{l,k}^*(M))= \G(U, M)\ .
\end{equation}
Since the space $\G(U, M)$ is a $\DD_U$-module, the relation \rf{gjM} describes
the structure of  a $\DD_U$-module on the space $\G(U,j_{l,k}^*(M))$.

Let now $N$  be a coherent $\D_{X^k}$-module and $U\subset X^l$
an open subset. Then $j_{l,k}^{-1}(U)= U\cap X^k\subset U$
is an open subset of $X^k$ and
\begin{equation}\label{gjN}
\G(U,j_{l,k,*}(N))= \G(j_{l,k}^{-1}(U), N)\ .
\end{equation}
The structure of a $\D_X$-module on the sheaf $j_{l,k,*}(N)$ is described
as follows.

Let $n \in \G(U,j_{l,k,*}(N))$
be a section
corresponding to the section $\widetilde{n}\in \G(j_{l,k}^{-1}(U), N)$ via
identification \rf{gjN}. Let
$d\in \DD_{U}$ be a differential operator on $U$. Then the product
$dn\in  \G(U,j_{l,k,*}(N)), N)$ corresponds via \rf{gjN} to the section
$\widetilde{d}\ \widetilde{n}\in \G(j_{l,k}^{-1}(U), N)$, where
$\widetilde{d}=d|_{U\cap X^k}$.

\medskip

The functor $j_{l,k}^*$ is exact and maps the subcategory
$\MD_{X^l}^{hol}\subset
\MD_{X^l}$ to the subcategory $\MD_{X^k}^{hol}\subset \MD_{X^k}$.
 {}For any $0\leq k_1<k_2<k_3\leq N$, we
have
$$
 j_{k_3,k_1}^*\ =\ j_{k_2,k_1}^*\cdot j_{k_3,k_2}^*\ .
$$
We use the notation $j_0^*$ for the functor\
$j_{N,0}^*\ :\ \MD_{X}=\MD_{X^N}\ \to\ \MD_{X^0}$.

 The functor $j_{l,k,*}: \MD_{X^k}\to \MD_{X^l}$
is right adjoint to the functor of inverse image
$j_{l,k}^*: \MD_{X^l}\to \MD_{X^k}$
and for any $0\leq k_1<k_2<k_3\leq N$:
$$
 j_{k_3,k_1,*}= j_{k_3,k_2,*}\cdot j_{k_2,k_1,*}\ .
$$
We use the notation $j_{0,*}$ for the functor\
$j_{N,0,*}\ :\  \MD_{X^0}\ \to\ \MD_{X^N}=\MD_{X}$.
It is exact since the map $j_0: X^0\to X$ is affine.

The following theorem is an immediate consequence of the constructions
 in the previous section.

\begin{theorem}
\label{proposition2} The  diagram below is
commutative\ :
$$
\begin{array}{ccc}
\Qui_{\G^l}&\stackrel{\j_{l,k}^*}{\longrightarrow}
& \Qui_{\G^k}\\
\downarrow\lefteqn{E^l}&&\downarrow\lefteqn{E^k}\\
\MD_{X^l}^{hol}&\stackrel{j_{l,k}^*}{\longrightarrow}&\MD_{X^k}^{hol}\ \ .
\end{array}
$$
\hfill $\square$
\end{theorem}
\medskip

 Let $0\leq k<l\leq N$. For any $\Ga^k$-quiver $\VA_{\G^k}$ and
$\Ga^l$-quiver ${\mathcal U}_{\G^l}$ denote by
$\a_{{\mathcal U}_{\G^l}, \VA_{\G^k}}:$
 $\Hom_{\Qui_{\G^l}}({\mathcal U}_{\G^l},\j_{l,k,*}(\VA_{\G^k}))$ $\to$
$\Hom_{\Qui_{\G^k}}(\j_{l,k}^*({\mathcal U}_{\G^l}),\VA_{\G^k})$
  the isomorphism of adjunction of the functors
$\j_{l,k}^*$ and $\j_{l,k,*}$. For any $\D_{X^k}$-module $M_k$ and
any $\D_{X^l}$-module $M_l$ denote by  $\tilde{\a}_{M_l,M_k}:$
$\Hom_{{\mathcal M}_{X^l}}(M_l,j_{l,k,*}(M_k))$ $\to$
$\Hom_{{\mathcal M}_{X^k}}(j_{l,k}^{*}(M_l),M_k)$
 the isomorphism of adjunction of the functors
$j_{l,k}^*$ and $j_{l,k,*}$.

\begin{theorem}
\label{prop4} 
Let $\VA_{\G^k}$ be a $\G^k$ quiver
such that its quiver direct image $\j_{l,k,*}(\VA_{\G^k})$ is strongly
 non-resonant.

\begin{enumerate}
\item[(i)]  Then the direct image of the $\D_{X^k}$-module $E^k(\VA_{\G^k})$ with respect to
$j_{l,k}$ is isomorphic to the quiver
$\D_{X^l}$-module, associated with
 $\j_{l,k,*}(\VA_{\G^k})$\ :
\begin{equation}\label{jE}
 j_{l,k,*}(E^k (\VA_{\G^k}))\ \simeq\
 E^l (\j_{l,k,*}(\VA_{\G^k}))\ .
\end{equation}
\item[(ii)] Identify $\D_{X^l}$-modules $j_{l,k,*}(E^k (\VA_{\G^k}))$ and
$ E^l (\j_{l,k,*}(\VA_{\G^k}))$ by isomorphism \rf{jE}. Then for any
$\G^l$-quiver ${\mathcal U}_{\G^l}$ and any morphism $\phi:$
${\mathcal U}_{\G^l}$ $\to$ $\j_{l,k,*}(\VA_{\G^k})$ of $\G^l$-quivers
we have
$E^k\a_{{\mathcal U}_{\G^l}, \VA_{\G^k}}=
\tilde{\a}_{E^l({\mathcal U}_{\G^l}), E^k(\VA_{\G^k})}(E^l\phi)$.
\end{enumerate}
\end{theorem}

The theorem is proved  in Section \ref{proofs}.

\medskip

\begin{corollary}\label{cor2}
Let $\VA_{\G^0} = \{V, A^i\}$ be a level zero quiver.
Suppose that each operator $A^i$ has a single eigenvalue.
Suppose that $\VA_{\G^0}$ has a non-resonant spectrum, see Section
\ref{Spectrum}. Then
$$
 j_{0,*} (E^0 (\VA_{\G^0}))\ \simeq\
 E (\j_{0,*}(\VA_{\G^0}))\ .
$$
\end{corollary}
\medskip

\subsection{Specialization of a quiver $\D_X$-module to a stratum}
\label{section4.6}
There is a {\it specialization} construction in the theory of $\D$-modules due to
Kashi\-wa\-ra \cite{Kash2}. This construction gives
a collection of functors between quiver $\D$-modules.

Let $Y$ be a smooth complex algebraic variety and $Z\subset Y$ a
smooth complex algebraic subvariety. Let $T_ZY \to Z$ be the normal bundle
of the subvariety $Z$.
Let $I\subset {\cal O}_Y$ be the sheaf of ideals
of functions vanishing on  $Z$.

Define a decreasing filtration $F(\D_{Y})$ of the sheaf $\D_Y$ by the
formula
\begin{equation}\label{FkD}
F^k(\D_{Y})\ =\ \{P\in\D_{Y}|\ P(I_\b^j)\subset I_\b^{j+k}\quad
\mbox{\rm for any }\ j\}\ .
\end{equation}
The associated graded quotient is isomorphic to the sheaf of differential
operators $\D_{T_ZY}$ on the total space $T_ZY$ of the normal bundle to $Z$.

Let $M$ be a $\D_{Y}$-module.  A decreasing filtration $F(M)$ of
$M$ is called {\it good} with respect to $F(\D_{Y})$, if
\begin{enumerate}
\item[(i)] $F^k(\D_{Y})F^j(M)\subset F^{k+j}(M)\ $
for any $k$ and $j$,
\item[(ii)]
$F^k(\D_{Y})F^j(M)= F^{k+j}(M)\ $
if $j\gg 0$ and $k\geq 0\ $ or  $\ j\ll 0$ and $k\leq 0$,
\item[(iii)] $ F^{j}(M)\ $  is a coherent $F^0(\D_{Y})$-module,
\item[(iv)]
$M=\cup_jF^j(M)$.
\end{enumerate}

Let $\theta$ be a vector field on $Y$ tangent to $Z$ and acting on
$I/I^2$ as the identity. Let $\overline{G}\subset\CC$ be a {\it
non-resonant section} of $\CC$, that is $\overline{G}$ contains zero
and for any $a\in\CC$ the intersection $\overline{G}\cap(a+\ZZ)$
consists of a single point.

Kashiwara's theorem \cite{Kash2} states  that for any
holonomic $\D_Y$-module $M$ with regular singularities there exists a
unique filtration of $M$ which is good with respect to $F(\D_{Y})$ and
such that
\begin{enumerate}
\item[(v)] there exists a polynomial $b(x)$ with zeros in $\CG$ such that
$$
b(\theta-k)\,F^k(M)\ \subset\ F^{k+1}(M)
$$
for any $k$.
\end{enumerate}

The associated graded quotient $\gr M$ has a structure of
a $\D_{T_ZY}$-module and is called {\it the specialization} of $M$ to $Z$.
We denote the specialization by $Sp_Z(M)$.

Thus, having $Y, Z, \theta, M, \overline{G}$, Kashiwara's theorem gives the
specialization $\D_{T_ZY}$-module $Sp_Z(M)$.

\medskip

In this section, we consider the case in which
\begin{enumerate}
\item[$\bullet$]
$Y$ is a vector space $X=\CC^N$ with a central arrangement $\A$,
\item[$\bullet$]
$Z$ is the closure $\X_\a$ of a stratum  of that arrangement,
so $Z$ is a vector subspace,
\item[$\bullet$]
$\theta
=  z_1\frac{\partial}{\partial{z_1}} + \ldots +
z_k\frac{\partial}{\partial{z_k}}$,
where $z_1,...,z_N$ are linear coordinates on
$X$, such that $Z$ is defined by equations $z_1=\cdots = z_k=0$,
\item[$\bullet$] $M$ is a quiver $\D_Y$-module associated with that
arrangement.
\end{enumerate}
Then it turns out that the specialization is also a quiver
$\D_{T_ZY}$-module for a suitable arrangement  of hyperplanes
in $T_ZY$ and a suitable
quiver associated with that arrangement in $T_ZY$.

\medskip

More precisely, let $\A$ be a central arrangement of hyperplanes in a
vector space
$X=\CC^N$. Let $\Ga$ be the graph of the arrangement $\A$.
Fix $\a\in\v(\Ga)$. The arrangement $\A$ induces an
arrangement $\A_\a$ in the total space $T_{\X_\a}X$
of the normal bundle $T_{\X_\a}X \to \bar X_\a$.

It is defined
as follows. Identify $T_{\X_\a}X$ with the direct sum
$\X_\a\oplus X/\X_\a$ and let $\pi_\a:X\to X/\X_\a$ be the natural
projection.  Assign to a hyperplane $H$ of $\A$ the hyperplane
$$
\tilde{H}\ =\ H\cap \X_\a\oplus \pi_\a(H)\subset \X_\a\oplus X/\X_\a\ .
$$
Then
the arrangement $\A_\a$ consists of all distinct hyperplanes among the
hyperplanes $\tilde{H}_j, j\in J(\A)$.

The graph $\Ga_\a$ of the arrangement $\A_\a$ can be described
as follows. Define an equivalence relation $\equiv_\a$
on points of the set $\v(\Ga)$\,:\  we say that
$\b\equiv_\a\!\gamma$ if
\begin{enumerate}
\item[(i)]
$\bar X_\a\cap\bar X_\b\,=\,\bar X_\a\cap\bar X_\gamma$\,,
\item[(ii)] the set of all $\delta\in\v(\Ga)$,
such that $\d\geq \a$ and $\d\geq\b$ coincides with the set of all
$\delta\in\v(\Ga)$, such that $\d\geq \a$ and $\d\geq\g$.
\end{enumerate}
Then the vertices of $\Ga_\a$ are the equivalence classes $\bar{\b}$
relative to
the relation $\equiv_\a$\,.

Two vertices $\bar{\b}$, $\bar{\gamma}\in
\v(\Ga_\a)$ are connected by an arrow in $\Ga_\a$ if there exist representatives
$\b\in \bar{\b}$ and $\gamma\in \bar{\gamma}$, connected by an arrow
in $\Ga$.

The length function $l : I(\Gamma_\a) \to \ZZ_{\geq0}$ is given by the rule:
the length $l(\bar{\b})$ of any equivalence class
$\bar{\b}\in\v(\Ga_\a)$ is the length of any representative
$\b'\in\bar{\b}\subset \v(\Ga)$.

\medskip

Let $\VA=\{V_\b, A_{\b,\g}\}$ be a quiver of $\Ga$. We define the new quiver
$Sp_\a(\VA)=\{{ W}_{\bar{\b}}, C_{\bar{\b},\bar{\gamma}}\}$ of
the graph $\Ga_\a$
which will be called {\it the specialization} of
the  quiver $\VA$ at $\a \in I(\Ga)$.
Namely, we set
$$
{ W}_{\bar{\b}}\ =\ \oplus_{\b'} \,V_{\b'}\ ,
$$
where the sum is  over
all $\b'$ from the equivalence class $\bar{\b}$. We set
$$
C_{\bar{\b},\bar{\gamma}}\ =\ \sum_{\b',\gamma'}\, A_{\b',\gamma'}\ ,
$$
where the sum is
over all  $\b'$ from the equivalence class $\bar{\b}$ and
$\g'$ from the equivalence class $\bar{\g}$. It is easy to check that
$Sp_\a(\VA)$  is a quiver of $\Ga_\alpha$.

\medskip

Let $F(\D_X)$ be the decreasing filtration \rf{FkD} of the sheaf
$\D_X$, related to the ideal of functions vanishing on
$\X_\a$. Define the decreasing filtration $F(\EVA)$ of the sheaf
$\EVA$ as follows. For any $j=0,1,...,l=l(\a)$, denote by
$\overline{V}_j$ the vector space
$$
\overline{V}_j\ =\ \oplus_\gamma\
\overline{\Omega}_\gamma\ot V_\gamma\ ,
$$
where the sum is taken over all
$\gamma$ with $\codim_{\X_\gamma} \X_\gamma\cap \X_\a=j$. We set
\begin{equation}\label{FEVA}
F^k(\EVA)=\left\{
\begin{array}{ll}
\sum_{j=0}^l F^{k+l-j}(\D_X)\overline{V}_j\,, & k\leq -l\,,
\\ [0.2em]
\sum_{j=k+l}^l F^{k+l-j}(\D_X)\overline{V}_j\,,  & -l\leq k\leq 0\,,
\\ [0.2em]
F^{k}(\D_X)\overline{V}_l\,,&  0\leq k\, .\\
\end{array}\right.
\end{equation}

{}Fix a nondegenerate bilinear form $(\,,\,)$ on $X=\C^N$.
For any vertex $\gamma\in \v(\Ga)$, the bilinear form determines an
isomorphism $\nu_\gamma : \X_\gamma \to \X_{\overline{\gamma}}$ of
vector spaces where
$\X_\gamma\subset X$ and $\X_{\overline{\gamma}}\subset
T_{\X_\a}X\equiv \X_\a\oplus X/\X_\a$. Namely, let $\X'_\gamma$ be the
orthogonal complement to $\X_\gamma\cap\X_\a$ in $\X_\gamma$.  We set
$\nu_\gamma(x + x')=x+\pi_\a(x')$ for any $x\in \X_\gamma\cap\X_\a$ and
$x'\in X'_\gamma$, where $\pi_\a: X\to X/\X_\a$ is the natural
projection. Let $\overline{\nu}_\gamma=\left(\nu_\gamma^*\right)^{-1}:
\overline{\Omega}_\gamma\to \overline{\Omega}_{\overline{\gamma}}$ be
the associated isomorphism of the spaces of top degree holomorphic
differential  forms  invariant with respect to translations.

{}For any vertex $\b\in\v(\Ga)$ define the linear operator $S_\b^{(\a)}:
V_\b\to V_\b$ by the formula
$$
S_\b^{(\a)}\ =\ \sum_\gamma A_{\b,\gamma}A_{\gamma,\b}\ ,
$$
where the sum is over
all $\gamma\in\v(\Ga)$ such that
$\X_{\a}\cap\X_\gamma=\X_\b\cap\X_\gamma$. We define
the operator
$$
S^{(\a)}\ :\ \oplus_{\b\in\v(\Ga)}V_\b\ \to\ \oplus_{\b\in\v(\Ga)}V_\b
$$
by the formula
 $$
S^{(\a)}\ =\ \sum_{\beta\in\v(\Ga)}S_\b^{(\a)}\,.
$$
\begin{lemma}${}$
\label{newlemma}
\begin{enumerate}
\item[ (i)]

The filtration $F(\EVA)$ is good with respect to $F(\D_X)$.

\item [(ii)]
The assignment $\gr( \overline{\omega}_\gamma\ot
v_\gamma)\,\mapsto \,\overline{\nu}_\gamma(\overline{\omega}_\gamma)\ot
v_\gamma$, with $\gamma\in\v(\Ga)$, $v_\gamma\in V_\gamma$,
$\overline{\omega}_\gamma\in \overline{\Omega}_\gamma$, establishes an
isomorphism of the  $\D_{T_{\X_\a}X}$-modules $\gr \EVA$ and $E
Sp_\a(\VA)$.
\item[(iii)] For  $\b\in\v(\Ga)$ let $j  = \codim_{\X_\b} \X_\b\cap \X_\a$.
Then  for any $\om_\b\ot v_\b\in \overline{\Omega}_\b\ot V_\b$ we have
$$
\left( \theta-j+l\right)(\om_\b\ot v_\b)-1\ot S_\b^{(\a)}(\om_\b\ot v_\b)
\quad \in \quad
 F^{j-l+1}(\EVA)\,.
$$
 \end{enumerate}
\hfill{$\square$}
\end{lemma}

The proof is by direct calculation like in the proof of
Proposition \ref{theorem1}.

\medskip

\begin{lemma}
\label{newcor}
Let ${G}^{(\a)}$ be a non-resonant set and
$\overline{G}^{(\a)}$ a non-resonant section such that
$\overline{G}^{(\a)}\supset {G}^{(\a)}$.
Assume that all
eigenvalues of  $S^{(\a)}$
belong to ${G}^{(\a)}$.  Then the
$\D_{T_{\X_\a}X}$-module $\gr \EVA$ from Lemma \ref{newlemma} is the
specialization module $Sp_{\X_\a}(\EVA)$ relative to the non-resonant
section  $\overline{G}^{(\a)}$ and  vector field $\theta$
defined above.
\end{lemma}

{\it Proof}.  Due to statements (i) and (ii) of Lemma \ref{newlemma},
it is sufficient to verify condition (v) for the filtration
$F(\EVA)$. But condition (v)
follows from the definition of the filtration in \rf{FEVA},
statement (iii) of Lemma \ref{newlemma}, and Proposition \ref{prop1}.
\hfill{$\square$}

\begin{corollary}
\label{COR}
Under the above assumptions
the $\D_{T_{\X_\a}X}$-modules
$E Sp_\a(\VA)$ and
\linebreak
$Sp_{\X_\a}(\EVA)$ are isomorphic with
respect to the isomorphism from
Lemma \ref{newlemma}.
\end{corollary}

\medskip
\subsection{Free resolution of quiver $\D_{X^0}$-modules}
 Let $\VA_{\G^0} = \{ V, A^j\ |\ j\in J(\A)\}$ be a
quiver in $\Qui_{\G^0}$ and $E^0\VA_{\G^0}$ the corresponding
$\D_{X^0}$-module.

Set
$$
\R^rE^0\VA_{\G^0}=\D_{X^0}\ot \bigwedge\nolimits^{-r} T_\emptyset\ot
\overline{\Omega}_\emptyset\ot V\ ,
\qquad
\R E^0\VA_{\G^0}=\oplus_{r\in\ZZ}\R^rE^0\VA_{\G^0}\ .
$$
 These are free
left $\D_{X^0}$-modules.

Set $\deg \R^rE^0\VA_{\G^0}=r$. Define $\D_{X^0}$-linear maps
$d_{i}:\R E^0\VA_{\G^0}\to \R
E^0\VA_{\G^0}$, $i=0,1$,\ of degree $1$ and  a $\D_{X^0}$-linear map
$\nu: \R^0E^0\VA_{\G^0}\to
E^0\VA_{\G^0}$ by the formulas:
 $$
\begin{array}{rcl}
d_{1}\left(1\ot \xi_1\wedge...\wedge \xi_r\ot \om_\emptyset\ot
v\right)&=&\!\!\!\sum\limits_{i=1}^r (-1)^{i+1}\xi_i\ot
\xi_1\wedge... \widehat{\xi_{i}}...\wedge \xi_r\ot
\om_\emptyset\ot v,
\\
d_{0}\left(1\ot \xi_1\wedge...\wedge \xi_r\ot \om_\emptyset\ot
v\right)&= &\!\!\!\sum\limits_{j\in J(\A)} f_j^{-1}\ot
i_{f_j} ( \xi_1\wedge...\wedge \xi_r)\ot \om_\emptyset\ot
A^j(v),\\
\nu(1\ot \om_\emptyset\ot v)&=&\om_\emptyset\ot
v.\end{array}
$$
Here $f_j=0$ is an equation of the hyperplane $H_j$. Set also
$d=d_1-d_0$. We have
$$d_{i}^2=0,\quad {\rm for}\quad\, i=0,1, \qquad d_{0}d_{1}+d_{1}d_{0}=0.$$

\begin{theorem}\label{prop11}
The complex
$$
0\to\R^{-N}E^0\VA_{\G^0}\stackrel{d}{\to}\R^{-N+1}E^0\VA_{\G^0}
\stackrel{d^{}}{\to}\cdots \stackrel{d}{\to}\R^0E^0\VA_{\G^0}
\stackrel{\nu}{\to}E^0\VA_{\G^0}\to 0
$$
is acyclic and presents a free $\D_{X^0}$-module resolution of
the $\D_{X^0}$-module $E^0\VA_{\G^0}$.
\end{theorem}

The  theorem is proved in Section \ref{proofs}.

\subsection{Free resolution of quiver $\D_X$-modules}

 In this section we generalize the result of the previous section and
attach to  a quiver $\VA\in\Qui_\Ga$ an acyclic complex
$$
\Rlong\EVA:\qquad 0\to \R^{-N}\EVA\to\R^{-N+1}\EVA\to\ldots\to
\R^0\EVA\to \EVA\to 0
$$
of left $\D_X$-modules, where
all   $\R^i\EVA$ are finitely generated  free
 left $\D_X$-modules. This complex gives  a free resolution  of the quiver
$\D_X$-module $\EVA$.

\bigskip

In the construction of the resolution we will use the following objects introduced in
Section \ref{Framing}.

For $\a \in I(\G)$, we will use the associated spaces
$\TTT{\a} = \Ta \oplus \Fa$, and  $\Oa$.
For strata  $X_\a$ and $X_\b$,\
$\X_\a\supset X_\b$,\ we will use the
  inclusions
$$
\mu_{\b,\a}:\ \bigwedge^p\Fa\hookrightarrow\bigwedge^p\Fb\ ,
\qquad
\overline{\mu}_{\a,\b}:\ \bigwedge^p\Tb\hookrightarrow\bigwedge^p\Ta\ .
$$
 For $\a\succ\b$, $f\in \Fb$,  $\xi\in\Ta$, \ we will use the linear maps
\begin{equation}
i_f^{\b,\a}\ :\ \bigwedge^p\Ta\ \to\
\bigwedge^{p-1}\Tb\ ,
\qquad
\ii_\xi^{\a,\b}\ :\ \bigwedge^p \Fb\ \to\ \bigwedge^{p-1}\Fa\ .
\notag
\end{equation}
\bigskip

{} For $r=0, \dots ,-N$, set
\begin{equation}\label{Rr}
 \R^r\EVA = \ooplus\limits_{\a\in\v(\Ga)} \D_X\ot
\bigwedge\nolimits^{-r}\TTT{\a}\ot \Oa\ot V_\a\ .
\end{equation}
This is a free left $\D_X$-module. Set  $\R^r\EVA=0$ for $r>0$
and $r<-N$. Set
$$
\R\EVA\ =\ \oplus_{r\in\ZZ}\R^r\EVA\ .
$$

\bigskip

For an edge   $(\a,\b)\in E(\Ga)$, and
 $r,\ 1\leq r\leq N$, \ we define a linear map
$$
A_{\b,\a}^r\, : \bigwedge\nolimits^r\TTT{\a}\ot \Oa\ot V_\a\to
\bigwedge\nolimits^{r-1}\TTT{\b}\ot \Ob\ot V_\b \ .
$$
The construction of the map uses an edge framing
$\{ f_{\a,\b},\xi_{\a,\b}\}$, but the map does not depend on the choice
of the framing.

The map
$$
A_{\b,\a}^1\, : \bigwedge\nolimits^r\TTT{\a}\ot \Oa\ot V_\a\to
\bigwedge\nolimits^{r-1}\TTT{\b}\ot \Ob\ot V_\b \
$$
is defined as follows.
For
$\xi\in\Ta$, $f\in\Fa$, $\om_\a\ot v_\a\in \Oa\ot V_\a$ we
set
$$
A_{\b,\a}^1(\xi\ot \om_\a\ot
v_\a)=\mylangle \xi,f_{\b,\a}\myrangle\ \pi_{\b,\a}(\om_\a)\ot
A_{\b,\a}(v_\a),\qquad A_{\b,\a}^1(f\ot \om_\a\ot v_\a)=0
$$
if $\a\succ\b$ and
$$
A_{\b,\a}^1(f\ot \om_\a \ot
v_\a)=\mylangle \xi_{\b,\a},f\myrangle\ \pi_{\b,\a}(\om_\a)\ot
A_{\b,\a}(v_\a),\qquad A_{\b,\a}^1(\xi\ot \om_\a\ot v_\a)=0
$$
if $\b\succ\a$.
\medskip

Define the map $A_{\b,\a}^r$ for an arbitrary $r$.
Let $\xxi{}\wedge\f{}\in
\bigwedge^p\Ta\wedge\bigwedge^{r-p}\Fa\subset\bigwedge^{r}\TTT{\a}$
and $\om_\a\ot v_\a\in \Oa\ot V_\a$. If $\a\succ\b$, set
\begin{equation}
\label{Apq+} A^{r}_{\b,\a}(\xxi{}\wedge\f{}\otimes \om_\a\ot
 v_\a)=i_{f_{\b,\a}}^{\b,\a}(\xxi{})\wedge\f{} \otimes
\pi_{\b,\a}(\om_\a)\ot A_{\b,\a}(v_\a) ,
\end{equation}
where the polyvector $\f{}$ in the right hand side of \rf{Apq+} is considered
as an element of $\bigwedge^{r-p}\Fb$ by means of the tautological embedding
$\mu_{\a,\b}$.
If $\b\succ\a$,  set
\begin{equation}
\label{Apq-} A^{r}_{\b,\a}(\xxi{}\wedge\f{}\otimes \om_\a\ot
v_\a)=(-1)^p \xxi{}\wedge
i_{\xi_{\b,\a}}^{\b,\a}(\f{})
 \otimes\pi_{\b,\a}(\om_\a)\ot A_{\b,\a}(v_\a) ,
\end{equation}
where the polivector $\xxi{}$ in the right hand side of \rf{Apq-} is considered
as an element of $\bigwedge^{p}\Tb$ by means of the tautological embedding
$\overline{\mu}_{\a,\b}$.

\bigskip

 For
 $i=0, 1\ $, define  $\D_X$-linear maps $\ d_{i}:\R\EVA\to \R\EVA\ $
by the formulas:
\begin{equation}
\begin{array}{ccc}
d_{1}\left(1\ot t_1\wedge...\wedge t_r\ot \om_\a\ot
v_\a\right)&=&\!\!\!\sum\limits_{i=1}^r (-1)^{i+1}t_i\ot
t_1\wedge... \widehat{t_{i}}...\wedge t_r\ot \om_\a\ot v_\a,
\label{dr}
\\
d_{0}\left(1\ot t_1\wedge...\wedge t_r\ot \om_\a\ot v_\a\right)&=
&\!\!\!\sum\limits_{\b,\ (\a,\b)\in E(\Ga)} 1\ot A^r_{\b,\a}
(t_1\wedge...\wedge t_r\ot \om_\a\ot v_\a),\end{array}
\end{equation}
where $t_i\in\TTT{\a}$, $\om_\a\ot v_\a\in \Oa\ot V_\a$.
Set
$$
d_i(1\ot\om_\a\ot v_\a)\ =\ 0 \ .
$$
It is clear that  $d_i(  \R^r\EVA) \subset \R^{r+1}\EVA$ for $i = 0, 1$.

Define a morphism $\nu: \R^0\EVA \to \EVA$
 of $\D_X$-modules by the formula
$$
\nu(1\ot \om_\a\ot v_\a)\ =\ \om_\a\ot v_\a\ .
$$

\begin{proposition} ${}$\\
{\em(i)} The maps $d_{i}$ are  commuting differentials,
 $$
d_{i}^2=0,\quad {\rm for}\quad\, i=0,1, \qquad d_{0}d_{1}+d_{1}d_{0}=0;
$$
 {\em (ii)} Set $d=d_{1}-d_{0}$. Then
$$
\nu  d\ =\ 0\ .
$$
\label{proposition5}
\hfill
$\square$
\end{proposition}

\medskip

 Denote by ${\EVA}^{(\cdot)}$ the $\D_X$-module $\EVA$
 considered as a complex with the zero differential,
 concentrated in degree zero.
Then, by Proposition \ref{proposition5},
 we can introduce two complexes,
\begin{gather}\nonumber
\R\EVA:\quad 0\to\R^{-N}\EVA\stackrel{d^{}}{\to}
\R^{-N+1}\EVA\stackrel{d^{}}{\to}\cdots
\stackrel{d^{}}{\to}\R^0\EVA\to 0\ ,
\\ \nonumber
\Rlong\EVA:\quad 0\to\R^{-N}\EVA\stackrel{d^{}}{\to}
\R^{-N+1}\EVA\stackrel{d^{}}{\to}\cdots
\stackrel{d^{}}{\to}\R^0\EVA \stackrel{\nu}{\to}\EVA\to 0
\end{gather}
with $\deg \R^k\EVA =k$.  The second complex can be  considered as
 a morphism of complexes $\nu: \R\EVA\to{\EVA}^{(\cdot)}$\ .

The following statement is one of the main results of the paper.
\begin{theorem}
\label{theorem3} For any quiver  $\VA\in\Qui_\Ga$,
the complex $\Rlong\EVA$ is acyclic.
\end{theorem}

The theorem is proved in Section \ref{proofs}.

\medskip

 By  Theorem \ref{theorem3}, the complex $\Rlong\EVA$ gives a free
$\D_X$-module resolution  of the $\D_X$-module $\EVA$. Equivalently, the
 complex $\R\EVA$ of free $\D_X$ modules
is quasi-isomorphic to ${\EVA}^{(\cdot)}$,
 that is, the morphism  $\nu: \R\EVA\to{\EVA}^{(\cdot)}$ induces an isomorphism
of cohomology groups of complexes.

\medskip
\noindent
{\it Example}.
Let $X=\CC^1$, $\A=\{z=0\}$. The graph $\G$ of the arrangement has two
vertices, $\ee$ and $\a$, $l(\a)=1$, connected by an edge. Let $\VA=
\{V_\ee, V_\a, A_{\ee,\a},A_{\a,\ee}\}$ be a quiver of $\G$.
 The free $\D_X$-module resolution $\R E\VA$
of the quiver $\D_X$-module $E\VA$ is the complex
$$\R E\VA:\  0\to\R^{-1} E\VA\,\stackrel{d}{\to}\, \R^{0}
 E\VA\,\stackrel{v}{\to}\, E\VA\to 0.$$
 Here
$\R^{0} E\VA$ is the free $\D_X$-module
 $\D_X\ot \overline{\Omega}_\ee \ot V_\ee \oplus
\D_X\ot  \overline{\Omega}_\a\ot V_\a$. The term $\R^{-1} E\VA$
is the free $\D_X$-module
$\D_X\ot T_\ee\ot  \overline{\Omega}_\ee\ot V_\ee \oplus \D_X\ot F_\a\ot
 \overline{\Omega}_\a\ot V_\a$, where $T_\ee$ is the
one-dimensional space $\CC\cdot \frac d{dz} $ of vector fields
with constant coefficients, parallel to $X_\ee$, $F_\a$ is the
one-dimensional space $\CC\cdot  z$ of affine functions, vanishing at
$X_\a$. The differentials $d$, $\nu$ are given by the relations:
\begin{align*}
d(D\ot , \frac d{dz}\ot dz \ot v_\ee)&=D\cdot \frac d{dz}
 \ot dz\ot v_\ee -D\ot 1\ot A_{\a,\ee}(v_\ee) ,\\
d(D\ot   z\ot 1 \ot v_\a)&=D\cdot  z
 \ot 1\ot v_\a -D\ot dz\ot A_{\ee,\a}(v_\a),\\
\nu(D\ot 1 \ot v_\a)=D\cdot(1&\ot v_\a), \qquad
\nu(D\ot dz \ot v_\ee)=D\cdot(dz\ot v_\ee)\ ,
\end{align*}
 for any $D\in D_X$, $v_\ee\in V_\ee$, $v_\a\in V_\a$. Here $D\cdot(1\ot v_\a)$
and $D\cdot(dz\ot v_\ee)$ denote the application of $D$ to elements of
$E\VA$.

\subsection{Derived categories and derived functors}\label{derived}
In the following we use the language of derived categories.
Let us briefly recall necessary definitions and constructions,
see e.g. \cite{GM}.

Let $\Abb1$ be an abelian category.
Denote by ${\cal Kom}(\Abb1)$ (resp.\  ${\cal Kom}^b(\Abb1)$) the category
of  complexes (resp.\ bounded complexes) over $\Abb1$ with differential of degree one.
 Morphisms in both categories
are  morphisms of complexes, commuting with differentials:
$fd_{ A}=d_{ B}f$.

{} For any complex $A^\mycircle$ $\in$
 ${\cal Kom}(\Abb1)$ and $n\in\ZZ$ denote by $A[n]^\mycircle$ the same complex
 with the
shifted gradation, $A[n]^k=A^{n+k}$. For any object $A\in\Abb1$
denote by $A^{(\cdot)}$ the complex with zero differential, defined by the
condition $(A^{(\cdot)})^0=A$ and $(A^{(\cdot)})^k=0$ for $k\not=0$.

A morphism of complexes
 $f:  A^\mycircle $ $\to$
$ { B}^\mycircle$ is called a quasi-isomorphism, if it induces
an isomorphism  of  cohomology groups ${ H}({ A}^\mycircle)$ and
${ H}({ B}^\mycircle)$.
It is  denoted by
$f$: $ A^\mycircle $ $\stackrel{quis}{\to}$ $ { B}^\mycircle$.
  Morphisms $f,f'$: ${ A}^\mycircle $ $\to$
$ { B}^\mycircle$ are called homotopic, if there exists a map of complexes
  $h$: ${ A}^\mycircle $
$\to$
$ { B}^\mycircle$ of degree $-1$, (that is, a morphism $h$: ${ A}^\mycircle $
$\to$
$  B[1]^\mycircle$), such that $f-f'=hd_{ A}+d_{ B}h$.

Denote
 by $\cal {K}(\Abb1)$ (resp.\ $\cal {K}^b(\Abb1)$) the homotopy (respectively,
bounded homotopy) category of complexes over $\Abb1$.
  Objects of $\cal {K}(\Abb1)$ (resp.\  $\cal {K}^b(\Abb1)$)
are  objects of
${\cal Kom}(\Abb1)$, (resp.\ ${\cal Kom}^b(\Abb1)$). Morphisms are  classes of morphisms of complexes
 modulo homotopy equivalence.

The derived category  $D(\Abb1)$ is defined as the localization
of the category $\cal {K}(\Abb1)$ over the set of all classes of
quasi-isomorphisms.
 The bounded derived category  $D^b(\Abb1)$ is defined as the localization
of the category $\cal {K}^b(\Abb1)$ over the set of all classes of
quasi-isomorphisms.

Objects of $D(\Abb1)$ are   complexes
over $\Abb1$, and every morphism
$f$: ${ A}^\mycircle $ $\to$
$ { B}^\mycircle$ can be presented as a
triple $\{{ C}^\mycircle$ $\in$ $\cal {K}(\Abb1)$,
$f'$ $\in$ $Hom_{\cal {K}(\Abb1)}({ A}^\mycircle, {C}^\mycircle)$,
$f''$ $\in$ $Hom_{\cal {K}(\Abb1)}({ B}^\mycircle, { C}^\mycircle)\}$,
where $f'$ is a class of quasi-isomorphisms. An analogous description holds for
the bounded derived category $D^b(\Abb1)$.

Let $\Abb1$ and $\Ab2$ be abelian categories and $F$: $ \Abb1\to\Ab2$
 a left (resp. right) exact additive functor. It admits a natural extension
to a functor
$F^\mycircle : K(\Abb1)\to K(\Ab2)$  of homotopy categories. Usually
it is denoted by the same symbol $F$.

A class $\Abb1_F$ of objects
 in $Ob(\Abb1)$
is called adjusted to $F$ if
\begin{itemize}
\item[(i)] It is closed under direct sums;
\item[(ii)] The functor $F^\mycircle$ transforms acyclic bounded from below
(resp. from above) complexes with graded
components in $\Abb1_F$ to acyclic complexes;
\item[(iii)] For any object $A\in\Abb1$, there exists an object $A'\in \Abb1_F$
 and an inclusion $\a:A\to A'$ (resp. an epimorphism $A'\to A$).
\end{itemize}

If a class $\Abb1_F$ of objects adjusted to $F$ is chosen, then
the right  derived functor $RF  : D(\Abb1)\to D(\Ab2)$
\ (resp. the left derived functor $LF : D(\Abb1)\to D(\Ab2)$ )\
is defined as follows. By (iii), for any complex $A^\mycircle$ $\in$ $D(\Abb1)$,
there exists a complex ${A'}^\mycircle$ over $\Abb1_F$ and a quasi-isomorphism
$\phi : A^\mycircle \to {A'}^\mycircle$ (resp. a complex ${A''}^\mycircle$ and a
quasi-isomorphism $\psi: {A''}^\mycircle \to A^\mycircle$).
By definition, $RF(A^\mycircle)=F^\mycircle({A'}^\mycircle)$ (resp.
$LF(A^\mycircle)=F^\mycircle({A''}^\mycircle)$).

It can be proved that the derived functor does not depend on the choice
of the class of adjusted objects.

We keep the natural notation $RF(A) = RF(A^{(\cdot)})$ and
$LF(A) = LF(A^{(\cdot)})$ for any object $A\in\Abb1$.

In most cases   left and
right derived functors admit well defined restrictions to
 bounded derived categories.
 It takes place, in particular,
 if any object $A$ in $\Abb1$ admits a finite resolution
$\phi : A^{(\cdot)}$  $\stackrel{quis}{\to}$ ${A'}^\mycircle$ over $\Abb1_F$
for the case of a left exact functor $F$ (resp. admits a finite resolution
$\phi:{A'}^\mycircle$  $\stackrel{quis}{\to}$ $A^{(\cdot)}$ over $\Abb1_F$
for the case of a right exact functor $F$).
It is always so in our considerations.

{} For a smooth algebraic variety $Y$ we denote by $D^b(\MD_Y)$ the
bounded derived category of the category  of
coherent $\D_Y$-modules. For a topological
space $Y$ we denote by $D^b(Y,\CC)$ the bounded derived category of the
category  $Sh(Y,\CC)$ of sheaves of complex vector spaces on $Y$.

\subsection{Duality of quiver $\D_X$-modules}\label{dualsection}
Let $Y$ be a smooth algebraic variety over $\CC$ and $\Omega_Y$  the
sheaf
of germs of top exterior differential
forms on $Y$. The sheaf $\Omega_Y$ is an
invertible locally free
 sheaf of ${\cal O}_Y$-modules of rank $1$.
Let $\Omega_Y^{-1}$ be the corresponding inverse sheaf.
Denote by $\widetilde{\D}_Y$ the sheaf $\D_Y\ot_{{\cal O}_Y}\Omega_Y^{-1}$.
It has two structures of a (locally free) left $\D_Y$-module.

The first structure  is given locally by the  relations:
\begin{equation}\label{structure1}
f\cdot (d\ot \tilde{\omega})= df\ot \tilde{\omega},
\qquad \xi \cdot(d\ot \tilde{\omega})=-d\xi \ot \tilde{\omega}+
d\ot Lie_\xi(\tilde{\omega})\ ,
\end{equation}
where $f$, $d$, $\tilde{\omega}$ are local sections of the sheaves ${\cal O}_Y$,
$\D_Y$, and $\Omega_Y^{-1}$, respectively, $\xi$ is a vector field over $Y$
and $Lie_\xi$ is the
Lie derivative in the direction of $\xi$.

The second structure
of a left $\D_Y$ module on $\widetilde{\D}_Y$
is given by the left multiplication of
$\D_Y$ on itself:
\begin{equation}\nonumber
d_1\star (d_2\ot \tilde{\omega})= d_1d_2\ot \tilde{\omega}.
\end{equation}
Let $M^\mycircle$ be a complex of $\D_Y$-modules.  We  set
$$
\Delta_Y M^\mycircle\ =\ { R}\, Hom_{\D_Y}(M^\mycircle  ,  \widetilde{\D}_Y)[\dim Y]\ .
$$
This formula has the following meaning.
 The complex $\Delta_Y M^\mycircle$ of left $\D_Y$-modules is
 quasi-isomorphic to the complex $Q^\mycircle$, where $Q^i$ is the sheaf
$$
Q^i\ =\
Hom_{\D_Y}(P^{-i-\dim Y},\widetilde{\D_Y})\ .
$$
The complex $P^\mycircle$ is a complex of locally free left $\D_Y$-modules,
quasi-isomorphic to $M^\mycircle$. The local $Hom$ functor
 is taken between the left $\D_Y$-module $P^{-i-\dim Y}$ and the
left $\D_Y$-module $\widetilde{\D_Y}$,  considered with
respect to the first structure.
The structure of a left $\D_Y$-module on $Q^i$ is induced from the
second module structure on  $\widetilde{\D}_Y$ as follows.
\ Let $\phi$ be a local homomorphism $\phi: P^{-i-\dim Y}$ $\to$
 $\widetilde{\D_Y}$, satisfying the condition $\phi(d'\cdot p)=d'\cdot
\phi(p)$ for any $p\in P^{-i-\dim Y}$ and  $d'\in \D_Y$. Then for any
$d\in \D_Y$,  the local $\D_Y$-linear
map $d\cdot\phi$ is defined by the rule
$(d\cdot\phi)(p)$ $=$ $d\star \phi(p)$.

It is known that if $M$ is a holonomic $\D_Y$-module, then
$H^i(\Delta_Y M^{(\cdot)})=0$ for $i\not= 0$.
The map $M \mapsto H^0(\Delta_Y M^{(\cdot)})$ defines
an exact contravariant functor $\MD_Y^{hol}\to \MD_Y^{hol}$ which we denote
by the same symbol $\Delta_Y$. For any $M \in \MD_Y^{hol}$, the $\D_Y$-module
$\Delta_Y( \Delta_Y(M))$ is isomorphic to $M$.

Let $j : Z \to Y$ be an open embedding of smooth algebraic varieties.  Then, besides
 the functors $j^*: \MD_Y\to \MD_Z$ and $j_*: \MD_Z\to \MD_Y$ of the inverse
and direct images of $\D$-modules, the second direct image functor
$j_\my!^{(D)}: \MD_Z^{hol}\to \MD_Y^{hol}$ can be defined by the rule:
\begin{equation}\label{j!D}
j_\my!^{(D)}(M)= \Delta_Yj_*\Delta_Z(M), \qquad M\in\MD_Z^{hol}\ .
\end{equation}
It is the left adjoint functor to the functor  $j^*: \MD_Y^{hol}\to
\MD_Z^{hol}\ .$
{}For any $M\in  \MD_{Z}^{hol}$, we denote by $j_{!*}^{(D)}(M)$ the
 image of $j_{!}^{(D)}(M)$ in $j_{*}(M)$ with respect to the natural map
$s_j^{(D)}$: $j_{!}^{(D)}(M)$ $\to$ $j_{*}(M)$,
induced by  canonical transformations between the pairs of adjoint functors
($j^*$,  $j_*$) and  ($j_!^{(D)}$, $j^*$), see Section
\ref{Rdiq}.

Denote by  $j_{l,k,!}: \MD_{X^k}^{hol}\to \MD_{X^l}^{hol}$ the direct
image functor $(j_{l,k})_\my!^{(D)}$. It
is left adjoint to the functor of  inverse image
$j_{l,k}^*: \MD_{X^l}^{hol}\to \MD_{X^k}^{hol}$. For
$0\leq k_1<k_2<k_3\leq N$, we have
$$
 j_{k_3,k_1,!}= j_{k_3,k_2,!}\cdot j_{k_2,k_1,!}\ .
$$
We use the notation $j_{0,!}$ for the functor
$j_{N,0,!}:  \MD_{X^0}^{hol}\to \MD_{X^N}^{hol}=\MD_{X}^{hol}$.
It is exact since the map $ X^0$ is an affine variety \cite{B}.
For any $M\in  \MD_{X^0}^{hol}$,  we denote by   $j_{0,!*}(M)$
the $\D_X$-module $(j_{N,0})_{!*}^{(D)}(M)$.

\begin{theorem}\label{propdual} ${}$\\
{\em (i)} Let $\VA_{\G^0}$ be a level zero quiver. Then the
following $\D_{X^0}$-modules are isomorphic:
\begin{equation}\nonumber
\Delta_{X^0} (  E^0(\VA_{\G^0}))\ \approx \ E^0(\tau (\VA_{\G^0}))\ .
\end{equation}
{\em (ii)} Let $\VA$ be a  quiver of $\Ga$.
 Then the
following $\D_{X}$-modules are isomorphic:
\begin{equation}\nonumber
\Delta_{X} ( E(\VA))\ \approx\ E(\tau (\VA))\ .
\end{equation}
\end{theorem}
The theorem is proved  in Section \ref{proofs}.

\medskip

\begin{corollary}\label{corj!}
Let $\VA_{\G^0} = \{V, A^i\}$ be a level zero quiver.
Suppose that each operator $A^i$ has a single eigenvalue.
Suppose that $\VA_{\G^0}$ has a non-resonant spectrum, see Section
\ref{Spectrum}. Then the
following $\D_{X}$-modules are isomorphic:
\begin{eqnarray}
j_{0,!} ( E^0(\VA_{\G^0}))\
&\approx\ &  E ( \j_{0,!}(\VA_{\G^0}))\ ,\nonumber\\
j_{0,!*} ( E^0(\VA_{\G^0}))
\ &\approx\ & E ( \j_{0,!*}(\VA_{\G^0}))\ .\nonumber
\end{eqnarray}
\end{corollary}

{\it Proof.} The first statement follows from the definitions of the
 functors $j_{0,!}=(j_{N,0})_!$ and $\j_{0,!}=\j_{N,0,!}$, see
\rf{compj2}, \rf{j!D},   and Theorem \ref{propdual}.

The canonical morphism
$s_0:$ $\j_{0,!}(\VA_{\G^0})$ $\to$ $\j_{0,*}(\VA_{\G^0})$, see
\rf{phikl}, is determined by  the collection of isomorphisms of adjunction
$$
\a_{{\mathcal U}_{\G}, \VA_{\G^0}}\ : \
\Hom_{\Qui_{\G}}({\mathcal U}_{\G},\j_{0,*}(\VA_{\G^0})) \ \to \
\Hom_{\Qui_{\G^0}}(\j_{0}^*({\mathcal U}_{\G}),\VA_{\G^0})\ ,
$$
where
${\mathcal U}_{\G}\in \Qui_\G$,  and by the automorphisms $\tau_0$ and $\tau_N$.
Analogously, for any $M_0\in \MD_{X^0}$, the canonical morphism
$s^{(D)}_0:$ $j_{0,!}(M_0)$ $\to$ $\j_{0,*}(M_0)$,
 is determined by  the collection of isomorphisms of adjunction
 $$
\tilde{\a}_{M,M_0}\ :\
\Hom_{{\mathcal M}_{X}}(M,j_{0,*}(M_0)) \ \to\
\Hom_{{\mathcal M}_{X^0}}(j_{0}^{*}(M),M_0)\ ,
$$
where
$M\in \MD_X$, and by the automorphisms $\Delta_{X^0}$ and $\Delta_X$.
 By  Theorem \ref{propdual} and  Theorem \ref{prop4}, statement (ii),
we have an equality $Es_0(\j_{0,!}(\VA_{\G^0}))$ $=$
$ s_0^{(D)}(E\j_{0,!}(\VA_{\G^0}))$, which implies the second statement of
Corollary \ref{corj!}. \hfill{$\square$}
\medskip

\subsection{Fourier transform of  quiver $\D_X$-modules}
\label{Fourier transform of  quiver}

\subsubsection{}
Let $z_1,...,z_N$ be linear coordinates on $X=\CC^N$ and $\xi_1,...,
\xi_N$ the dual coordinates on the dual space $X^*$.
The assignment
$$
\xi_i \,\mapsto\, -\frac{\partial}{\partial{z_i}}\ ,
\qquad
\frac{\partial}{\partial{\xi_i}}\mapsto  z_i\
$$
for  $i=1,\dots,N$,  defines an isomorphism
 of the rings $D_{X^*}$ and $D_{X}$ called {\it the Fourier transform.}
The Fourier transform defines a functor from the category $\MD_X$ to the category
$\MD_{X^*}$.

\bigskip

Let $\A=\{H_j\}, j\in J(\A)$, be a central arrangement of
hyperplanes in the vector space $X=\CC^N$. Let $\Ga$ be
the graph of the arrangmeent, $\VA$ a quiver of $\Ga$, and
$E\VA$ the associated quiver $\D_X$-module. In this section  we
shall discuss
how to describe combinatorially  the Fourier transform of $E\VA$.
It turns out that the Fourier transform of  $E\VA$ is the  quiver  $\D_{X^*}$-module
associated with a suitable line arrangment in $X^*$.

First we shall describe the line arrangment.
 Then we shall introduce quivers  associated with the
line arrangment. Then we shall describe the Fourier transform
of $E\VA$.

\subsubsection{}
Let $\A=\{H_j\}, j\in J(\A)$, be a central arrangement of
hyperplanes in $X=\CC^N$.
Then the one-dimensional subspaces
$H_j^\perp\in X^*$, $j\in J(\A)$, form an arrangement  of lines in
 $X^*$, denoted by $\A^\perp$.

The arrangement $\A^\perp$ defines on $X^*$
the structure of a stratified space.  The strata of $X^*$ are labeled by
vertices of the graph $\Ga$ of the arrangement $\A$.
For $\a\in I(\Ga)$, the closed stratum
$\X_\a^* \subset X^*$
is defined to be the subspace $\X_\a^\perp$. Equivalently, if
$\X_\a \subset X$ is the intersection of hyperplanes $H_j, j\in J_\a$, then
$\X_\a^*$ is the sum of the one-dimensional subspaces
$H_j^\perp, j\in J_\a$.

Clearly, the adjacencies of strata in $X^*$ are described by edges of
$\Ga$. Thus the adjacency graph of the stratified space $X^*$ is
naturally isomorphic to the adjacency graph $\Ga$ of the stratified
space $X$.

\subsubsection{}

Having a quiver $\VA=\{V_\a,A_{\a,\b}\}$ of the graph $\Ga$ of the
arrangement $\A$, we shall construct  a $\D_{X^*}$-module denoted by
$\tE\VA$ and called {\it the quiver $\D_{X^*}$-module associated with
the line arrangement $\A^\perp$ and quiver $\VA$}.

\medskip
\noindent
{\bf Remark.}
All our previous constructions were associated with $\D$ modules
on $X$. This is the first construction which considers
$\D$-modules on the dual space $X^*$.

\medskip

The construction of the    $\D_{X^*}$-module $\tE\VA$ is analogous to the
construction of the $\D_X$ module $E\VA$ in
 Section \ref{Sec.3.2}. Namely,
the    $\D_{X^*}$-module $\tE\VA$ is
 the quotient
of the free  sheaf $\R^0 \tE\VA$ $=$
$\D_{X^*}\otimes_\C(\oplus_{\a\in \v(\G)} \Oa\ot V_\a)$ of
  $\D_{X^*}$-modules
over its subsheaf ${\R'} \tE\VA$
 of coherent  $\D_{X^*}$-modules,
 whose global sections are the following global sections of  $\R^0 \tE\VA$:
\begin{enumerate}
\item[(i)] the sections
 \begin{equation*}
u\,\xi \ot \om_\a\ot v_\a \ -\ \sum_{\b,\a\succ\b}
\mylangle \xi,f_{\b,\a}\myrangle \,u\ot\pi_{\b,\a}(\om_\a) \ot
A_{\b,\a}(v_\a)\ ,
\end{equation*}
where  $\a\in\v(\Ga)$,  $u\ot \om_\a\ot v_\a\in \DD_{X^*}\otimes_\C\Oa\ot V_\a$
and  $\xi\in\Ta$,
\item[(ii)] the sections
\begin{equation*}
u\, f  \ot \om_\a\ot v_\a\ -\ \sum_{\b,\b\succ\a}
 \mylangle \xi_{\b,\a},f\myrangle\, u\ot
\pi_{\b,\a}(\om_\a)\ot
 A_{\b, \a}(v_\a)\ ,
\end{equation*}
where $\a\in\v(\Ga)$,  $u \ot\om_\a\ot v_\a\in \DD_{X^*}\ot_\C\Oa\ot V_\a$ and
$f\in\Fa$.
\end{enumerate}
Here
$\Oa$ denotes the one-dimensional space of the top degree
holomorphic differential forms on $X_\a^*$, invariant with respect to
translations along $X_\a^*$, and
 $\{f_{\b,\a}, \xi_{\a,\b}\ |\ \a,\b\in\v(\Ga),\a\succ\b\}$
is an edge framing of $\A^\perp$.

It is easy to see that $\tE\VA$  does not depend on
the choice of the edge framing.

\subsubsection{}

Given a quiver $\VA$, denote
by $\tVA$ the quiver  $\tVA=\{W_\a, B_{\a,\b}\}$, where $W_\a=V_\a$ and
$B_{\a,\b}=\ve(\b,\a)A_{\a,\b}$.

\begin{proposition}
The Fourier transform of the quiver $\D_X$-module $E\VA$
is isomorphic to the quiver $\D_{X^*}$-module
$\tE\tVA$.
\end{proposition}

The proof is by
easy calculations. The choice of the isomorphism depends
on the choice of a top degree differential
holomorphic form on $X^*$, invariant with respect to translations.

\medskip
\noindent
{\bf Remark.}
The constructions of Section \ref{Fourier transform of  quiver}
 show that quiver $\D$-modules can be defined
not only for hyperplane arrangements but for  central line
arrangements as well. It is easy to see that the quiver $\D$-modules
of line arrangements are all holonomic with regular
singularities. They admit an explicit free resolution and have all
nice properties of quiver $\D$-modules of hyperplane arrangements.
\setcounter{equation}{0}
\section{Quiver perverse sheaves}

\medskip

\subsection{Perverse sheaves over analytic varieties}
Let $Y$ be a smooth analytic variety over $\CC$ of dimension $N$.
Let  $\mathcal S$ be a stratification of $Y$
by smooth subvarieties, satisfying
equisingularity conditions. It means that $Y$ is presented as a
finite disjoint union of strata, such that the closure of each stratum
 is a union of strata, and for any stratum $S$ and any points
$x,y\in S$, there exists a diffeomorphism $\phi$ of $Y$,
preserving the stratification and such that $\phi(x)=y$.

Denote by  $Y^k$, $k=0,1,...$,
 the closure of the union of all strata of codimension $k$.
The subspaces $Y^k$ form a filtration:
$$Y=Y^0\supset Y^1\supset \dots \supset Y^N\supset Y^{N+1}=Y^{N+2}=...=
\emptyset\ .$$
A complex $\mathcal F$ of sheaves of vector spaces with differential of degree
one is called a perverse sheaf
with respect to the stratification
if it satisfies the following conditions on supports of its cohomology
groups
 \cite{GMc} :
\begin{equation}\label{perverse}
{\rm supp}\ H^{-N\pm k} (\mathcal F)\ \subset\ Y^k
\end{equation}
{}for all $k=0,1, \dots $. A complex  of sheaves $\mathcal F\in D^b(Y, \C)$
 of vector spaces over $Y$ is called a perverse sheaf if it is perverse with
respect to some stratification of $Y$, satisfying the conditions above.

Any continuous map $j: Y_1\to Y_2$ of topological spaces determines
the functor $j^\mycircle:  Sh(Y_2,\CC)$ $\to$ $Sh(Y_1,\CC)$
of inverse image
 and two functors $j_\mycircle: Sh(Y_1,\CC)$ $\to$ $Sh(Y_2,\CC)$ and
$j_\my!: Sh(Y_1,\CC)$ $\to$ $Sh(Y_2,\CC)$ of direct image.

 For any sheaf ${\mathcal G}\in Sh(Y_2,\CC)$
and a point $y_1\in Y_1$,\ the stalk $j^\mycircle{\mathcal G}_{y_1}$ of
  $j^\mycircle{\mathcal G}_{}$ at  $y_1\in Y_1$ is the
stalk  ${\mathcal G}_{j(y_1)}$.
For any sheaf ${\mathcal F}\in Sh(Y_1,\CC)$ and an
 open  $U\subset Y_2$,\ the sections $\G(U, j_\mycircle{\cal F})$ of
  $j_\mycircle{\cal F}$ are
\begin{equation}\label{dirim}
\G(U, j_\mycircle{\cal F})\ =\ \G(j^{-1}(U), {\cal F})\ .
\end{equation}
 The sections
$\G(U, j_\my!{\cal F})$ are those sections
$s\in \G(j^{-1}(U), {\cal F})$, for which
 the map $j: supp(s)\to U$ is proper.

Both functors
$j_\mycircle$ and $j_\my!$  are left exact. The functor
 $j^\mycircle$ is exact. Denote by \linebreak
$Rj_!$ : $D^b(Y_1,\CC)\to
D^b(Y_2,\CC)$ and
 $Rj_\mycircle$ : $D^b(Y_1,\CC)$ $\to$ $ D^b(Y_2,\CC)$  the
corresponding derived functors.

{}For a topological space $Y$ denote by
 $\overline{\Delta}_{Y^{}}: D^b(Y^{},\CC)\to D^b(Y^{},\CC)$
 the Poincare-Verdier duality,
$$
\overline{\Delta}_{Y^{}}({\mathcal F})\ =\
{\rm R}\, Hom({\mathcal F},
{\mathbf D}_{Y^{}})\ ,
$$
where ${\mathbf D}_{Y^{}}$ is  {\it the dualizing complex}. For a smooth
complex analytic variety $Y$  the dualizing complex is isomorphic
in $D^b(Y,\CC)$
to the shifted constant sheaf:
$${\mathbf D}_{Y^{}}=\CC_{Y^{}}[2\,  \dim_{\CC} Y]\ .$$
For any ${\mathcal F}\in D^b({Y_1},\CC)$ there is an isomorphism in
$D^b({Y_2},\CC)$:
$$
R  j_\my!
\ \approx\
\overline{\Delta}_{Y_2}\, R j_\mycircle\, \overline{\Delta}_{Y_1}
\ .$$

Equalities \rf{dirim} determine isomorphisms of adjunction
$$
\a_{{\mathcal G},{\mathcal F}}\ :\
{\rm Hom}_{Sh(Y_2,\CC)}({\mathcal G},
j_{\mycircle}({\mathcal F}))\ \risom  \
{\rm Hom}_{Sh(Y_1,\CC)}(j^\mycircle({\mathcal G}),
{\mathcal F})\ ,$$
which extend to isomorphisms of adjunction in derived categories,
$$
\overline{\a}_{{\mathcal G},{\mathcal F}}\ :\
 {\rm Hom}_{D^b(Y_2,\CC)}({\mathcal G},\,
R\, j_{\mycircle}({\mathcal F}))\risom
{\rm Hom}_{D^b(Y_1,\CC)}(j^\mycircle({\mathcal G}),
\, {\mathcal F})\ ,
$$
 with respect to which
the functor $Rj_\mycircle$ : $D^b(Y_1,\CC)\to D^b(Y_2,\CC)$ is right adjoint
to the functor
$j^\mycircle$ : $D^b(Y_2,\CC)\to D^b(Y_1,\CC)$.

If the map $j: Y_1\to Y_2$ is an open embedding,
then the functor $j_\my!$ is exact, and its extension to derived categories:
$j_\my!$ : $D^b(Y_1,\CC)\to D^b(Y_2,\CC)$  is left
adjoint to the functor of inverse image $j^\mycircle$ :
$D^b(Y_2,\CC) \to D^b(Y_1,\CC)$ .

Let $j: Y_1\to Y_2$ be an open embedding of a
smooth affine variety $Y_1$ to a smooth algebraic variety $Y_2$. Denote by
the same symbol $j: Y_1^{an}\to Y_2^{an}$ the corresponding map of complex
analytic varieties. Then
the functors $Rj_\mycircle$ and $j_!$ map perverse sheaves
to perverse sheaves.
 In this situation  {\it the MacPherson extension} $j_{\my!*}({\mathcal
F})$ of a perverse sheaf ${\mathcal F}$ over $Y_1^{an}$ is defined
as the image of $j_{\my!}{\mathcal F}$ in $ Rj_\mycircle{\mathcal F}$ with
respect to the canonical map \linebreak
$s_j$ : $j_{\my!}({\mathcal F})$
 $\to$ $ Rj_\mycircle({\mathcal F})$, determined by the two pairs
$(j^{\mycircle}, Rj_{\mycircle})$ and $(j_\my!,j^{\mycircle})$
of adjoint functors, see Section \ref{Rdiq}.
The MacPherson extension is a perverse sheaf over
$Y_2^{an}$ whose restriction to $Y_1^{an}$ is isomorphic to
${\mathcal F}$ in the category $D^b(Y_1^{an},\CC)$.

\medskip

\subsection{De Rham functor}

{} Let $Y$ be a smooth algebraic variety  over $\C$. Denote by $Y^{an}$ the same
variety considered as an analytic variety. Denote by
$\mathcal{O}^{an}_Y$ the sheaf of analytic
functions on $Y^{an}$ and by $\D_Y^{an}$ the sheaf of rings of differential
operators on $Y^{an}$ with analytic coefficients, \
$\D_Y^{an}=\mathcal{O}^{an}_Y\ot_{\mathcal O_Y}\D_Y$. We  assign to any
$\D_Y$-module $M$ the corresponding $\D_Y^{an}$-module $M^{an}$,
$$
M^{an}\ =\ \mathcal{O}^{an}_Y\ot_{\mathcal O_Y}M\ .
$$

The sheaf $\Om_Y$  of top exterior forms on $Y$
 has a structure of the right $\D_Y$-module.   Differential
operators of degree zero (functions) act on the top forms by means of point-wise
multiplication and vector fields   act by means of minus Lie
derivative:
\begin{equation}\label{right}
\o\ \cdot\ f\ =\ f\o\ ,
\qquad
\o\ \cdot\ \xi\ =\ - Lie_\xi(\o)\ .
\end{equation}
Here $f\in\mathcal{O}_Y$, $\omega\in\Om_Y$, $\xi$ is a vector field
on $Y$.

Let $M^{}$ be a coherent $\D_Y$ module. Then the  assignment $M$ $\to$
 $\Oman_Y{\ot}_{\D_Y^{an}}M^{an}$ defines a functor from the category
$\MD_Y$ of coherent $\D_Y$-modules to the category $Sh({Y^{an},\CC})$
 of sheaves of complex vector spaces over $Y^{an}$.
This functor is right exact.

The left derived functor of this functor is called the de Rham functor:
\begin{equation}
\label{DR}
DR:\ D^b(\MD_Y)\to D^b(Y^{an},\CC),\qquad
DR(M^\mycircle)=\Oman_Y\stackrel{{\rm L}}{\ot}_{\D_Y^{an}}(M^\mycircle)^{an}
\end{equation}

To calculate the de Rham functor $DR(M^\mycircle)$,
the class of free $\D_Y$-modules can be used: if
$\tilde{M}^\mycircle$ is a complex of free $\D_Y$-modules quasi-isomorphic to
$M^\mycircle$,  then
$$
DR(M^\mycircle)\ =\ \Oman_Y{\ot}_{\D_Y^{an}}(\tilde{M}^\mycircle)^{an}\ .
$$
Another possibility is to replace the right $\D_Y$-module
$\Omega_Y$ by its free $\D_Y$-resolution $\tilde{\Omega}_Y^\mycircle$,
then
$$
DR(M^\mycircle)\ =\ (\tilde{\Omega}_Y^\mycircle)^{an}
{\ot}_{\D_Y^{an}}({M}^\mycircle)^{an}\ .
$$
A free $\D_Y$-resolution $\tilde{\Omega}_Y^\mycircle$ of the right
$\D_Y$-module $\Omega_Y$ looks as follows:
\begin{equation}\label{oresol}
0\to \D_Y\to\Om^1_Y\ot_{\mathcal{O}_Y}\D_Y\to\cdots\to
\Om^{N-1}_Y\!\!\ot_{\mathcal{O}_Y}\D_Y\to\Om^N_Y\ot_{\mathcal{O}_Y}
\D_Y\to\Om^N_Y\to 0
\end{equation}
where $N=\dim\, Y$, $\Om^k_Y$ is the sheaf of exterior $k$-forms over $Y$.
 In local coordinates $z_1, \dots ,z_N$,  the differential
is
\begin{equation}\label{differential}
d(\o\ot D)\ =\ d\o\ot D\ +\
\sum_{i=1}^N \left(dz_i\wedge \o\right)\ot \frac{\partial}
{\partial z_i}D\ .
\end{equation}
Here $\o\in\Om^k_Y$, $D\in \D_Y$.

 Let $M$ be a coherent $\D_Y$-module. Then,
due to \rf{oresol}, the complex of sheaves
\linebreak
$DR(M)=DR(M^{(\cdot)})$
is isomorphic in $\D^b(Y^{an},\CC)$  to the complex
\begin{equation}
\label{oresolM}
0\to M^{an}\to\left(\Om_Y^1\right)^{an}
\!\!\ot_{\Oan_Y}M^{an}\to\cdots\to
\left(\Om_Y^N\right)^{an}\!\!\ot_{\Oan_Y}M^{an}\to 0
\end{equation}
where $\deg\left(\left(\Om_Y^{N-k}\right)^{an}\!\!\ot_{\Oan_Y}M^{an}\right)=
-k$ and the
 differential is  induced from \rf{differential}.
\medskip

{\em Example.\ }
 Let $Y=X^0$ be the complement to the arrangement $\A$ in $X=C^N$.
  For each quiver
 $\VA_{\Ga^0}=\{V, A^i\}\in \Qui_{\G^0}$,
 the  $\D_{X^0}$-module $E^0 \VA_{\Ga^0}$ is a finitely generated free
 ${\cal O}_{X^0}$- module. By the Poincare lemma   the
  cohomology groups of the complex \rf{oresolM} in this case
are  nonzero only in  degree $-N$.
Thus we have an identification in $D^b(Y^{an},\CC)$:
\begin{equation}
DR(E^0\VA_{\Ga^0})
\approx{\mathcal L}\VA_{\Ga^0} [N]^{(\cdot)}
\label{locsys}
\end{equation}
where ${\cal L}\VA_{\G^0}$ is the locally constant sheaf (local system)
of  flat sections of the corresponding  connection \rf{3.12}, and
 $C[N]^\mycircle$ means the shift of degree by
$N$,  $C[N]^k=C^{N+k}$.

\medskip

\subsection{De Rham complex of a quiver $\D_X$-module}

{}For any quiver $\VA\in\Qui$, denote by $\QVA$ the tensor product over
${\Dan_X}$ of the right ${\Dan_X}$-module $ \Om_X^{an}$ and the complex
 of left  $\Dan_X$-modules $\R \EVA^{an}$:
 \begin{equation}\label{QVA1}
\QVA\ =\ \Om_X^{an}\ot_{\Dan_X}\R \EVA^{an}\ .
\end{equation}
 Theorem \ref{theorem3} implies
\begin{corollary}
\label{proposition7}
Let $\VA\in \Qui_\Ga$ be a  quiver of $\Ga$. Then the de Rham complex of the
quiver $\D_X$-module $\EVA$ is isomorphic in $D^b(X^{an},\CC)$ to the complex of sheaves
$\QVA$.
\end{corollary}
Due to the description in \rf{Rr} and \rf{dr} of the resolution
$\Rlong\EVA$, we get the following description of the complex
$\QVA:\  0\to{\mathcal Q}^{-N}\VA\stackrel{d}{\to}
\dots \stackrel{d}{\to}{\mathcal Q}^{0}\VA^0\to 0$
with $\deg {\mathcal Q}^{r}\VA = r$:
\begin{equation}\label{DREVA}
{\mathcal Q}^{r}\VA=\ooplus\limits_{\a\in\v(\Ga)} \Oman_X\ot \left(
\bigwedge\nolimits^{-r}\TTT{\a}\ot\Oa\ot V_\a\right)
\end{equation}
with differential $d=d_1-d_0$, where
\begin{equation}
\begin{array}{ccc}
d_1\left({\o}\ot t_1\wedge \dots \wedge t_r\ot \om_\a\ot
v_\a\right)&=&\!\!\!\!\sum\limits_{i=1}^r (-1)^{i+1} (\o)\cdot
t_i\ot t_1\wedge \dots \widehat{t_{i}} \dots \wedge t_r\ot \om_\a\ot v_\a,
\label{drom}
\\
d_0\left({\o}\ot t_1\wedge \dots \wedge t_r\ot \om_\a\ot
v_\a\right)&=&\!\!\sum\limits_{\b,\ (\a,\b)\in E(\G)} \o\ot
A^r_{\b,\a} (t_1\wedge \dots \wedge t_r\ot \om_\a\ot v_\a),\end{array}
\end{equation}
for $\o\in\Oman_X$, $t_i\in\TTT{\a}$,
 $\om_\a\ot v_\a\in \Oa\ot V_\a$.

Analogously, Theorem \ref{prop11} implies

\begin{corollary}\label{new corollary}
For any
 $\Ga^0$ quiver $\VA_{\Ga^0}=\{V, A^i\}$,   the  de Rham complex
$DR(E^0\VA_{\Ga^0})$ is isomorphic in $D^b(X^{0, an},\CC)$ to the
complex
${\QVA}_{\Ga^0}= \Omega_{X^0}^{an}\ot_{\Dan_{X^0}}\R E^0\VA^{an}_{\G^0}$,
where
$$
{\mathcal Q}^{r}\VA_{\Ga^0}\ =\ \Oman_{X^0}\ot (
\bigwedge\nolimits^{-r} T_{\emptyset}\ot\overline{\Omega}_\emptyset\ot
V)\ ,
$$
and  the  differential is $d=d_1-d_0$,
 $$
\begin{array}{rcl}
d_{1}\left(\omega\ot \xi_1\wedge \dots \wedge \xi_r\ot \om_\emptyset\ot
v\right)&=&\!\!\!\sum\limits_{i=1}^r (-1)^{i+1}(\omega)\cdot\xi_i\ot
\xi_1\wedge\dots \widehat{\xi_{i}} \dots \wedge \xi_r\ot
\om_\emptyset\ot v\ ,
\\
d_{0}\left(\omega\ot \xi_1\wedge \dots \wedge \xi_r\ot \om_\emptyset\ot
v\right)&= &\!\!\!\sum\limits_{j\in J(\A)}\omega \cdot f_j^{-1}\ot
i_{f_j} (\xi_1\wedge\dots\wedge \xi_r)\ot \om_\emptyset\ot
A_j(v)\ .
\end{array}
$$
\end{corollary}

Notice that this complex  is isomorphic in $D^b(X^{0, an},\CC)$
to the complex \rf{oresolM}.

\medskip

\subsection{Quiver perverse sheaves}

It is known \cite{Kash,Meb,B} that for any smooth algebraic variety $Y$, the
 de Rham functor
transforms holonomic $\D_Y$-modules to perverse sheaves. Moreover, if $M$
is  a holonomic $\D_Y$-module such that its  singular support $ss.M$ belongs
to the union of conormal bundles of nonsingular subvarieties $Z_i\subset Y$,
$$ss.M\subset \bigcap_i T^*_{Z_i}Y\ ,$$
then the complex of sheaves $DR(M)$ is a perverse sheaf with respect to any
equisingular stratification of $Y^{an}$, such that each $Z_i^{an}$ is a
union of strata.

Corollary \ref{proposition7} and Proposition \ref{prop1} imply
\begin{proposition}
\label{proposition8}
Let $\VA\in\Qui$ be a quiver of $\,\Ga$. Then the complex of sheaves
 $\QVA$ defined by
 \rf{QVA1}  is a perverse sheaf on $X^{an}$ with respect to the
stratification $\{ X_\a^{an}, \a\in\v(\Ga)\}$.

Analogously, for any ${\Ga^0}$-quiver $\VA_{\Ga^0}\in\Qui_{\Ga^0}$,
 the complex of sheaves $\QVA_{\Ga^0}$ is a perverse sheaf over $X^{0,
 an}$, isomorphic in $D^b(X^{0, an},\CC)$ to the shifted by $N$
 locally constant sheaf ${\cal L}\VA_{\G^0}[N]$, see \rf{3.12}.
\end{proposition}
 \medskip

The de Rham functor respects the duality and all of the functors of direct and
inverse image. In particular,
for a smooth complex algebraic variety $Y$ and
 $M \in \MD^{hol}_Y$,  we have \cite{Kash3, B}
\begin{equation}\label{Ddelta}
DR(\Delta_Y(M))\approx \overline{\Delta}_{Y^{an}}( DR(M))\ ,
\end{equation}
{}For  an  open embedding  $j: Y_1\to Y_2$
of a smooth affine
variety $Y_1$ into a smooth algebraic variety $Y_2$,
 any $M_1\in\MD^{hol}_{Y_1}$ and $M_2\in\MD_{Y_2}$,  we have
\begin{equation} \begin{split}
j^\mycircle (DR(M_2))&\approx DR(j^*(M_2))\ ,\\
Rj_\mycircle (DR(M_1))& \approx DR(j_*(M_1)\ ,\label{Djcom}\\
j_\my! (DR(M_1))&\approx DR(j_\my!^{(D)}(M_1)) \ ,\\
j_{\my!*} (DR(M_1))&\approx DR(j_{\my!*}^{(D)}(M_1))
\end{split}\end{equation}

\begin{corollary}
\label{propdual1} ${}$\\
{\em (i)} Let $\VA_{\G^0}$ be a level zero quiver. Then
we have in $D^b(X^{0,  an},\CC)$:
\begin{equation}\nonumber
\overline{\Delta}_{X^{0,an}}(DR (E^0(\VA_{\G^0})))\approx DR
 (E^0(\tau\VA_{\G^0}))\ .
\end{equation}
{\em (ii)} Let $\VA=\{V_\a, A_{\a,\b}\}\in \Qui$ be a  quiver of $\Ga$.
 Then                we have
in $D^b(X^{an},\CC)$:
\begin{equation}\nonumber
\overline{\Delta}_{X^{an}}(DR (E(\VA)))\approx DR (E(\tau\VA))\ .
\end{equation}
\end{corollary}

{\it Proof.}
 Both statements follow from Theorem \ref{propdual} and
relation \rf{Ddelta}. \hfill{$\square$}

\begin{corollary}\label{corjj!}
 Let $\VA_{\G^0}=\{V, A^i\}$ be a level zero quiver. Suppose that
each operator $A^i$ has a single eigenvalue and
 $\VA_{\G^0}$ has a non-resonant
spectrum, or suppose that
all operators $A^i$ are close to zero.  Then the
following perverse sheaves are isomorphic:
\begin{eqnarray}
R(j_{0})_{\mycircle}DR (E^0(\VA_{\G^0}))
\approx DR  (E\left(\j_{0,*}(\VA_{\G^0})\right))
\ ,\nonumber\\
(j_{0})_\my!DR
(E^0(\VA_{\G^0}))
\approx DR  (E\left(\j_{0,!}(\VA_{\G^0})\right))\ ,
\nonumber\\
(j_{0})_{\my!*}DR (E^0(\VA_{\G^0}))
\approx DR  (E\left(\j_{0,!*}(\VA_{\G^0})\right))\ .\nonumber
\end{eqnarray}
\end{corollary}
{\it Proof.} The first statement follows from
Corollary \ref{cor2} and  relation \rf{Djcom}.
The second and third statements follow from Corollary \ref{corj!}
and relation \rf{Djcom}.
\hfill{$\square$}

\medskip

\subsection {Cohomology groups of $\CC^N$ with coefficients in quiver
 perverse sheaves} \label{section4.5} ${}$
We consider $X=\CC^N$ as an  analytic variety.
For a quiver $\VA$,
 the complex $\QVA$ is a complex of free ${\cal
O}_X$-modules. The affine space $X=\C^N$ is a Stein space,
$$
H^k(X, \Oan)\ =\ 0, \qquad {\rm for \ all}\quad k\not= 0\ .
$$
Hence  the (hyper)cohomology groups of $X$ with coefficients
 in the complex $\QVA$ of sheaves
are  isomorphic to the cohomology groups of the complex $\overline{\QVA}=
\G(X, \QVA)$
 of global sections  of  $\QVA$:
\begin{equation}\label{stein}
H^k(X, \QVA)\ \simeq \ H^k(\overline{\QVA}) \ .
\end{equation}

\begin{theorem}
\label{theorem4}
Let $\A$ be a central arrangement
(that is all hyperplanes of $\A$ contain  \linebreak
$0\in \CC^N$).
Let $\VA\in \Qui_\G$ be a  quiver with  all  linear maps $A_{\a,\b}$ of the quiver
close to zero.
Then the complex $\overline{\QVA}$ of global sections of the complex
\rf{DREVA}-\rf{drom} is quasi-isomorphic to the shifted by $N$
 quiver complex $C_+(\VA)$,  introduced in
   \rf{complex}-\rf{differentials}:
\begin{equation}
\nonumber
\overline{\QVA}\  {} \,
\stackrel{quis}{\simeq}\ {}\, C_+(\VA)[N]    \ .
\end{equation}
\end{theorem}

Note that  $C_+(\VA)$ is a finite dimensional complex.

The theorem is proved in Section \ref{proofs}.

Theorem \ref{theorem4} and the isomorphism in  \rf{stein} imply the following
statement.
\begin{corollary}
\label{cor3} Let $\A$ be a central arrangement and $\VA\in \Qui_\G$ a
quiver with all linear maps $A_{\a, \b}$ close to zero. Then
for any $k\in \ZZ$,
 \begin{equation}\nonumber
 H^k(X, DR(E\VA))\ \simeq \ H^k(X, \QVA)\ \simeq\ H^{k+N}(C_+(\VA))\ .
 \end{equation}
\end{corollary}

Let $X^0$ denote  the complement to  $\A$ in $\C^N$.
For a local system ${\mathcal L}$ on $X^0$
introduce the intersection cohomology groups $I\!H^k(X^{},\mathcal{L})$ as
the cohomology groups $H^{k-N}(X,(j_{0})_{\my!*}\mathcal{L}[N])$ of the MacPherson
extension of the perverse sheaf $\mathcal{L}[N]\in D^b(X^0,\CC)$,
\begin{equation}\label{IH}
I\!H^k(X^{},\mathcal{L})\stackrel{def}{=}
H^{k-N}(X,(j_{0})_{\my!*}\mathcal{L}[N])\ .
\end{equation}

 Let ${\cal L}\VA$ denote  the local system of  flat
sections of the connection \rf{3.12}.

\begin{corollary}
\label{theorem5} Let $\VA=\{V,A^i\}$ be a level zero quiver of
 a central arrangement $\A$.
 Assume that all operators $A^i$ are close to zero.
 Then for any $k\in\ZZ$,
\begin{eqnarray}\nonumber
H^k(X^0,{\cal L}\VA)&\simeq& H^k(C_+(\j_{0,*}\VA))\ ,\\
\nonumber
I\!H^k_{}(X^{},{\mathcal L}\VA)&\simeq&
H^{k}(C_+(\j_{0,!*}\VA))\ .
\end{eqnarray}
\end{corollary}
{\it Proof.} The conditions on operators $A^i$ imply
 that the quiver $\VA$ has a non-resonant spectrum and Corollary
\ref{corjj!} can be applied to the quiver $\j_{0,!*}\VA$. By \rf{locsys},
${\mathcal L}\VA\approx DR(E^0\VA)[-N]$ and
\begin{equation}
\begin{split}\notag
H^k(X^0,{\mathcal L}\VA)& \simeq H^k(X, R(j_0)_{\mycircle}{\mathcal L}\VA)
\ \simeq\
H^k(X, DR(\j_{0,*}\VA)[-N])\ \simeq\
H^k(C_+(\j_{0,*}\VA))\ ,\\ \notag
I\!H^k(X^{},{\mathcal L}\VA)&\simeq
H^{k-N}(X, (j_0)_{\my!*}{\mathcal L}\VA[N])
\ \simeq\
H^{k-N}(X, DR(\j_{0,!*}\VA)) \ \simeq\
H^k(C_+(\j_{0,!*}\VA))\ .
\end{split}
\end{equation}
The second isomorphism in both lines is due to {\rm Corollary} \ref{corjj!},
the third isomorphism in both lines is due to {\rm Corollary} \ref{cor3}.
\hfill{$\square$}

\medskip

{\bf Example.} Let $\A$ be a central arrangement in $X$.
Let $a_j$, $ j\in \v_1(\A)$, and $\kappa$ be complex numbers.
Consider the trivial bundle  $\C\times X^0 \to X^0$ with the connection
given by the differential one-form
$
 \frac 1 \kappa \,  \sum_i\ a_i\ \frac{d\,  f_i}{f_i}.
$
Here $f_i=0$ is the equation of the hyperplane $H_i$ of  $\A$.
Let ${\cal L}\VA$ be  the local system of  flat
sections of the connection.

Let $\VA=\{\CC,A^i\}$ be the one-dimensional level zero quiver of
 $\Ga^0$, where $A^i:\C\to\CC$ is the operator of multiplication by
$a_i/\kappa$.

Let  $(\OS^\upcircle(\A), d_a)$ be the Aomoto complex, defined for $\A$
in Section \ref{Orlik}. Let $(S_a(\FF^\upcircle(\A)), d_a)$ be
its subcomplex of flag forms.

By Corollary \ref{theorem5}, for $|\kappa|\gg 1$, we have
$$
H^\upcircle(\OS^\upcircle(\A), d_a)\ \simeq \ H^\upcircle(X^0,{\cal L}\VA)\ ,
\qquad
H^{\upcircle} (S_a(\FF^\upcircle(\A)), d_a)\ \simeq\
I\!H^{\upcircle}_{}(X^{},{\cal L}) \ .
$$
Here $H^\upcircle(X^0,{\cal L}\VA)$ are the cohomology groups of
$X^0$ with coefficients in the local system ${\cal L}\VA$ and
$I\!H^\upcircle_{}(X^{},{\cal L}\VA)$ are the
intersection cohomology groups of
$X^0$ with coefficients in the local system ${\cal L}\VA$.

\medskip
We see that the calculation of cohomology groups of $X$ with
coefficients in the perverse sheaf associated with  the quiver MacPherson
extension of a level zero quiver $\VA_{\G^0}$ gives a topological meaning to
the complex of the flag forms. One of the aims of this paper was to give
a topological meaning to the complex of flag forms.

\medskip

\setcounter{equation}{0}
\section{Equivariant structures}
\subsection{Arrangement with a group action}
Suppose a finite group $G$ acts on  $X = \C^N$ by  affine
transformations $\TG{X}$ : $X \to X$ preserving an arrangement $\A$.
This means that $G$ acts on the set $J(\A)$ so that
$\TG{X}(H_j) = H_{g(j)}$  for $g\in G, \ j\in J(\A)$.

The action of $G$  preserves  the principal open subsets $X^n$.

The group $G$ naturally acts on the
graph $\G$ of the arrangement and on the truncated
graphs $\G^n$.

The $G$-action on $X$ induces a $G$-action on the ring $D_X$. For an affine
transformation $L : X \to X$, the associated
automorphism $L_{D_X}$ : ${D_X}\ \to \ D_X$ is defined by
 the formula  $L_{D_X}(f)(z) = f(L^{-1}z)$
for a function $f$ on $X$, and by the formula
$L_{D_X}(\eta)(z) = {\rm d} L \cdot \eta(L^{-1}z)$
for a vector field $\eta$. Here ${\rm d} L$ is the derivative of $L$.

For $g\in G$, the transformation  $\TG{D_X} : D_X \to D_X$
preserves the space $F$ of affine functions
on $X$ and for  $\a\in\v(\G)$, maps the subspace $\Fa\subset F$ to
 $F_{g(\a)}$.

The transformation  $\TG{D_X}$ preserves the space $T$ of vector fields
on $X$ with constant coefficients
 and for any $\a\in\v(\G)$,\ maps the subspace $\Ta\subset T$ to $T_{g(\a)}$.

The $G$-action on $X$ induces a $G$-action on the algebra
$\Omega_X^*$ of differential forms.
For an affine transformation $L : X \to X$, the associated
automorphism $L_{\Omega^*_X} : \Omega_X^* \to \Omega_X^*$
is defined by  the formula $L_{\Omega^*_X}(\omega)(z) =
(({\rm d} L)^t)^{-1}(\omega)(L^{-1}(z))$ for $w\in \Omega^1_X$.

For $g\in G$ and $\a \in I(\G)$, the transformation
$L_{\Omega^*_X}(g)$ maps ${\overline{\Omega}}_\a$ to
${\overline{\Omega}}_{g^{}(\a)}$.

\subsection{Equivariant quivers}

An equivariant quiver (or a quiver with a group action)
\linebreak
 $\VA_G=\{V_\a, A_{\a,\b},\TG{V_\a}\}$  is a quiver
$\VA=\{V_\a, A_{\a,\b}\}$ together with a collection of linear maps
${\TG{V_\a}}:V_\a\to V_{g(\a)}$, for $\a\in\v(\Ga)$ and $g\in G$, such
that  for  $\a,\b\in\v(\Ga)$, we have
\begin{enumerate}
\item[$\bullet$]
$\TXG{V_{g''(\a)}}{g'}\,\TXG{V_\a}{g''}\ =\ \TXG{V_\a}{g'g''}$\,,
\item[$\bullet$]
$\TXG{V_\a}{g^{-1}}\ =\ (\TXG{V_{g^{-1}(\a)}} {g} )^{-1}$\,,
 \item[$\bullet$]
${\TG{V_\a}}\,A_{\a,\b}\ =\ A_{g(\a),g(\b)}\,{\TG{V_\b}}$\,.
\end{enumerate}

In the same manner we define level $l$ $G$-equivariant quivers.

Denote by $Qui_{G}=Qui_{\G,G}$ the category of $G$-equivariant
quivers, by $Qui_{\G^l,G}$ the category of level $l$ $G$-equivariant
quivers and by $\forget^l: \Qui_{\G^l,G}\to \Qui_{\G^l}$ the forgetful
functor, which assigns to an equivariant quiver $\VA_{\G^l,G}$ the
same quiver with no morphisms of group action. In this case we
will use notation
$\VA_{\G^l}$ for the quiver $\forget^l(\VA_{\G^l,G})$.

Let  $V_{\G^l,G}=\{V_\a, A_{\a,\b}, A_\a^\b, \TG{V_\a}\}$ be a level $l$
equivariant quiver. Then for  $k < l$, the
 collection of maps $\TG{V_\a}: V_\a\to V_\a$,
 for $g\in G$ and $\a\in\v(\G^k)$, determines a group action on the
level $k$ quiver $\j_{l,k}^*(V_{\G^l})$.

This assignment determines a functor  $\j_{l,k,G}^{*}$ :
$\Qui_{\G^l,G}$ $\to$ $\Qui_{\G^k,G}$. This functor
commutes with forgetful functors,
$$
\forget^k\ \j_{l,k,G}^{*}\ =\ \j_{l,k}^{*}\ \forget^l\ .
$$

The functor $\j_{l,k,G}^{*}:\ \Qui_{\G^l,G}\to \Qui_{\G^k,G}$ admits left
and right adjoint functors. Let
$\j_{l,k,*,G}^{}:\ \Qui_{G^k,G}\to \Qui_{G^k,G}$ be the right adjoint
functor to  $\j_{l,k,G}^{*}$. Let
$\j_{l,k,!,G}^{}:\ \Qui_{G^k,G}\to \Qui_{G^k,G}$  be the left adjoint
functor to $\j_{l,k,G}^{*}$. We have
$$
\forget^l\ \j_{l,k,*,G}^{}\ =\ \j_{l,k,*}^{}\ \forget^k\ ,
\qquad
\forget^l\ \j_{l,k,!,G}^{}\ =\ \j_{l,k,!}^{}\ \forget^k\ .
$$

\subsection{Equivariant quiver $\D$-modules} Let $X$ be a topological space.
Any homeomorphism $a: X\to X$ determines an exact functor
  $a^{\inverse}:Sh(X,\CC)\to Sh(X,\CC)$
of inverse image of sheaves of vector spaces,
 which can be described as follows.
  The fiber $a^{\inverse}{\cal F}_x$ of the sheaf
$a^{\inverse}{\cal F}$ at a point $x\in X$ is the fiber ${\cal F}_{a(x)}$
of ${\mathcal F}$ at the point $a(x)$, and for any morphism \linebreak
$\varphi$ :
${\mathcal F}$ $\to$ ${\mathcal G }$ of sheaves over $X$ the map of fibers
$a^{\inverse}\varphi(x)$ $:$ $ a^{\inverse}{\mathcal F}_x$ $\to$ $
a^{\inverse}{\cal G}_x$ coincides with the map $ {\varphi}(a(x))$ $:$
 $ {\mathcal F}_{a(x)}$ $\to$ $ {\cal G}_{a(x)}$.

Suppose a finite group $G$ acts on $X$ by homeomorphisms $L_X(g)$.
 For  ${\mathcal F}\in
Sh(X,\CC)$ and  $g\in G$, denote by  ${\mathcal F}^g$ the sheaf
$(L_X(g)^{-1})^\inverse {{\mathcal F}}$.
 A sheaf
${\mathcal F} \in Sh(X, \C)$  is called $G$-equivariant, if
for $g \in G$, a  morphism
$L_{{\mathcal F}}(g) : {{\mathcal F}}\to {\mathcal F}^g$
is given such that
\begin{enumerate}
\item[$\bullet$]
$\TXG{{{\mathcal F}^{g_2}}}{g_1}\  \TXG{\mathcal F}{g_2}\ =\
\TXG{{\mathcal F}}{g_1g_2}$ ,
\item[$\bullet$]
$\TXG{{\mathcal F}^g}{g^{-1}} \ \TXG{{\mathcal F}}{g}\ =\
Id_{\mathcal F}.$
\end{enumerate}
Here the morphism of sheaves $\TXG{{{\mathcal F}^{g_2}}}{g_1}$ : ${\mathcal F}^{g_2}$
$\to$ $({\mathcal F}^{g_2})^{g_1}$ is defined by the relation
$\TXG{{{\mathcal F}^{g_2}}}{g_1} =$ $ \TXG{X}{g_2^{-1}}^{\inverse}
(\TXG{{\mathcal F}}{g_2^{-1}g_1g_2})\ .$

In other words, an equivariant sheaf  is a sheaf  ${\mathcal F}\in
Sh(X,\CC)$ and a collection of linear maps
 $\TXG{\{{\mathcal F},\ U\}}{g}$ :
$\Ga(U,\ {\mathcal F})$ $\to$ $\Ga(\TXG{X}{g}(U),\ {\mathcal F})$,
attached to any open set $U\subset X$ and an element  $g\in G$, such
that
\begin{equation}
\begin{split}\nonumber
\TXG{\{{\mathcal F,\ \TXG{X}{g_2}(U)}\}}{g_1}\ \TXG{\{{\mathcal F},\ U\}}{g_2}&=
\TXG{\{{\mathcal F},\ U\}}{g_1g_2}\ ,
\\
\TXG{\{{\mathcal F},\ \TXG{X}{g}(U)\}}{g^{-1}}\ \TXG{\{{\mathcal F},\ U\}}{g}&=
Id_{\G(U,\ \mathcal{F})}\ .
\end{split}
\end{equation}

One defines similarly a $G$-equivariant complex of sheaves of vector spaces.

Denote by $Sh(X,G,\CC)$ the category of equivariant sheaves over $X$ and by
$\DB(X,G,\CC)$ the corresponding bounded derived category.
\bigskip

Assume that a finite group $G$ acts on
a smooth complex algebraic variety $X$ by  algebraic automorphisms $L_X(g)$.
Then the sheaves ${\mathcal O}_X$, $\Omega_X^*$, and $\D_X$ are  equivariant sheaves.
The morphisms   $\TXG{\{{\D_X},\ U\}}{g}$ and $\TXG{\{{\mathcal O}_X,\ U\}}{g}$
are defined by the relations
\begin{equation}
\begin{split}\nonumber
\TXG{\{{\D_X},\ U\}}{g}(f)(x)&=
\TXG{\{{{\mathcal O}_X},\ U\}}{g}(f)(x)\  =\ f(\TXG{X}{g^{-1}}(x))\ ,
\\
\TXG{\{{\D_X},\ U\}}{g}(\eta)(x) & = {\rm d} \TXG{X}{g}(
\eta)(\TXG{X}{g^{-1}}(x))\ ,
\end{split}
\end{equation}
where $f$ is a  function and $\eta$ is a vector field on $U$.
The morphisms $\TXG{\{{\Omega}_X^*,\ U\}}{g}$ are defined by the relations
\begin{equation}\label{gOmega}
\TXG{\{{\Omega}_X^*,\ U\}}{g}(\omega)(x) \ =
\ (({\rm d} \TXG{X}{g})^t)^{-1}
(\omega)(\TXG{X}{g^{-1}}(x))
\end{equation}
for one-forms $\omega$.

A sheaf $M$ of $\D_X$-modules is called $G$-equivariant
if it is equipped with a structure of an equivariant sheaf, such that the maps
$\TXG{\{{M},\ U\}}{g}$ respect the structure of the $\D_X$-module. This means that
$$
\TXG{\{{M},\ U\}}{g}(d\cdot m)\ = \ \TXG{\{{\D_X},\ U\}}{g}(d)\cdot
\TXG{\{{M},\ U\}}{g}(m)\ ,
$$
for any $m\in \G(U, M)$ and $d\in\G(U, \D_X)$.
In particular, the space $\G(X,M)$ of global sections is an
equivariant module over the ring $D_X = \G(X,\D_X)$. Denote by
 $\MD_{X,G}$ the
category of equivariant coherent $\D_Y$-modules.

If $X$ is an affine space and $M$ is a coherent $\D_X$-module, then the
equivariant $D_X$-module of its global sections determines $M$ as an
equivariant $\D_X$-module.

The sheaf ${\mathcal O}_X$
is an equivariant left $\D_X$-module. The sheaf $\Omega_X$
of top exterior forms is an equivariant sheaf of right $\D_X$-modules.

\bigskip

Let $\VA_G=\{V_\a, A_{\a,\b},\TG{\VA}\}\in \Qui_{\G,G}$
be a $G$-equivariant quiver of $\G$. Denote by $\VA=\forget^N(\VA_G)$ the
same quiver with no group action.
We equip the $\D_X$-module $E\VA$ with a structure of a $G$-equivariant
$\D_X$-module in the following way.

{}Choose an edge framing $\{f_{\a,\b}, \  \xi_{\a,\b}\}$
of the arrangement $\A$.
Then the
$\D_X$-module $E\VA$ is defined as the quotient of
 the left free $\D_X$-module
$\R^0E\VA=\oplus_{\a\in\v(\G)}\D_X\ot\Oa\otimes V_\a$ by the
submodule  $\R^0E\VA'$, whose global sections are given by
 \rf{2.3} and \rf{2.4}.
 For $g\in G$, define a map $\TG{\R^0E\VA} : \R^0E\VA\to
(\R^0E\VA)^g$ by the formula
\begin{equation}\label{TGE}
\TXG{\{\R^0E\VA,X\}}{g}(d \ot \om_\a\ot v_\a)=\TG{D_X}(d)\ot
\TG{\overline{\Omega}_\a}(\om_\a)\ot\TG{V_\a}(v_a)
\end{equation}
for   $\a\in\v(\G)$   and
 $d\ot\om_\a\ot v_\a\in D_X\ot \Oa\ot V_\a$.

\begin{proposition}\label{Gprop1}
${}$

\begin{enumerate}
\item[ (i)]
 The maps (\ref{TGE}) determine a $G$-equivariant structure on the free
$\D_X$-module $\R^0E\VA$.
\item[ (ii)]
 The subsheaf $\R^0E\VA'$ is
invariant with respect to the maps $\TG{\R^0E\VA}$.
\item[ (iii)]
 Denote by $\TG{E\VA}$ the maps, induced by the action of the
maps  $\TG{\R^0E\VA}$ in the quotient space $E\VA$. Then the maps
$\TG{E\VA} : E\VA\to E\VA^g$ determine a $G$-equivariant structure on the quiver
$\D_X$-module $E\VA$.

\end{enumerate}

\hfill{$\square$}
\end{proposition}

By Proposition \ref{Gprop1}, we have a functor $E_G : \Qui_G\to
\MD_{X,G}$ from the category of equivariant quivers to the category of
equivariant $\D_X$-modules.

\bigskip

 For a nonnegative $n$,
let $\VA_{\G^n,G}=\{V_\a, A_{\a,\b}, A_\a^\b, \TG{\VA}\}\in \Qui_{\G,G}$
be a level $n$ $G$-equivariant quiver. Denote by
$\VA_{\G^n}=\forget^n(\VA_{\G^n,G})$ the
same quiver with no  group action.
We equip the sheaf $E^n\VA_{\G^n}$ of $\D_{X^n}$-modules
 with a structure of a $G$-equivariant sheaf of $\D_{X^n}$-modules
as follows.

Let  ${\mathcal F}_n =\{{f}_\a\ | \ \a\in\v_{n+1}(\G)\}$
be a level $n$ vertex framing and $U_{\FF_n}$
the corresponding open subset of $X^n$, see Section
\ref{Sec.Xn}. The   $\D_{X^n}$-module  $E^n\VA_{\G^n}$ is defined
as the quotient
of the free  sheaf $\R^0 E^n\VA_{\G^n}$ $=$
$\D_{X^n}\otimes_\C(\oplus_{\a\in \v(\G^n)} \Oa\ot V_\a)$ of
 $\D_{X^n}$-modules
by its subsheaf $\R' E^n\VA_{\G^n}$ of coherent  $\D_{X^n}$-modules,
where the sections $\G(U_{\FF_n}, \R' E^n\VA_{\G^n})$
of the sheaf $\R' E^n\VA_{\G^n}$ over  $U_{\FF_n}$ are given by
formulas \rf{3.7}--\rf{3.9}.

For $g\in G$,  define a morphism  $\TG{\R^0E^n\VA_{\G^n}}$ :
$\R^0E^n\VA_{\G^n}$ $\to$ $(\R^0E^n\VA_{\G^n})^g$
by the formula
\begin{equation}
\TG{\{\R^0E^n\VA_{\G^n},X^n\}}
(d \ot \om_\a\ot v_\a)\ =\ \TG{D_{X^n}}(d)\ot
\TG{{\Omega}_X^*}(\om_\a)\ot\TG{V_\a}(v_a) \
\label{TGEn}
\end{equation}
for any  $\a\in\v(\G^n)$,  and
 $d\ot\om_\a\ot v_\a\in D_{X^n}\ot \Oa\ot V_\a$.

The transformation $\TXG{X}{g}$ transforms  a level $n$ vertex framing
$\FF_n$ to the level $n$ vertex framing $\FF_n^g=\{f_{\a}^g\}$, where
$f_{\a}^g(x) = f_{g(\a)}(g^{-1}x)$.
The transformation $\TXG{X}{g}$ maps
 the subset $U_{\FF_n}$
to the subset $U_{\FF_n^g}$. Formula \rf{TGEn}
induces the maps
\begin{equation}\label{TGEn-}
\TG{\{\R^0E^n\VA_{\G^n},\ U_{\FF_n}\}} : \G(U_{\FF_n},\R^0E^n\VA_{\G^n})
\to \G(U_{\FF_n}^g,\R^0E^n\VA_{\G^n})\ .
\end{equation}

\medskip
{\bf Example.}

Let $\VA_{\G^0,G}=\{V,A^j,\TG{V}\}\in\Qui_{\G^0,G}$ be an equivariant
 quiver of level 0.
It means, in particular, that linear maps $A^j : V \to V$, as well as linear maps
$A^j_g= \TG{V}A^j\TG{V}^{-1}$
for each $\a\in\v_2(\G)$, $j\succ\a$,
satisfy the relations
$$
[ A^j , \sum_{i,\ \, i\succ\a} A^{i}]\ =\  0\ ,\qquad
[ A^j_g , \sum_{i,\ \, i\succ\a} A^{i}_g]\ =\  0\ .
$$
The associated equivariant $\D_{X_0}$-module   $E^{0}_G\VA_{\Ga^0,G}$
 is generated by the vectors $\o_\ee\ot v$, where
$\o_\ee = {\rm d}z_1\wedge
\ldots \wedge {\rm d}z_N$ and  $v\in V$,
 subject to the relations \rf{3.133}. The group action is defined by the
relation
$$\TG{E^{0}_G\VA_{\Ga^0,G}}(\o\ot v)=\TG{\Omega_{X^0}^*}(\o)\ot \TG{V}(v)
\ .$$
This module is isomorphic to the
 $\D_{X^0}$-module, associated to the  flat connection \rf{3.12}
in the trivial bundle on $X^0$ with fiber $V$.
The group $G$ acts in the space of the global sections of the
trivial bundle ${\rm L}$ by the relation
\begin{equation}
\label{5.5}
\TG{{\rm L}}(f\ot v)=\TG{{\mathcal O}_{X^0}}(f)\ot {\det}_G(g)\cdot \TG{V}( v)\ .
\end{equation}
Here $f\in {\mathcal O}_{X^0}$, $v\in V$ and
$\det_G$ is the one--dimensional representation
 of the group $G$ in the space
$\overline{\Omega}_{\ee}$.
\medskip

\begin{proposition}
\label{Gprop2}${}$
\begin{enumerate}
\item[ (i)] The maps \rf{TGEn-} determine an equivariant
structure on  the sheaf $\R^0E^n\VA_{\G^n}$ of
\linebreak
$\D_{X^n}$-modules.

\item[(ii)]
For $g\in G$, the map $\ \TG{\{\R^0E^n\VA_{\G^n},\ U_{\FF_n}\}}\ $
maps the subspace $\ \G(U_{\FF_n},\R' E^n\VA_{\G^n})$ $\subset$
$\G(U_{\FF_n},\R^0E^n\VA_{\G^n})$ to the subspace
$\G(U_{\FF_n}^g,\R' E^n\VA_{\G^n})$ $\subset$
 $\G(U_{\FF_n}^g,\R^0E^n\VA_{\G^n})$.

\item [ (iii)] Denote by $\TG{E^n\VA_{\G^n}}$ the maps, induced by the action
of   $\TG{\R^0E^n\VA_{\G^n}}$ on the quotient sheaf $E^n\VA_{\G^n}$.
 Then the maps $\TG{E^n\VA_{\G^n}}$
 determine  an equivariant structure on  $E^n\VA_{\G^n}$.

\end{enumerate}
${}$
\hfill{$\square$}
\end{proposition}

By Proposition \ref{Gprop2}, we have a functor $E^{n}_G : \Qui_{\G^n,G}
\to\MD_{X^n,G}$ from the category of level $n$ equivariant quivers to the
 category of equivariant sheaves of $\D_{X^n}$-modules.

Since the map $j_{l,k} : X^k\to X^l $ is an embedding of $G$-invariant
open subsets of $X$, the functor $j^*_{l,k}$ of inverse image of $D$-modules
admits the canonical extension to the functor $j^{*}_{l,k,G} : \MD_{X^l,G}\to
\MD_{X^k,G} $ of inverse image of equivariant $\D$-modules. We
define the functors $j_{l,k,*,G} : \MD_{X^k,G}\to
\MD_{X^l,G}$ and  $j_{l,k,*,G} : \MD_{X^l,G}\to
\MD_{X^l,G}$ as right and left adjoint respectively
to the functor $j^{*}_{l,k,G}$.
 We have forgetful functors $\forget^n :
\MD_{X^n,G}\to \MD_{X^n}$ and the functors of direct
and inverse image of equivariant $\D$-modules commute with the forgetful
functors.

Equivariant versions of Theorems  \ref{proposition2} and \ref{prop4}
are valid. For any $k<l$ and an equivariant level $l$ quiver
$\VA_{G^k,G}$, we have an isomorphism
$$
E^{k}_G ( \j_{l,k,G}^{*}(\VA_{G^k,G}))\ \simeq\ j_{l,k,G}^{*}
( E^{l}_G(\VA_{G^k,G}))\ .
$$
{}For any $k<l$ and an equivariant level $k$ quiver
$\VA_{G^k,G}$, such that the quiver direct image
$\j_{l,k,*}(\forget^k(\VA_{G^k,G}))$ is strongly non-resonant,
we have an isomorphism
$$j_{l,k,*}^G ( E^{k}_G (\VA_{G^k,G}))\ \simeq
\  E^{l}_G ( \j_{l,k,*,G}(\VA_{G^k,G}))\ .
$$

Let $\VA_G=\{V_\a, A_{\a,\b},\TG{\VA}\}\in \Qui_{\G,G}$
be an equivariant quiver  and $\VA=\forget^N(\VA_G)$.
Define an equivariant structure on the free $\D_X$-module
$\R\EVA=\oplus_{r=0}^N\R^r\EVA$, where
$$
 \R^r\EVA = \ooplus\limits_{\a\in\v(\Ga)} \D_X\ot
\bigwedge\nolimits^{-r}\TTT{\a}\ot \Oa\ot V_\a\ ,
$$
see \rf{Rr}. We define the map
 $\TG{\R^rE\VA} : \R^rE\VA\to
(\R^rE\VA)^g$ by the formula
\begin{equation*}
\begin{split}
\TG{\{\R^rE\VA,\  X\}}
(d\ot\overline{\tau} \ot \om_\a\ot v_\a)=
\hspace{7cm}
\\
=\TG{D_X}(d)\ot
\bigwedge\nolimits^{-r}\!\!\TG{D_X}
(\overline{\tau})\ot
\TG{{\Omega}_X^*}(\om_\a)\ot\TG{V_\a}(v_a) \ ,
\end{split}
\end{equation*}
for   $\a\in\v(\G^n)$,  and
 $d\ot\overline{\tau}\ot\om_\a\ot v_\a\in D_X\ot
 \bigwedge\nolimits^{-r}\TTT{\a}\ot \Oa\ot V_\a$.

\begin{proposition}\label{Gprop3}
${}$
\begin{enumerate}
\item[ (i)]
 The maps $\TG{\R^rE\VA}$ determine an equivariant structure on
the $\D_X$-module $\R^rE\VA$ and commute with the differentials $d_0$ and
$d_1$ defined in  \rf{dr},
$$
d_i\ \TG{\R^rE\VA}\ =\ \TG{\R^{r+1}E\VA}\ d_i\ ,
\qquad i = 1 , 2\ .
$$
\item[ (ii)]
 The maps $\TG{\R^rE\VA}$ and $\TG{E\VA}$,
determine an equivariant structure on the complexes
$\R E\VA$ and $\overline{\R}E\VA$ of sheaves of $\D_X$-modules.

\end{enumerate}
\hfill{$\square$}
\end{proposition}

Denote by $\R E\VA_G$ and $\overline{\R}E\VA_G$
the complexes $\R E\VA$ and $\overline{\R}E\VA$ equipped with the
equivariant structures described above.

Proposition \ref{Gprop3} implies that  for any equivariant quiver
$\VA_G$, the complex  $\overline{\R}E\VA_G$
 presents a free resolution of the equivariant
$\D_X$-module $E_G (\VA_G)$ in the category of equivariant $\D_X$-modules.

In particular,
equivariant counterparts of Theorem \ref{propdual} and Corollary
\ref{corj!} are valid. Namely, for any level zero equivariant
quiver $\VA_{G^0,G} \,=\,\{V, A^i,\TG{V}\}$ with non-resonant
spectrum the
following equivariant $\D_{X}$-modules are isomorphic :
\begin{gather}\nonumber
{}j_{0,!,G} ( E^{0}_{G}(\VA_{\G^0,G}))\ \approx \
E_G ( \j_{0,!,G}(\VA_{\G^0,G}))\ ,
\\
j_{0,!*,G} ( E^{0}_{G}(\VA_{\G^0,G}))\ \approx\
 E_G ( \j_{0,!*,G}(\VA_{\G^0,G})) \ . \nonumber
\end{gather}

\subsection{Equivariant quiver perverse sheaves}
Let $X$ be a complex algebraic variety equipped with an algebraic action of
a finite group $G$. Then the sheaf $\Omega_X^{an}$  of
top exterior analytic forms on $X$ is the right equivariant
$\D_X$-module with respect to the action  \rf{right} and equivariant structure
\rf{gOmega}. For any $M\in \MD_{X,G}$ the diagonal action of $G$ on
 sections of the sheaf
 $\ \Oman_X\stackrel{}{\ot}_{\D_X^{an}}(M)^{an}$
equips this sheaf with an equivariant structure. The assignment
$M\mapsto  \ \Oman_X\stackrel{}{\ot}_{\D_X^{an}}(M)^{an}$ defines
a functor from the category $\MD_{Y,G}$ to the category
$Sh(X^{an},G,\CC)$.
The left derived functor to this functor is called the equivariant de Rham
functor $ DR_G:\ D^b(\MD_{X,G})\to \DBG(X^{an},G,\CC)$.

In the remaining part of this subsection, $X$ is the affine space $\CC^N$, considered
 as an  analytic variety equipped with an affine
action of a finite group $G$.

Let  $\VA_G$ be an equivariant quiver of $\Ga$.  The diagonal action
 of $G$ equips the complex
 $$
\QVA\ =\ \Om_X^{an}\ot_{\Dan_X}\R \EVA^{an}
$$
with a structure of an equivariant complex of sheaves  over $X^{an}$.
We describe this equivariant structure for global sections of the sheaf
$\QVA$. The graded component $\QVA^r$ of the complex $\QVA$
was described in \rf{DREVA} as
$$
\QVA^r=\ooplus\limits_{\a\in\v(\Ga)} \Oman_X\ot \left(
\bigwedge\nolimits^{-r}\TTT{\a}
\ot\Oa\ot V_\a\right)
$$
and the  equivariant structure is determined by the maps
\begin{equation}
\begin{split}
\TG{\{ \QVA^r, \ X \}}
(\omega\ot\overline{\tau} \ot \om_\a\ot v_\a)\ =
\hspace{7cm}
\\
=\ \TG{\{ \Omega_X, \ X\}}(\omega)\ot
\bigwedge\nolimits^{-r}\!\!\TG{D_X}
(\overline{\tau})\ot
\TG{{\Omega}_X^*}(\om_\a)\ot\TG{V_\a}(v_a) \ ,
\label{TGEo}
\end{split}
\end{equation}
for any  $\a\in\v(\G^n)$,  and
 $\omega\ot\overline{\tau}\ot\om_\a\ot v_\a\in \G(X,\Omega_X)\ot
 \bigwedge\nolimits^{-r}\TTT{\a}\ot \Oa\ot V_\a$.
We have
\begin{theorem}
\label{Gprop4}
Let $\VA_G$ be an equivariant quiver. Then the equivariant de Rham complex
of the equivariant quiver $\D_X$-module $E_G\VA$ is isomorphic
in $\DBG(X^{an},G,\CC)$
 to the perverse sheaf
$\QVA$ equipped with the equivariant structure \rf{TGEo}.
\hfill{$\square$}
\end{theorem}

Denote by $Y$ the quotient space $Y=X/G$ and by $q:X\to Y$ the quotient
map. We regard $Y$ as a space with the trivial action of $G$. We also use
the notation $Y^0$ for the quotient space $Y^0=X^0/G$ with the quotient map
$q_0:X^0\to Y^0$. The open embedding $Y^0\to Y$ is denoted by $j_0:Y^0\to
Y$.

Let ${\mathcal F}\in Sh(X,G,\CC)$ be a $G$-equivariant sheaf of complex
vector spaces over $X$. The direct image $q_{\mycircle}{\mathcal F}$
of the sheaf ${\mathcal F}$ inherits from ${\mathcal F}$ the $G$-equivariant
structure. We regard the sheaf $q_{\mycircle}{\mathcal F}$ as an element
of $Sh(Y,G,\CC)$.

Denote by $q^G_{\mycircle}{\mathcal F} \in Sh(Y,\CC)$
the sheaf of $G$-invariant sections of $q_{\mycircle}{\mathcal F}$.

The functor $q_\mycircle^G$: $Sh(X,G,\CC)\to Sh(Y,\CC)$ is the composition
of two exact functors: direct image functor $q_\mycircle$ 
 and the functor of invariants of a finite group $G$. It sends
\cite{BL} equivariant
perverse sheaves on $X$ to perverse sheaves on $Y$. Thus we have
\begin{corollary}
Let $\VA_G$ be an equivariant quiver. Then the complex of shaves
$q_\mycircle^G\QVA$ is a perverse sheaf on the quotient space $Y$.
\end{corollary}

Denote by $\overline{\QVA}^{G_{}}$ the complex
$\G^{G_{}}(X,\QVA)=\G(Y,q_\mycircle^G\QVA)$
 of $G$-invariant global sections of the
equivariant complex of sheaves $\QVA$.
Denote by $\det_G$ the one--dimensional representation
 of the group $G$ in the space
$\overline{\Omega}_\emptyset$ of top exterior analytic forms on $X$,
invariant with respect to affine translations.

We have the following equivariant version of Theorem \ref{theorem4}.
\begin{theorem}
\label{Gprop5}
Let $\A$ be a  central arrangement of hyperplanes in $\CC^N$. Let
$G$ be a finite group of linear transformations of\
$\CC^N$ preserving \ $\A$. Let
$\VA_G$ be an  equivariant quiver of the  graph $\Ga$ of the arrangement, such
that all linear maps $A_{\a,\b}$ of the quiver
are close to zero.
Then the complex  $\overline{\QVA}^{G_{}}$ is quasi-isomorphic to the
  subcomplex of $G$-invariants $((C_+(\VA)\ot\det_G)[N])^G$ of the complex
 $(C_+(\VA)\ot\det_G)[N]$:
$$
 \overline{\QVA}^{G_{}}\stackrel{quis}{\simeq}
\left((C_+(\VA)\ot\det\nolimits_G)[N]\right)^{G_{}}\ .
$$
\end{theorem}
The theorem is proved in Section \ref{proofs}.

\medskip
\subsection {Cohomology  of $\CC^N$ with coefficients in
equivariant quiver perverse sheaves} ${}$ \label{section5.5} The
category $Sh(X,G,\CC) $ of equivariant sheaves has  enough injectives
\cite{Gr}. Let $\G^G(X,\cdot) :$  $Sh(X,G,\CC)\to Vect_\CC$ be the functor of
invariant global sections. Following \cite{Gr}, define
$H^k(X,G, \cdot)$ to be the $k$-th right derived functor of the functor
$\G^G(X,\cdot)$, such that
$H^k(X,G, {\mathcal F})= R^k\G^G(X,{\mathcal F})$ for any complex
of equivariant sheaves ${\mathcal F}\in D^b(Sh(X,G,\CC))$.
Since $G$ is a finite group,
for any ${\mathcal F}\in D^b(Sh(X,G,\CC))$ the
 graded space $H^\upcircle(X,G, {\mathcal F})$ coincides with
 the 
 equivariant cohomology group of ${\mathcal F}$ \cite{BL}.

Since the functor $q_\mycircle^G$ : $Sh(X,G,\CC)\to Sh(Y,\CC)$ is the composition
of two exact functors, we have an equality
\cite{Gr} :
$$H^k(X,G, {\mathcal F})=H^k(Y,q_\mycircle^G{\mathcal F})\ .$$
The same statement is valid for the space $X^0$ :
$H^k(X^0,G, {\mathcal F})=H^k(Y^0,(q_0)_\mycircle^G{\mathcal F})\ .$

 The affine space $X$ satisfies the equivariant Stein
condition, that is, for a finite group $G$
\begin{equation}
H^k(X,G, \Oan)=0, \qquad {\rm for \ all}\quad k\not= 0.
\label{stein1}
\end{equation}

Statement \rf{stein1} can be proved as follows.
In the non-equivariant case to show that $H^k(X, \Oan)$ is trivial
we replace the sheaf $\Oan$ (modulo constants) by its Dolbeault resolution
$(\Omega^{0,q},\partial_{\overline{z}})$.
 Then the
$\partial_{\overline z}$-Poincare lemma
says that any $\partial_{\overline z}$-closed  $(0,q)$-form
is a $\partial_{\overline z}$-coboundary.
 In the equivariant case we can use the same
Dolbeault resolution, since it consists of equivariant soft sheaves.
Statement \rf{stein1} requires that an invariant
 $\partial_{\overline z}$-closed $(0,q)$-form $\omega$
has to be the $\partial_{\overline z}$-coboundary of an invariant form. By
the $\partial_{\overline z}$-Poincare lemma, we may conclude that $w$ is the
$\partial_{\overline z}$-coboundary of a possibly non-invariant form $y$.
 Then the non-invariance of $y$
is corrected by averaging over   $G$.

\begin{corollary}
\label{Gprop6} Let $\A$ be a central arrangement of hyperplanes
 in $\CC^N$. Let $G$ be a finite group of linear transformations of
$\CC^N$ preserving $\A$.  Let $\VA_G$ be an equivariant
quiver with all linear maps $A_{\a,\b}$ close to zero. Then
for any $k\in \ZZ$
 \begin{equation}\nonumber
 H^k(X,G, DR(E\VA))\simeq  H^k(X,G, \QVA)\simeq
H^{k+N}\left((C_+(\VA)\ot\det\nolimits_G)^{G}\right)\ .
\end{equation}
\hfill{$\square$}
\end{corollary}

For an equivariant local system ${\mathcal L}$ on $X^0$
introduce equivariant intersection cohomology groups
$I\!H^k_G(X,G,\mathcal{L})$ as
the equivariant cohomology groups $H^{k-N}(X,G,(j_{0})_{\my!*,G}\mathcal{L}[N])$
of the MacPherson extension of the equivariant perverse sheaf
$\mathcal{L}[N]$,
$$
I\!H^k(X,G,\mathcal{L})\stackrel{def}{=}
H^{k-N}(X,G,(j_{0})_{\my!*,G}\mathcal{L}[N])\ .
$$
The functors $q_\mycircle$,  $(q^0)_\mycircle$,  $q_\mycircle^G$, and
$(q^0)_\mycircle^G$ transform perverse sheaves to
perverse sheaves \cite{BBD}, so that we have an equality
$$
I\!H^k(X,G,\mathcal{L}) = I\!H^k(Y, (q^0)_\mycircle^G \mathcal{L})\,,
 $$
where $I\!H^k(Y,  \mathcal{L}')$ means the same as in Section \ref{section4.5},
$$
I\!H^k(Y^{},\mathcal{L'})\stackrel{def}{=}
H^{k-N}(Y,(j_{0})_{\my!*}\mathcal{L'}[N])\ .
$$
for any local system  ${\cal L'}$ over $Y^0$.

 Let ${\cal L}\VA$ be  the local system of  flat
sections of the connection \rf{3.12} equipped with the group action \rf{5.5}.

\begin{corollary} \label{cor5.7}
Let $\VA_G=\{V,A^i,\TG{V}\}$ be an equivariant
 level zero quiver of a central arrangement $\A$.
 Assume that    all
 operators $A^i$ are close to zero.
 Then for any $k\in\ZZ$,
\begin{eqnarray}\nonumber
H^k(Y^0,(q_0)^G_\mycircle{\cal L}\VA)\ \simeq \H^k(X^0,G,{\cal L}\VA)&\simeq&
H^k\left((C_+(\j_{0,*}\VA)\ot\det\nolimits_G)^{G_{}}\right),\\
\nonumber
I\!H^k(Y,(q_0)^G_\mycircle{\cal L}\VA)\ \simeq\
I\!H^k(X,G,{\cal L}\VA)&\simeq&
H^{k}(C_+(\j_{0,!*}\VA)\ot\det\nolimits_G)^{G_{}})\ .
\end{eqnarray}
\end{corollary}

\subsection{Application to discriminantal arrangements} \label{section5.6}
In this section we consider the discriminantal arrangements and local
systems appearing in hypergeometric solutions of the KZ equations in
\cite{SV, V}. We calculate the equivariant intersection cohomology
groups of the discriminantal arrangements in terms of representation
theory.

In this section we follow notations of \cite{V}, chapter 11.

We fix the following data:
\begin{itemize}
\item[(a)] a finite-dimensional complex space $\hh$;
\item[(b)] a non-degenerate bilinear form $(\,,)$ on $\hh$;
\item[(c)] linearly independent covectors ('simple roots') $\a_1,...\a_r \in
\hh^*$.
\end{itemize}
Denote by $b : \hh\to\hh^*$ the isomorphism induced by $(\,,)$. Transform the
form $(\,,\,)$ to $\hh^*$, using $b$. Set  $b_{i,j}=(\a_i,\a_j)$,
$h_i=b^{-1}(\a_i)\in\hh$.

Denote by $\ggg$ the Lie algebra, generated by $e_i,f_i$, and $\hh$, subject to
the relations
\begin{equation}\notag
\begin{array}{cc}\ds [h,e_i]=\< h,\a_i\> e_i\ ,\qquad
&[h,f_i]=-\< h,\a_i\> f_i\ ,\\
\ds [e_i,f_j]=\delta_{i,j}h_i\ ,& [h,h']=0\ ,
\end{array}
\end{equation}
for all $i,j=1,...,r$ and $h,h'\in\hh$. Denote by $\nn_+$ and $\nn_-$
the subalgebras of
$\ggg$ , generated by $e_i$ and $f_i$, respectively.

There is a unique symmetric bilinear form $S(\,,)$ on $\ggg$, satisfying the
conditions
\begin{itemize}
\item[(i)] the subspaces $\nn_+$, $\hh$, and $\nn_-$ are mutually orthogonal
with respect to $S$;
\item[(ii)] the restriction of $S$ to $\hh$ coincides with $(\,,)$,
\item[(iii)] $S(e_i,e_j) = S(f_i,f_j)=\delta_{i,j}$; \
$S([f_i,x],y) = S(x,[e_i,y])$.
\end{itemize}
Set $\gggg=\ggg/\Ker S$, $\nnn_\pm=\nn_\pm/(\nn_\pm\cap \Ker S)$.  These
 are Lie algebras. If $b_{i,i}\not=0$ for all $i$,
then $\gggg$ is the Kac-Moody algebra associated with the matrix
$a_{i,j}=2b_{i,j}/b_{i,i}$.

 For $\Lambda\in\hh^*$, we
denote by $U_\Lambda$ the irreducible
$\gggg$-module with highest weight $\Lambda\in\h^*$.

{}For a highest weight $\gggg$-module
$U$, we denote by $C_{\mycircle}(\nnn_-,U)$ the standard chain complex of
$\nnn_-$ with coefficients in $U$.  The chain complex
$C_{\mycircle}(\nnn_-,U)$ has the canonical structure of an $\hh$-module.
For  $\mu\in\hh^*$, we denote by $C_{\mycircle}(\nnn_-,U)_\mu$ the subcomplex of
chains of weight $\mu$.

\bigskip

 Let
$\A_{1,N}$ be the arrangement  in $X=\CC^N$ of hyperplanes
$$
\begin{array}{crl}\displaystyle
H_i\ {} :& t_i=0,& \qquad 1\leq  i \leq N;\\
H_{i,j}:& t_i-t_j=0,& \qquad 1\leq  i <j\leq N\ .
\end{array}
$$
Fix an $r$-tuple $\l=(k_1,...,k_r)$ $\in$ $\ZZ_{\geq 0}^r$, $k_1+ \dots + k_r = N$,
 and a map   $\pi_\l:\{1,...,N\}$
 $\to$ $\{1,...,r\}$, such that $\# \pi_\l^{-1}(i)=k_i$.
 We denote by  $\Sigma_\l$ the subgroup of the symmetric group $\Sigma_N$, which
consists  of permutations $\sigma$ of
$\{1,...,N\}$, satisfying the condition $\pi_\l(\sigma(i))=\pi_\l(i)$ for all $i$.
We have an isomorphism
$\Sigma_\l\approx \Sigma_{k_1}\times...\times \Sigma_{k_r}$.
The group $\Sigma_\l$, as well as the symmetric group $\Sigma_N$,
acts on $X$ by permutations of coordinates. This action preserves
 $\A_{1,N}$.

{}For  $\Lambda\in\hh^*$ and  $\kappa\in\CC$, we denote
by  $a_{{\Lambda},\pi_\l,\kappa}$ the collection of exponents
$a_{{\Lambda},\pi_\l,\kappa}$: $\A_{1,N} \to \CC$, where
$$
a_{{\Lambda},\pi_\l,\kappa}(H_{i,j})=(\a_{\pi_\l(i)},\a_{\pi_\l(j)})/\kappa,
\qquad
a_{{\Lambda},\pi_\l,\kappa}(H_{i})=-(\a_{\pi_\l(i)},\Lambda)/\kappa\ .$$

Let  $\ (\OS^{\mycircle}(\A_{1,N}), d_{a_{{\Lambda},\pi_\l,\kappa}})\ $
be the Aomoto complex of $\ \A_{1,N}$, see Section \ref{Orlik}.
Let

\noindent
$(\overline{\mathcal F}^{\mycircle}_{a_{{\Lambda},\pi_\l,\kappa}},
d_{a_{{\Lambda},\pi_\l,\kappa}})$  be its subcomplex of flag forms.

The action of $\Sigma_\l$ on $\A_{1,N}$ induces the action of $\Sigma_\l$  on
$(\OS^\mycircle(\A_{1,N}), d_{a_{{\Lambda},\pi_\l,\kappa}})$. This action preserves
$(\overline{\mathcal F}^{\mycircle}_{a_{{\Lambda},\pi_\l,\kappa}},
d_{a_{{\Lambda},\pi_\l,\kappa}})$.

The following isomorphism of complexes was established in
\cite{SV}, Corollary 6.13.
\medskip

\begin{theorem}[\cite{SV}] \label{thm from[SV]}
We have
\begin{equation}\label{dis1}
C_{\mycircle}(\nnn, U_{\Lambda})_{\Lambda-{\overline{\l}}}\approx
\left(\overline{\mathcal F}^{N-\mycircle}_{a_{{\Lambda},\pi_\l,\kappa}}
 (\A_{1,N})\ot {\rm sgn}\right)^{\Sigma_\l}\ ,
\end{equation}
Here $\bar \l=k_1\a_1+\ldots+k_r\a_r$ and
${\rm sgn}$ is the restriction to $\Sigma_\l$ of the one-dimensional
sign representation of  $\Sigma_N$.
\end{theorem}

We rewrite \rf{dis1} in a quiver language. The collection of exponents
$a_{{\Lambda},\pi_\l,\kappa}$ determines the level zero
 equivariant quiver
$\VA_{\Sigma_\l}({\Lambda},\pi_\l,\kappa)$ of the arrangement
$\A_{1,N}$ by the rule
 $\VA_{\Sigma_\l}({\Lambda},\pi_\l,\kappa)$ $=$ $\{V, A^l, \TG{V}\}$.
Here $V=\CC$. The operator $A^{l}:\CC\to \CC$ is the operator of multiplication
 by the exponent $a_{{\Lambda},\pi_\l,\kappa}(H(l))$ of the
 corresponding hyperplane. All operators $\TG{V}$ are equal to $Id_V$.

 In these notations the flag complex
$\overline{\mathcal F}^{\mycircle}_{a_{{\Lambda},\pi_\l,\kappa}}
 (\A_{1,N})$ coincides with the complex  of the direct image
$\j_{0,!*}\VA^{}({\Lambda},\pi_\l,\kappa)$ of the quiver
$\VA({\Lambda},\pi_\l,\kappa)$ $=$
${\mathcal N}^0\VA_{\Sigma_\l}({\Lambda},\pi_\l,\kappa)$,
$$
\overline{\mathcal F}^{\mycircle}_{a_{{\Lambda},\pi_\l,\kappa}}
 (\A_{1,N})\ =\ C_+\left(\j_{0,!*}\VA^{}({\Lambda},\pi_\l,\kappa)\right)
\ .
$$
Statement \rf{dis1} is equivalent to the statement
\begin{equation}\label{dis2}
C_{\mycircle}(\nnn, L_{\Lambda})_{\Lambda-{\overline{\l}}}\approx
\left(C_+^{N-\mycircle}\left(\j_{0,!*,\Sigma_\l}\VA_{\Sigma_\l}({\Lambda},\pi_\l,\kappa)\right)
\ot {\rm sgn}\right)^{\Sigma_\l}\ .
\end{equation}
Let ${\mathcal L}_{{\Lambda},\pi_\l,\kappa}$ be the local system on $X^0$,
associated to the connection in the trivial line bundle
${\rm L}$ over
$X^0$ with the connection-form
$$
\sum_{1\leq i<j\leq N}
\frac{(\a_{\pi_\l(i)},\a_{\pi_\l(j)})}{\kappa}\,d\log(t_i-t_j)\ -\
\sum_{1\leq i\leq N}
\frac{(\a_{\pi(i)},\Lambda_j)}{\kappa}\,d\log t_i\ .
$$
We identify sections of ${\rm L}$ with scalar functions.
We define the action of  $\Sigma_\lambda$ on
sections by the formula
$$
\TG{{\rm L}}(f(x))\ =\ {\rm sgn}(g)\cdot f(\TG{X}^{-1}(x))\ .
$$

In notations of Section \ref{section5.5}, we denote by
$\overline{{\mathcal L}}_{{\Lambda},\overline{\l},\kappa}$ the local
system on $Y^0$ given by the formula
$$
\overline{{\mathcal L}}_{{\Lambda},\overline{\l},\kappa}\ = \
(q_0)_\mycircle^{\Sigma_\l}
{\mathcal L}_{\Lambda,\pi_\l,\kappa}\ .
$$

Combining \rf{dis2} and Corollary \ref{cor5.7}, we get
\begin{corollary}\label{cor5.8}
Let   $X^0$ be the complement to
 $\A_{1,N}$ in $\CC^N$.
 Then for  $|\kappa|\gg 1$ and  $k\in\ZZ$, we have
$$
I\!H^k(Y,\overline{{\mathcal L}}_{{\Lambda},\overline{\l},\kappa})\ =\
I\!H^k(X,{\Sigma_\l},
{\mathcal L}_{\Lambda,\pi_\l,\kappa})\simeq
H_{N-k}(\nnn, L_{\Lambda})_{\Lambda-\overline{\l}}\ .$$
\end{corollary}

\bigskip

Let  $\gggg$
be a  simple Lie algebra with Weyl group $W$
and half sum of positive roots $\rho$. Let $U_\Lambda$
be the irreducible  finite-dimensional $\gggg$-module with highest weight $\Lambda$.
By the Borel--Weil--Bott theorem,
we have
\begin{equation}\notag
H_k(\nnn,U_\Lambda)=\ooplus_{w\in W,\ l(w)=k}\CC_{w(\Lambda+\rho)-\rho}
\end{equation}
where $\CC_\mu$ is a one-dimensional space of weight $\mu$ and $l(w)$ is the length of
$w$. The dominant weight $\Lambda+\rho$ is regular, that is, the equality
$w(\Lambda+\rho)= \Lambda+\rho$ implies  $w = id$.  Hence
each weight subspace of $H_k(\nnn,U_\Lambda)$ is at most one-dimensional.

\bigskip

To a simple Lie algebra $\gggg$ and an element $\overline{\l}=
k_1\a_1+\ldots +k_r\a_r$, $k_i\in\ZZ_{\geq 0}$,
 of the root lattice we  assign the vector
$\l=(k_1,...,k_r)$, the number $N = k_1+\ldots+ k_r$, and the arrangement $\A_{1,N}$
in  $X=\CC^N$, equipped with the action of
the product of symmetric groups $\Sigma_\l$ as at the beginning of this section.
 An integral dominant weight $\Lambda\in\h^*$ determines the $\Sigma_\l$-equivariant
 local system ${\mathcal L}_{\Lambda,\pi_\l,\kappa}$ on the
 complement $X^0$ to the arrangement $\A_{1,N}$ and the local system
$\overline{{\mathcal L}}_{{\Lambda},\overline{\l},\kappa}$ on the
quotient space $Y^0=X^0/\Sigma_\l$, which is an open subset of the quotient
space $Y=X/\Sigma_\l$.

\begin{corollary}\label{cor5.9} For   $|\kappa|\gg 1$ and  $k\in\ZZ$,
the  intersection cohomology
$I\!H^k(Y, \overline{{\mathcal L}}_{{\Lambda},\overline{\l},\kappa})$
is at most one-dimensional and is nonzero if and only there exists $w\in W$, such
that $l(w)=N-k$ and $w(\Lambda+\rho)-\rho=\Lambda-\overline{\lambda}$.
\end{corollary}

\setcounter{equation}{0}

\section{Proofs}\label{Proofs} \label{proofs}
\subsection{Proof of Propositions \ref{lemmaTb} and \ref{propcc}}\label{pl1}
We start with the proof of Proposition \ref{lemmaTb}.  The operators
   $A_\b^\a$ and $A_\b^\g$ commute if $\a\succ\b$ and $\b\succ\g$.
   Indeed, using relations (b) and (c) in Section 2.1 we get
 $$  
 \begin{array}{c} A_\b^\a
 A_\b^\g=A_{\b,\a}A_{\a,\b}A_{\b\g}A_{\g\b}= -\sum\limits_{
 \begin{subarray}{l}\b',\, \b'\not=\b
 \\  \a\succ\b'\succ\g\end{subarray}} A_{\b,\a}A_{\a\b'}A_{\b'\g}A_{\g\b}=\\
 =-\sum\limits_{\begin{subarray}{l}
\b',\,\b'\not=\b \\ \a\succ\b'\succ\g\end{subarray}}
 A_{\b\g}A_{\g\b'}A_{\b'\a}A_{\a,\b}=
 A_{\b\g}A_{\g\b}A_{\b,\a}A_{\a,\b}
= A_\b^\g A_\b^\a\ .
 \end{array}
 $$   
 To prove  { (i)} we have to check only that
 $[S_\b,A_\b^\a]=0$
 for $\a$ such that $\a\succ\b$. This is equivalent to the equality
 $[\sum\limits_
{\substack{\a',\;\,\a'\not=\a,\\ \a'\succ\b}}\! A_\b^{\a'},A_\b^\a]\ =\ 0$.
 Consider the product
 $A_\b^\a\cdot\!\!\!\!\!\sum\limits_
{\substack{\a_j,\;\,\a_j\not=\a,\\ \a_j\succ\b}}A_\b^{\a_j}$.
A nonzero contribution to the product is given only by the summands with $\a_j$ for which
there exists
 $\d_j$ with  $\d_j\succ\a_j$ and $\d_j\succ\a$. We have
\begin{equation}\begin{split}
 A_\b^\a\sum_{\substack{\a_j,\;\,\a_j\succ\b,\\ \a_j\not=\a}}A_\b^{\a_j} =
A_{\b,\a}A_{\a,\b} \sum\limits_{\d_i,\,
\d_i\succ\a}
\ \sum\limits_{\begin{subarray}{l}\a_j,\,\a_j\not=\a,\\
\,\d_i\succ\a_j\succ\b\end{subarray}}A_{\b\a_j}A_{\a_j\b}=\\
\label{commut} 
- \sum\limits_{\d_i,\,
\d_i\succ\a}\ \sum\limits_{\begin{subarray}{l}\a_j,\,\a_j\not=\a,\\
\,\d_i\succ\a_j\succ\b\end{subarray}}A_{\b,\a}A_{\a\d_i}A_{\d_i\a_j}A_{\a_j\b}=
 \sum\limits_{\d_i,\, \d_i\succ\a}
A_{\b,\a}A_{\a\d_i}A_{\d_i\a}A_{\a,\b} \ .
\end{split}\end{equation}
We can apply the same procedure to the product
 $\sum\limits_{\substack{\a_j,\,\a_j\not=\a,\\ \a_j\succ\b}}\!A_\b^{\a_j}A_\b^\a$ and get
 the same result \rf{commut}. This proves  statement  (i) of Proposition
\ref{lemmaTb}.

  To prove  (ii) we note
 that for any vertex $\delta>\b$ with $l(\d)=l(\b)-2$ and for any
  vector $e_\d\in V_\d$, the vector $\oplus_{\a,\,\d\succ\a\succ\b}
  A_{\a,\d}e_\d$ belongs
  to the kernel of the operator  $T_\b$ due to defining relation (b) of
 Section 2.1.
  {}  Fix a vertex $\a$, such that $\a\succ\b$, and a vector
  $e_\a\in V_\a$. We have
  $$   
  T_\b(e_\a)=A_{\a,\b}A_{\b,\a}(e_\a)+\sum\limits_{\begin{subarray}{r}
\a_j,\,\a_j\not=\a\
  \\ \,\a_j\succ\b\end{subarray}}A_{\a_j\b}A_{\b,\a}(e_\a).$$
   Due to relation (c) of Section 2.1,
   in the last sum we have a nonzero contribution only from
  the summands with  $\a_j$ such that there exists $\d_j$
  satisfying $\d_j\succ\a_j$ and $\d_j\succ\a$. Thus
  \begin{equation*}\begin{split}
  T_\b(e_\a)=A_{\a,\b}A_{\b,\a}(e_\a)-\sum\limits_{\d_i,\d_i\succ\a}
  \sum\limits_{\begin{subarray}{r}\a_j,\,\a_j\not=\a\\
  \,\d_i\succ\a_j\succ\b\end{subarray}}
  A_{\a_j\d_i}A_{\d_i\a}(e_\a)=\\
  =A_{\a,\b}A_{\b,\a}(e_\a)+\sum\limits_{\d_i,\d_i\succ\a} A_{\a\d_i}
  A_{\d_i\a}(e_\a)-\sum\limits_{\d_i,\d_i\succ\a}
  \sum\limits_{\a_j,\,\d_i\succ\a_j\succ\b}
  A_{\a_j\d_i}A_{\d_i\a}(e_\a).\end{split}\end{equation*}
  The last term belongs to the kernel of $T_\b$. So the action of
   $T_\b$ on  $e_\a$ in the quotient space
  $\oplus_{\a:\a\succ\b}V_\a/{\rm Ker}\, T_\b$ is the same as the action
  of  $A_\a^\b+T_\a^0$ on  $e_\a$. This   proves statement
(ii) of Proposition \ref{lemmaTb}.
\hfill{$\square$}

\smallskip

Let us prove statement (i) of Proposition \ref{propcc}.
Since $A_\b^\a$ and $A_\b^\g$
   commute if $\a\succ\b$ and $\b\succ\g$, we have to prove that
$[\widetilde{S}_\b,A_\b^\a]=0$ for any $\a\in\v(\G)$, $\b\succ\a$.
Fix a vertex $\a\in\v(\G)$, such that $\b\succ\a$.
We have $\widetilde{S}_\b=\sum_{\gamma,\,\b\succ\gamma}A_\b^\gamma$.
Since the arrangement is central, for any $\gamma,\,\b\succ\gamma$,
the planes $\X_\a$ and $\X_\g$ intersect in codimension one and we can
rewrite $\widetilde{S}_\b$ as
 $$\widetilde{S}_\b=\sum_{\delta, \, \b>\delta, l(\delta)=l(\b)+2}
\widetilde{S}_{\b,\delta}\ ,$$
where
\begin{equation}
\label{Sb}
\widetilde{S}_{\b,\delta}=\sum_{\gamma,\,
\b\succ\gamma\succ\delta}A_\b^\gamma\ .
\end{equation}
The defining relations of quivers (see the relation (v) of Section
\ref{Seclevn}) imply that
$[A_\b^\a, \widetilde{S}_{\b,\delta}]=0 $
for any such $\delta$. Thus $[A_\b^\a, \widetilde{S}_{\b}]=0 $ and
statement (i) is proved.
\smallskip

Let us prove statement (ii) of Proposition \ref{propcc}.
Fix $\a,\b\in\v(\G)$, such that $\a\succ\b$. For any $\delta$,
$\b\succ\d$ we have by the second  equality of Section \ref{Seclevn} (iv)
\begin{equation}\label{Sba}
 A_{\a,\b}A^\delta_\b=(\widetilde{S}_{\a,\delta}-A_\a^\b)A_{\a,\b}\ .
\end{equation}
The summation of \rf{Sba} over all  $\delta$, $\b\succ\delta$ gives
\begin{equation}\label{Sbb}
 A_{\a,\b}\widetilde{S}_\b=(\widetilde{S}_{\a}-A_\a^\b)A_{\a,\b}\ .
\end{equation}
On the other hand, we have, analogously to \rf{commut},
\begin{equation*}\begin{split}
 A_{\a,\b}\sum_{\substack{\a_j,\;\,\a_j\succ\b,\\ \a_j\not=\a}}A_\b^{\a_j} =
A_{\a,\b} \sum\limits_{\d_i,\,
\d_i\succ\a}
\ \sum\limits_{\begin{subarray}{l}\a_j,\,\a_j\not=\a,\\
\,\d_i\succ\a_j\succ\b\end{subarray}}A_{\b\a_j}A_{\a_j\b}=\\
\label{commuta} 
- \sum\limits_{\d_i,\,
\d_i\succ\a}\ \sum\limits_{\begin{subarray}{l}\a_j,\,\a_j\not=\a,\\
\,\d_i\succ\a_j\succ\b\end{subarray}}A_{\a\d_i}A_{\d_i\a_j}A_{\a_j\b}=
 \sum\limits_{\d_i,\, \d_i\succ\a}
A_{\a\d_i}A_{\d_i\a}A_{\a,\b} \ ,
\end{split}\end{equation*}
that is,
\begin{equation}
\label{commutb}
 A_{\a,\b}(S_\b-A_\b^\a)=S_\a A_{\a,\b} \ .
\end{equation}
Combining \rf{Sbb} and \rf{commutb}, we get
$$A_{\a,\b}(S_\b+\widetilde{S}_\b-A_\b^\a)=
(S_\a+\widetilde{S}_\a-A_\a^\b) A_{\a,\b}\ ,$$
which implies \rf{commutc}, since $A_{\a,\b}A_\b^\a=A_\a^\b A_{\a,\b}$
by definition of operators $A_\a^\b$ and $A_\b^\a$, see \rf{Aab}.
 The proof of \rf{commutc} for $\b\succ\a$ is analogous. The relations
 \rf{commutc} imply the centrality of the operator $S$. \hfill{$\square$}

\subsection{Proof of Propositions \ref{lemma2}--\ref{propcmon}}
First we prove  Proposition \ref{lemma2}.

Consider the direct image
$\j_{0,*}\UG_{\G^0}$ of the level zero quiver  $\UG_{\G^0} =\{W, B^i\}$.
Consider a vertex $\a\in\v(\Ga)$ and a flag
 $F_{\a_,...,\a_m=\a}$ which ends at $\a$.
Let $\b\in\v(\Ga)$ be a vertex such
that $\b\succ\a$. According to \rf{3.13} and \rf{flag-}, the
vector $A_{\a,\b}A_{\b,\a}\left(F_{\a_,...,\a_m}\ot u\right)$ in
$\left(\j_{0,*}\UG_{\G^0}\right)_\a$ is equal to zero if
$F^\b_{\a_0,...,\a_m}=0$ or is equal to
\begin{equation*}
(-1)^{((\b;\a_0,...,\a_m)+m-1)}F_{\b_0,...,\b_{m-1}\a_m}\ot
\!\!\!\!\!\! \sum\limits_ {\substack{i\in\v_1(\Ga),\\ i\geq\a_{k+1},\,
 i\not\geq \a_{k}} }\!\!\!\!\!
 B^i(u),
\end{equation*}
if $F^\b_{\a_,...,\a_m}\not=0$. Here $F_{\b_0,...,\b_{m-1}}=
F^\b_{\a_,...,\a_m}$. Denote by $I_k(\a_0,\a_1,...\a_m)$ the set
of all $\b\in\v(\G)$ such that 
$F^\b_{\a_,...,\a_m}\not=0$ and $(\b;\a_0,...,\a_m)=k$. Applying
relations \rf{flagrelation},  we get
\begin{equation}
\label{Ik} \sum_{\b\in
I_k(\a_0,...,\a_m)}A_{\a,\b}A_{\b,\a}(F_{\a_,...,\a_m=\a}\ot u)=
F_{\a_,...,\a_m}\ot\sum\limits_ {\substack{i\in\v_1(\Ga),\\ i\geq\a_{k+1},\,
 i\not\geq \a_{k}} }\!\!\!\!\!  B^i(u)\ .
\end{equation}
Relations \rf{Ik} imply that for any $\a\in\v(\Ga)$,
\begin{equation}
S_\a=Id_{\FF_\a}\ot\sum_{i\in\v_1(\Ga), \,i\geq \a}B^i \ ,\label{Sa}
\end{equation}
as an operator in $\FF_\a\ot U$. Denote by $\widetilde{B}^\a$ the
operator
 $\widetilde{B}^\a=\sum_{i\in\v_1(\Ga),\ i\geq \a}B^i: W\to W$.
 According to Proposition \ref{lemmaTb},
the operator $S_\a$ commutes with any product $A_{\a,\b}A_{\b,\a}$,
if $\b$ is adjacent to $\a$. Applying this fact to relations \rf{Ik}
and \rf{Sa} we get
$$
[\widetilde{B}^\a,\sum\limits_{\substack{i\in\v_1
(\Ga),\\ i\geq\a_{k+1},\,
 i\not\geq \a_{k}} }\!\!\!\!\!
 B^i]=0\ ,
$$
and, in particular, for $k=0$, we get  $[\widetilde{B}^\a,B^{\a_1} ]=0$.
We can repeat these arguments for any flag ending at $\a$ and
we conclude that the operator  $\widetilde{B}^\a$ commutes with
all its own summands $B^i$. In this situation all eigenvalues of
the operator
$\widetilde{B}^\a$ are   equal to $\l_\a$.
Indeed, for any eigenvalue $\l$  of $\widetilde{B}^\a$ the corresponding subspace
of eigenvectors is invariant with respect
to any $B^i$, $i\geq \a$. Taking the trace over this subspace of
the equality
 $\widetilde{B}^\a=\sum_{i\in\v_1(\Ga):\, i\geq \a}B^i$, we get
$\l=\l_\a$. Due to \rf{Sa} we get statement (i) of the proposition.

To prove (iii), note that due to  \rf{3.13} and \rf{flag-} for any
$\a,\b\in \v(\Ga)$, such that $\a\succ\b$, the operator
$A_\a^\b=A_{\a,\b}A_{\b,\a}:\FF_\a\ot U\to \FF_\a\ot U$ is equal to
\begin{equation}\label{Sa1}
A_\a^\b=Id_{\FF_\a}\ot \sum\limits_{\substack{i\in\v_1(\Ga),\\ i\geq\b,\,
 i\not\geq \a}}B^i\ .
\end{equation}
Statement {(iii)} follows now from { (i)} and
{ (ii)} and Proposition
\ref{lemmaTb}.
\hfill{$\square$}

\bigskip

Consider the quiver $\j_{0,!}\UG_{\G^0} =
\tau_N^{-1}\j_{0,*}\tau_0\UG_{\G^0}$.  The quiver
$\tau_0\UG_{\G^0}=\{U,C^i\}$ consists of the vector space $U=W^*$ and
operators $C^i=(B^i)^*$, respectively dual to $B^i$. Each of $C^i$ has
a single eigenvalue $\l_i$. By Proposition \ref{lemma2}, the operators
$S_\a$, $A_\a^\b$, $T_\a$, $T_\a'$, attached to the quiver
$\j_{0,*}\tau_0\UG_{\G^0}$, have eigenvalues, prescribed by that
proposition.  The operators $S_\a$, $A_\a^\b$, $T_\a$, $T_\a'$,
attached to the quiver $\tau_N^{-1}\j_{0,*}\tau_0\UG_{\G^0}$, are dual
to the corresponding operators for the quiver
$\j_{0,*}\tau_0\UG_{\G^0}$ and thus have the same eigenvalues. This
proves Proposition \ref{8}.  \hfill{$\square$}
\smallskip

In order to prove Proposition \ref{propcmon}, it is  sufficient to prove it
for the quiver $\j_{0,*}\UG_{\G^0}$ and then apply the duality arguments.

{}For $\a\in\v(\Ga)$, formula \rf{Sa1} implies that
$$
\widetilde{S}_\a=Id_{\FF_\a}\ot \sum B^i\ ,
$$
where the sum 
 is taken over all $i\in\v_1(\Ga)$, such that
the hyperplane $H_i$ does not contain $X_\a$, but its
intersection with $X_\a$ is non-empty. For a central arrangement, the
last condition is satisfied always. Thus
\begin{equation}\label{Sa2}
\widetilde{S}_\a=Id_{\FF_\a}\ot \sum\limits_{{i\in\v_1(\Ga), \,
 i\not\geq \a}}B^i\ .
\end{equation}
Combining \rf{Sa2} with \rf{Sa}, we get
\begin{equation}
\label{Sa3}
S=Id_{\FF}\ot \sum\limits_{{i\in\v_1(\Ga), \,}}B^i\ .
\end{equation}
The trace arguments as above show that the operator $S$ has a single
eigenvalue $\l_\infty=\sum_{j\in J(\A)}\l_j$, \ if every operator $B^i$
has a single eigenvalue $\l_i$. \hfill{$\square$}
\subsection{Proof of Proposition \ref{lemma1}}
 Choose a scalar product on $X$ invariant with respect to parallel translations.
 Without loss of generality we  assume that for any
  pair of adjacent strata $\X_\a\supset X_\b$, the
 function $f_{\b,\a}\in \Fb\setminus \Fa$ is chosen to be orthogonal
 to $\Fa$  and the vector field
 $\xi_{\a,\b}\in\Ta\setminus\Tb$ belongs to
 $\TTb$.
 Set $M=\oplus_\a M_\a$.
 Due to
relations \rf{2.3} and \rf{2.4}, the space of global sections of
the $\D_X$-module $\EVA$ is
 isomorphic to a quotient space of $M$. To prove
 \rf{2.6}, we define inductively the structure of a $D_X$-module
on $M$.

\medskip

Define the grading on $M$.
It is given by a grading on each summand
$M_\a$, where we set $\deg (v_\a\otimes \om_\a)=0$,
$\deg\xi=\deg f=1$ for  $\xi \in \TTa$ and $f \in \FFa$.

We will use the increasing Bernstein filtration
 $BD^\bullet$ of the ring
 $D_X$, in which the subspace $BD^{\leq k}\subset D_X$
 is spanned by the elements $z_{i_1}\cdots z_{i_l}\der{z_{i_{l+1}}}\cdots
 \der{z_{i_m}}$\ ,
 $$
BD^{\leq k}\ =\ \CC\left\langle z_{i_1}\cdots z_{i_l}\der{z_{i_{l+1}}}\cdots
 \der{z_{i_m}}\ {}\ \vert \ {}\ m\leq k\right\rangle\ .
$$

Define the action of $BD^{\leq 1}$ on elements
 of zero degree\ $M^0\ \subset\ M$. Take a stratum $X_\a$.
 An element  $d \in BD^{\leq 1}$ can be presented as a sum
  \begin{equation}
  d=c+d_\a+d'_\a,
  \label{2.7}
  \end{equation}
  where $c$ is a constant, $d_\a \in \Ta+\Fa$,
  $d'_\a \in \TTa+\FFa$.
 In the ring $D_X$, we have
\begin{equation}
 [d_1,d_2]=0\qquad{\rm if} \qquad d_1,d_2\in\Ta+\Fa\quad
 {\rm or}\quad d_1,d_2\in\TTa+\FFa\ .
 \label{2.0}
 \end{equation}
 We use  relations \rf{2.3} and \rf{2.4}
  to define the
  action. We set
  \begin{equation}
 d_\a(\om_\a\ot v_\a)\ =\ \sum_{\b, \ \a\succ\b}
\mylangle d_\a,f_{\b,\a}\myrangle
\pi_{\b,\a}(\om_a)\ot A_{\b,\a}(v_\a) \in \sum_{\b,\ \a\succ\b} M_\b^0\ ,
\label{2.10}
\end{equation}
if $d_\a \in \Ta$, and
\begin{equation}\label{2.11}
 d_\a(\om_\a\ot v_\a)\ =\  \sum_{\b,\ \b\succ\a}
 \mylangle \xi_{\b,\a},d_\a\myrangle
\pi_{\b,\a}(\om_\a)A_{\b, \a}(v_\a)\in \sum_{\b,\ \b\succ\a}M_\b^0\ ,
\end{equation}
if $d_\a \in \Fa$.
 We set
 \begin{equation}\label{2.12}
 d'_\a(\om_\a\ot v_a)\ =\ d'_\a\ot(\om_\a\ot v_a)\in M_\a^1\ ,
 \end{equation}
if    $d'_\a \in \TTa+\FFa$.

 As the next step we define the action
 of $BD^{\leq 1}$ on the subspace of degree one
elements $M^1\subset M$.
 Again we   use decomposition \rf{2.7}.
Define the action of $d' \in \TTa+\FFa$ on
 $d'\otimes (\om_\a\otimes v_\a)\in M_\a^1$ by the formula
\begin{equation}\label{2.8}
 d'_\a \cdot (d'\otimes (\om_\a\otimes v_\a))=
 (d'_\a  d')\otimes (\om_\a\otimes v_\a)\in M_\a^2\ .
 \end{equation}
Define the action of  $d_\a \in \Ta+\Fa\subset BD^{\leq 1}$ by the formula
 \begin{equation}\label{2.9}
 d_\a \cdot (d'\otimes (\om_\a\otimes v_\a))=
 [d_\a,d']\otimes (\om_\a\otimes v_\a)+
  d'\cdot d_\a (\om_\a\otimes v_\a)\ ,
  \end{equation}
where the first term in the right hand side is in $M_\a^{0}$, while
the action of $BD^{\leq 1}$ for the
second  term was  already defined.

We claim that formulas \rf{2.10}--\rf{2.9}
define also an action of $BD^{\leq 2}$ on $M^0$.

To prove this we note first that  relations
\rf{2.10}--\rf{2.9} define an action of all second degree
monomials on $M^0$ and the only thing to check is
to verify that this action respects quadratic relations in
$BD^{\leq 2}$. The action of the basic relation
\begin{equation}\label{2.13}
d_1d_2-d_2d_1\ =\ [d_1,d_2]\
\end{equation}
on $M^0_\a$ can be split into several types depending
on the types in \rf{2.7} of the first degree generators $d_i$.

The case, when both $d_i$ belong to $\TTa+\FFa$, is trivial. In that case
$d_i$ commute and their actions commute due to \rf{2.12}.

If $d_1 \in \TTa+\FFa$ and  $d_2\in \Ta+\Fa$, then relation
\rf{2.13} is preserved due to \rf{2.9}.

The nontrivial case is
when both $d_i$ belong to $\Ta+\Fa$ and should commute.

Suppose first that both $d_i\in\Ta$. Then
$$
d_1d_2(\om_\a\ot v_a)\ =\
 \sum_{\b,\ \a\succ\b}
\mylangle d_2,f_{\b,\a}\myrangle\
 d_1\left(\pi_{\b,\a}(\om_\a)\ot A_{\b,\a}(v_\a) \right)\ .
$$
We can present each of $d_i$ as
$$
d_i
\ =\ c_i^\b\xi_{\a,\b}+d'_i\ ,
$$
where $\xi_{\a,\b} \in\Ta\cap\TTb$,
 $f_{\b,\a} \in\Fb\cap T_\a^\perp$,  $d'_i$ belongs to $\Tb$, and the constants
  $c_i^\b$ are given by the formula $c_i^\b=\tfrac{\mylangle d_i,f_{\b,\a}\myrangle}
  {\mylangle \xi_{\a,\b},f_{\b,\a}\myrangle}$. Then
\begin{equation}\begin{split}
(d_1d_2-d_2d_1)(\om_\a\ot v_\a)\ =\ \sum_{\b,\ \a\succ\b}\left(c_1^\b
\mylangle d_2,f_{\b,\a}\myrangle -c_2^\b\mylangle d_1,f_{\b,\a}\myrangle\right)
\pi_{\b,\a}(\om_\a)\ot A_{\b,\a}v_\a+
\\
\label{2.14}
+\sum_{\b,\ \a\succ\b}
\left(\mylangle d_2,f_{\b,\a}\myrangle d'_1-
\mylangle d_1,f_{\b,\a}\myrangle d'_2\right)\cdot
\left(\pi_{\b,\a}(\om_\a)\ot
A_{\b,\a}v_\a\right) .\end{split}
\end{equation}
The first summand in the right hand side of \rf{2.14} is zero, since
$$
\left(c_1^\b
\mylangle d_2,f_{\b,\a}\myrangle -c_2^\b\mylangle d_1,f_{\b,\a}\myrangle
\right) = 0
$$
for each
$\b\in\v(\Ga)$, such that $\a\succ\b$.
The second summand can be rewritten as a double sum
\begin{equation}\label{2.15}
\sum_{\substack{\b,\ \a\succ\b,\\ \g,\ \b\succ\g}}
C_{\a,\b,\g}\ot A_{\g\b}A_{\b,\a}v_\a=\!\!
\sum_{\substack{\b,\ \a\succ\b,\\ \g,\ \b\succ\g}}
i_{d_1\wedge d_2}\left(df_{\b,\a}\wedge df_{\g\b}\right)
\pi_{\g,\b}\pi_{\b,\a}(\om_\a)\ot
A_{\g\b}A_{\b,\a}v_\a \ .
\end{equation}
{}For a fixed $\g$ the coefficient  $C_{\a,\b,\g}$ in
\rf{2.15} does not depend on
index $\b$. This  can be seen from the equality
$$C_{\a,\b,\g}=\frac{\om_\a}{df_{\b,\a}\wedge df_{\g\b}}
i_{d_1\wedge d_2}\left(df_{\b,\a}\wedge df_{\g\b}\right)\ .$$
 Then the sum in \rf{2.15} equals zero due to relation (b) of Section 2.1.

 The case when both $d_i$ belong to \Fa is treated in
 analogous manner. Let now $d_1\in\Fa$, $d_2\in\Ta$.
 Then
\begin{equation}\begin{split}
d_1d_2(\om_\a\ot v_\a)=
 \sum_{\b,\ \a\succ\b}
\mylangle d_2,f_{\b,\a}\myrangle
 d_1\left(\pi_{\b,\a}(\om_\a)\ot A_{\b,
\a}v_\a \right)=\\
\label{2.16}
=\sum_{\substack{\b,\ \a\succ\b\\ \g,\ \g\succ\b}}
\mylangle d_2,f_{\b,\a}\myrangle \mylangle \xi_{\g\b},d_1\myrangle\pi_{\g,\b}\pi_{\b,\a}(\om_\a)\ot
A_{\g\b}A_{\b,\a}v_\a\ ,
\end{split}\end{equation}
since $d_1$ belongs to $\Fb\subset \Fa$ for each summand in the first line of
\rf{2.16}, and
\begin{equation}\begin{split}
d_2d_1(\om_\a\ot v_\a)=
 \sum_{\d,\ \d\succ\a}
\mylangle \xi_{\d\a},d_1\myrangle
 d_2\left(\pi_{\d,\a}(\om_\a\ot A_{\d,
\a}v_\a\right)=\\
\label{2.17}
=\sum_{\substack{\d,\ \d\succ\a\\ \tilde{\g},\ \d\succ\tilde{\g}}}
\mylangle d_2,f_{\tilde{\g}\d}\myrangle
\mylangle \xi_{\d\a},d_1\myrangle\pi_{\tilde{\g},\d}\pi_{\d,\a}(\om_\a)
\ot A_{\tilde{\g}\d}A_{\d\a}v_\a\ ,
\end{split} \end{equation}
since $d_2$ belongs to $\T{\d}\subset\T{\a}$ for each summand
 in the first line of \rf{2.17}.

 We have the equalities of codimensions $l(\g) = l(\tilde{\g}) = l(\a)$,
and the summation in the last lines of \rf{2.16} and \rf{2.17} goes
 over one index $\g$ or $\tilde{\g}$ while $\X_\b=\X_\a\cap \X_\g$ in
 \rf{2.16} and  $X_\d$ is the unique hyperplane which contains both
 $X_\a$ and $X_{\tilde{\g}}$ in codimension 1 for \rf{2.17}.

Due to \rf{2.0}, we should prove that the difference of the
right hand sides
of \rf{2.16} and \rf{2.17} is zero. First of all, there are no
terms with $\g=\a$
in the right hand side of \rf{2.16}, since in this case
the  coefficient $\mylangle \xi_{\g\b},d_1\myrangle$  is zero. Analogously, there
are no terms in the right hand side of \rf{2.17} with $\tilde{\g}=\a$ since in
this case the coefficient $(d_2, f_{\tilde{\g},\delta})$ vanishes.

 {}For a  vertex $\g'$ of the graph $\Ga$  there are three
different possibilities:
\begin{itemize}
\item[(1)] 
there is a summand in \rf{2.16} with index $\g=\g'$
and  a summand in \rf{2.17} with index $\tilde{\g}=\g'$;
\item[(2)] 
there is a summand in \rf{2.16} with index $\g=\g'$ and
 there is no summand in \rf{2.17} with $\tilde{\g}=\g'$;

\item[(3)] 
there is no summand in \rf{2.16} with index $\g=\g'$ and
 there is a summand in \rf{2.17} with $\tilde{\g}=\g'$.
 \end{itemize}
In  case (2) the corresponding term vanishes due to
relation (c) in Section 2.11. In case  (3), the strata $X_\a$ and
$X_{\tilde{\g'}}$ are parallel inside $X_\d$, so the coefficient
$\mylangle d_2,f_{\g'\d}\myrangle$  vanishes as well as the corresponding
summand. In case  (1) we can apply  relation
(c) of Section 2.1 if the sum of the coefficients in \rf{2.16} and
\rf{2.17} is zero. So the only thing to prove is the equality
\begin{equation}\label{2.18}
\mylangle d_2,f_{\b,\a}\myrangle
\mylangle \xi_{\g\b},d_1\myrangle\pi_{{\g},\b}\pi_{\b,\a}(\om_\a)
+
(d_2,f_{{\g}\d})\mylangle \xi_{\d\a},d_1\myrangle
\pi_{{\g},\d}\pi_{\d,\a}(\om_\a)
=0 ,
\end{equation}
when $X_\a$ and $X_\g$ are strata of equal dimension,
 $\X_\b=\X_\a\cap\X_\g$, $\X_\d=\X_\a+\X_\g$ if an origin
  of $\C^N$ is taken inside $X_\b$, $d_1\in\Fa$ and $d_2\in\Ta$.
 Then we can choose $\xi_{\d\a}=\xi_{\g\b}$ and $f_{\b,\a}=f_{\g\d}$.
 With such an identification the left hand side of \rf{2.18} is equal to
$$(d_2,f_{\g,\d})\mylangle \xi_{\g\b},d_1\myrangle\left(\pi_{{\g},\b}\pi_{\b,\a}(\om_\a)
+\pi_{{\g},\d}\pi_{\d,\a}(\om_\a)\right).$$
The equality $\om_\g=\pi_{{\g},\b}\pi_{\b,\a}(\om_\a)$ is equivalent to
the equality
$$\om_\a=df_{\b,\a}\wedge i_{\xi_{\g,\b}}\om_\g=
df_{\g,\d}\wedge i_{\xi_{\g,\b}}\om_\g, $$
and the equality
 $\om_\g=\pi_{{\g},\d}\pi_{\d,\a}(\om_\a)$ is equivalent to
$$\om_\a= i_{\xi_{\d,\a}}(df_{\g,\d}\wedge\om_\g)=
 i_{\xi_{\g,\b}}(df_{\g,\d}\wedge\om_\g). $$
 The elements  $\xi_{\g\b}$ and  $f_{\g\d}$ commute in $D_X$ by \rf{2.0}.
This implies the equality
$df_{\g,\d}\wedge i_{\xi_{\g,\b}}\om_\g+
 i_{\xi_{\g,\b}}(df_{\g,\d}\wedge\om_\g)=0$ which proves \rf{2.18}.

Suppose $n\geq 0$ and  $k, 0\leq k\leq n$ are given. Suppose
 an action of $B^{\leq k}$ on $M^{l}$ is defined, if $l\leq n-k$.
Let $d \in BD^{\leq 1}$.  For a given stratum $X_\a$ and
 a monomial$D$  of degree $n$ over $\C[\TTa+\FFa]$ define
\begin{align}
\label{2.19}
&d\cdot D\otimes \om_\a\ot v_\a= [d,D]\ot \om_\a\ot v_\a +
D\cdot d(\om_\a\ot v_\a) & &{\rm if} && d\in\Ta+\Fa,\\
&d\cdot D\otimes \om_\a\ot v_\a= dD\ot
\om_\a\ot v_\a&&
  {\rm if} && d\in\TTa+\FFa\ .
\label{2.20}
\end{align}
Then the action of all  monomials
 from $B^{\leq k}$ on $M^{n-k+1}$ is defined.
We should check that this action preserves the commutation relations
 in the ring $D_X$, more precisely, the relations of the type
 \begin{equation}
 \label{2.21}
 d_1\cdot d_2\cdot (D\ot \om_\a\ot v_\a)-
 d_2\cdot d_1\cdot (D\ot \om_\a\ot v_\a)
 =[d_1, d_2] (D\ot \om_\a\ot v_\a)
 \end{equation}
where $d_i \in BD^{\leq 1}$ and
 $D$ is a monomial of degree $n-1$ over $\C[\TTa+\FFa]$.
 If both $d_i \in  \TTa+\FFa$ then \rf{2.21} is satisfied due to
 \rf{2.20}. Let now $d_1\in \Ta+\Fa$, $d_2\in \TTa+\FFa$. Then
 by \rf{2.19} and by induction assumptions
 \begin{equation}\notag
 \begin{split}
 d_1d_2(D\ot \om_\a\ot v_\a)&=[d_1,d_2D]\ot \om_\a\ot v_\a +
 d_2Dd_1( \om_\a\ot v_\a)\   ,\\
 d_2d_1(D\ot \om_\a\ot v_\a)&=d_2[d_1,D]\ot \om_\a\ot v_\a+
 d_2(Dd_1( \om_\a\ot v_\a))\ ,
 \end{split}
 \end{equation}
 so  we get
 $$(d_1d_2-d_2d_1)\cdot(D\ot \om_\a\ot v_\a)=[d_1,d_2]D\ot
  \om_\a\ot v_\a). $$
Let now $d_1,d_2\in\Ta+\Fa$. Then by \rf{2.19}
$$
 d_1d_2(D\ot \om_\a\ot v_\a)=
 [d_1,[d_2,D]]\ot \om_\a\ot v_\a +
 [d_2,D](d_1(\om_\a\ot v_\a))+d_1(Dd_2(\om_\a\ot v_\a))
$$
and
 $$(d_1d_2-d_2d_1)\cdot(D\ot \om_\a\ot v_\a)=
 D(d_1d_2-d_2d_1)(\om_\a\ot v_\a)=0$$
 by induction assumptions. This ends the induction step and the proof of
the proposition.\hfill{$\square$}
\subsection{Proof of Proposition \ref{prop1}}

The description \rf{2.3} and \rf{2.4} of a quiver $\D_X$-module
imply that the principal filtration \rf{pfil} is good and the map
$$\phi:
\ooplus\limits_{\a\in\v(\Ga)}\left(S(\TTa)\otimes \Oa\ot V_\a\right)
 \otimes \CC[\X_\a]\to \overline{E}\VA\ ,$$
where $\overline{E}\VA = {\rm gr}\,{E}\VA$, is a well defined map of
${\mathcal O}_X$-modules. Proposition \ref{lemma1}
implies that this map is an isomorphism.

Consider the product
 $$
\prod_{\a\in  \v( \Ga)} \overline{d}_\a\ ,
$$
 where $d_\a\in \Ta+\Fa$ and $\overline{d}_\a={\rm gr}\,(d_a)$.
Relations \rf{2.3}, \rf{2.4}, \rf{2.19} and \rf{2.20}
 show that  this product
 annihilates $ \overline{E}\VA$.
 It means that $E\VA$ is a holonomic $\D_X$ module
 with singular support contained in $\cup_\a T^*_{X_\a}\C^N$.
\hfill{$\square$}

\subsection{Proof of Proposition \ref{theorem1a}}\label{prth1a}
{}Fix a nonzero exterior form $\om_\ee\in\overline{\Om}_\ee$.
Identify vectors $v\in V$ with
$\om_\ee\ot v\in \overline{\Om}_\ee\ot V$, and  $w\in W$ with
$\om_\ee\ot w\in \overline{\Om}_\ee\ot W$.

Let $\Phi: E^0 \VA \to E^0 {\mathcal W}$ be a morphism of $\D_{X^0}$-
modules. Since $X_0$ is an affine variety, the morphism is determined
by the corresponding map of global sections,  denoted by $\Phi$
also.  Since the $D_{X^0}$-module $E^0\VA$ is generated by the space
$V$, the map of global sections is determined by its restriction to
$V$, denote $\varphi=\Phi|_V$.  The map $\varphi$ is a regular rational
function on $X^0$ with values in $\Hom (V,W)$. We need to prove that
$\varphi$ is constant.

{}For $j\in J(\A)$, choose affine coordinates $z_1,...,z_N$ on $X$,
such that the hyperplane $H_j$ has equation $z_1=0$. Take a generic
point $(0,\varepsilon_2,...,\varepsilon_N)$ in
$H_j$.  Expand $\varphi$ in the Laurent series,
\begin{equation}\label{Loran}
\varphi=\frac{\varphi_k}{z_1^k}+\frac{\varphi_{k-1}}{z_1^{k-1}}+\cdots,
\end{equation}
where $\varphi_k$ are functions in $z_2,...,z_N$ holomorphic at
$(\varepsilon_2,...,\varepsilon_N)$. We assume that
\linebreak
 $\varphi_k
(\varepsilon_2,...,\varepsilon_N)\neq 0$. For $v\in V$, consider the Laurent
expansion for both sides of the equality $\Phi
(\partial_{z_1}v)=\partial_{z_1}(\Phi(v))$.  We have
$$
\Phi(\partial_{z_1}v)=\Phi\left(\frac{A^j(v)}{z_1}+\cdots\right)=
\frac{1}{z_1}\varphi(A^j(v))+...=\frac{\varphi_k(A^j(v))}{z_1^{k+1}}+...
\ ,
$$
and
$$
\partial_{z_1}(\Phi(v))=\partial_{z_1}(\varphi(v))=
\frac{(B^j-k)\left(\varphi_k(v)\right)}{z_1^{k+1}}+\cdots\ .
$$
We get $\varphi_k(A^j(v))\,=\,(B^j-k)\left(\varphi_k(v)\right)$.
Since eigenvalues of $A^j$ and of $B^j$ are in the same weakly
non-resonant set $G^j$, we conclude that $k=0$.  Thus the rational
function $\varphi$ has no poles and is a polynomial function on $X$.

Let $\C[X]$ be the ring of polynomials on $X$ with complex coefficients.
Denote $\widetilde{V}=\CC[X]\cdot V$, $\widetilde{W}=\CC[X]\cdot W$.
We proved  that $\Phi(\widetilde{V})\subset \widetilde{W}$.

Let $O\in X$ be the center of the arrangement. Let $z_1,\dots,z_n$ be
affine coordinates on $X$ with center at $O$. Let $\eta$ be the
associated Euler vector field, $\eta\,=\,z_1\partial_{z_1}+
... +z_N\partial_{z_N}$. For any $v\in V$, we have\ $\eta(v) =
\sum_{j\in J(\A)} A^j(v)$.  Let $\widetilde{V} = \oplus_\mu
\widetilde{V}_\mu,$ \ $\widetilde{W} = \oplus_\mu \widetilde{W}_\mu$
be the decompositions into eigenspaces of $\eta$.
It is easy to see that the summations are over
$\mu\in\widetilde{G}_\ee+\ZZ_{\geq 0}$.
We have $\Phi(\widetilde{V}_\mu) \subset \widetilde{W}_\mu$ for any $\mu$.
Since the set $\widetilde{G}_\ee$ is weakly non-resonant, we have equalities
$$
V\ =\ \ooplus_{\mu\in\widetilde{G}_\ee}\widetilde{V}_\mu\,,
\qquad
W\ =\ \ooplus_{\mu\in\widetilde{G}_\ee}\widetilde{W}_\mu\ ,
$$
and $\Phi(V)\subset W$. Thus $\varphi$ is constant.
The condition $[\Phi,\xi](v)=0$ for constant vector fields $\xi$ imply
 that $\varphi$ is a morphism of quivers and $\Phi\, =\, E^0(\varphi)$.
\hfill{$\square$}

\subsection{Proof of Proposition \ref{theorem1}}\label{prth1}
Part (i) of the proposition follows from the description in Section
 \ref{section3.3} of global sections of a quiver $\D_{X^n}$-module.
 Namely, the quiver $\D_{X^n}$-module $E^n\VA_{\G^n}$ has a
 finite-dimensional subspace
 $\overline{V}=\oplus_{\a\in\v(\G^n)}\Oa\ot V_\a$ of the space of its
 global sections. If $\varphi:\VA_{\G^n}\to{\mathcal W}_{\G^n}$ is a
 nonzero morphism of quivers, then the associated morphism
 $E^n\varphi:E^n\VA_{\G^n}\to E^n{\mathcal W}_{\G^n}$ of
 $\D_{X^n}$-modules is nonzero on the space $\overline{V}$. Thus the
 functor $E^n$ is faithful.
\smallskip

The proof of part (ii) is by induction on $n$.  The proof is based on
 the Beilinson-Kashiwara-Malgrange gluing construction.

First we will recall the gluing construction, following
\cite{Kash2}. Then we will introduce special open subsets
$X^{n,\b}\subset X^n$, on which the gluing procedure will be
performed.  For a level $n$ quiver $\VA_{\Ga^{n}}$, we will restrict
the quiver $\D_{X^{n}}$-module $E^{n}\VA_{\Ga^{n}}$ to an open set
$X^{n,\b}$ and describe, how this $\D_{X^{n,\b}}$-module can be
obtained by gluing quiver $\D$-modules of level $n-1$.  We will note
that the sheaf axioms imply that the category $\MD_{X^n}$ of
$\D_{X^n}$-modules is equivalent to a category of fibered products of
$\MD_{X^{n,\b}}$-modules.  In this way we will describe the
$\D_{X^{n}}$-module $E^{n}\VA_{\Ga^{n}}$ by means of quiver
$\D$-modules of level $n-1$ quivers. With such a description we will
perform the induction step.

\subsubsection{Kashiwara-Malgrange filtrations and the gluing construction}
\label{Kashiwara}
Recall that a subset $\CG\subset \CC$ is called a {\it non-resonant section} of
$\CC$, if $\CG$ contains zero and for any $a\in\CC$, the intersection
$\CG\cap (a+\ZZ)$ consists of a single point.

Let $Y$ be a smooth complex algebraic variety and $Z\subset Y$ a
smooth complex algebraic subvariety. We assume that $Z$ is a
principal divisor, defined by the equation $f=0$ for some regular function $f$.
Let $I\subset {\cal O}_Y$ be the sheaf of ideals
 of functions vanishing on  $Z$ and
  $F(\D_{Y})$ a corresponding (see Section \ref{section4.6}
decreasing filtration of the sheaf $\D_Y$:
$$
F^k(\D_{Y})\ =\ \{P\in\D_{Y}|\ P(I_\b^j)\subset I_\b^{j+k}\quad
\mbox{\rm for any }\ j\}\ .
$$
 Let
$\theta$ be a vector field in $Y$ tangent to $Z$ and acting on $I/I^2$
as the identity.  We suppose that  $\theta$ is presented as a product
 $\theta=f\xi$, where $\xi$ is
a vector field on $Y$ such that $\xi(df)=1$.

Let $M$ be a $\D_{Y}$-module.
{}For any good filtration $F(M)$ of a $\D_Y$-module $M$, the quotient
 spaces $\quad {\rm gr}\,F^k(M) = F^k(M) /F^{k+1}(M)$ have the
 structure of $\D_Z$-modules. and the function $f$ and the vector field
 $\xi$
define the $\D_Z$-module maps
 $$
{\rm gr}\,f: {\rm gr}\,F^k(M)\to {\rm gr}\,F^{k+1}(M)\qquad {\rm and}\qquad
{\rm gr}\,\xi : {\rm gr}\,F^{k}(M)\to {\rm gr}\,F^{k-1}(M)\ ,
$$
such that ${\rm gr}\, f\, \cdot \,{\rm gr}\, \xi\,= \,{\rm gr} \, \theta:\, {\rm
gr}\,F^k(M)\to {\rm gr}\,F^{k}(M)$.

Suppose a non-resonant section
$\CG\subset\CC$ is chosen. Let $F_{\CG}(M)$ be the unique good
filtration of $M$, satisfying the condition {(v)} of Section \ref{section4.6}.
 Denote $\Psi^{(k)}_{\CG}(M)= {\rm
gr}\,F_{\CG}^k(M)$.   Then for all $k\geq 0$, the $\D_Z$-modules
$\Psi^{(k)}_{\CG}(M)$ and the maps ${\rm gr}\,\theta_{}:
\Psi^{(k)}_{\CG}(M)\to \Psi^{(k)}_{\CG}(M)$ depend on the restriction
$M|_{Y\setminus Z}$ only.  We use this remark in notations and set
$\Psi^{(0)}_{\CG}(M|_{Y\setminus Z})=\Psi^{(0)}_{\CG}(M)$.  Denote by
$a$ the map $a={\rm gr}\,\xi
|_{\Psi^{(0)}_{\CG}(M)}:\Psi^{(0)}_{\CG}(M)\to \Psi^{(-1)}_{\CG}(M)$
and by $b$ the map $b={\rm gr}\,f|_{\Psi^{(-1)}_{\CG}(M)}:
\Psi^{(-1)}_{\CG}(M)\to \Psi^{(0)}_{\CG}(M)$.

The Beilinson-Malgrange-Kashiwara gluing theorem \cite{Be1}, \cite{Kash2},
\cite{M} states that for a
 given non-resonant section $\CG\subset \CC$  the
assignment
$$M\to (M|_{Y\setminus Z}, \Psi^{(-1)}_{\CG}(M), a,b)$$
establishes an equivalence of the category
${\cal M}_Y^{hrs}\subset\MD_Y^{hol}$ of holonomic
$\D_Y$- modules with regular singularities and  the category
 of quadruples $ (\tilde{M},   N, a, b)$, where
$\tilde{M}\in{\cal M}_{Y\setminus Z}^{hrs}$, $N\in{\cal M}_{ Z}^{hrs}$,
$a:\Psi^{(0)}_{\CG}(\tilde{M})\to N$ and $b:N\to \Psi^{(0)}_{\CG}(\tilde{M})$
 are $\D_Z$-module maps such
that $ba=\gr \theta:\Psi^{(0)}_{\CG}(j_*\tilde{M})\to
\Psi^{(0)}_{\CG}(j_*\tilde{M})$.
\medskip

\subsubsection{The spaces $X^{n,\b}$}
Take a vertex $\b\in\v(\G^n)$, such that $l(\b)=n$.
Let $X^{n,\b}\subset X^n$ be the following open subset of $X^n$,
$$
X^{n,\b}=X^{n-1}\ \bigcup \ X_\b=X^n\setminus
\left(\bigcup\limits_{\a\not=\b,\ l(\a)=n}X_\a\right)\, .
$$
The sets $X^{n,\b}$  cover  $X^n$,
$$
X^n=\bigcup\limits_{\b,\ l(\b)=n}X^{n,\b}\, ,
$$
and for any $\b_1,\b_2\in\v_n(\G)$ we have
$$
X^{n,\b_1}\bigcap X^{n,\b_2}=X^{n-1}\, .
$$
Thus the category $\MD_{X^n}$ of $\D_{X^n}$-modules is the fibered product
of the categories $\MD_{X^{n,\b}}$:
\begin{equation}\label{dfiber}
\MD_{X^n}\approx \MD_{X^{n,\b_1}}\times_{{\mathcal M}_{X^{n-1}}}
\MD_{X^{n,\b_1}}
\times_{{\mathcal M}_{X^{n-1}}}\cdots
\times_{{\mathcal M}_{X^{n-1}}}\MD_{X^{n,\b_m}}\, ,
\end{equation}
where the product is taken over all $\b_k\in\v_n(\G)$. An object
of the category
 $\MD_{X^{n,\b_1}}\times_{{\mathcal M}_{X^{n-1}}}\cdots
\times_{{\mathcal M}_{X^{n-1}}}\MD_{X^{n,\b_m}}$ is a collection
 $M_1\in \MD_{X^{n,\b_1}}$, $...$,
$M_m\in \MD_{X^{n,\b_m}}$, such that
 $M_1|_{X^{n-1}}=M_2|_{X^{n-1}}= ... =M_m|_{X^{n-1}}$.

\medskip

 Let  ${\mathcal F}_{n-1}=\{\ff{\a}\,|\ \a\in\v_{n}(\Ga)\}$ be a
  level $n-1$ vertex framing of the arrangement $\A$. The framing  determines an
 affine open subset $U_{\FF_{n-1}}^\b\subset X^{n,\b}$ of  $X^{n,\b}$:
\begin{equation}\label{UFnb}
U_{\FF_{n-1}}^{\b}=
X\setminus \bigcup\limits_{\a\in\v_{n}(\Ga),\ \a\not=\b} \{\ff{\a}=0\}\, .
\end{equation}
Subsets $U_{\FF_{n-1}}^{\b}$ form a covering of $X^{n,\b}$ when
$\FF_{n-1}$ varies.

\medskip

\subsubsection{Gluing data for $\Ga^n$-quivers}
{}Fix a vertex $\b\in\v_n(\Ga)$. Let
 $ \ff{\b}\in F_\b$ be a generic affine function vanishing on $X_\b$.
In the notations of Section \ref{Kashiwara}, put
  $Y=X^{n,\b}$ and $Z=Y\cap H$, where $H=\{\ff{\b}=0\}$.

Intersections $H_j\cap H$, $j\in J(\A),$ determine in $H$ an
arrangement, denoted by $\A_{H}$, which  equips $H$ with a structure
 of a stratified
space. The space $Z$ coincides with the principal
 level $n-1$ open subset of $H$:
$$Z=H^{n-1}= H\setminus\bigcup_{{\rm codim}_H H_\a \geq n}H_\a \, .$$

The strata of $Z$ are  the intersections $H\cap X_\a$, where
$\a\in\v(\Ga)$, such that $l(\a)<n$.

 The graph  of adjacencies $\G^{n,\b}$ of the stratified
space $Z$ is  a certain contraction of
the  graph $\G^{n-1}$.

 More precisely, let $\v_{n-1}^\b(\Ga)$ be the set of classes of isomorphisms
 of  elements of  $\v_{n-1}(\Ga)$ with respect to the following equivalence
relation: $\a_1\approx \a_2$ if $\a_1\succ\b$ and $\a_2\succ\b$, or
 $\a_1=\a_2$. Thus, the set
$\v_{n-1}^\b(\Ga)$ consists of the vertices $\a\in\v_{n-1}(\Ga)$, such that
$\a\not\succ\b$ and of the equivalence class $\overline{\b}$ of the
 vertices $\a\in\v_{n-1}(\Ga)$, such that
$\a\succ\b$.

Let  $\tilde{\v}_n^\b(\Ga)$ be the subset of ${\v}_n(\Ga)$, which
consists of the elements $\b'\in {\v}_n(\Ga)$, such that the stratum
 $X_{\b'}$ is not equal to the stratum $X_\b$ and  is not parallel to it
(not obtained by a parallel translation). Let
 $\v_{n}^\b(\Ga)$  be the set of classes of isomorphisms
 of  elements of  $\tilde{\v}_{n}^\b(\Ga)$ with respect to the
following equivalence
relation: $\a_1\approx \a_2$ if $\a_1=\a_2$ or
the  intersections
$\X_{\a_1}\cap \X_\b$ and $\X_{\a_2}\cap \X_\b$ coincide and have
codimension one in $\X_\b$:
$$\X_{\a_1}\cap \X_\b=\X_{\a_1}\cap \X_\b\qquad {\rm and} \qquad
\codim_{\X_\b}\X_{\a_1}\cap \X_\b=1.$$
In these notations the set of  vertices of  $\G^{n,\b}$ is the union
$$\v(\Ga^{n,\b})=\v(\Ga^{n-2})\cup I_{n-1}^\b(\Ga).$$
Define the edges of $\G^{n,\b}$. The vertices $\a_1,\a_2\in\v(\Ga^{n,\b})$,
 such that $\a_1\not=\overline{\b}$,
 $\a_2\not=\overline{\b}$ and $\a_1\not=\a_2$, are connected by an edge
 in $\G^{n,\b}$ if and only if they are connected by an edge in $\Ga^{n-1}$.
The vertices
$\overline{\b}$ and  $\a\in\v(\G^{n,\b})$, $\a\not=\b$  define an edge in
 $\G^{n,\b}$
if there exists $\a_1\succ\b$, such that
$\a_1$ and $\a$ are connected by an edge in $\G^{n-1}$.

The graph  $\G^{n,\b}$ has loops
 $(\overline{\a}_1,\overline{\a}_1)^{\overline{\a}_2}$,
 where $\overline{\a}_1$ $\in $
$\v_{n-1}^\b(\Ga)$, $\overline{\a}_2$ $\in$
$\v_{n}^\b(\Ga)$ and there exist representatives
${\a}_1$ of $\overline{\a}_1$, ${\a}_1\in \v_{n-1}(\Ga)$, and
$\a_2$ of $\overline{\a_2}$,  ${\a_2}\in \tilde{\v}_{n}^\beta(\Ga)$,
such that ${\a}_1\succ{\a_2}$.

Let $\VA_{\Ga^{n-1}}=\{ V_\a, A_{\a_1,\a_2},
 A^{\a_1}_{\a_2}\}\in \Qui_{\G^{n-1}}$ be a quiver
 of $\G^{n-1}$. Define a new $\Ga^{n,\b}$-quiver
$\overline{\psi}^{(0)}(\VA_{\Ga^{n-1}}) =
\{ W_\a, B_{\a_1,\a_2},B_{\a_1}^{\a_2} \}$
 in the following way.

Set $W_\a=V_\a$ for all $\a\in\v_{n-1}(\Ga)$, and
$W_{\overline{\b}}=\ooplus\limits_{\a,\ \a\succ\b}V_\a$.

 Set
$B_{\a_1,\a_2}=A_{\a_1,\a_2} :\, W_{\a_2}\to W_{\a_1}$
if both
$\a_1$ and $\a_2$ are different from $\overline{\b}$.
Set $B_{\overline{\b},\a_2} = \oplus_{\a,\  \a\succ\b}A_{\a,\a_2}:$
$ W_{\a_2}\to W_{\overline{\b}}$ and
$B_{\a_1,\overline{\b}} = \oplus_{\a,\  \a\succ\b}A_{\a_1,\a}:$
$ W_{\overline{\b}}\to W_{\a_1}$.

Let $({\a}_1,{\a}_1)^{\overline{\a}_2}$ be a loop of  $\Ga^{n,\b}$ where
$\a_1\in \v_{n-1}^\b(\Ga)$,  $\overline{\a}_2\in \v_{n}^\b(\Ga)$ and
 $\a_1\not=\overline{\b}$. Define the map $B_{\a_1}^{\overline{\a}_2}:
W_{\a_1}\to W_{\a_1}$ by the formula
$$B_{\a_1}^{\overline{\a}_2}=\sum_{j=1}^{m_{\overline{\a}_2}}
A_{\a_1}^{\tilde{\a}_j}\ ,$$
where $\tilde{\a}_1,...,\tilde{\a}_{m_{{\overline{\a}}_2}} \in
\tilde{\v}_{n}^\beta(\Ga)$
are all representatives
of the equivalence class $\bar {\a}_2\in  \v_{n}^\b(\Ga)$.

Let $(\overline{\b},\overline{\b})^{\overline{\a}_2}$ be a loop of
$\G^{n,\b}$. Let $\tilde{\a}_1,...,\tilde{\a}_{m_{\overline{\a}}} \in
\tilde{\v}_{n}^\beta(\Ga)$
be all representatives
of the equivalence class $\bar{\a}_2\in  \v_{n}^\b(\Ga)$, and
${\a}_1,...,{\a}_{m} \in
{\v}_{n-1}(\Ga)$
 all representatives
of the equivalence class $\overline{\b}$.
Define the map
$B_{\overline{\b}}^{\overline{\a}}:
W_{\overline{\b}}\to W_{\overline{\b}}$ by the formula
$$B_{\overline{\b}}^{\overline{\a}_2}=
\sum_{j_1=1}^m\sum_{j_2=1}^{m_{\overline{\a}}}
A_{\a_{j_1}}^{\tilde{\a}_{j_2}}-
\sum_{j_1,j_2=1,...,m,\, j_1\not=j_2}A_{\a_{j_1},
\delta(\a_{j_1}, \a_{j_2})}A_{\delta(\a_{j_1}, \a_{j_2}),\a_{j_2}} \, .$$
 The vertex $\delta(\a_{j_1}, \a_{j_2})$
is defined by the following rule. Suppose there exists a vertex
$\gamma\in\{\tilde{\a}_1,...,\tilde{\a}_{m_{\overline{\a}}}\}$, such that
$\a_{j_1}\succ\gamma$ and $\a_{j_2}\succ\gamma$. Then we set
$\delta(\a_{j_1}, \a_{j_2})$ to be the unique vertex  of $\G^{n-1}$, such
that $\delta(\a_{j_1}, \a_{j_2})\succ\a_{j_1}$ and
$\delta(\a_{j_1}, \a_{j_2})\succ\a_{j_2}$. Otherwise the corresponding
term is supposed to be zero.

 We define a quiver morphism $\overline{T}_\b$:
$\overline{\psi}^{(0)}(\VA_{\Ga^{n-1}}) \to
\overline{\psi}^{(0)}(\VA_{\Ga^{n-1}})$. Morphism $\overline{T}_\b$
is a collection of linear operators
$\overline{T}_\b|_{W_\a}: W_\a\to W_\a$, $\a\in\v(\G^{n,\b})$, commuting
with the maps $B_{\a_1,\a_2}$ and $B_{\a_1}^{\a_2}$. The operators
$\overline{T}_\b|_{W_\a}$ are zero for all $\a$ except $\a=\overline{\b}$,
 and
\begin{equation}\label{mm}
\overline{T}_\b|_{W_{\overline{\b}}}=\sum_{j=1}^m
A_{\a_j}^\b-
\sum_{j_1,j_2=1,...,m,\, j_1\not=j_2}A_{\a_{j_1},
\delta(\a_{j_1}, \a_{j_2})}A_{\delta(\a_{j_1}, \a_{j_2}),\a_{j_2}} \, .
\end{equation}
 Here ${\a}_1,...,{\a}_{m} \in
{\v}_{n-1}(\Ga)$
are all representatives
of the equivalence class $\overline{\b}$.
The vertex $\delta(\a_{j_1}, \a_{j_2})$ is the
  unique vertex (if exists) of $\G^{n-1}$, such
that $\delta(\a_{j_1}, \a_{j_2})\succ \a_{j_1}$ and
$\delta(\a_{j_1}, \a_{j_2})\succ\a_{j_2}$.

Let $\VA_{\Ga^{n}}$
 be a quiver of $\G^{n}$.
We set ${\psi}^{(0)}(\VA_{\Ga^n})$ $=$
$\overline{\psi}^{(0)}(\j_{n,n-1}^*\VA_{\Ga^n}).$

\begin{lemma}\label{lemma10}${}$
For any level $n$ quiver $\VA_{\G^n}$ the linear map
$\overline{T}_\b|_{W_{\overline{\b}}}$
 coincides with
the operator
$T_\b$, defined in \rf{2.22}.
\end{lemma}
{\it Proof}. The statement follows from relations (ii),
Section \ref{Seclevn}.\hfill{$\square$}
\medskip

The closure $X'=\X_\b$ of the stratum $X_\b$ is an affine plane.
The arrangement $\A$ induces in $X'$ an arrangement denoted by $\A_\b$. The stratum
$X_\b$ is the complement in $X'$ to the union of hyperplanes of the arrangement
$\A_\b$. The hyperplanes $\X'_\g$ of the arrangement $\A_\b$ are
closures of the strata $X_\g\subset\X_\b$ of codimension one in $\X_\b$.
They are in one-to-one
correspondence with elements of the set $I^\b_n$: to each equivalence
class $\{\a_i\}$ of $I^\b_n$ we attach a hyperplane $\X_{\a_i}\cap \X_\b$ in
$X'=\X_\b$.

 Denote by $\G^\b$ the level zero graph of the arrangement $\A_\b$.
We have a map  of graphs $i_\b: \G^\b\to \G^{n,\b}$, which sends the vertex
$\ee$ of $\G^\b$ to the vertex $\overline{\b}$ of $\G^{n,\b}$ and the loop
$(\ee,\ee)^\gamma$ of $\G^\b$ to the  loop
 $(\overline{\b},\overline{\b})^{i_\b(\gamma)}$
of $\G^{n,\b}$, where $i_\b(\gamma)\in I_n^\b$ is characterized by the condition that
$\X_\a\cap \X_\b=\X'_\gamma$ for any representative
$\a$ of $i_\b(\gamma)$.

{} To any
 $\G^\b$-quiver $\VA'=\{V'_\ee, {A'}_{\ee}^\gamma\}$ we attach
 the following $\G^{n,\b}$-quiver $i_\b(\VA)$ $=$
 $\mathcal{U}=\{U_\a, C_{\a_1,\a_2},{C}^{\a_2}_{\a_1}\}$.
The spaces  $U_\a$ of the quiver ${\mathcal U}$  are all zero-dimensional except
for
$\a=\overline{\b}$. Set $U_{\overline{\b}}$ $=$ $V'_\ee$. All operators
$C_{\a_1,\a_2},{C}^{\a_2}_{\a_1}$ are the zero operators except for
operators
${C}_{\overline{\b}}^{i_\b(\gamma)}$. We set
${C}_{\overline{\b}}^{i_\b(\gamma)}= {A'}_{\ee}^\gamma: U_{\overline{\b}}
\to U_{\overline{\b}}$.

Let $\VA_{\Ga^{n}}=\{ V_\a, A_{\a_1,\a_2},
 A^{\a_1}_{\a_2}\}\in \Qui_{\G^{n}}$ be a quiver
 of $\G^{n}$. Define a quiver
$\overline{\psi}_\b^{(-1)}(\VA_{\Ga^n})=\{V'_\ee, {A'}_\ee^\gamma\}$
 of $\G^\b$ as follows.
We set $V'_\ee=V_\b$ and ${A'}_\ee^\gamma=A_\b^\gamma$.
Define a  quiver
${\psi}_\b^{(-1)}(\VA_{\Ga^n})=\{V'_\ee, {A'}_\ee^\gamma\}$
 of $\G^{n,\b}$ as $i_\b(\overline{\psi}_\b^{(-1)}(\VA_{\Ga^n}))$.

Define two morphisms of the constructed quivers:
 $\overline{a}_\b: \psi_\b^{(0)}(\VA_{\Ga^n})\to
\psi_\b^{(-1)}(\VA_{\Ga^n})$ and $\overline{b}_\b:
\psi_\b^{(-1)}(\VA_{\Ga^n})\to
\psi_\b^{(0)}(\VA_{\Ga^n})$. To define these morphisms,
it is sufficient to define linear maps
$\,\overline{a}_\b|_{W_{\overline{\b}}}$: $W_{\overline{\b}}\to
U_{\overline{\b}}$ and $\,\overline{b}_\b|_{U_{\overline{\b}}}$:
$U_{\overline{\b}}\to W_{\overline{\b}}$. We set
$$\overline{a}_\b(v_\a)=A_{\b,\a}(v_\a)\, ,\qquad
\overline{b}_\b(v_\b)=\sum_{\a\in\v_{n-1}(\Ga),\ \a\succ\b}
A_{\a,\b}(v_\b)\, ,$$
for any $v_\b\in V_\b$ and
$v_\a\in V_\a$.
We have the equality
 $$\overline{b}_\b\overline{a}_\b= \overline{T}_\b\ .$$

Let $Q^{n,\b}$ be the category of quadruples $(\VA_{\G^{n-1}},{\mathcal
U}_{\G^\b}, a_\b,b_\b)$, where $\VA_{\G^{n-1}}$ is a quiver of $\G^{n-1}$,
${\mathcal U}_{\G^\b}$ is a quiver of $\G^\b$, $a_\b$ and $b_\b$
are morphisms of $\G^{n,\b}$-quivers,
${a}_\b: \overline{\psi}_\b^{(0)}(\VA_{\Ga^{n-1}})\to
i_\b({\mathcal U}_{\G^\b})$ and ${b}_\b:i_\b({\mathcal U}_{\G^\b})
\to \overline{\psi}_\b^{(0)}(\VA_{\Ga^{n-1}})$, such that the composition
$ba:$
$\overline{\psi}^{(0)}(\VA_{\Ga^{n-1}}) \to
\overline{\psi}^{(0)}(\VA_{\Ga^{n-1}})$ is equal to
$\overline{T}_\b$, see \rf{mm}. A morphism between objects
$(\VA_{\G^{n-1}},{\mathcal U}_{\G^\b}, a_\b,b_\b)$ and
$(\VA'_{\G^{n-1}},{\mathcal U}'_{\G^\b}, a'_\b,b'_\b)$ of $Q^{n,\beta}$
is a pair $(\xi, \eta_\b)$, where
 $\xi\in \Hom_{Qui_{\G^{n-1}}}(\VA_{\G^{n-1}}, \VA'_{\G^{n-1}})$,
 $\eta_\b\in\Hom_{Qui_{\G^\b}}({\mathcal U}_{\G^\b},{\mathcal U}'_{\G^\b})$,
such that $a'_\b\overline{\psi}^{(0)}(\xi)=i_\b(\eta_\b) a_\b$, and
$\overline{\psi}^{(0)}(\xi) b_\b=b'_\b i_\b(\eta_\b)$.

\medskip
Now we assume that the index $\b$ is not fixed, but range over the set
 $I_n(\G)$.

Let $Q^{n}$ be the category of collections $\{\VA_{\G^{n-1}},{\mathcal
U}_{\G^\b}, a_\b,b_\b\ |\
\b\in\v_n(\G)\}$, where $\VA_{\G^{n-1}}$ is a quiver of $\G^{n-1}$,
${\mathcal U}_{\G^\b}$ is a quiver of $\G^\b$, $a_\b$ and $b_\b$
are morphisms of $\G^{n,\b}$-quivers,
${a}_\b: \overline{\psi}_\b^{(0)}(\VA_{\Ga^{n-1}})\to
i_\b({\mathcal U}_{\G^\b})$ and ${b}_\b:i_\b({\mathcal U}_{\G^\b})
\to \overline{\psi}_\b^{(0)}(\VA_{\Ga^{n-1}})$, such that the composition
$ba:$
$\overline{\psi}^{(0)}(\VA_{\Ga^{n-1}}) \to
\overline{\psi}^{(0)}(\VA_{\Ga^{n-1}})$ is equal to
$\overline{T}_\b$.
 A morphism between objects
%
$\{\VA_{\G^{n-1}}, {\mathcal U}_{\G^\b}, a_\b, b_\b\ |$\ $\b \in \v_n(\G)\}$ and
$\{\VA'_{\G^{n-1}},{\mathcal U}'_{\G^\b}, a'_\b,b'_\b\ |\ \b\in\v_n(\G)\}$
is a collection $\{\xi, \eta_\b\ |\ \b\in\v_n(\G)\}$, where
 $\xi\in \Hom_{Qui_{\G^{n-1}}}(\VA_{\G^{n-1}}, \VA'_{\G^{n-1}})$, and
 $\eta_\b\in\Hom_{Qui_{\G^\b}}({\mathcal U}_{\G^\b},{\mathcal U}'_{\G^\b})$,
such that $a'_\b\overline{\psi}^{(0)}(\xi)=i_\b(\eta_\b) a_\b$, and
$\overline{\psi}^{(0)}(\xi) b_\b=b'_\b i_\b(\eta_\b)$  for any
$\b\in\v_n(\G)$.

{}For any $\b\in\v_n(\G)$, $A\in Q^{n,\b}$,  $A=(\VA_{\G^{n-1}},{\mathcal
U}_{\G^\b}, a_\b,b_\b)$, denote by $\nu(A)$ the
$\G^{n-1}$-quiver $\VA_{\G^{n-1}}$. Define the fibered product
$\widetilde{Q}^n$ of the categories $Q^{n,\b}$,
\begin{equation}
\label{fiber}
\widetilde{Q}^n=Q^{n,\b_1}\times_{Qui_{\G^{n-1}}}Q^{n,\b_2}
\times_{Qui_{\G^{n-1}}}\cdots
 \times_{Qui_{\G^{n-1}}}Q^{n,\b_m}\, ,
\end{equation}
where the product is taken over all $\b_k\in\v_n(\G)$. An object
of the category
$\widetilde{Q}^n$ is the collection of quadruples
 $A_1\in Q^{n,\b_1},A_2\in Q^{n,\b_2}$, $...,$ $A_m\in Q^{n,\b_m}$, such that
$\nu(A_1)=\nu(A_2)=\cdots =\nu(A_m)$.
\begin{lemma}\label{lemmafiber}
The categories $Q^n$ and $\widetilde{Q}^n$ are equivalent.
\end{lemma}\hfill{$\square$}

Define a functor
$\mu: Qui_{\G^n}\to Q^{n}$ as follows. For any $\VA_{\G^n}\in Qui_{\G^n}$
we set
$$\mu(\VA_{\G^n})=\{\j_{n,n-1}^*\VA_{G^n},
\overline{\psi}_\b^{(-1)}(\VA_{\Ga^n}),\overline{a}_\b, \overline{b}_\b\
|\ \b\in\v(\G^n)\}\ .$$

\begin{lemma}\label{lemma10a}
The functor $\mu: Qui_{\G^n}\to Q^{n}$ is an
equivalence of categories.
%
\end{lemma}
{\it Proof}. In order to prove that $\mu$ is a functor, one should check
the commutativity of certain diagrams of maps between quivers.
 A direct check shows that the relations
 (ii)-(v) of Section \ref{Seclevn} imply the required equalities.
The functor $\mu$ is an equivalence of categories, since the spaces
$V_\a$ and the maps $A_{\a,\b}$ of the quiver $\VA_{\G^n}$
can be reconstructed from $\mu(\VA_{\G^n})$.


\medskip

\subsubsection{Gluing data for  quiver $\D_{X^{n,\b}}$-modules}
\label{glue}
Fix a vertex $\b\in\v_n(\G)$.
 Choose
 a generic constant vector field $\xi_\b$ such that $\mylangle \xi_\b,
 \ff{\b}\myrangle=1$
and set, following the notations of the previous section,
$f=\ff{\b}, \xi=\xi_\b, \theta=\ff{\b}\xi_\b$.
We apply the construction of  Section \ref{Kashiwara}
to  $Y=X^{n,\b}, Z=Y\cap H,\theta$,
and the $D_{Y}$-module
$E^n_Y\VA_{\Ga^n}=E^n\VA_{\Ga^n}(Y)$, which is the restriction of the $D_{X^n}$-module
$E^n\VA_{\Ga^n}$ to the open subset $Y=X^{n,\b}$.
We denote by
$i_\b$ the closed embedding of $X_\b$ to $Z$.

Denote by $\V_\pm$ the following finite-dimensional subspaces of the space
of global sections of  the $D_{Y}$-module $E^n_Y\VA_{\Ga^n}$:
$$\V_+=\ooplus\limits_{\a\in\v(\Ga^{n-1})}\Oa\ot V_\a,\qquad
 \V_-=\Ob\ot V_\b.$$
Let $F^k(\D_Y)$ be the filtration of $\D_Y$,
associated to the sheaf of ideals, generated by $f=\ff{\b}$ as in
Section \ref{Kashiwara}. Set
\begin{equation}
F^k(E^n_Y \VA_{\Ga^n})=\left\{
\begin{array}{ll}\label{filtration1}
F^k(\D_Y)\cdot\V_+,& k\geq 0\\
F^{k+1}(\D_Y)\V_-+F^k(\D_Y)\cdot \V_+,& k<0\, .
\end{array}
\right.
\end{equation}
Formula \rf{filtration1} determines a filtration of the sheaf
$E^n  \VA_{\Ga^n}$.
\begin{lemma}\label{lemma11}${}$
\begin{enumerate}
\item[ (i)] The filtration \rf{filtration1} is good with respect to
 $F(\D_Y)$.
\item[ (ii)] The sheaf
$F^0(E^n_Y \VA_{\Ga^n})/F^1(E^n_Y \VA_{\Ga^n})$ of $\D_Z$-modules
is isomorphic to the $\D_Z$-module $E^{n-1} \psi_\b^{(0)} ( \VA_{\Ga^n})$.
\item[ (iii)] The sheaf
$F^{-1}(E^n_Y \VA_{\Ga^n})/F^0(E^n_Y \VA_{\Ga^n})$ of $\D_Z$-modules
is isomorphic to the quiver $\D_Z$-module $E^{n-1} \psi_\b^{(-1)}
( \VA_{\Ga^n})$.
\item[ (iv)]  $\D_Z$-module maps $\ {\rm gr}\,\ff{\b}:$ $\
F^{-1}(E^n_Y \VA_{\Ga^n})/F^0(E^n_Y \VA_{\Ga^n})$ $\to$
 $F^0(E^n_Y \VA_{\Ga^n})/F^1(E^n_Y\VA_{\Ga^n})\ $ and
${\rm gr}\,\xi_{\b}:$ $
F^{0}(E^n_Y \VA_{\Ga^n})/F^1(E^n_Y \VA_{\Ga^n})$ $\to$
 $F^{-1}(E^n_Y \VA_{\Ga^n})/F^0(E^n_Y \VA_{\Ga^n})$
are isomorphic to corresponding $ \D_Z$-module maps
$\ E^{n-1}\overline{b}_\b:\ E^{n-1} \psi_\b^{(-1)} ( \VA_{\Ga^n})\ \to \
E^{n-1} \psi_\b^{(0)} ( \VA_{\Ga^n}) $ and
$ E^{n-1}\overline{a}_\b:$\ $\ E^{n-1}  \psi_\b^{(0)}  ( \VA_{\Ga^n})$ $\to$
$E^{n-1} \psi_\b^{(-1)} ( \VA_{\Ga^n})$.
\item[ (v)] The minimal polynomial of the map ${\rm gr}\,\theta_\b$
coincides
with minimal polynomial of the map
$\overline{b}_\b\overline{a}_\b$: $\oplus_{\a\in\v_{n-1}(\Ga)}V_\a$ $\to$
$\oplus_{\a\in\v_{n-1}(\Ga)}V_\a$.
\end{enumerate}
\end{lemma}

The proof of statements { (ii)-(v)} is based on Proposition \ref{prop1}
and  consists of direct
 calculations with quiver $\D_Y$-modules on affine open sets
$U_{\FF_{n-1}}^{\b}$.
  \hfill{$\square$}

Let $\CG_\b$ be a non-resonant section of $\CC$.

Denote by $\MMD_{X^n}^{hrs}= \MMD_{X^n}^{hol}\cap\MD_{X^n}^{hrs}$
 the category of holomomic $\D_{X^n}$-modules with regular singularities,
 smooth along the strata of $X^n$.

Denote by  ${\widetilde{\mathcal N}}^{n,\b}(\CG_\b)$ the category
of quadruples  $ ({M},   N, c_\b, d_\b)$
 where $M\in \MMD_{X^{n-1}}^{hrs}$,
$N\in\MMD_{X_{\b}}^{hrs}$,\,
$c_\b:\Psi^{(0)}_{\CG_\b}({M})\to (i_\b)_* N$ and
$d_\b:(i_\b)_*N\to \Psi^{(0)}_{\CG_\b}({M})$
 are $\D_Z$-module morphisms such
that $d_\b c_\b={\rm gr}\, \theta_\b:\Psi^{(0)}_{\CG_\b}({M})\to
\Psi^{(0)}_{\CG_\b}({M})$ and all eigenvalues of $d_\b c_\b$ are in the set
 $\CG_\b$.
 The gluing construction, see Section \ref{Kashiwara} says
 that the functor
$H(\CG_\b): \MMD^{hrs}_{X^{n,\b}}\to {\widetilde{\mathcal N}}^{n,\b}(\CG_\b)$,  attaching to
every $M\in\MMD^{hrs}_{X^{n,\b}}$ the quadruple
$(M|_{X^{n-1}}, (i_\b)^*\Psi^{(-1)}_{\CG_\b}(M),
 c_\b= {\rm gr}\, \xi_\b, d_\b= {\rm gr}\, \overline{f}_\b)$,
establishes an equivalence of the two categories. For any object
${\widetilde N}\in {\widetilde{\mathcal N}}^{n,\b}(\CG_\b)$,
${\widetilde N}=({M},   N, c_\b, d_\b)$ denote by $\tilde{\nu}_\b({\widetilde N})$
the $\D_{X^{n-1}}$-module $M$.

Let $G_\b\subset \CG_\b$ be a non-resonant set.

Consider the category $Q^{n,\b}(G_\b)$ of the quadruples
 $Q=(\VA_{\G^{n-1}},{\mathcal
U}_{\G^\b}, a_\b,b_\b)$ in $Q^{n,\b}$, such that the eigenvalues of the
operator \rf{mm} are contained in $G_\b$. Then an assignment
$$F^{n,\b}:\,Q\mapsto (E^{n-1}(\VA_{\G^{n-1}}), E^0({\mathcal U}_{\G^\b}),
 E^{n-1}(a_\b), E^{n-1}(b_\b))\, ,$$
 defines a functor
$F^{n,\b}:Q^{n,\b}(G_\b)\to {\widetilde{\mathcal N}}^{n,\b}(G_\b)$.

Consider the category
 $Qui_{\G^n}(G_\b)\subset Qui_{\Ga^{n}}$ of level $n$  quivers
$\VA_{\Ga^n}$,
such that the eigenvalues of operator $T_\b$, defined in \rf{2.22}, are
contained in the set $G_\b$. By Lemma \ref{lemma11}, we have a
commutative diagram of functors:
\begin{equation}\label{CD}
\begin{CD}Qui_{\G^n}(G_\b)@>E^n(X^{n,\b})>>
\MMD_{X^{n,\b}}^{hrs}\\
@VV\mu_\b V @VV H(\CG_\b)V\\
Q^{n,\b}(G_\b)@>F^{n,\b}>>\widetilde{{\mathcal
N}}^{n,\b}(\CG_\b)\, .
\end{CD}
\end{equation}
The vertical lines in the diagram \rf{CD} represent equivalences of
categories.
\medskip

\subsubsection{The induction}

Suppose that for every vertex $\b\in\v_n(\G)$ there are chosen:
(a) a non-resonant section $\CG_\b$ of $\CC$; (b) a non-resonant  set
 $G_\b\subset\CG_\b$; (c)  a generic affine function
$\overline{f}_{\b}\in\Fb$ and a generic constant vector field $\xi_\b$ such that $\mylangle \xi_\b,
 \ff{\b}\myrangle=1$.  For every $\b\in\v_n(\G)$ denote by $Z_\b$
the intersection $X^{n,\b}\bigcap\{\overline{f}_{\b}= 0\}$, by $i_\b$
the closed embedding of $X_\b$ to $Z_\b$.

 Consider the category $\widetilde{{\mathcal
N}}^n(\{\CG_\b\})$, which consists of the collections
 $ \{{M},   N_\b, a_\b, b_\b\ |\ \b\in\v_n(\G)\}$,
 where $M$ is a holonomic
$\D_{X^{n-1}}$-module with regular singularities,
$N_\b$ is a holonomic
$\D_{X_\b}$-module with regular singularities,
$a:\Psi^{(0)}_{\CG_\b}({M})\to (i_\b)_*(N)$ and
$b:(i_\b)_*(N)\to \Psi^{(0)}_{\CG_\b}({M})$
 are $\D_{Z_\b}$-module maps such
that $ba={\rm gr}\, \theta_\b:\Psi^{(0)}_{\CG_\b}({M})\to
\Psi^{(0)}_{\CG_\b}({M})$ and all eigenvalues of $ba$
are in the set $\CG_\b$.
 The category $\widetilde{{\mathcal
N}}^n(\{\CG_\b\})$ is equivalent to the fibered product of the categories
$\widetilde{{\mathcal N}}^{n,\b}(\CG_\b)$. This fibered product consists
of the collections $({\widetilde N}_{\b_1},...,{\widetilde N}_{\b_m})$,
where
${\widetilde N}_{\b_k}\in \widetilde{{\mathcal N}}^{n,\b_k}(\CG_{\b_k})$
and $\widetilde{\nu}_{\b_1}({\widetilde N}_{\b_1})=...=
\widetilde{\nu}_{\b_m}({\widetilde N}_{\b_m})$.

Denote by  $Qui_{\G^n}(\{G_\b\})\subset Qui_{\Ga^{n}}$
 the category of level $n$  quivers
$\VA_{\Ga^n}$,
such that the eigenvalues of operator $T_\b$, defined in \rf{2.22}, are
contained in the set $G_\b$ for every $\b\in\v_n(\G)$.

Denote by $Q^{n}(\{G_\b\})$ the category  of  collections
 $Q=\{\VA_{\G^{n-1}},{\mathcal
U}_{\G^\b}, a_\b,b_\b\ |\ \b\in \v_n(\G)\}$ in $Q^{n}$, such that
 the eigenvalues of the
operator $\overline{T}_\b$ are contained in $G_\b$ for every
$\b\in\v_n(\G)$.

The assignment
$$F^{n}:\,Q\mapsto \{E^{n-1}(\VA_{\G^{n-1}}), E^0({\mathcal U}_{\G^\b}),
 E^{n-1}(a_\b), E^{n-1}(b_\b)\ |\ \b\in\v_n(\G)\}\, ,$$
 defines a functor
$F^{n}:Q^{n}(\{G_\b\})\to {\widetilde{\mathcal N}}^{n}(\{G_\b\})$.
We have a commutative diagram
\begin{equation}\label{CD2}
\begin{CD}Qui_{\G^n}(\{G_\b\})@>E^n>>\MMD_{X^{n}}^{hrs}\\
@VV\mu V @VV H(\{\CG_\b\})V\\
Q^{n}(\{G_\b\})@>F^{n}>>{\widetilde{\mathcal N}}^{n}(\{\CG_\b\})\, .
\end{CD}
\end{equation}
The functor $H(\{\CG_\b\})$ attaches to any $M\in \MMD_{X^{n}}^{hrs}$
the collection $(M_{X^{n-1}},\, (i_\b)^*\Psi^{(-1)}_{\CG_\b}(M),$
 $ c_\b= {\rm gr}\, \xi_\b, d_\b= {\rm gr}\, \overline{f}_\b\ |\
\b\in\v_n(\G))$. This functor is equivalent to the fibered product of functors
$H(\CG_\b)$ and thus is an equivalence of categories.
 We know also that the functor $\mu$ is an equivalence of categories.
 Thus if the functor $F^n$ is full, then
the functor $E^n$ is full as well.
\medskip

Consider the category $Qui_{\G^{n},\{G_\a,\widetilde{G}_\b,G_{\g,\d}\}}$
 associated with a collection of non-resonant sets $\{G_\a\}$,
$\a\in\v(\G^n)\setminus\ee$,  and
a collection of weakly non-resonant sets
$\widetilde{G}_\b$, $\b\in\v(\G^n)\setminus\ee$, and
 $\{G_{\g,\d}\}$, where $(\g,\d)\in E(\G^{n+1})$, $\g\succ\d$.
{}For the same collection of sets $\{G_\a,\widetilde{G}_\b,G_{\g,\d}\}$
 denote by
$Q^n_{\{G_\a,\widetilde{G}_\b,G_{\g,\d}\}}$ the category of sequences
$Q=\{\VA_{\G^{n-1}},{\mathcal
U}_{\G^\b}, a_\b,b_\b\ |\ \b\in \v_n(\G)\}$
in $Q^{n}$, such that
\begin{itemize}
\item[(a)]
for any $\b\in\v_n(\G)$ the eigenvalues of the operator $T_\b$,
 defined in \rf{2.22}, are contained in the set $G_\b$,
\item[(b)] the level $n-1$ quiver $\VA_{\G^{n-1}}$ is in the category
$Qui_{\G^{n-1},\{G_\a,\widetilde{G}_\b,G_{\g,\d}\}}$,
\item[(c)] for any $\b\in\v_n(\G)$ the (level zero) $\G^\b$-quiver
 ${\mathcal U}_{\G^\b}=\{U_\ee, C_\ee^\g\}$
is in the category $Qui_{\G^{\b},\{\widetilde{G}_\b,G_{\b,\g}\}}$, that is,
the eigenvalues of the operator $C_\ee^\g$ lie in the set $G_{\b,i_\b(\g)}$
and the eigenvalues of the operator $\sum_\g A_\ee^\g$ lie in the set
$\widetilde{G}_\b$.
\end{itemize}
The functor $\mu$ establishes an equivalence of categories
$Qui_{\G^{n},\{G_\a,\widetilde{G}_\b,G_{\g,\d}\}}$ and
$Q^n_{\{G_\a,\widetilde{G}_\b,G_{\g,\d}\}}$.
\medskip

Suppose that the statement (ii) of Proposition \ref{theorem1} is valid for
quivers of level less than $n$. Then the restriction of the
functor $E^{n-1}$ to the category
$Qui_{\G^{n-1},\{G_\a,\widetilde{G}_\b,G_{\g,\d}\}}$ is full and faithful.
Proposition \ref{theorem1a} implies that the restriction of the functor
$E^0$ to the category $Qui_{\G^{\b},\{\widetilde{G}_\b,G_{\b,\g}\}}$ is
 full and faithful. Thus the restriction of the functor $F^n$
to the category $Q^n_{\{G_\a,\widetilde{G}_\b,G_{\g,\d}\}}$  is full. This
proves that the restriction of the functor $E^n$ to the category
$Qui_{\G^{n},\{G_\a,\widetilde{G}_\b,G_{\g,\d}\}}$ is full. \hfill{$\square$}
\medskip

\bigskip

\subsection{Proof of Theorem \ref{prop4}} ${}$ \label{section6.6}
In this section we use the following notion.
Let $\widetilde{Y}$ be an affine nonsingular variety. Let $z_1,...,z_k$
be regular functions on $\widetilde{Y}$. Assume that
 $\widetilde{Z}\subset \widetilde{Y}$ defined by equations $z_1=0, \ldots
, z_k=0$ is nonsingular. Let $N$ be a $\D_{\widetilde{Y}}$-module.
Set $N^{\widetilde{Z}}=\{ n\in N\ |\ z_j^pn=0$ for $j=1,...,k$ and  $p\in\NN,
p\gg 1\}$. Then $N^{\widetilde{Z}}$ is
the maximal submodule of $M$ with  support in
$\widetilde{Z}$.
\medskip

Keep the notations of  Section \ref{prth1}.
Let  $\VA_{\G^{n-1}} = \{V_\a, A_{\a_1,\a_2},A_{\a_1}^{\a_2}\}$ be
a $\Ga^{n-1}$-quiver.
Fix a vertex $\b\in\v_n(\G)$, a level $n-1$ vertex framing $\FF_{n-1}$, and a
generic affine function $\ff{\b}\in F_\b$.
They define open affine subsets
$U_{\FF_{n-1}}\subset X^{n-1}$, see \rf{UFn}, and
$U=U_{\FF_{n-1}^\b}\subset X^{n,\b}$, see \rf{UFnb}.
Then $X_\b\subset U$ is a closed nonsingular subvariety of $U$.
 Denote by $M$ the restriction to $X^{n,\b}$ of
the sheaf $ j_{n,n-1,*}(E^{n-1}\VA_{\G^{n-1}})$
of $\D_{X^n}$-modules  and by
$M_U$ the restriction of $M$ to $U$.

Let $H=\{\ff{\b}=0\}$,
$H_U=U\cap H\subset Z=X^{n,\b}\cap H$, and $\widetilde{U}=U\setminus H_U $.
 Let  $\xi$ be a constant vector field on $X$, such that $\xi(d\ff{\b})=1$.
Denote by  $\theta $ the vector field $\theta= f\xi$.
Let $\CG_\b$ be a non-resonant section of $\CC$ and $G_\b\subset\CG_\b$
 a non-resonant set.
 Let $\Psi^{(k)}_{\CG_\b}(M)$ be the graded factors of the
Kashiwara-Malgrange
filtration of $M$, associated with the ideal, generated by the
 function $\ff{\b}$, and the set $\CG_\b$.

\medskip

\begin{lemma} \label{lemma1pr10}${}$
\begin{enumerate}
\item[ (i)] The image of the map $ba:\  \Psi^{(0)}_{\CG_\b}(M)\to
\Psi^{(0)}_{\CG_\b}(M)$
has  support in $X_\b\subset H_\b$,
$${\rm Im}\ ba\subset  \left(\Psi^{(0)}_{\CG_\b}(M)\right)^{X_\b}.$$
\item[ (ii)] The $\D_{Z}$-module $\Psi^{(-1)}_{\CG_\b}(M)$ has support in
$X_\b$,
$$\Psi^{(-1)}_{\CG_\b}(M)=\left(\Psi^{(-1)}_{\CG_\b}(M)\right)^{X_\b}.$$
\item[ (iii)] The gluing data for $\D_U$-module $M$ are isomorphic to
$$\left(M|_{\widetilde{U}},
\left(\Psi^{(0)}_{\CG_\b}(M)\right)^{X_\b},a,b\right),$$
where $b:\left(\Psi^{(0)}_{\CG_\b}(M)\right)^{X_\b}\to\Psi^{(0)}_{\CG_\b}(M)$ is the
canonical inclusion and $a:\Psi^{(0)}_{\CG_\b}(M)\to
\left(\Psi^{(0)}_{\CG_\b}(M)\right)^{X_\b}$ is the restriction of the monodromy
map $ba$ to $\left(\Psi^{(0)}_{\CG_\b}(M)\right)^{X_\b}$.
\end{enumerate}
\end{lemma}
{\bf Proof}. The statement { (i)} follows from statements {(ii)} and {
(iv)} of Lemma \ref{lemma11}, Section \ref{glue}.  Let us prove statement {
(ii)}. Note that since $\ff{\b}$ is generic, for any
$\a\in\v(\G^{n-1})$ there exists a vector field $\eta\in\Ta$, such
that $\mylangle\eta,\ff{\b}\myrangle\not= 0$. This implies that in the
open set $U$ any element in $ F^{-1}_{\CG_\b}(M_U)$ can be presented as a
fraction $\frac{m}{\ff{\b}}$, where $m\in F^{0}_{\CG_\b}(M_U)$.  Indeed, for
any element $n\in M_U$, we have the inclusion $\ff{\b}^kn\in
F^{0}_{\CG_\b}(M_U)$ for sufficiently big $k$, so the elements of
$F^0_{\CG_\b}(M_U)$ generate $M_U$ as a $D_U[\ff{\b}^{-1}]$-module, that is,
generate $M_U$ as a $\D_U$ module in  the open set
$\ff{\b}\not=0$. Then any element $n\in F^{-1}_{\CG_\b}(M_U)$ can be written
as $n=\ff{\b}^{-1}m'$, where $m'\in F^{0}_{\CG_\b}(M_U)$, or as $n=\nu m'$,
where $m'\in F^{0}_{\CG_\b}(M_U)$, and $\nu$ is a vector field with constant
coefficients, such that $\mylangle\nu,\ff{\b}\myrangle\not=0$. In the
last case we can find a constant field $\nu'$, such that
$\mylangle\nu',\ff{\b}\myrangle=\mylangle\nu,\ff{\b}\myrangle$ and the
element $n'=\nu'm'\in F^{-1}_G(M_U)$ satisfies the following
properties. It contains a smaller number of derivatives in its
expression in terms of vectors $\om_\a\ot v_\a$ and there is an
inclusion $n-n'\in F^{0}_{\CG_\b}(M_U)$. By induction we prove that there
exists $m\in F^{0}_{\CG_\b}(M_U)$, such that $n=f^{-1}m$. Let $f'$ be another
generic affine function, vanishing on $X_\b$. By the same arguments we
show that there exists an element $\tilde{m}\in F^{0}_{\CG_\b}(M_U)$, such
that $n={f'}^{-k}\tilde{m}$ for some integer $k>0$. Put
$x=\ff{\b}-f'$.  The equality
$$n=\frac{m}{\ff{\b}}=\frac{\tilde{m}}{(\ff{\b}+x)^k}$$
implies that $x^kn\in F^{0}_{\CG_\b}(M_U)$, that is
$$\overline{x}^kn=0\qquad {\rm in}\qquad \Psi^{-1}_{\CG_\b}(M_U)\, ,$$
where $\overline{x}=x|_{H_U}$. This proves  statement { (ii)} of
the lemma.

Denote by $\widetilde{j}$ the inclusion $\widetilde{j}: \widetilde{U}
\hookrightarrow U$ and by $\widetilde{M}$ the $\D_{\widetilde{U}}$-module
$\widetilde{M} = M|_{\widetilde{U}} = E^{n-1}\VA_{\G^{n-1}}(\widetilde{U})$.
Let $M'$ be the $\D_U$-module, given by the gluing data
$$\left(\widetilde{M},
\left(\Psi^{(0)}_{\CG_\b}(M)\right)^{X_\b},a,b\right),$$
where the maps $a$ and $b$ are described in statement { (iii)} of the
lemma. For any  $\D_{U}$-module $N$, whose gluing data satisfy the condition
that $\Psi^{(-1)}_{\CG_\b}(N)$ is supported at $X_\b$, we have the equality
\begin{equation}\label{tilde}
\Hom_{\D_{\widetilde{U}}}(\widetilde{j}^*(N),\widetilde{M})\simeq
\Hom_{\D_{U}}((N),M')\ .\end{equation}
Moreover, $M'$ is the unique $\D_U$-module, satisfying conditions
 \rf{tilde}.
Denote be  $j_\b$ the inclusion $j_\b:U_{\FF_{n-1}}\hookrightarrow U$.
 Since  $\ff{\b}$ is generic, we have the equality
$$\Hom_{\D_{\widetilde{U}}}(\widetilde{j}^*(N),\widetilde{M})\simeq
\Hom_{\D_{U_{\FF_{n-1}}}}(j_\b^*(N),
E^{n-1}\VA_{\G^{n-1}}|_{U_{\FF_{n-1}}}) .$$
This equality, relation \rf{tilde} and
statement {(ii)} of lemma imply statement { (ii)}.
\hfill{$\square$}
\begin{lemma}\label{lemmaker}
Let $E^{n}\VA_{\Ga^{n}}$ be the $\D_{X^n}$-module, attached to a $\Ga^n$-quiver
$\VA_{\Ga^{n}}$. Let $\b\in\v(\G^n)$, $l(\b)=n$, be a vertex  such that
 the operator $S_\b$, defined in \rf{Tb0} has no positive integer
eigenvalues. Then the $\D_{X^n}$-module $\left(E^{n}\VA_{\Ga^{n}}\right)^{X_\b}$
is generated by  vectors of the vector space $V'_\b\ot \Ob$, where
\begin{equation}\label{v1b}
V'_\b=\bigcap_{\a, \ \a\succ\b}{\rm Ker}\
A_{\a,\b}\subset V_\b\ .
\end{equation}
\end{lemma}
The proof of the lemma goes in two steps.
 First, we show that the lemma holds in a
special case of a quiver $\D$-module supported at a flag of a normal crossing
arrangement and then show how to reduce the general case
to the special.

Denote by $\A_{n,N}$ the arrangement of $n$ hyperplanes $\{z_j=0\}$,
 $j=1,..,n$,
 in the affine space $\CC^N$ with coordinates $z_k$, $ k=1,...,N$.
Let $\overline{\G}_{n}$ be the graph of the arrangement $\A_{n,N}$.
The vertices  of $\overline{\G}_{n}$ are labeled by subsets
$J\subset\{1,...,n\}$. We say that a $\overline{\G}_{n}$-quiver $\VA$ is
supported at the flag $\emptyset,\{1\}, \{1,2\},...,\{1,2,...,n\}$ of vertices
of $\overline{\G}_{n}$ if the vector spaces
$V_J$ of the quiver are nonzero only for $J$ equal to one of the sets
$J_0=\emptyset$, $J_1=\{1\}$, $J_2= \{1,2\}$, ... , $J_n=\{1,2,...,n\}$.

Let $\VA$ be a $\overline{\G}_{n}$-quiver supported at the flag
$\{1\}, \{1,2\},...,\{1,2,...,n\}$.
Choose basic exterior forms $\om_{J_k}\in\overline{\Om}_{J_k}$,
$\om_{J_k}=dz_{k+1}\wedge...\wedge dz_N$ and denote
$\overline{v}_{J_k}=\om_{J_k}\ot v_{J_k}$. Then the defining relations
\rf{2.3} and \rf{2.4} for the quiver $\D_X$-module $E\VA$ look as
follows:
\begin{equation}\label{normal}
\left\{\begin{array}{ll}
z_i\overline{v}_{J_k}=0,& i<k,\\
z_i\overline{v}_{J_k}=\overline{A_{J_{k-1},J_k}v}_{J_k},& i=k,
\end{array}\right.
\qquad
\left\{\begin{array}{ll}
\partial_{z_i}\overline{v}_{J_k}=0,& i>k+1,\\
\partial_{z_i}\overline{v}_{J_k}=\overline{A_{J_{k+1},J_k}v}_{J_k},& i=k+1.
\end{array}\right.
\end{equation}
 Denote by $\hat{Z}$ the stratum
$X_{J_n}$, given by  equations $z_1=...=z_n=0$.
Suppose $x\in (E\VA)^{\hat{Z}}$. We can assume that
\begin{equation}\label{x}
x=P(z_{n+1},...,z_N)\partial_{z_1}^{k_1}\cdots \partial_{z_n}^{k_n}
 (\overline{v}_{J_n}),
\end{equation}
where  $v_{J_n}\in V_{\{1,2,...,n\}}$.
We claim that  $x\in (E\VA)^{\hat{Z}}\Rightarrow v_{J_n}\in {\rm Ker}\ A_{J_{n-1,n}}$
if the operator $S_{J_n}=A_{J_{n,n-1}}A_{J_{n-1,n}}: $
$V_{J_{n}}\to V_{J_{n}}$ has no positive integer
eigenvalues.

Suppose we prove this claim for $x$, such that
$k_1+...+k_n<k$, see \rf{x}. Let us prove the claim for
$k_1+...+k_n=k$. Since $z_i( \overline{v}_{J_n})=0\,$ for $\,i=1,...,n-1$,
the condition $z_i^Mx=0$ for $M>0$ implies $\,k_i=0\,$ for $\,i=1,...,n-1$.
 Let $x=P(z_{n+1},...,z_N)
 \partial_{z_n}^{k_n}(\overline{v}_{J_n})$. If $k_n=0$, then
the condition $z_n^Mx=0$ implies $x\in {\rm Ker}\ A_{J_{n-1,n}}$ by
\rf{normal}. Let $k_n\not=0$. Then
$$z_nx=P(z_{n+1},...,z_N)\partial_{z_n}^{k_n-1}\overline{
\left( A_{J_{n,n-1}} A_{J_{n-1,n}}
-k_n\right)v}_{J_n}.$$
The condition $z_n^Mx=0$ for some $M>0$ implies by assumption of induction that
$\left(A_{J_{n,n-1}} A_{J_{n-1,n}}-k_n\right){v}_{J_n}\in
{\rm Ker}\ A_{J_{n-1,n}}$, that is
\begin{equation}\label{Ker}
 A_{J_{n-1,n}}\left(A_{J_{n,n-1}} A_{J_{n-1,n}}-k_n\right){v}_{J_n}=
\left(A_{J_{n-1,n}} A_{J_{n,n-1}}-k_n\right)A_{J_{n-1,n}}{v}_{J_n}=0.
\end{equation}
The operator $A_{J_{n-1,n}} A_{J_{n,n-1}}$, as well as the operator
$S_{J_n}=A_{J_{n,n-1}}A_{J_{n-1,n}}$ has no positive integer
eigenvalues, so the operator $\left(A_{J_{n-1,n}} A_{J_{n,n-1}}-k_n\right)$ is
invertible and \rf{Ker} implies $A_{J_{n-1,n}}{v}_{J_n}=0$.

The passage from the general case to the special case of normal crossing,
described above, is based on the following observation.

Keep the notations of Section \ref{Kashiwara}.
Let $Y$ be a smooth algebraic complex variety, let $Z\subset Y$ be a
smooth subvariety of $Y$ and let $I\subset {\cal O}_Y$ be a sheaf of ideals
 of the functions vanishing on  $Z$. Let $M$ be a holonomic $\D_Y$-module
and let $F(M)$ be a good filtration with respect to $I$.
Note that  in general $\cap_n F^n(M)\not=0$.
\begin{lemma} \label{lemmamon}
Let $m$ be a nonzero element of a $\D_Y$-module $M$ annihilated by some power of the
 ideal $I$, $I^km=0$. Then there exists $n\in \ZZ$ such that
 $m\not\in F^n(M)$. Thus the symbol $\,\gr (m)$ is
well defined and is a
nonzero element of the associated graded quotient $\,\gr (M)$.
\end{lemma}
The claim follows from  condition (ii), Section \ref{Kashiwara},
 for good filtrations and
the fact that $F^k(\D_Y)\subset F^0(\D_Y)\cdot I^k$ for all $k>0$.
\hfill{$\square$}

\medskip
Return to the proof of Lemma \ref{lemmaker} for a general arrangement
$\A$.
Let $E^{n}\VA_{\Ga^{n}}$ be the $\D_{X^n}$-module, attached to the
$\Ga^n$-quiver
$\VA_{\Ga^{n}}$. Let $\b\in\v(\G^n)$, $l(\b)=n$, be a vertex such that
 the operator $S_\b$ has no positive integer
eigenvalues. Let  $z_1,...,z_N$ be affine coordinates in $X=\CC^N$ such
that  the closure
$\Xb$ of the stratum $X_\b$ is given by the equations
 $z_1=0,\ldots, z_n=0$, and the functions $z_1,...,z_n$ are generic affine
functions in $\Fb$.

Denote by $\left(E^{n}\VA_{\Ga^{n}}\right)^0$ the $\D_{X^n}$-module,
 generated
 by vectors of the vector space $V'_\b\ot \Ob$, where the space $V'_\b$ is defined in
\rf{v1b}. Clearly, $\left(E^{n}\VA_{\Ga^{n}}\right)^0\subset
 \left(E^{n}\VA_{\Ga^{n}}\right)^{X_\b}$. We will prove that
\begin{equation}
\label{kereq}
\left(E^{n}\VA_{\Ga^{n}}\right)^0=
 \left(E^{n}\VA_{\Ga^{n}}\right)^{X_\b}
\end{equation}
under the conditions of Lemma \ref{lemmaker}.

For $k=1,...,n$ let $I_k$ be the left ideal of ${\mathcal O}_X$, generated
by $z_1,z_2,...,z_k$.
Let ${\rm gr}_k$ denote the graded quotient space with respect to a good
 filtration, determined by the ideal $I_k$. Then
${\rm gr}_{1,...,n}(E^{n}\VA_{\Ga^{n}}):=
{\rm gr}_n {\rm gr}_{n-1}\cdots {\rm gr}_1(E^{n}\VA_{\Ga^{n}})$ is
a quiver
 $\D_Y$-module, where $Y$ denotes the  affine space $\CC^N$ stratified by intersections
of coordinate hyperplanes $\{z_j=0\}$, $j=1,...,n$,
and the corresponding quiver is
supported at the flag $\emptyset,\{1\}, \{1,2\},...,\{1,2,...,n\}$ of vertices
of $\overline{\G}_{n}$, described above.
More precisely, ${\rm gr}_{1,...,n}(E^{n}\VA_{\Ga^{n}})$ is a
$\overline{\G}^n$-quiver $\D_Y$-module $E\tilde{\VA}_{\overline{\G}^n}$,
such that the space $\tilde{V}_{J_k}$ of the quiver
$\tilde{\VA}_{\overline{\G}^n}$ is
 isomorphic to the direct sum
$$\tilde{V}_{J_k}\simeq \ooplus\limits_{\a,\ \a>\b, l(\a)=k}V_\a\, ,$$
and the maps $A_{J_k,J_{k\pm 1}}$ are the sums of the corresponding maps
$A_{\a_i,\a_j}$. In particular, $\tilde{V}_{J_k}\approx V_\b$
(and $Y_{J_n}\simeq \Xb$), and
 the map ${S}_{J_n}$ is isomorphic to the map $S_\b$. Thus we conclude that
 the conditions of Lemma \ref{lemmaker} imply the equality
\begin{equation}
\label{grkereq}
\left({\rm gr}_{1,...,n}\left(E^{n}\VA_{\Ga^{n}}\right)\right)^0=
 \left({\rm gr}_{1,...,n}\left(E^{n}\VA_{\Ga^{n}}\right)\right)^{\Xb}\ .
\end{equation}
By Lemma \ref{lemmamon}, any element $m\in \left( E^{n}\VA_{\Ga^{n}}\right)^{\Xb}$
has a well defined nonzero symbol
${\rm gr}_{1,...,n}(m)$. Thus
relation
\rf{grkereq} implies relation \rf{kereq}. This proves Lemma
\ref{lemmaker}.
\hfill{$\square$}.
\medskip

Return to the proof of Theorem  \ref{prop4}.
 Since the operator $\Tprime_\b$ for the quiver $E^{n-1}\VA_{\Ga^{n-1}}$ is
isomorphic to the operator $S_\b$ for the quiver
$\psi^{(0)}_\b(\VA_{\Ga^n})$, the $\D_{H_\b}$-module
$\left(\Psi^{(0)}_G(M)\right)^{X_\b}$ is generated by the space
 $W_\b\subset \oplus_{\a\in\v_{n-1}(\Ga),\ \a\succ\b}V_\a$,
which  consists of the vectors
$w=\sum_{\a\in\v_{n-1}(\Ga),\ \a\succ\b}v_\a$. Here $v_\a$  are
 vectors in $V_\a$, such that
for any $\d\in\v_{n-2}(\Ga)$, satisfying the condition $\d>\b$, the sum
$\sum_{\a\in\v_{n-1}(\Ga),\ \a\succ\b}A_{\d,\a}v_\a$ equals zero.

Then statement {(iii)} of Lemma \ref{lemma1pr10}  implies that the gluing
 data for the $\D_{X^n}$-module $M=j_{n,n-1,*}E^{n-1}\VA_{\Ga^{n-1}}$
coincide with the gluing data of the $\D_{X^n}$-module
$E^{n}\j_{n,n-1,*}\VA_{\Ga^{n-1}}$. This proves statement (i) of
 Theorem \ref{prop4}.
\hfill{$\square$}.

\medskip

Let us prove statement (ii) of Theorem \ref{prop4}.

 We need
the following general property of adjoint functors.
Let ${\mathcal A}$ and ${\mathcal B}$ be  additive categories. Let
$F:{\mathcal A}\to {\mathcal B}$   and
$G:{\mathcal B}\to {\mathcal A}$  form a pair of left and right
adjoint functors. The adjunction is established by a collection of
isomorphisms of adjunction,
$$
\a_{X,Y}:\,\Hom_{{\mathcal A}}\left(X,G (Y)\right)\risom
\Hom_{{\mathcal B}}\left(F(X), Y\right),\qquad X\in {\mathcal A},\ {} \ Y\in
{\mathcal B} .
$$
Let
$\sigma_{F,G}(X):\, X\to GF(X),\quad$
$\tau_{F,G}(Y):\,  FG(Y)\to Y,$ 
be the corresponding adjunction morphisms, see \rf{stau}.
One can see that
the adjunction morphisms completely determine the isomorphisms of adjunction.
Namely, for any $X\in{\mathcal A}$, $Y\in {\mathcal B}$ and  morphisms
$\varphi\in \Hom_{{\mathcal A}}(X,G(Y))$  and
$\psi\in \Hom_{{\mathcal B}}(F(X),Y)$
we have the equalities
\begin{equation}\label{6.26}
\a_{X,Y}(\varphi)=\tau_{F,G}(Y)\cdot F(\varphi)\, ,\qquad
\a_{X,Y}^{-1}(\psi)=G(\psi)\cdot \sigma_{F,G}(X)\, .
\end{equation}

 Statement (ii) of Theorem \ref{prop4} relates the isomorphisms of
adjunction of two pairs of adjoint functors.
The first pair consists of the functors
 $\j_{l,k}^{*}$: $Qui_{\G^l}\to Qui_{\G^k}$ and
$\j_{l,k,*}$: $Qui_{\G^k}\to Qui_{\G^l}$, where $0\leq k<l\leq N$.
The second pair consists of the functors
$j_{l,k}^{*}$: $\MD^{hol}_{X^l}\to \MD^{hol}_{X^k}$ and
$j_{l,k,*}$: $\MD^{hol}_{X^k}\to \MD^{hol}_{X^l}$.

Relation \rf{6.26} says that it is sufficient to prove that the
corresponding adjunction morphisms $\tau$ commute with the functor
$E^k$,
 \begin{equation}\label{6.27}
E^k\tau _{\j_{l,k}^{*},\j_{l,k,*}}(\VA_{\G^k})=
\tau _{j_{l,k}^{*},j_{l,k,*}}(E^k\VA_{\G^k})
\end{equation}
for any strongly non-resonant level $k$ quiver
 $\VA_{\G^k}$.
By the constructions of Section \ref{section2.4}, the morphism
$\tau _{\j_{l,k}^{*},\j_{l,k,*}}:
 \j_{l,k}^{*}\j_{l,k,*}(\VA_{\G^k}) \to (\VA_{\G^k})$ is realized as the
identity map, $\tau _{\j_{l,k}^{*},\j_{l,k,*}}|_{V_\a}= {\rm Id}_{V_\a}$,
 $l(\a)\leq k$. On the other hand, for any $M\in \MD^{hol}_{X^k}$  the map
$\tau _{j_{l,k}^{*},j_{l,k,*}}:
 j_{l,k}^{*}j_{l,k,*}(M) \to M$ is the identification  with $M$ of the restriction
to the open set $X^k$ of the direct image $j_{l,k,*}(M)$.
 Relation \rf{6.27} means that the restriction to $X^k$ of both sides of
 relation \rf{jE} is equal to $E^k(\VA_{\G^k})$. This is true since all
constructions of Section \ref{section6.6} do not change quiver
$\D$-modules on the open set $X^k$.\hfill{$\square$}.
\medskip
\subsection{Proof of Theorem \ref{prop11}}

 Introduce a filtration on the complex $(\R
E^0\VA_{\G^0},d)$ as follows. Set $\deg (D\ot
\xi_1\wedge...\wedge\xi_r\ot\om_\emptyset\ot v_\emptyset)=k$ if
$D$ is a differential operator of degree $k-r$ and denote by
 $\R E^0\VA_{\G^0}^{(k)}$ the vector subspace of $\R
E^0\VA_{\G^0}$, generated by  elements of degree less than or equal
to $k$. We have $d(\R E^0\VA_{\G^0}^{(k)})\subset \R
E^0\VA_{\G^0}^{(k)}$,
\begin{equation}
0\subset \R E^0\VA_{\G^0}^{(0)}\subset...\subset
 \R E^0\VA_{\G^0}^{(k)}
 \subset ...\subset\R E^0\VA_{\G^0},
\label{filtered}
\end{equation}
 and $\bigcup_{k\geq 0}\R E^0\VA_{\G^0}^{(k)}=\R E^0\VA_{\G^0}$.
The image of the differential $d$ in the quotient space $\R
E^0\VA_{\G^0}^{(k)}/\R E^0\VA_{\G^0}^{(k+1)}$ coincides with the
image of the differential $d_1$.  The quotient complex $\gr\R
E^0\VA_{\G^0}^{(k)} = \R E^0\VA_{\G^0}^{(k)}/\R
E^0\VA_{\G^0}^{(k+1)}$ is isomorphic to the tensor product of the
 free ${\cal O}(X^0)$-module ${\cal O}(X^0)\ot\overline{\Omega}_\emptyset
\ot V_\emptyset$ and the $k$-th graded component of the Koszul
 complex $\left(S(\C^N)\ot\wedge(\C^N),\ d_{Kos}\right)$, where
$d_{Kos}=\sum_{i=1}^N e_i\ot \dfrac{\partial}{\partial e_i}$, for a given basis
 $e_1,
 ,...,e_N$ of $\C^N$:
 $$\Gr\R
E^0\VA_{\G^0}^{(k)}\approx {\cal O}(X^0) \ot\overline{\Omega}_\emptyset
\ot V_\emptyset\ot
\left(S(\C^N)\ot\wedge(\C^N)\right)^{(k)}.$$

The complex $\left((S(\C^N)\ot\wedge(\C^N))^{(k)},\ d_{Kos}\right)$ is
 acyclic for $k>0$ and consists of a one dimensional space
$S^0(\C^N)\ot\wedge^0(\C^N)$ with the zero differential for $k=0$.
 As a consequence, the
 complex  $\Gr\R E^0\VA_{\G^0}^{(k)}$ is acyclic for $k>0$, while
  $ RE^0\VA_{\G^0}^{(0)}$ consists of the vector space  ${\cal O}(X^0) \ot
\overline{\Omega}_\emptyset
\ot V_\emptyset \ot S^0(\C^N)\ot\wedge^0(\C^N)$ with the zero
differential. This implies that the inclusion
$i_0:  RE^0\VA_{\G^0}^{(0)}\hookrightarrow RE^0\VA_{\G^0}$ is a
quasi-isomorphism. On the other hand, by Proposition \ref{prop1},
{ (i)}, the composition $\nu
i_0$ is an isomorphism of the vector spaces ${\cal O}(X^0) \ot
\overline{\Omega}_\emptyset
\ot V_\emptyset$ and $E^0\VA_{\G^0}$. Thus the map $\nu$ is also a
quasi-isomorphism of the complex $(\R
E^0\VA_{\G^0},d)$ and of $E^0\VA_{\G^0}$, treated as the complex with
zero differential. \hfill{$\square$}

\subsection{
 Proof of Theorem \ref{theorem3}.}

Let $D^k_X\subset D_X$ denote the subspace of differential operators of
order less than or equal to $k$.

{}Fix a quiver $\VA=\{ V_\a, A_{\a,\b}\}$.
Let $\a\in\v(\Ga)$. For any $v_\a\in V_\a$ and
$\om_\a\in \Oa$,  set $\deg(\om_\a\ot v_\a)=l(\a)$. Set also
$\deg \xi =1$ for any $\xi\in \Ta$ and $\deg f =0$ for any $f\in \Fa$.
 Hence $\deg x =p+l(\a)$ for any $x\in \bigwedge^p\Ta\wedge \bigwedge^q\Fa
 \ot \Oa \ot V_\a$. Define $\Rlong\EVA^k$ as the subspace of
 $\Rlong\EVA$, spanned
 by  linear combinations
$$\sum\limits_{\a\in\v(\Ga)} y_\a x_\a,$$
where $x_\a\in\bigwedge^*\TTT{\a}\ot\Oa\ot V_\a$, $y_\a\in \D_X$ and
 $\deg x_\a +\deg y_\a\leq k$.
\begin{lemma}\label{lemma5}${}$

\noindent
{\em (i)} For any $k\geq 0$, $\Rlong\EVA^k$ is a subcomplex of  $\Rlong\EVA$.
\\
{\em (ii)} For any $k,l\geq 0$, there are inclusions
 $\D^k_X \Rlong\EVA^l\subset \Rlong\EVA^{k+l}$.
\end{lemma} \hfill{$\square$}

By Lemma \ref{lemma5} we have a filtration
$$ 0\subset \Rlong\EVA^0\subset\Rlong\EVA^1\subset ...\subset \Rlong\EVA$$
of the complex $\Rlong\EVA$ by  subcomplexes $\Rlong\EVA^k$.
We can view this filtration as a pair:
a filtration
$$ 0\subset \R\EVA^0\subset\R\EVA^1\subset ...\subset \R\EVA$$
of the complex $\R\EVA$, and a filtration
$$
 0\subset \EVA^0\subset\EVA^1\subset ...\subset \EVA
$$
of the $\D_X$-module $\EVA$. This filtration coincides with the
principal filtration $F_p\EVA$ of $E\VA$, see Section \ref{Sec.3.2}. We
also have
  the map $\nu$ of filtered complexes,
$\nu: \R\EVA\to {\EVA}_{(\cdot)}$.
Let  $\GR\EVA^k=\R\EVA^k/\R\EVA^{k-1}$, and $\GM{E}\VA^k=\EVA^k/\EVA^{k-1}$.
 The space $\GR\EVA^k$ is spanned by
 linear combinations
$\sum\limits_{\a\in\v{\Ga}} y_\a x_\a$,
where $x_\a\in\bigwedge^*\TTT{\a}\ot\Oa\ot V_\a$ is an element
of some degree $j$,
$y$ is the symbol of a differential operator of degree $k-j$ in
$\D_X^{k-j}/\D_X^{k-j-1}$.

 The differential $\GM{d}$   in the
 graded quotient complex  $\GR\EVA=\oplus_{k\geq 0}\GR\EVA^k$
can be presented as the difference
 $\GM{d}=\GM{d}_1-\GM{d}_0$, where
$$\begin{array}{ccc}
\GM{d}_{1}\left(y_\a\ot t_1\wedge...\wedge t_r\ot \om_\a\ot
v_\a\right)&=&\!\!\!\sum\limits_{i=1}^r (-1)^{i+1}y_\a t_i\ot
t_1\wedge... \widehat{t_{i}}...\wedge t_r\ot \om_\a\ot v_\a,
 \\
\GM{d}_{0}\left(y_\a\ot t_1\wedge...\wedge t_r\ot
\om_\a\ot v_\a\right)&=
&\!\!\!\sum\limits_{\b,\ \a\succ \b} y_\a\ot A^r_{\b,\a}
(t_1\wedge...\wedge t_r\ot \om_\a\ot v_\a).\end{array}
$$
By Proposition \ref{lemma1}, the space
 $\EVA$ can be identified with the space $M$,
$$M=\ooplus\limits_{\a\in\v(\Ga)}S(\TTa)\ot S(\FFa)\ot \Oa\ot V_\a.$$
Here $S(\TTa)=\oplus_{j\geq 0}S^j(\TTa)$ and
$ S(\FFa)\oplus_{j\geq 0}S^j(\FFa)$ are the symmetric algebras of the spaces
$\TTa$ and $\FFa$ introduced in Section \ref{Framing}.
The identification map $\varphi:M\to\EVA$,
\begin{equation}\label{phi}
\varphi(\xi_{i_1}\cdots \xi_{i_n}\otimes f_{j_1}\cdots f_{j_m}\ot\om_\a\ot
 v_\a)=\xi_{i_1}\cdots \xi_{i_n} f_{j_1}\cdots f_{j_m}(\om_\a\ot
 v_\a)
\end{equation}
is an application of commuting elements $\xi_{i_k}\in\TTa\subset\D_X$ and
$f_{j_l}\in\FFa\subset\D_X$ to the vectors $\om_\a\ot v_\a\in\EVA$.

For any $\a\in \v(\Ga)$, $\xi\in\TTa$, and $f\in\FFa$, set $\deg \xi=1$,
 $\deg f=0$, and $\deg(\om_\a\ot v_\a)=l(\a)$. Thus, $M$ is a graded space
with graded components
$$M^{(k)}=\ooplus\limits_{\a\in \v(\Ga)}\ooplus\limits_{j, j+l(\a)=k }
S^j(\TTa)\ot S(\FFa)\ot \Oa\ot V_\a.$$

Set $\,M^k=\ooplus\limits_{j\leq k}M^{(j)}$. The spaces $M^k$, $k\geq0$,\,
 give an
 increasing filtration on $M$:\, $0\subset M^0\subset...\subset M^k\subset...
\subset M$.

\begin{lemma}
\label{lemma7}
The map \rf{phi} establishes an isomorphism of the filtered spaces $M$ and
$\EVA$.
\end{lemma}

{\it Proof}. We know by  Proposition
\ref{lemma1} that   the linear map \rf{phi} is an isomorphism.
So it is sufficient to establish an equality of vector spaces
\begin{equation}\label{mk}
\varphi(M^k)=\EVA^k.
\end{equation}
for any $k\geq 0$.
 We prove \rf{mk} by induction over $k$. By definition of $M^k$ and
$\EVA^k$, we have $\varphi(M^k)\subset \EVA^k$ and $\varphi(M^0)=\EVA^0$.
Suppose we know that $\varphi(M^{k-1})=\EVA^{k-1}$ for some $k\geq 1$.
 Consider an element
$x=y \cdot \om_\a\ot v_\a\in \EVA^k$, where $y\in\D_X^{k-l(\a)}$. We prove by
 induction on $\dim X_\a=N-l(\a)$ that $x\in\varphi(M^k)$.
  Present $x$ as
\begin{equation}\label{xif}
x=\sum\nolimits_i\xi'_if'_i\xi_i f_i(\om_\a\ot v_\a),
\end{equation}
where $\xi'_i,\ f'_i,\  \xi_i,\ f_i\in \D_X$
are homogeneous elements of the symmetric algebras over $\TTa,\ \FFa,\ \Ta, \
\Fa$, respectively.
  The summands in \rf{xif}
with $f_i\in S^+(\Fa)=\oplus_{j>0}S^j(\Fa)$  belong to $\EVA^{k-1}$.
 The summands in \rf{xif} with $f_i=1$ but
$\xi_i\in S^+(\Ta)$ can be presented in the form
$\sum_{\b, l(\b)>l(\a)}y_\b(\om_\b\ot v_\b)$ and belong to $\EVA^k$
 according
to the assumption of the induction on $N-l(\a)$.
For vertices $\a$,  $l(\a)=N$,
the summands in \rf{xif}  with $f_i=1$ but
$\xi_i\in S^+(\Ta)$ are missing since $\Ta=0$ in this case.
 The rest of the summands
 is of the form $\xi'_if'_i(\om_\a\ot v_\a)$ with
$\xi'_i\in S(\TTa)$, $f'_i\in S(\FFa)$, and belong to
$\varphi(M^k)$ by the definition of  $\varphi$.
\hfill{$\square$}
\medskip

Let $\widetilde{\nu}:\GR\EVA_0\to
\GM{E}\VA=\oplus_{k\geq 0}\GM{E}\VA^k$ be the map, induced by
$\nu$.  Identify $\GM{E}\VA^k$ with $M^{(k)}$.
The proof of Lemma \ref{lemma7} shows that for any
$y\in \D_X^l/\D_X^{l-1}$ and $\a\in\v(\Ga)$ we can write
the image $\widetilde{\nu}(y\ot\om_\a\ot v_\a)$ as
$$\widetilde{\nu}(y\ot\om_\a\ot v_\a)=\ooplus\limits_{\b,\  \b\geq\a}
\widetilde{\nu}_{\b,\a}(y\ot\om_\a\ot v_\a),$$
where $\widetilde{\nu}_{\b,\a}(y\ot\om_\a\ot v_\a)\in
 S(\TTb)\ot S(\FFb)\ot \Ob\ot V_\b$ and for any
 $\xi'\in S(\TTa)$,  $f'\in S(\FFa)$,  $\xi\in S(\Ta)$,
 $f\in S(\Fa)$,
\begin{equation}\nonumber
\widetilde{\nu}_{\a,\a} (\xi'f'\xi f\ot\om_\a\ot v_\a)= \varepsilon(\xi)
\varepsilon(f)\,
\xi'\ot f' \ot\om_\a\ot v_\a.
\end{equation}
Here
\begin{equation}\label{augmentation}
\varepsilon: S(V)\to\CC
\end{equation}
 is the augmentation map, equal to the projection of the symmetric algebra
of a vector space $V$ to its zero graded component
$S^0(V)\approx\CC$ parallel to $S^+(V)=\oplus_{k>0}S^k(V)$.

{}For a fixed $k\geq 0$ and $l=0,...,N=\dim X$, define
$\GR\EVA^{k,l}$ as the  vector subspace of the space $\GR\EVA^{k}$,
 spanned by the linear combinations
$$ \sum\limits_{\a ,\ l(\a)\geq l} y_\a x_\a,$$
where $x_\a\in\bigwedge^*\TTT{\a}\ot\Oa\ot V_\a$, $y_\a$ is the symbol of
a homogeneous differential operator in $\D_X$, such that
 $\deg x_\a +\deg y_\a= k$. Define also $\GM{E}\VA^{k,l}$ to be the subspace
of the space  $\GM{E}\VA^{k}$ generated by the elements
$\xi'\ot f' \ot\om_\a\ot v_\a$ with $l(\a)\geq l$.
\begin{lemma}
\label{lemma6} Let $k>0$ and $l=0,...,N=\dim X$. Then

\noindent
{\em (i)}  $\,\GR\EVA^{k,l}$
is a subcomplex of $\GR\EVA^{k}$,

\noindent
{\em (ii)} $\,\GR\EVA^{k,l}\subset \GR\EVA^{k,l'}$ and $\GM{E}\VA^{k,l}
\subset \GM{E}\VA^{k,l'}$ if
$l>l'$,

\noindent
{\em (iii)} $\,\widetilde{\nu}(\GR\EVA^{k,l}_0)\subset
\GM{E}\VA^{k,l}$.
\end{lemma} \hfill{$\square$}

By Lemma \ref{lemma6}, for any $k>0$, we have the filtration of complexes,
$$0\subset \GR\EVA^{k,N}\subset ... \subset \GR\EVA^{k,0}=\GR\EVA^{k},$$
 the filtration of the space $\GM{E}\VA^k$,
$$
 0\subset \GM{E}\VA^{k,N}\subset ...\subset \GM{E}\VA^{k,0}= \GM{E}\VA^{k}
$$
 and the map $\GM{\nu}$ of filtered complexes
$\GM{\nu}: \GR\EVA^k\to \GM{E}\VA^k_{(\cdot)}$.
{}For $l=0,...,N$ let
 $\GGR\EVA^{k,l}=\GR\EVA^{k,l}/\GR\EVA^{k,l+1}$ be
the quotient complex of $\GR\EVA^k$ and
$\GGM{E}\VA^{k,l}=\GM{E}\VA^{k,l}/\GM{E}\VA^{k,l+1}$
 the  quotient module of $\GM{E}\VA^k$.
 The space $\GGR\EVA^{k,l}$ is spanned by
 linear combinations
$\sum_{\a,\ l(\a)=l} y_\a x_\a$,
where $x_\a\in\bigwedge^*\TTT{\a}\ot\Oa\ot V_\a$,
$y_\a$ is the symbol of a homogeneous differential operator in $\D_X$
 such that $\deg x_\a+\deg y_\a=k$. The space $\GGM{E}\VA^{k,l}$
 is generated by  linear combinations
$\sum_{\a,\ l(\a)=l} \xi'_\a\ot f'_a\ot\om_\a\ot v_\a$,
where $\xi'_\a\in S^{k-l}(\TTa)$, $f'_\a\in S(\FFa)$.

Let $\GGM{d}$
be the differential in the quotient complex
$\GGR\EVA=\oplus_{k,l}\GGR\EVA^{k,l}$ and $\GGM{\nu}$ the map
 $\GGR\EVA_0\to \GGM{E}\VA=\oplus_{k,l}\GGM{E}\VA^{k,l}$ induced by
$\GM{\nu}$. We have
$\GGM{d}=\GGM{d}_1$, where
\begin{equation}
\GGM{d}_{1}\left(y_\a\ot t_1\wedge...\wedge t_r\ot \om_\a\ot
v_\a\right)=
\!\!\!\sum\limits_{i=1}^r (-1)^{i+1}y_\a t_i\ot
t_1\wedge... \widehat{t_{i}}...\wedge t_r\ot \om_\a\ot v_\a,
\label{ddr}
\end{equation}
and
\begin{equation}
\GGM{\nu} (\xi'f'\xi f\ot\om_\a\ot v_\a)= \varepsilon(\xi)
\varepsilon(f)\,
\xi'\ot f' \ot\om_\a\ot v_\a.
\label{nnu}
\end{equation}
 Relations \rf{ddr} and \rf{nnu} show that the complex $\GGR\EVA$ is the
direct sum of its subcomplexes $\GGR\EVA_\a$,  generated by the elements of
the form $y_\a x_\a$ with a fixed $\a\in\v(\G)$. The vector space
$\GGM{E}$ is the direct sum of its subspaces $\GGM{E}_\a$, generated by the
elements $\xi'\ot f' \ot\om_\a\ot v_\a$. The complex $\GGR\EVA_\a$, the
space $\GGM{E}_\a$ and the map $\GGM{\nu}$ form the complex
$$
\overline{\GGR}\EVA:\qquad 0\to{\GGR}^{-N}\EVA\stackrel{\GGM{d}}{\to}
\GGR^{-N+1}\EVA\stackrel{\GGM{d}^{}}{\to}\cdots
\stackrel{\GGM{d}}{\to}\GGR^0\EVA \stackrel{\GGM{\nu}}{\to}\GGM{E}\VA\to 0 .$$
This complex is isomorphic to  the tensor product of
Koszul complexes related  to symmetric algebras $S(\Ta)$ and $S(\Fa)$ with
coefficients in the rings $S(\TTa)\ot S(\FFa)$:
$$\overline{\GGR}\EVA\quad\approx\
\ooplus\limits_{\a\in\v(\Ga)}S(\TTa)\ot S(\FFa)\ot Kos_{T_\a}\ot
Kos_{F_\a}.$$
Recall that for a vector space $V$ of the dimension $N$, the Koszul complex
 $Kos_V$is  the total complex of free $S(V)$-module resolution
of the trivial $S(V)$-module $\CC$:
$$Kos_V:\qquad 0\to S(V)\ot\bigwedge\nolimits^N(V^*)\stackrel{{d_{Kos}}^{}}{\to}
\cdots
\stackrel{{d}_{Kos}}{\to}S(V)\ot\bigwedge\nolimits^0(V^*)
\stackrel{{\varepsilon}^{}}{\to}
\CC\to 0,$$
where $d_{Kos}=\sum_{i=1}^N e_i\ot \frac{\partial}{\partial e^i}$ with
 $\{e_i\}$ and $\{e^i\}$, $i=1,...,N$, being dual bases of $V$ and $V^*$,
and
$\varepsilon$ is the augmentation map
$\varepsilon:S(V)\to\CC$, see  \rf{augmentation}.

 The complex  $\overline{\GGR}\EVA$ is acyclic, since Koszul complexes
are acyclic. Then
the complex  $\overline{\GR}\EVA$, as well as
 ${\Rlong}\EVA $ is acyclic, since all its graded factors are acyclic.
  Theorem \ref{theorem3} is proved.
\hfill $\square$
\medskip

\subsection{Proof of Theorem \ref{propdual}}
We prove  statement { (ii)} of the theorem. Statement
 {(i)} is a specialization
of {(ii)} to the open set  $X^0$.

First we reformulate statement (ii). Let $\VA$ be a quiver of $\G$.
Consider the complex
$${\rm R}Hom_{\D_X}(\R E\VA,\widetilde{\D}_X)[N]\ ,$$
and denote it by $\widetilde{\R} E\VA$. The complex has $N+1$ terms.
Write it as
$$\widetilde{\R} E\VA:\ 0\to \widetilde{\R}_{-N} E\VA\to \cdots \to
\widetilde{\R}_0 E\VA\to 0\ .$$
Then $\deg \widetilde{\R}_k E\VA =k$, and
\begin{equation}
\begin{split}\nonumber
\widetilde{\R}_{-r} E\VA= Hom_{\D_X}\left(\ooplus\limits_{\a\in\v(\Ga)}
\D_X\ot \left(\bigwedge\nolimits^{N-r}\TTT{\a}\ot \Oa\ot V_\a\right),
\widetilde{\D}_X\right)\approx\\ \nonumber
\approx\ooplus\limits_{\a\in\v(\Ga)}\widetilde{\D}_X\ot
\left(\bigwedge\nolimits^{N-r}\TTT{\a}\ot \Oa\ot V_\a\right)^*\ .
\end{split}
\end{equation}
We have to prove that $\widetilde{\R} E\VA$ is quasi-isomorphic
 to the complex $\R E(\tau \VA)$.

Choose a vector $\eta\in\bigwedge^NT_\emptyset$,
$\eta=c\frac{\partial}{\partial{z_1}}\wedge\cdots\wedge
\frac{\partial}{\partial{z_N}}$ in local affine coordinates $z_1,...,z_N$
of $\CC^N$, and define a pairing
$$\left(\bigwedge\nolimits^{r}\TTT{\a}\ot \Oa\ot V_\a\right)\ot
\left(\bigwedge\nolimits^{N-r}\TTT{\a}\ot \Oa\ot V^*_\a\right)\to\CC$$
in the following way. Let $\a\in\v_l(\Ga)$,
$\xxi{p}, \vec{\xi}'_{p'}\in\bigwedge^{p}\Ta$,
$\f{q}, \f{q'}'\in \bigwedge^q \Fa$, where $p+q=r$, $\om_a,\om'_a\in\Oa$,
 $v_\a\in V_\a$,
$v'_a\in V^*_\a$. Set
\begin{equation}\begin{split}\label{523}
\mylangle\xxi{p}\wedge\f{q}\ot\om_\a\ot v_\a,\ \vec{\xi}'_{p'}\wedge\f{q'}'
\ot\om'_\a\ot v'_\a\myrangle&=\\ =
(-1)^\frac{p(p-1)+q(q-1)}{2}\delta_{p', N-l-p}\delta_{q',l-q}
\mylangle\xxi{p}\wedge\vec{\xi}'_{p'},\
\om_\a\myrangle&\cdot\mylangle\eta,\ d\f{q}
\wedge\om'_\a\wedge d\f{q'}'\myrangle\cdot \mylangle v_\a,v'_\a\myrangle\ .
\end{split}
\end{equation}
Due to \rf{523} we have an isomorphism
$$
\widetilde{\R}_{-r} E\VA
\approx\ooplus\limits_{\a\in\v(\Ga)}\widetilde{\D}_X\ot
\left(\bigwedge\nolimits^{r}\TTT{\a}\ot \Oa\ot V_\a\right)\, .
$$
 The differential $\widetilde{d}: \widetilde{\R}_{-r} E\VA \to
\widetilde{\R}_{-r+1} E\VA$, adjoint to the differential \rf{dr},
takes the form $\widetilde{d}=\widetilde{d}_1-\widetilde{d}_0$,
\begin{equation}
\begin{array}{ccc}\nonumber
\widetilde{d}_{1}\left(\widetilde{D}\ot t_1\wedge...\wedge t_r\ot \om_\a\ot
v_\a\right)&=&\!\!\!\sum\limits_{i=1}^r (-1)^{i+1}t_i^{(I)}(\widetilde{D})\ot
t_1\wedge... \widehat{t_{i}}...\wedge t_r\ot \om_\a\ot v_\a,
\\
\widetilde{d}_{0}\left(\widetilde{D}\ot t_1\wedge...\wedge t_r\ot \om_\a\ot v_\a\right)&=
&\!\!\!\sum\limits_{\b,\ (\a,\b)\in E(\Ga)} \widetilde{D}\ot A^r_{\b,\a}
(t_1\wedge...\wedge t_r\ot \om_\a\ot v_\a),\end{array}
\end{equation}
where $t_i\in\TTT{\a}$, $\om_\a\ot v_\a\in \Oa\ot V_\a$, $\widetilde{D}\in
\widetilde{\D}_X$ and $x^{(I)}(\widetilde{D})$ means the action of
  $x\in\D_X$ on the element $\widetilde{D}\in \widetilde{\D}_X$ with
respect to the first $\D_X$-module structure on $\widetilde{\D}_X$, defined
in Section \ref{dualsection}. We note that
 $\widetilde{\D}_X$ is the free $\D_X$-module of rank one with respect
 to the second $\D_X$-module structure.  A global section
$\tilde{\omega}\in \Omega_X^{-1}$, say $\tilde{\omega} =
\left(dz_1\wedge...\wedge dz_N\right)^{-1}$, defines an
isomorphism of free left $\D_X$-modules
$\phi_{\tilde{\omega}}: \widetilde{\D}_X\to \D_X$,
$$\phi_{\tilde{\omega}}(D\ot{\tilde{\omega}})=D\cdot 1.$$
Applying the isomorphism $\phi_{\tilde{\omega}}$, we get the
 isomorphism of $\D_X$-modules,
$$
\widetilde{\R}_{-r} E\VA
\approx\ooplus\limits_{\a\in\v(\Ga)}{\D}_X\ot
\left(\bigwedge\nolimits^{r}\TTT{\a}\ot \Oa\ot V_\a\right)\ .
$$
Under this identification  the differential in
$\widetilde{\R} E\VA$ becomes $d'=d'_1-d'_0$,
\begin{equation}\nonumber
\begin{array}{ccc}
{d}'_{D}\left(D\ot t_1\wedge...\wedge t_r\ot \om_\a\ot
v_\a\right)&=&\!\!\!\sum\limits_{i=1}^r (-1)^{i+1}D\cdot t_i^\tau\ot
t_1\wedge... \widehat{t_{i}}...\wedge t_r\ot \om_\a\ot v_\a,
\\
{d}'_{0}\left(D\ot t_1\wedge...\wedge t_r\ot \om_\a\ot v_\a\right)&=
&\!\!\!\sum\limits_{\b,\ (\a,\b)\in E(\Ga)} D\ot A^r_{\b,\a}
(t_1\wedge...\wedge t_r\ot \om_\a\ot v_\a),\end{array}
\end{equation}
where $D\in\D_X$ and ${}^\tau$ is an anti-isomorphism of $\D_X$, such that
$t_i^\tau= t_i$ and $\left(\frac{\partial}{\partial{t_i}}\right)^\tau=
-\frac{\partial}{\partial{t_i}}$.

Define a linear map
$$\nu: \ooplus\limits_{\a\in\v(\Ga)}{D}_X\ot
\left(\bigwedge\nolimits^{r}\TTT{\a}\ot \Oa\ot V_\a\right) \to
\ooplus\limits_{\a\in\v(\Ga)}{D}_X\ot
\left(\bigwedge\nolimits^{r}\TTT{\a}\ot \Oa\ot V_\a\right)$$
as multiplication by $(-1)^p$ on each subspace of the form
${D}_X\ot{\bigwedge\nolimits^{p}\Ta\wedge\bigwedge\nolimits^{q}\Fa
\ot \Oa\ot V_\a}$. Conjugation of $d'$ with $\nu$  transforms
 the complex $\widetilde{\R} E\VA$ to the complex
${\R} E(\tau\VA)$. The theorem is proved. \hfill{$\square$}
\medskip

\subsection{Proof of Theorem \ref{theorem4}}
Let $\A$ be a central arrangement of hyperplanes in $X=\CC^N$. Let
 $\VA\in\Qui$, $\VA=\{ V_\a, A_{\a,\b}\}$,
 be a quiver of the arrangement $\A$  and ${\Gamma}\QVA$  the
complex of  global sections of  the de Rham complex
$\QVA$. The  component ${\Gamma}\QVA^r$, $-n\leq r\leq 0$, of this
complex is isomorphic to
$$
\G\QVA^r=\ooplus\limits_{\a\in\Ga} \G\Oman_X\ot \left(
\bigwedge\nolimits^{-r}\TTT{\a}\ot\Oa\ot V_\a\right),$$
 where $\G\Oman_X$ denotes the space of global holomorphic
 $N$-forms $f(z_1,...,z_N)dz_1\wedge...\wedge dz_N$ on $\C^N$
 with analytic coefficients. The differential $d$,
$d( \G\QVA_r)\subset \G\QVA_{r+1}$ in the complex
${\Gamma}\QVA$ is given by formula \rf{drom}.

 We construct below a new differential $\dd: \G\QVA\to
\G\QVA$,
 such that $\dd( \G\QVA_r)\subset \G\QVA_{r-1}$.
We construct the Laplace operator $\Delta=d\dd +\dd d$.
Using the Laplace operator, we prove  that the
complex ${\Gamma}\QVA$ is quasi-isomorphic to the
quiver complex $C_+(\VA)[N]$
 under the assumption that the matrix coefficients of the
 linear maps $A_{\a,\b}$ are small enough.

The definition of the differential $\dd$ depends on the choice of a
non-degenerate symmetric bilinear form $\lform\,,\rform$ on the vector
 space $\C^N$.

{} For any $\a\in\v(\G)$,
 let $\{\xi_i\}$ be an orthonormal basis of the space $\Ta$,
 $\lform\xi_i,\xi_j\rform=\delta_{i,j}$, and let $\{\xi^k\}$ be an
 orthonormal basis of the space $\TTa$ (see Section \ref{Framing}),
 $\lform\xi^k,\xi^l\rform=\delta_{k,l}$. We have
 $\lform\xi_i,\xi^k\rform=0$, since the
spaces $\Ta$ and $\TTa$ are orthogonal.
 We denote by $\ita$ the set of the indices
of the basis $\{\xi_i\}$ of the space $\Ta$ . It
contains $\dim \Ta$ elements.  We denote by $\ifa$ the set of
the indices of the basis
$\{\xi^k\}$ of the space $\TTa$. It contains
$\dim \Fa$ elements.

 Let $f_l$,
$l\in \ifa$, be the elements of $\Fa$, such that $\mylangle \xi^k,
f_l\myrangle=\delta_{k,l}$. Let  $f^j$, $j\in \ita$, be the elements of
$\FFa$, such that
 $\mylangle \xi_i, f^j\myrangle=\delta_{i,j}$.  They form an orthonormal basis
 of the space $\tilde{F}\simeq \Fa\oplus \FFa$ of all linear function
on $\CC^N$.

Let $d_{1}^*: \G\QVA \to \G\QVA$ be the linear map,
defined by the formula
$$d_{1}^*\left({\o}\ot \overline{t}\ot \om_\a\ot
v_\a\right)=\sum\limits_{i\in\ita}\o\cdot {f}^i\ot \xi_i\wedge
\overline{t}\ot \om_\a\ot v_\a + \!\sum\limits_{k\in\ifa}\o\cdot
{\xi}^k\ot f_k\wedge \overline{t}\ot \om_\a\ot v_\a,$$
 where $\a\in\v(\G)$, $\o\in \G\Oman_X$,
$\overline{t}\in\bigwedge^*\TTT{\a}$,
 $\om_\a\in\Oa$, $v_\a\in V_\a$.

Choose an edge framing $\{f_{\b,\a},
 \xi_{\a,\b} \}$ in such a way
 that for any edge $(\a,\b)\in E(\G), \a\succ\b$, we have
 $\lform\xi_{\a,\b},\xi_\b\rform=0$ for any $\xi_\b\in T_\b$,
 $\mylangle \xi',f_{\b,\a}\myrangle=0$ for any
 $\xi'\in\TTa$ and $\mylangle \xi_{\a,\b},f_{\b,\a}\myrangle=1$. In particular, the
 last equality implies that the maps $\pi_{\a,\b}$ and
 $\pi_{\b,\a}$ are inverse to each other. We call such a framing
 orthogonal.

The bilinear form $(,)$ on $\CC^N$ defines for any $\a\succ\b$
 orthogonal projections $T_\a\to T_\b$ and $F_\b\to F_\a$, which extend to
projections
 $P_{\b,\a}:\bigwedge^*\Ta\to\bigwedge^*\Tb$ and
 $P_{\a,\b}:\bigwedge^*\Fb\to\bigwedge^*\Fa$. For an orthogonal
 edge framing $\{f_{\b,\a},
 \xi_{\a,\b} \}$ we have
 $$\begin{array}{ccc}
 P_{\b,\a}(\xi_\a\wedge\overline{\xi}_\b)&=&(\xi_\a-\mylangle \xi_\a,
 f_{\b,\a}\myrangle
 \xi_{\a,\b})\wedge\overline{\xi}_\b,\\
 P_{\a,\b}(f_\b\wedge\overline{f}_\a)&=&(f_\b-\mylangle \xi_{\a,\b},
 f_{\b}\myrangle
 f_{\b,\a})\wedge\overline{f}_\a,
 \end{array} $$
if $\xi_\a\in\Ta$, $\overline{\xi}_\b\in\bigwedge^*\Tb\subset
\bigwedge^*\Ta$, $f_\b\in\Fb$, $\overline{f}_\a\in
\bigwedge^*\Fb\subset\bigwedge^*\Fa$.

 {}For each edge $(\a,\b)\in E(\G)$ define a
linear map $$d^{*}_{0,(\b,\a)}:\G\Oman_X\ot
\bigwedge\nolimits^*\TTT{\a}\ot\Oa\ot V_\a\to \G\Oman_X\ot
\bigwedge\nolimits^*\TTT{\b}\ot\Ob\ot V_\b\ ,$$
 as follows. Let
 $\o\in \G\Oman_X$,
$\overline{\xi}_\a\in\bigwedge^*\Ta$,
$\overline{f}_\a\in\bigwedge^*\Fa$, $\om_\a\in\Oa$, $v_\a\in
V_\a$.
  We set
$$d^{*}_{0,(\b,\a)}(\o\ot\overline{\xi}_\a\wedge
\overline{f}_\a\ot\om_a\ot v_\a)= \o\ot f_{\b,\a}\wedge
P_{\b,\a}(\overline{\xi_\a})\wedge
\overline{f}_\a\ot\pi_{\b,\a}(\om_\a)\ot A_{\b,\a}(v_\a)$$ if
$\a\succ\b$, and
$$d^*_{0,(\b,\a)}(\o\ot \overline{\xi}_\a\wedge
\overline{f}_\a\ot\om_a\ot v_\a)= \o\ot\xi_{\b,\a}\wedge
\overline{\xi_\a}\wedge
P_{\b,\a}(\overline{f}_\a)\ot\pi_{\b,\a}(\om_\a)\ot
A_{\b,\a}(v_\a)$$ if $\b\succ\a$. Let
$$
d^*_{0}\ =\ \sum_{(\a,\b)\in E(\G)}d^*_{0,(\b,\a)}\ .
$$
\begin{lemma}\label{lemma*}
 The maps $d_{i}^*$, $i=0,1$, form two commuting
differentials,
 $$(d_{i}^*)^2=0\quad {\rm for}\quad
 i=0,1,\qquad d_{0}^*d_{1}^*+d_{1}^*d_{0}^*=0.
 $$
\end{lemma} \hfill{$\square$}

 Define $d^*=d^*_{1}+d^*_{0}$ and $\Delta=dd^*+d^*d$. We
decompose $\Delta$ as
$$\Delta=\Delta^{(2)}+\Delta^{(1)}+\Delta^{(0)},$$
where $\Delta^{(2)} =d_1^*d_1+d_1d_1^*$, $\Delta^{(1)}
=d_0^*d_1+d_1d_0^*-d_1^*d_0-d_0d_1^*$, $\Delta^{(0)}
=-d_0^*d_0-d_0d_0^*$.

 Denote by $\V_\a$ the following finite-dimensional space:
$$\V_\a=\bigwedge\nolimits^*\TTT{\a}\ot\Oa\ot V_\a\\ ,$$
and let $\ \V=\ooplus\limits_{\a\in\v(\Ga)}\V_\a$.
\begin{lemma}\label{laplace}
 For any $\o\in \G\Oman_X$, $\overline{v}_\a\in\V_\a$, and
$\overline{v}\in\V$,
 we have
\begin{align*}
&\mbox{(i)}&&\Delta^{(2)}\left(\o\ot
 \overline{v}_a\right)&&=&&
 \o\cdot\left(\sum\limits_{i\in\ita}f^i\xi_i+
 \sum\limits_{k\in\ifa}\xi^k f_k\right)\ot\overline{v}_\a,\hspace{4cm}\\
&\mbox{(ii)}&&\Delta^{(1)}\left(\o\ot\overline{v}\right)&&=&&0,\\
&\mbox{(iii)}&&
\Delta^{(0)}\left(\o\ot\overline{v}\right)&&=&&\o\ot
 \overline{\Delta}^{(0)}(\overline{v})\quad \text{{\rm for some linear
operator $\overline{\Delta}^{(0)}:\V\to \V$.}}
\end{align*}
\end{lemma}
\hfill{$\square$}

Let $\a\in\v(\Ga)$. For any $v_\a\in V_\a$ and
$\om_\a\in \Oa$,  set $\deg(\om_\a\ot v_\a)=\dim X_\a= N-l(\a)$. Set
$\deg \xi =-1$ for any $\xi\in \Ta$ and $\deg f =1$ for any $f\in \Fa$.
 Thus for any
$x\in \bigwedge^p\Ta\wedge \bigwedge^q\Fa
 \ot \Oa \ot V_\a$, we have
 $\deg x =p-q+ N-l(\a)$.  This grading gives the decomposition
$\V=\ooplus\limits_{k=0}^N\V^k$,
where $\V^k$ consists of  elements $v\in\V$ of degree $k$.
 We say that a polynomial differential form
 $\omega=f(z_1,...,z_N)dz_1\wedge...\wedge dz_N$ has degree
$k\geq 0$ and write $\deg \omega= k$, if
$f(z_1,...,z_N)$ is a homogeneous polynomial of degree $k$
 variables $z_1,...,z_N$.
We say that a holomorphic form
 $\omega=f(z_1,...,z_N)dz_1\wedge...\wedge dz_N$ has degree not less than
$k\geq 0$ and write $\deg \omega\geq k$, if the holomorphic function
$f(z_1,...,z_N)$ has zero of order not less than $k$ at the origin
$z_1=z_2=...=z_N=0$.

Define $\G\QVA^k$ as the subspace of the vector space $\G\QVA$,
 spanned by linear combinations
$$
\sum\limits_{j=0}^N \omega^{k-j} v^j\ ,
$$
where $v^j\in\V^j$ and
$\omega^{k-j}\in \G\Oman_X$ is a holomorphic form of degree not less than
 $k-j$.

\begin{lemma}\label{lemma8}
 For a central arrangement $\A$   and for any $k\geq 0$, $\G\QVA^{k}$
 is a  subcomplex of  $\G\QVA$.
\end{lemma} \hfill{$\square$}

By Lemma \ref{lemma8} we have a decreasing filtration
$$ \G\QVA=\G\QVA^0\supset \G\QVA^1
\supset...\supset  \G\QVA^k\supset ...\supset 0$$
of the complex $\G\QVA$ by  subcomplexes $\G\QVA^k$, such that
$\cap_{k\geq 0}\G\QVA^k=0$.
Let us prove that the subcomplex $\G\QVA^1$
is acyclic.

Take an element $y\in\G\QVA^1$, such that $dy=0$.
It can be presented as a convergent series
\begin{equation}\label{series}
y=\sum\limits_{j=0}^N\sum\limits_{k\geq 1}\omega^{j,k-j}(z)\ot v^j,
\end{equation}
where $v^j\in \V^j$ and  $\ \omega^{j,k-j}(z)$ is a homogeneous polynomial
 form  of degree $k-j$.
The convergence means that for every $j=0,...,N$ and for any
$\eta_j\in (\V^j)^*$ the series
$$\sum\limits_{k\geq 1}\omega^{j,k-j}(z)\mylangle v^j,\eta_j\myrangle$$
converges as a power series in  $z$.

Since the differential $d$ is homogeneous, we have $\sum_{j=0}^Nd(
\omega^{j,k-j}(z)\ot v^j)=0$
for any $k\geq 1$.
 Lemma  \ref{laplace} says that
for any $v^j\in\V^j$ and for any homogeneous  polynomial
form $\omega^k$, $\deg \omega^k=k$, we have
\begin{equation}\label{Laction}
\Delta(\omega^k\ot v^l)=\omega^k\ot (\overline{\Delta}^{(0)}+k+l) v^l
\end{equation}
Since the Laplace operator $\Delta $ commutes with $d$ and the eigenvalues of
$\overline{\Delta}^{(0)}$ are close to zero, \rf{Laction} implies also that
$ d(\omega^{j,k-j}(z)\ot v^j)=0$
for all $j$ and $k\geq 1$.

 The space $\V^j$ is invariant with respect to $\overline{\Delta}^{(0)}$,
$\overline{\Delta}^{(0)}(\V^j)\subset \V^j$. Let $\V^j=\oplus_\l \V^j_\l$
be the decomposition of  $\V^j$ in the generalized
 eigenspaces of the operator $\overline{\Delta}^{(0)}$,
 $$
\V^j_\l\ =\ \bigcup\limits_{k=0}^\infty
{\rm Ker}(\overline{\Delta}^{(0)}-\l)^k{}|_{\V^j}\ .
$$
Present the element $\omega^{j,k-j}(z)\ot v^j$ as
\begin{equation}\label{eigendec1}
\omega^{j,k-j}(z)\ot v^j=\omega^{j,k-j}(z)\ot\sum_\l v^j_\l,
\end{equation}
where $v^j_\l\in\V^j_\l$.
 Relation \rf{Laction} implies that $\omega^{j,k-j}(z)\ot v^j_\l\in
 \bigcup\limits_{m=0}^\infty
{\rm Ker}({\Delta}-\l-k)^m$ and $d(\omega^{j,k-j}(z)\ot v^j_\l)=0$.

Standard arguments with the Laplace operator say the following.
Let an element $x$ be such that $dx=0$ and let
$x$ belong to a finite-dimensional
 generalized eigenspace of the Laplace operator $\Delta=dd^*+d^*d$ with eigenvalue
$\mu\not=0$. Then $x$ can be presented as the boundary
$$
x=dd^*\left(\frac{x}{\mu}-\frac{(\Delta-\mu)x}{\mu^2}+
\frac{(\Delta-\mu)^2x}{\mu^3}-...\right)\ ,
$$
which is a finite sum.
This implies that for any $j=0,...,N$, $k\geq 1$, and eigenvalue $\l$ of
the operator $\overline{\Delta}^{(0)}|_{\V^j}$, we have
\begin{equation}\label{boundary}\begin{split}
\omega^{j,k-j}(z)\ot v^j_\l=dd^*\left(
\frac{\omega^{j,k-j}(z)\ot v^j}{\l+k}-
\frac{\omega^{j,k-j}(z)\ot(\overline{\Delta}^{(0)}-\l-k)
 v^j}{(\l+k)^2}+\right.\\ \left.
+\frac{\omega^{j,k-j}(z)\ot(\overline{\Delta}^{(0)}-\l-k)^2
 v^j}{(\l+k)^3}-
...\right)
\end{split}
\end{equation}
The sum in \rf{boundary} is finite and the number of terms does not depend on
$k$. The substitution of  \rf{boundary} to \rf{eigendec1} and then to
\rf{series} gives the presentation of  the element $y$ as the boundary
 of a convergent series of homogeneous elements. Convergence follows from the
 structure of \rf{boundary}.
This proves that the subcomplex  $\G\QVA^1\subset \G\QVA$ is acyclic and
the complex $\G\QVA$ is quasi-isomorphic to the quotient complex
$\G\QVA/\G\QVA^1$,  which is isomorphic to the quiver
 complex  $C_+(\VA)[N]$.
  This proves  Theorem \ref{theorem4}.
 \hfill{$\square$}

\subsection{Proof of Theorem \ref{Gprop5}}
The proof of Theorem \ref{Gprop5} repeats the arguments of the
proof of Theorem \ref{theorem4}. For the construction of the corresponding
Laplace operator we use the symmetric bilinear form on $X$, invariant
with respect to
the  $G$-action. Such a form exists, since $G$ is finite.
 Then the Laplace operator commutes with the
 $G$-action and establishes a quasi-isomorphism of $G$ -invariants of
the complexes. \hfill{$\square$}

\end{document}